\documentclass[11pt]{article}
\usepackage{amsmath}
\usepackage{graphicx}
\usepackage{enumerate}
\usepackage{natbib}
\usepackage{url} 
\usepackage[T1]{fontenc}
\usepackage{a4wide, url}
\usepackage[english]{babel}
\usepackage{graphicx}
\graphicspath{{}}
\usepackage{array}
\usepackage{setspace}
\usepackage{multirow}
\usepackage{amsthm}
\usepackage{amsmath}
\usepackage{amsfonts}
\usepackage{amssymb}
\usepackage{pgf}
\usepackage{bm}
\usepackage{paralist}
\usepackage{pdflscape}
\usepackage{xfrac}
\newcolumntype{P}[1]{>{\centering\arraybackslash}p{#1}}
\usepackage{bbm}
\usepackage{enumitem,kantlipsum}

\usepackage{booktabs}
\usepackage{float}

\usepackage{mwe} 
\usepackage{caption}

\usepackage{xr}
\DeclareMathAlphabet\mathbfcal{OMS}{cmsy}{b}{n}

\newcommand{\beq}{\begin{equation}}
\newcommand{\eeq}{\end{equation}}
\newcommand{\bea}{\begin{eqnarray}}
\newcommand{\eea}{\end{eqnarray}}
\newcommand{\ba}{\begin{array}}
\newcommand{\ea}{\end{array}}
\newcommand{\bit}{\begin{itemize}}
\newcommand{\eit}{\end{itemize}}
\newcommand{\ben}{\begin{enumerate}}
\newcommand{\een}{\end{enumerate}}
\newcommand{\bpm}{\begin{pmatrix}}
\newcommand{\epm}{\end{pmatrix}}
\newcommand{\bbm}{\begin{bmatrix}}
\newcommand{\ebm}{\end{bmatrix}}

\renewcommand{\l}{\left}
\renewcommand{\r}{\right}
\newcommand{\E}[0]{\mathrm{E}}

\newcommand{\nn}{\nonumber}
\newcommand{\wh}{\widehat}
\newcommand{\wt}{\widetilde}

\newcommand{\mc}{\mathcal}

\newcommand{\numberset}{\mathbb}
\newcommand{\N}{\numberset{N}}
\newcommand{\Z}{\numberset{Z}}

\newcommand{\R}{\numberset{R}}
\newcommand{\C}{\numberset{C}}

\newtheorem{thm}{Theorem}[section]
\newtheorem{coroll}{Corollary}[section]
\newtheorem{Lemma}{Lemma}[section]

\newtheorem{Ass}{Assumption}[section]
\newtheorem{Rem}{Remark}[section]

\newcommand{\blind}{1}

\begin{document}

\def\spacingset#1{\renewcommand{\baselinestretch}%
{#1}\small\normalsize} \spacingset{1}


\if1\blind
{
  \title{\textsc {An algebraic estimator\\ for large spectral density matrices}
  \footnote{This work was supported by the \emph{Accademia Nazionale dei Lincei} 
   under Grant \emph{British Academy 2018.}}}
  \author{Matteo Barigozzi\hspace{.2cm}\\
    Department of Economics, University of Bologna\\
    and \\
    Matteo Farn\`{e}\thanks{Corresponding author: \texttt{matteo.farne@unibo.it}} \\
    Department of Statistics, University of Bologna}
  \maketitle
} \fi

\if0\blind
{
  \bigskip
  \bigskip
  \bigskip
  \begin{center}
    {\Large\textsc {An algebraic estimator\\ for large spectral density matrices}}
\end{center}
  \medskip
} \fi

\bigskip
\begin{abstract}
We propose a new estimator of high-dimensional spectral density matrices, called UNshrunk ALgebraic Spectral Estimator (UNALSE), under the assumption of an underlying low rank plus sparse structure, as typically assumed in dynamic factor models. The UNALSE is computed by minimizing a quadratic loss under a nuclear norm plus $l_1$ norm constraint to control the latent rank and the residual sparsity pattern. The loss function requires as input the classical smoothed periodogram estimator and two threshold parameters, the choice of which is thoroughly discussed. We prove consistency of UNALSE as both the dimension $p$ and the sample size $T$ diverge to infinity, as well as algebraic consistency, i.e., the recovery of latent rank and residual sparsity pattern with probability one. The finite sample properties of UNALSE are studied by means of an extended simulation exercise as well as an empirical analysis of US macroeconomic data.
\end{abstract}

\noindent%
\vfill

\newpage

\section{Introduction}
An appealing, natural, and classical way to model  time series data is through spectral analysis \citep{brillinger2001time}. Given a $p$-dimensional vector stochastic process, its $p\times p$ spectral density matrix characterizes all second order dependencies. Moreover, conditional second order dependencies can also be extracted starting from the inverse of the spectral density matrix. The spectral approach is appealing since, once we move from the time domain to the frequency domain, data become asymptotically independent, as the sample size $n$ grows to infinity.

Statistical methods for the study of time series based on spectral analysis include: pseudo-maximum likelihood estimation \citep{dahlhaus2000likelihood,velasco2000whittle},
linear regression \citep{harvey1978linear}, cointegration tests or information criteria based on  the zero-frequency spectral density matrix of a vector of time series \citep{stock1988testing,barigozzi2021large}, and similarly seasonal cointegration tests based on the spectral density matrix at selected frequencies \citep{joyeux1992tests}, de-trending methods \citep{corbae2002band}, Granger causality tests \citep{breitung2006testing,farne2018bootstrap}, and the analysis of low frequency co-movements \citep{muller2018long}. Finally, the inverse spectral density matrix is at the basis of graphical models and dynamic network analysis \citep{granger1969investigating,dahlhaus2000graphical,eichler2007granger,davis2016sparse,barigozzi2019nets}.

The use of spectral analysis is widespread in many applied fields. Examples are the construction of business cycle indicators in macroeconomics \citep{sargent1977business,altissimo2010new}, portfolio optimization at different horizons in finance \citep{chaudhuri2015spectral}, and the study of brain activity in biostatistics \citep{ombao2001automatic,ombao2005slex,fiecas2011generalized,fiecas2016modeling}.

All above methods and applications require as input an estimator of the spectral density matrix or of its inverse. Just like for the covariance matrix estimation in time domain,  estimation of a spectral density matrix is a hard problem when the dimension of the process $p$ is comparable, or even larger, than the sample size $T$. In this case, the classical smoothed periodogram estimator is not positive definite simply due to lack of degrees of freedom. Given the increased availability of large datasets in the recent years, this issue becomes of fundamental importance. \citet{wu2018asymptotic} provide  consistency results for the smoothed periodogram estimator in high-dimension, which hold uniformly over all frequencies (see also \citealp{zhang2021convergence}).

To solve the problem of the curse of dimensionality, here, we start from observing that the second moments of most high-dimensional time series tend to have both a low rank and a sparse component. Indeed,
on the one hand, most economic datasets are known to be mainly ``dense'' rather than sparse \citep{giannone2017economic}. Moreover, there exist mathematical results proving that large dimensional panels of time series can always be represented as having a factor structure \citep{forni2001generalized,hallin2013factor}. On the other hand, once the common factors are controlled for, there is evidence of sparseness in the second order structure of the residuals \citep{barigozzi2017network}.

In this paper, we assume that the spectral density matrix of $p$-dimensional time series has the low rank plus sparse structure:
\[
\Sigma(\theta) = L(\theta)+S(\theta), \qquad \theta\in[-\pi,\pi],
\]
where $L(\theta)$ has rank $r$ independent of $p$ and such that $r<p$, and $S(\theta)$ is a sparse matrix. Based on this assumption, our estimators $\wh{L}(\theta)$ and $\wh{S}(\theta)$ of the two components of the spectral density matrix are obtained by regularizing the smoothed periodogram estimator, $\wt{\Sigma}(\theta)$, by means of  a nuclear norm plus $l_1$ norm penalization. Specifically, at each given frequency our estimators are defined as
\begin{equation}
\l(\wh{L}(\theta),\wh{S}(\theta)\r)=\arg\!\!\! \min_{\underline L(\theta),\underline S(\theta)} \frac{1}{2}\Vert \wt{\Sigma}(\theta)-(\underline L(\theta)+\underline S(\theta))\Vert_{F}^2
+\psi\Vert\underline L(\theta)\Vert_{*}+\rho\Vert\underline S(\theta)\Vert_{1}, \nn
\end{equation}
where $\underline L(\theta)$ and $\underline S(\theta)$ indicate generic values of the matrices belonging to appropriate algebraic matrix varieties,
$\Vert\underline L(\theta)\Vert_{*}=\text{tr}(\underline L(\theta))$  and $\Vert\underline S(\theta)\Vert_{1}=\sum_{i,j=1}^p\vert \underline S_{ij}(\theta)\vert$, and $\psi$ and $\rho$ are threshold parameters.

The above optimization problem is solved by iterating between a singular value thresholding step \citep{cai2010singular}, giving $\wh{L}(\theta)$, and a soft-thresholding step \citep{daubechies2004iterative}, giving $\wh{S}(\theta)$. The algorithm we employ has also been described in \citet{luo2011recovering} for the case of covariance estimation. We also apply the un-shrinkage step of estimated latent eigenvalues by \cite{farne2020large}, that optimizes the finite sample Frobenius loss with respect to the smoothed periodogram while retaining algebraic consistency. We call the resulting estimator of the overall spectral density matrix, $\wh{\Sigma}(\theta)=\wh{L}(\theta)+\wh{S}(\theta)$, UNshrunk ALgebraic Spectral Estimator (UNALSE).

We prove the algebraic and parametric consistency of UNALSE uniformly over frequencies, as both the dimension $p$ and the sample size $T$ diverge. By algebraic consistency, we mean that, with probability tending to $1$: (i) the UNALSE low rank estimate $\wh{L}(\theta)$ is positive semidefinite with the true rank $r$, (ii) the UNALSE residual estimate $\wh{S}(\theta)$ is positive definite having the true sparsity pattern, and (iii) $\wh{\Sigma}(\theta)$ is positive definite. The usual parametric consistency holds because UNALSE estimates are close to their targets in spectral norm (rescaled by the dimension $p$) with probability approaching $1$.  Our consistency results are obtained by generalizing to our framework the results of \cite{wu2018asymptotic} for the smoothed periodogram. We also provide a thorough discussion on the selection of the threshold parameters $\psi$ and $\rho$.

Our approach is based on the fundamental identifiability assumptions we make on the behavior of the eigenvalues of the spectral density matrix. We assume the $r$ eigenvalues of the low rank component, $L(\theta)$, to be diverging at a rate $p^\alpha$ with $\alpha\in [0,1]$, possibly different across frequencies. In the language of factor models, this means we are allowing for the presence of factors with different degrees of pervasiveness across frequencies, i.e., both weak and strong factors. Moreover, we assume the sparse component, $S(\theta)$, to have eigenvalues diverging at a rate $p^\delta$ with $\delta\in [0,1/2]$ and $\delta<\alpha$, possibly different across frequencies. These assumptions imply the existence of an eigen-gap in the spectrum of the spectral density matrix, $\Sigma(\theta)$, which has varying width across frequencies.

There exist alternative approaches to the estimation of large spectral density matrices.
 \citet{forni2000generalized} propose principal component analysis in the frequency domain to recover the low rank component, \cite{bohm2008structural,bohm2009shrinkage} propose to shrink the smoothed periodogram towards either a reduced rank target or the identity, respectively, \cite{fiecas2014data} propose a penalized likelihood approach, and \citet{fiecas2019spectral} consider constrained $l_1$ minimization for estimating the inverse. While some of those works assume either a low rank or a sparsity structure, none of them considers both assumptions jointly. Our approach is thus encompassing all the estimators assuming one of the two settings.

Similar approaches based on a low rank plus sparse assumption exist also in time domain, i.e., for the estimation of the covariance matrix.
\cite{fan2013large} consider principal components to recover the low rank component and then, in a second step, apply soft or hard thresholding to the orthogonal complement to obtain a sparse and positive definite residual. Their resulting estimator is called POET.
\cite{farne2020large} adopt a minimization algorithm analogous to the one considered in this paper which recovers the covariance matrix consistently, both algebraically and parametrically. Their resulting estimator is called UNALCE and they show that it systematically outperforms POET both in terms of parametric consistency, and, more importantly, because it provides the algebraic recovery of latent rank and sparsity pattern. A similar approach was proposed by \citet{luo2011recovering}, however it is based on the assumption of bounded eigenvalues for the
covariance matrix, which does not allow for the joint identification of the two components.

Our assumption of a low rank plus sparse decomposition of the spectral density matrix is strictly related to, and inspired by, the Generalized Dynamic Factor Model (GDFM) representation of a large panel of time series, originally proved by \citet{forni2001generalized}. This is a very popular approach to dimension reduction \citep[see, e.g., the application in][]{altissimo2010new}. In the GDFM, $r$ latent factors are loaded by each series in a dynamic way, i.e., not only contemporaneously but also with lags. The key assumptions are: (i) pervasiveness of the  factors modeled via $r$ leading spiking spectral eigenvalues, and (ii) weak serial and cross-correlation in the residuals, modeled via boundedness of the spectral eigenvalues. These in our notation imply $\alpha = 1$.

\citet{forni2000generalized,forni2005generalized,forni2017dynamic} consider different estimators of the GDFM, which are all built starting from a consistent estimator of the spectral density matrix. In particular, in all those approaches the low rank component of the spectral density is estimated via the $r$ leading dynamic principal components, i.e., the principal components of the spectral density matrix across frequencies of the smoothed periodogram (see also \citealp{brillinger2001time}). The consistency of this method relies on the pervasiveness of spectral eigenvalues with respect to the dimension $p$. The spectral density of the residual component, called idiosyncratic component in the GDFM literature, is then estimated as the difference between the estimated spectral density of the observed data and its estimated low rank component. Hence, by construction, the spectral density of the idiosyncratic component has rank $p-r$, i.e., it is not positive definite, and, therefore, not invertible. There exist also few papers dealing with determining the dynamic rank, $r$: \cite{hallin2007determining} propose an information criterion, and \citet{onatski2009testing} proposes a test based on the asymptotic distribution of the spectral eigenvalues.

The above approaches to the estimation of the GDFM suffer from some drawbacks.  First, any estimator of the spectral density matrix based on the principal components of an input estimator, like the smoothed periodogram,  is likely to suffer from numerical instability, especially if $p$ is large, due to the \citet{marvcenko1967distribution} law. Second, the strict pervasiveness assumption of spectral eigenvalues ($\alpha=1$) is rarely satisfied in practice, since both the factor number and their strength might vary across frequencies, e.g., due to common, frequency specific, features. Third, the weak correlation assumption increases the number of parameters when $p$ is large, which prevents the residual component to be identified.

The estimator we propose in this paper can be used as input of all the estimators of the GDFM considered in the literature, and, given its algebraic consistency, it is also a consistent estimator of the latent rank $r$. {Moreover, as already noticed above, our assumptions generalize the GDFM setting in that they are compatible with frequency specific factor numbers and strengths, and goes beyond a factor structure in that if no factor is present, then our method would not return any low rank plus sparse decomposition.}

The paper is organized as follows. In Section \ref{sec:res} we present our main results using the GDFM setting as a guiding example.
In Sections \ref{sec:model} and \ref{sec:cons} we present the general framework, describe estimation, and prove consistency. In Section \ref{sec:thresholds} we discuss the choice of the threshold parameters. Sections \ref{sec:sim} and \ref{sec:real} present numerical results for simulated and real datasets.

\subsection*{Notation}

Let us define a $p\times p$ symmetric positive-definite complex matrix ${M}$, and denote its transposed complex conjugate as $M^\dag$.
We denote by $\lambda_i({M})$, $i=1,\ldots,p$, the eigenvalues of ${M}$ in descending order (note that they are all real numbers), and by $M_{ij}$ the $(i,j)$th entry of $M$. We also define $\overline M_{ij}$ as the complex conjugate of $M_{ij}$, thus the complex modulus is $\vert M_{ij}\vert=\sqrt{M_{ij}\overline{M}_{ij}}$, while the real and imaginary parts are indicated as $\text{Re}(M_{ij})$ and $\text{Im}(M_{ij})$, respectively. To indicate that $M$ is positive definite or semidefinite we use the notations: ${M} \succ 0$ or ${M} \succeq 0$, respectively.

Element-wise norms:
$l_0$ norm: $\Vert  M \Vert_0 =\sum_{i=1}^p \sum_{j=1}^p \mathbbm{1}(M_{ij}\ne 0)$, which is the total number of nonzeros;
$l_1$ norm: $\Vert M \Vert_1=\sum_{i=1}^p \sum_{j=1}^p \vert M_{ij}\vert $; 
Frobenius norm: $\Vert M \Vert_F=\sqrt{\sum_{i=1}^p \sum_{j=1}^p \vert M_{ij}\vert^2}=\sqrt{\text{tr}(MM^\dag)}$; 
maximum norm: $\Vert M \Vert_\infty= \max_{1\le i,j\leq p} \vert M_{ij}\vert $. 
Vector-induced norms:
 $\Vert M \Vert_{0,v}= \max_{1\le i \leq p} \sum_{j =1}^p \mathbbm{1}(M_{ij} \ne 0)$, which is the maximum number of nonzeros per column,
defined also as the maximum ``degree'' of ${M}$; 
$\Vert M \Vert_{1,v}= \max_{1\le j \leq p} \sum_{i =1}^p \vert M_{ij}\vert $; 
$\Vert M \Vert_{\infty,v}= \max_{1\le i \leq p} \sum_{j =1 }^p \vert M_{ij}\vert $;
spectral norm: $\Vert M \Vert_{2}=\sqrt{\lambda_1({MM^\dag})}=\lambda_1({M})$; 
the nuclear norm: $\Vert M \Vert_{*}=\mbox{tr}(M)=\sum_{i=1}^p \lambda_i({M})$.
The minimum nonzero off-diagonal element of $M$ in absolute value is denoted as $\Vert M \Vert_{min,off}=\min_{ \substack{1\le i,j \leq p\\
i \ne j , M_{ij} \ne 0}}{\vert M_{ij} \vert}$.


\section{Main results}\label{sec:res}

In this section, we present the main features of our estimator under the assumption that the data follow a GDFM as defined by \cite{forni2001generalized} and \cite{hallin2013factor}. The GDFM setting has to be considered just as a motivating example, which is well suited to allow the reader to immediately appreciate the contribution of this paper with respect to the state of art. In the following sections, we present our theory in more detail showing that the validity of our results is actually much broader than the case here considered.

Let $X=\{X_{it}, i=1,\ldots ,p\,, t \in \Z\}$ be a $p$-dimensional panel of time series. We assume that for any $p\in\mathbb N$ the process $X$ is second-order stationary, and, without loss of generality, we also assume that $\E[X_{it}]=0$ and  $V(X_{it})=\E[{X_{it}^2}]>0$, for any $i\in\mathbb N$. The set of all $L_2$-convergent linear combinations of
$X_{it}$'s and their limits, as $p\to\infty$, of $L_2$-convergent sequences thereof,
is a Hilbert space, denoted by $\mathcal{H}_{X}$. Hence, for all $t\in\mathbb Z$ and all $p\in\mathbb N$, any dynamic linear combination of $X_{it}$s, $y_t=\sum_{i=1}^p \sum_{k=-\infty}^{\infty} a_{ik} X_{i,t-k}$, such that $\sum_{i=1}^p \sum_{k=-\infty}^{\infty} a_{ik}^2=1$, belongs to $\mathcal{H}_{X}$.
Following Definitions 2.1 and 2.2 in \cite{hallin2013factor}, we define as \emph{common variable} the $L_2$-limit of any standardized dynamic linear combination of the $X$s, say {$y_t/V(y_t)^{1/2}$}, such that $V(y_t)\to\infty$, as $p\to\infty$.
The Hilbert space of all \emph{common variables} is denoted by $\mathcal{H}_{\textit{com}}$,
while its orthogonal complement with respect to $\mathcal{H}_{X}$, denoted as $\mathcal{H}_{\textit{idio}}$, contains
all the \emph{idiosyncratic variables}, i.e., all dynamic linear combinations $y_t$
with bounded variance $V(y_t)$ for all $p\in\mathbb N$. \cite{hallin2013factor} prove that there exist two unique stochastic processes $\chi_{it}\in \mathcal{H}_{\textit{com}}$ and $\epsilon_{it}\in \mathcal{H}_{\textit{idio}}$, mutually orthogonal at all leads and lags, such that
\begin{equation}
X_{it}=\chi_{it}+\epsilon_{it}
\label{gdfm}
\end{equation}
for all $i \in \N$ and $t \in \Z$. The process $\chi_{it}$ is called \emph{common component}, the process $\epsilon_{it}$ is called \emph{idiosyncratic component}. Representation \eqref{gdfm} is the GDFM.  In vector terms, we can write $X_t=\chi_t+\epsilon_t$, where $X_t$, $\chi_t$
and $\epsilon_t$ are $p$-dimensional random vectors. The GDFM encompasses the \emph{approximate static} factor models of \cite{fan2013large}, as well as the \emph{exact dynamic} factor models of \cite{sargent1977business}.

Let us define the spectral density matrices of $\{\chi_t\}$ and of $\{\epsilon_t\}$ as
$$
L(\theta)=\frac 1{2\pi}\sum_{k=-\infty}^\infty \Gamma_\chi(k) \mathrm{e}^{-{\rm i}\theta k} \;\mbox{ and }\; S(\theta)=\frac 1{2\pi}\sum_{k=-\infty}^\infty \Gamma_\epsilon(k) \mathrm{e}^{-{\rm i}\theta k},\qquad \theta\in[-\pi,\pi],
$$
respectively, where $\Gamma_\chi(k)=\E[\chi_{t+k}\chi_t']$ and $\Gamma_\epsilon(k)=\E[\epsilon_{t+k}\epsilon_t']$.
\cite{forni2001generalized} prove that: (i) the common component $\chi_{it}$ is driven by a $r$-tuple of mutually orthogonal white noises loaded by a linear time filter, and $\epsilon_{it}$ is orthogonal to those white noises at all leads and lags, and (ii) $X_{it}$ follows the GDFM representation (\ref{gdfm}) if and only if the $r$ eigenvalues of $L(\theta)$ diverge almost everywhere across $[-\pi,\pi]$ as $p$ diverges, while the eigenvalues of  $S(\theta)$ remain bounded for all $p$.

Hence, the spectral density matrix of $X_t$ is such that $\Sigma(\theta)=L(\theta)+S(\theta)$, for any $\theta\in[-\pi,\pi]$, and, as usual in the GDFM literature, in this section
we adopt the assumption (relaxed later on) that the $r$ eigenvalues of $p^{-1}L(\theta)$ are bounded away from $0$ for all $p$
almost everywhere across the frequency range $[-\pi,\pi]$. Similarly, the definition of idiosyncratic variable leads to the condition $\Vert S(\theta)\Vert_2<\infty$ almost everywhere across the frequency range $[-\pi,\pi]$ for all $p$. These assumptions on $L(\theta)$ and $S(\theta)$ cause the fact
that the gap between the $r$th and the $(r+1)$th eigenvalue of the spectral density matrix $\Sigma(\theta)=L(\theta)+S(\theta)$ increases at each $\theta \in [-\pi,\pi]$ as $p$ diverges, making the recovery of the low rank component easier.

In this paper, we further control the idiosyncratic spectral density matrix $S(\theta)$ at each $\theta$ by enforcing $\Vert S(\theta)\Vert_{0,v}$ to be bounded and finite for all $p$. Since $\Vert S(\theta)\Vert_{2} \leq \Vert S(\theta)\Vert_{0,v}$, the original assumption $\Vert S(\theta)\Vert_2<\infty$ still holds. This is done in order to enforce element-wise sparsity on $S(\theta)$ at each $\theta$.


Suppose now that we observe a sample of $p$-dimensional data vectors with size $T$.
A classical estimator of the spectral density matrix, which is our \emph{pre-estimator},
is the smoothed periodogram, defined as
\begin{equation}
\widetilde{\Sigma}(\theta_h)=\frac 1{2\pi}\sum_{k=-(T-1)}^{T-1} K\left(\frac{k}{M_T}\right) \mathrm{e}^{-{\rm i} \theta_h k}\, \wh{\Gamma}_X(k),\qquad \theta_h=\frac{h\pi}{M_T}, \quad \vert h\vert\le \lfloor M_T\rfloor,\label{kernel}
\end{equation}
where $\wh{\Gamma}_X(k)=T^{-1}\sum_{t=1}^{T-\vert k\vert } X_t X_{t+k}'$, and  $K(\cdot)$ is a suitable kernel function with $M_T$ being the associated smoothing parameter.
According to \cite{brillinger2001time}, for any given $\theta_h$, $\widetilde{\Sigma}(\theta_h)$ is
consistent if $\tfrac{M_T}{T}\rightarrow 0$ while $M_T\rightarrow\infty$ and $T\rightarrow\infty$. \citet{wu2018asymptotic} prove the consistency
of $\widetilde{\Sigma}(\theta_h)$ uniformly over the frequencies, under appropriate assumptions to be discussed later.

Under the GDFM setting described above, augmented with the sparsity assumption for $S(\theta)$, we define the UNshrunk ALgebraic Spectral Estimator (UNALSE) estimator of the spectral density matrix $\Sigma(\theta)$ as $\wh{\Sigma}(\theta)=\wh{L}(\theta)+\wh S(\theta)$, where
$\wh{L}(\theta)$ and $\wh S(\theta)$ are such that:
\begin{equation}
\l(\wh{L}(\theta),\wh{S}(\theta)\r)=\arg\!\!\! \min_{\underline L(\theta),\underline S(\theta)} \frac{1}{2}\Vert \wt{\Sigma}(\theta)-(\underline L(\theta)+\underline S(\theta))\Vert_{F}^2
+\psi\Vert\underline L(\theta)\Vert_{*}+\rho\Vert\underline S(\theta)\Vert_{1},  \label{func:ob_spec}
\end{equation}
where
$\underline L(\theta)$ and $\underline S(\theta)$ indicate generic values of the matrices belonging to appropriate algebraic matrix varieties (see Section \ref{sec:model} for the details),
while $\psi$ and $\rho$ are threshold parameters. The minimization problem \eqref{func:ob_spec} is a non-smooth convex optimization problem
which is the tightest convex relaxation of the following NP-hard problem:
\begin{equation}
\min_{\underline L(\theta),\underline S(\theta)} \frac{1}{2}\Vert \wt{\Sigma}(\theta)-(\underline L(\theta)+\underline S(\theta))\Vert_{F}^2
+\psi\, \text{\upshape{rk}}(\underline L(\theta)) + \rho \Vert\underline S(\theta)\Vert_{0},\label{prob:orig_spec}
\end{equation}
which would be the natural target under the low rank plus sparse assumption. Indeed, we know that:
\begin{inparaenum}[(i)]
\item
$\Vert\underline S(\theta)\Vert_{1}$ is the tightest convex relaxation of $\Vert\underline S(\theta)\Vert_{0}$ 
\citep{donoho2006most};
\item $\Vert\underline L(\theta)\Vert_{*}$ is the tightest convex relaxation of $\text{\upshape rk} (\underline L(\theta))$ \citep{fazel2001rank}.
\end{inparaenum}

In practice, the solution of \eqref{func:ob_spec} is computed as follows. For any given frequency $\theta_h= \frac{\pi h}{M_T}$, with $\vert h\vert\le \lfloor M_T\rfloor$, we apply the following iterative procedure:
\begin{enumerate}
\item set $\left(L_{0}(\theta_h),S_{0}(\theta_h)\right)=\l(\frac{\text{\it diag}(\wt{\Sigma}(\theta))}{2},\frac{\text{\it diag}(\wt{\Sigma}(\theta))}{2}\r)$, $\eta_0=1$, and  initialize $Y_{0}(\theta_h)=L_{0}(\theta_h)$ and $Z_{0}(\theta_h)=S_{0}(\theta_h)$;
\item for $k\ge 0$, repeat:
\begin{compactenum}
\item compute $\frac{\partial \frac{1}{2}\Vert Y_{k-1}(\theta_h)+Z_{k-1}(\theta_h)-\wt{\Sigma}(\theta_h)\Vert_{F}^2}{\partial Y_{k-1}(\theta_h)}=\frac{\partial \frac{1}{2}\Vert  Y_{k-1}(\theta_h)+Z_{k-1}(\theta_h)-\wt{\Sigma}(\theta_h)\Vert_{F}^2}{\partial Z_{k-1}(\theta_h)}=Y_{k-1}(\theta_h)+Z_{k-1}(\theta_h)-\wt{\Sigma}(\theta_h)$;
\item apply the \emph{singular value thresholding} operator of \cite{cai2010singular} $T_{\psi}(\cdot)$ to $\mathcal E_{Y,k}(\theta_h)=Y_{k-1}(\theta_h)- \frac{1}{2}(Y_{k-1}(\theta_h)+Z_{k-1}(\theta_h)-\wt{\Sigma}(\theta_h))$ and set $L_{k}(\theta_h)=T_{\psi}(\mathcal E_{Y,k}(\theta_h))$;
\item apply the \emph{soft-thresholding} operator of \cite{daubechies2004iterative} $T_\rho(\cdot)$ to $\mathcal E_{Z,k}(\theta_h)=Z_{k-1}(\theta_h)- \frac{1}{2}(Y_{k-1}(\theta_h)+Z_{k-1}(\theta_h)-\wt{\Sigma}(\theta_h))$ and set $S_{k}(\theta_h)=T_\rho({\mathcal E_{Z,k}(\theta_h)})$;
\item set $(Y_{k}(\theta_h),Z_{k}(\theta_h))=(L_{k}(\theta_h),S_{k}(\theta_h))+\frac{\eta_{k-1}-1}{\eta_{k}}[(L_{k}(\theta_h),S_{k}(\theta_h))-(L_{k-1}(\theta_h),S_{k-1}(\theta_h))]$
where $\eta_{k}=\frac{1+\sqrt{1+4 \eta_{k-1}^2}}{2}$;
\item stop if $\frac{\Vert L_{k}-L_{k-1}\Vert_{F}}{\Vert 1+L_{k-1}\Vert_{F}}+\frac{\Vert S_{k}-S_{k-1}\Vert_{F}}{\Vert 1+S_{k-1}\Vert_{F}} \leq \varsigma$, where $\varsigma$ is a prescribed precision level (we set $\varsigma=0.01$);
\end{compactenum}
\item set $\wh{L}(\theta_h)=Y_{k}(\theta_h)$ and $\wh{S}(\theta_h)=Z_{k}(\theta_h)$.
\end{enumerate}

The two thresholding operators introduced in the above algorithm are defined as follows.
\ben
\item [(I)] Singular value thresholding operator: let the singular value decomposition of a positive semi-definite complex symmetric matrix $M$ be $M=U_M \Lambda_M U_M^\dag$, then, define $T_{\psi}(M)$=$U_M \Lambda_{\psi,M} U_M^\dag$, where $\Lambda_{\psi,M}$ is a diagonal matrix with $i$th diagonal element $\Lambda_{\psi,ii,M}= \max{(\Lambda_{ii,M}-\psi,0)}$.
\item [(II)] Soft-thresholding operator: for a positive definite complex symmetric $M$ define\linebreak $T_\rho(M_{ij})=\frac{M_{ij}}{ \sqrt{M_{ij}\overline M_{ij}}}
\max \l(\sqrt{ M_{ij}\overline M_{ij}}-\rho,0\r)$.
\een

In this paper we prove the following results for the UNALSE estimator.

\begin{thm}\label{th:introTh} For all $p\in\mathbb N$, assume that: (i) the $r$ nonzero eigenvalues of $L(\theta)$ are such that $\frac{\lambda_j(L(\theta))}p$ is finite and bounded away from zero for all $j=1,\ldots, r$, and (ii) $\Vert S(\theta)\Vert_{0,v}$ is bounded. Then, under the regularity conditions in Section \ref{sec:cons}, there exist finite positive reals $G_1$, $G_2$, $G_3$, $G_4$, and $G_5$, independent of $p$ and $T$, such that, as $p,T\to\infty$, with probability approaching $1$, for $\theta_h=\tfrac{h\pi}{M_T}$:
\begin{compactenum}
\item $\text{\upshape {rk}}(\wh{L}(\theta_h))=\text{\upshape {rk}}({L}(\theta_h))=r$ and $\max_{\vert h\vert \le \lfloor M_T\rfloor} \frac{1}{p}\Vert \wh{L}(\theta_h)-{L}(\theta_h)\Vert_{2}\leq G_1 \sqrt{\frac{M_T \log(M_T)}{T}}$;
\item $\max_{\vert h\vert \le \lfloor M_T\rfloor} \frac{1}{p}\Vert \wh{S}(\theta_h)-{S}(\theta_h)\Vert_{\infty}\leq G_2 \sqrt{\frac{M_T \log(M_T)}{T}}$, and, consequently, \\
$\max_{\vert h\vert \le \lfloor M_T\rfloor} \frac{1}{p}\Vert \wh{S}(\theta_h)-{S}(\theta_h)\Vert_{2}\leq G_2 \Vert{S}(\theta_h)\Vert_{0,v} \sqrt{\frac{M_T \log(M_T)}{T}}$;
\item
$\max_{\vert h\vert \le \lfloor M_T\rfloor} \frac{1}{p}\Vert  \wh{\Sigma}(\theta_h)-\Sigma(\theta_h)\Vert_{2} \leq G_3 \sqrt{\frac{M_T \log(M_T)}{T}}$.
\end{compactenum}
Furthermore, if $\frac{\lambda_{p}(S(\theta))}{p}\ge 2 G_4 \sqrt{\frac{M_T \log(M_T)}{T}}$ and $\frac{\lambda_{p}(\Sigma(\theta))}{p}\ge 2 G_5 \sqrt{\frac{M_T \log(M_T)}{T}}$, then, as $p,T\to\infty$,
with probability approaching $1$, for $\theta_h=\tfrac{h\pi}{M_T}$:
\begin{compactenum}
\item [4.] $\wh{S}(\theta_h)$ is positive definite and $\max_{\vert h\vert \le \lfloor M_T\rfloor}\frac{1}{p}\Vert \wh{S}(\theta_h)^{-1}~-~ S(\theta_h)^{-1}\Vert_{2} \leq G_4 \sqrt{\frac{M_T \log(M_T)}{T}}$;
\item [5.] $\wh{\Sigma}(\theta_h)$ is positive definite and $\max_{\vert h\vert \le \lfloor M_T\rfloor}\frac{1}{p}\Vert \wh{\Sigma}(\theta_h)^{-1}-\Sigma(\theta_h)^{-1}\Vert_{2} \leq G_5 \sqrt{\frac{M_T \log(M_T)}{T}}$.
\end{compactenum}
\end{thm}

The convergence speed depends on the rapidity of decay of the physical dependence, on the finite moment of highest order, and on the smoothing parameter $M_T$ (see Section \ref{sec:cons} for further details and a more detailed exposition). We stress that the standard condition $\frac{M_T} {T} \rightarrow 0$ (as $M_T,T \rightarrow \infty$) ensures that the relative errors in the above Theorem vanish asymptotically uniformly over the frequency grid. Moreover, if $p\sqrt{\frac{M_T}{T}}\to 0$, the conditions for parts 4 and 5 to hold are certainly satisfied as $p$ diverges.

The results of Theorem \ref{th:introTh} contribute to the literature in three ways. First, the exact dynamic rank recovery in part 1 allows to bypass the use of existing criteria for determining the number of factors, like those by \cite{hallin2007determining} and \cite{onatski2009testing}. Second, assuming that $S(\theta)$ is full rank, we derive a consistency result also for the estimator of the idiosyncratic spectral density  $\wh{S}(\theta)$, which allows to obtain the same error bound also for the overall estimator $\wh{\Sigma}(\theta)$. Third, we obtain results also for the estimators of the inverse spectral densities.

The validity of our estimation framework lies well beyond the standard GDFM assumptions.
First, we can relax the strict pervasiveness assumption on latent dynamic factors, by allowing the $r$ eigenvalues of the matrix $p^{-\alpha}L(\theta)$, with $\alpha\le 1$, to be bounded away from $0$ almost everywhere across the frequency range $[-\pi,\pi]$.
Second, we allow for the maximum number of nonzeros per row in $S(\theta)$,
$\Vert S(\theta)\Vert_{0,v}$, to be at most proportional to $p^{\delta}$, with
$ \delta \in[0,1/2]$ and $\delta<\alpha$. This means that we allow the idiosyncratic spectrum
to be quite far from the diagonal matrix. Our setting reduces to the GDFM one
when $\alpha=1$.  Third, the rank $r$ can be in fact frequency dependent.

\section{Model setup}\label{sec:model}

The aim of this paper is estimating the spectral density matrix of a $p$-dimensional process $X=\{ X_{it}, i=1,\ldots, p,\, t \in \Z\}$. We consider the following data generating process for $X$:
\begin{align}
X_t&=\chi_t+\epsilon_t,\quad\quad\quad\;\;\, t\in\Z,\label{mod_base}\\
\chi_t&=\sum_{s=0}^{\infty}B_s u_{t-s},\qquad t\in\Z,\label{mod_chi}\\
\epsilon_t&=\sum_{s=0}^{\infty}C_s e_{t-s},\qquad t\in\Z,\label{mod_eps}
\end{align}
where $X_t$, $\chi_t$, $\epsilon_t$, and $e_t$ are $p$-dimensional, $u_t$ is $r$-dimensional, the $B_s$ are $p\times r$, and the $C_s$ are $p\times p$. Note that, differently from the original works on GDFM by \citet{forni2000generalized} and \citet{forni2001generalized} who derive \eqref{mod_chi} using two-sided filters, here we follow more the recent works by \citet{hallin2013factor} and \citet{forni2017dynamic} who show that also a one-sided representation is possible. We make the following assumptions on the processes $\{u_t\}$ and $\{e_t\}$ in \eqref{mod_chi} and \eqref{mod_eps}.

\begin{Ass}\label{assGDFM}
\begin{compactenum}
\item [(i)] $\{u_t, t\in\Z\}$ is a $r$-dimensional independent process with $r<p$ and independent of $p$, $\E(u_t)={0}_r$ and $\E[u_t u_t']={I_r}$;

\item [(ii)] there exists a $K_u$ independent of $j$ such that $\E[|u_{jt}|^{4+\epsilon}]\le K_u$ for some $\epsilon>0$ and $j=1,\ldots,r$;

\item [(iii)] $\{e_t, t\in\Z\}$ is a $p$-dimensional independent process with $\E(e_t)={0}_p$ and $\E[e_t e_t']={I_p}$;

\item [(iv)] there exists a $K_e$ independent of $j$ such that $\E[|e_{jt}|^{4+\epsilon}]\le K_e$ for some $\epsilon>0$ and $j=1,\ldots,p$;

\item [(v)] $\{u_{t}\}$ and $\{e_t\}$ are two mutually independent processes.

\end{compactenum}
\end{Ass}

Under Assumption \ref{assGDFM}, processes \eqref{mod_chi} and \eqref{mod_eps} are zero-mean linear and weakly stationary, and consequently process (\ref{mod_base}) also is.

We define for $z\in\mathbb C$ the filter of the common component  as $B(z)=\sum_{s=0}^{\infty}B_s z^s$, with $B_s$ being $p\times r$, and the filter of the idiosyncratic component as
$C(z)=\sum_{s=0}^{\infty}C_s z^s$, with $C_s$ being $p\times p$.
We develop our theory under two different settings, imposing a basic linear shape for the filters, where all matrices $B_s'B_s$ have the same condition number and all matrices $C_sC_s'$ have the same sparsity pattern,
and a general linear shape, where both features are allowed to vary across frequencies.

We start from the first setting.
\begin{Ass}[Basic linear filters]\label{basic}
\begin{inparaenum}
\item [(i)] $B_s=U_L\Lambda_{L,s}$, where $U_L$ is a $p\times r$ matrix such that $U_L'U_L=I_r$, and
$\Lambda_{L,s}=\lambda_s \sqrt{\Lambda_u}$ with $\lambda_s \in\mathbb R$ and such that $\sum_{s=0}^\infty \lambda_s^2=1$ and with $\Lambda_u$ a diagonal $r\times r$  positive definite matrix;
\item [(ii)] $C_s=U_S\Lambda_{S,s}$, where $U_S$ is a $p\times p$ matrix such that $U_S'U_S=I_p$ and $\Vert U_SU_S'\Vert_0=q$ with $q \ll p^2$, and $\Lambda_{S,s}=\lambda_s \sqrt{\Lambda_e}$ with $\lambda_s \in\mathbb R$ and such that $\sum_{s=0}^\infty \lambda_s^2=1$
and with $\Lambda_e$ a diagonal $p\times p$ positive definite matrix.
\end{inparaenum}
\end{Ass}

It immediately follows that $\E[\chi_t\chi_t']=\Gamma_\chi(0)=U_L \Lambda_u U_L'$ has rank $r$.
At the same time, $\E[\epsilon_t\epsilon_t']=\Gamma_\epsilon(0)=U_S \Lambda_e U_S'$
is sparse with $q$ nonzero elements. Notice that in this simple setting the scalar $\lambda_s$ is the same for both filters.
Moreover, since the filters are linear, for all $k\in\mathbb Z$, we have:
\begin{align}
\E[\chi_t\chi_{t+k}']&=\Gamma_\chi(k)=\sum_{s=0}^\infty B_sB_{s+k}'=U_L\l(\sum_{s=0}^\infty \lambda_s \lambda_{s+k}\r)\Lambda_u U_L'=a_k \Gamma_\chi(0),\nn\\
\E[\epsilon_t\epsilon_{t+k}']&=\Gamma_\epsilon(k)=\sum_{s=0}^\infty C_sC_{s+k}'=U_S\l(\sum_{s=0}^\infty \lambda_s \lambda_{s+k}\r)\Lambda_e U_S'=a_k \Gamma_\epsilon(0),\nn
\end{align}
where $a_k=\sum_{s=0}^\infty \lambda_s \lambda_{s+k}$, with $a_{-k}=a_k$ and $a_0=1$. Obviously, $\vert a_k\vert \leq a_0$, because $\Vert \Gamma_\chi(k)\Vert_{2} \leq \Vert \Gamma_\chi(0)\Vert_{2} $, for all $k$.
Since, by Assumption \ref{assGDFM}, $\{u_t\}$ and $\{e_t\}$ are uncorrelated processes,
we obtain $\E[X_tX_{t+k}']=\Gamma_X(k)=\Gamma_\chi(k)+\Gamma_\epsilon(k)=a_k(\Gamma_\chi(0)+\Gamma_\epsilon(0))$ for all $k$, which has a low rank plus sparse structure.

If we define $a(\theta)=\frac{1}{2\pi}\sum_{k=-\infty}^{\infty}a_k \mathrm{e}^{-\rm{i}\theta k}$, for $\theta\in[-\pi,\pi]$,  then the spectral density matrices of $\{\chi_t\}$ and $\{\epsilon_t\}$ are $L(\theta)= a(\theta)\Gamma_\chi(0)$ and $S(\theta)=a(\theta)\Gamma_\epsilon(0)$
respectively. Therefore, the spectral density matrix of $\{X_t\}$, which is $\Sigma(\theta)=L(\theta)+S(\theta)$,
has a low rank plus sparse structure at all frequencies $\theta\in[-\pi,\pi]$. Note that Assumption \ref{basic} describes the simplified case where the spectrum has no phase component, i.e., it is a real matrix, since $a(\theta)$ is real because $\mathrm{e}^{-\rm{i}\theta k}+\mathrm{e}^{i\theta k}$ is the sum of two waves of opposite argument.

Alternatively, we assume a more general structure for the filters.

\begin{Ass}[Generalized linear filters]\label{gen}
\begin{inparaenum}
\item [(i)] $B_s=U_L\Lambda_{L,s}$, where $U_L$ is a $p\times r$ matrix such that $U_L'U_L=I_r$, and $\sum_{s=0}^{\infty}{\Lambda_{L,s}^2}$ is a diagonal positive definite $r\times r$ matrix, $\Vert\sum_{s=0}^\infty B_s\Vert_2\le M_B(p)$;
\item [(ii)] $C_s=U_{S,s}\Lambda_{S,s}$, where $U_{S,s}$ is a $p\times p$ matrix such that $U_{S,s}'U_{S,s}=I_r$ and $\Vert U_{S,s}U_{S,s}'\Vert_0=q_s$ with $q_s \ll p^2$, and $\sum_{s=0}^{\infty}{\Lambda_{S,s}^2}$ is a diagonal positive definite $p\times p$ matrix, $\Vert\sum_{s=0}^\infty C_s\Vert_2\le M_C(p)$.
\end{inparaenum}
\end{Ass}

We refer to Assumption \ref{ass3} below for the definition of the terms $M_B(p)$ and $M_C(p)$, which for a fixed $p$ are positive constants, but might diverge as $p$ diverges.

Under Assumption  \ref{gen}, for all $k\in\mathbb Z$, we have
$\E[\chi_t\chi_{t+k}']=\Gamma_\chi(k)=\sum_{s=0}^\infty B_sB_{s+k}'=\sum_{s=0}^\infty U_L \Lambda_{L,s}\Lambda_{L,s+k} U_L'$, and
$\E[\epsilon_t\epsilon_{t+k}']=\Gamma_\epsilon(k)=\sum_{s=0}^\infty C_sC_{s+k}'=\sum_{s=0}^\infty U_{S,s}\Lambda_{S,s}\Lambda_{S,s+k} U_{S,s+k}'$,
thus leading to $\Gamma_\chi(0)=U_L (\sum_{s=0}^\infty \Lambda_{L,s}^2) U_L'$, which has rank $r$ by assumption, and $\Gamma_\epsilon(0)=\sum_{s=0}^\infty U_{S,s}\Lambda_{S,s}^2 U_{S,s}'$,
which is sparse with $q_{\Gamma_\epsilon(0)}=\sum_{s=0}^{\infty} q_s$ nonzero elements, while
$\Gamma_\epsilon(k)$ has $q_{\Gamma_\epsilon(k)}$ nonzero elements.


Finally,
the spectral density matrix of $\{\chi_t\}$ is:
$L(\theta)=\frac 1{2\pi}U_L\l(\sum_{k=-\infty}^{\infty}\l(\sum_{s=0}^\infty \Lambda_{L,s} \Lambda_{L,s+k}\r)\mathrm{e}^{-\rm{i}\theta k}\r)U_L'
=U_L A(\theta) U_L'$, where
$A(\theta)=\frac 1{2\pi}\sum_{k=-\infty}^{\infty}(\sum_{s=0}^\infty \Lambda_{L,s} \Lambda_{L,s+k}) \mathrm{e}^{-\rm{i}\theta k}$. Therefore, $L(\theta)$ has rank $r$ and has the same orthogonal base as $\Gamma_\chi(0)$ at all frequencies $\theta\in[-\pi,\pi]$.  Although at first sight this might be a restrictive design we notice that all results in the next section hold locally in the algebraic variety of rank $r$ matrices (see \citet{chandrasekaran2012} and Remark \ref{rem:chandra}, below).
At the same time, we obtain
$S(\theta)=\frac 1{2\pi}\sum_{k=-\infty}^\infty\l(\sum_{s=0}^\infty U_{S,s}\Lambda_{S,s}\Lambda_{S,s+k}' U_{S,s+k}' \mathrm{e}^{-\rm{i}\theta k} \r)$, which has $q$ nonzero elements, where $q\leq q_{\Gamma_{\epsilon}(0)}+2\sum_{k=1}^{\infty} q_{\Gamma_{\epsilon}(k)}$.
Therefore, the spectral density matrix of $\{X_t\}$, which is $\Sigma(\theta)=L(\theta)+S(\theta)$,
has a low rank plus sparse structure at all frequencies $\theta\in[-\pi,\pi]$.

%

As a consequence of Assumptions \ref{assGDFM} and \ref{basic} for the basic filter specification, or of Assumptions \ref{assGDFM} and \ref{gen} for the general filter specification, the spectral density matrix has a low rank plus sparse structure. To make this formal we introduce the following algebraic matrix varieties:
\begin{align}
\mathcal{L}(\mathsf r)&=\{L \in \R^{p \times p} \mid {L} \succeq 0, L=UDU^\dag, U\in \C^{p
\times \mathsf r}, U^\dag U=I_{\mathsf r}, D \in \R^{\mathsf r \times \mathsf r} \text{diagonal}, \mathsf r<p\},\label{var:L}\\
\mathcal{S}(\mathsf q)&=\{S\in\C^{p\times p} \mid {S} \succeq 0, \vert \text{supp}(S)\vert  \leq
\mathsf q, \mathsf q<p^2\}.\label{var:S}
\end{align}
In other words, $\mathcal{L}(\mathsf r)$ is the variety of Hermitian matrices with at most
rank $\mathsf r$ and $\mathcal{S}(\mathsf q)$ is the variety of Hermitian complex sparse matrices with at
most $\mathsf q$ nonzero elements ($\text{supp}(S)$ is the orthogonal complement of $\text{{ker}}(S)$). Therefore, under our assumptions $L(\theta)\in \mathcal{L}(r)$ and $S(\theta)\in \mathcal{S}(q)$, for all $\theta\in[-\pi,\pi]$, and where $r$ is defined in Assumptions \ref{basic}(i) or \ref{gen}(i), and $q$ is defined in Assumption \ref{basic}(ii) or under Assumption \ref{gen}(ii).

\begin{Rem}[Rank and sparsity] \upshape{Notice that while our model assumptions, which combine Assumptions \ref{assGDFM} and \ref{basic} for the basic filter specification, as well as Assumptions \ref{assGDFM} and \ref{gen} for the general filter specification, prescribe a fixed and common rank $r$ and sparsity degree $q$ for the spectral density matrices across frequencies, our methodology is still able to consistently recover latent structures with frequency-varying rank and sparsity. However, we avoid this distinction here to keep the notation simple. }
\end{Rem}


For simplicity, hereafter we adopt the following notation $L^*= \Gamma_\chi(0)$ and $S^*= \Gamma_\epsilon(0)$. Under Assumption \ref{basic}, $L(\theta)=a(\theta)L^*$, which has eigenvalues $\lambda_r(L(\theta))=\lambda_r(L^*)a (\theta)$. Since $L^*=U_L\Lambda_uU_L'$, in order to characterize the behavior of $\lambda_r(L(\theta))$ we just need to focus on the diagonal elements of $\Lambda_u$ (notice that $a(\theta)$ is a positive scalar).

\begin{Ass}\label{ass1}
For all $p\in\mathbb N$, the elements of the $r \times r$ diagonal matrix $\Lambda_u$,
$\Lambda_{u,j}$, are such that for all $j=1,\ldots, r$,
$M^{\min}_j \le \frac{\Lambda_{u,j}}{p^{\alpha}}\le M^{\max}_j$,
with $0 \leq \alpha \leq 1$, and $M^{\min}_j$ and $M^{\max}_j$ independent of $p$ and $T$ such that $M^{\min}_j \ge M^{\max}_{j+1}$, for $j=1,\ldots,(r-1)$.
In addition, there exists some $\kappa_a>0$ such that $a(\theta) \leq \kappa_a$ for all $\theta \in [-\pi,\pi]$. Moreover, $r$ is finite for all $p\in\mathbb N$.
\end{Ass}

Under Assumption \ref{gen}, $L(\theta)=U_L A(\theta) U_L'$, we characterize the eigenvalues of $L(\theta)$ by means of the following assumption (notice that $A(\theta)$ is diagonal with positive entries).

\begin{Ass}\label{ass1bis} For all $p\in\mathbb N$ and all $\theta\in[-\pi,\pi]$,
the elements, $A_j(\theta)$, of the $r \times r$ diagonal matrix $A(\theta)$, are such that for all $j=1,\ldots, r$,
$M^{\min}_j \le \frac{A_{j}(\theta)}{p^{\alpha}}\le M^{\max}_j$,
with $0 \leq \alpha \leq 1$, and $M^{\min}_j$ and $M^{\max}_j$ independent of $p$ and $T$ such that $M^{\min}_j \ge M^{\max}_{j+1}$, for $j=1,\ldots,(r-1)$. Moreover, $r$ is finite for all $p\in\mathbb N$.
\end{Ass}

In other words, under both assumptions, we are assuming that the eigenvalues of $L(\theta)$, which are all real, are of order $p^{\alpha}$, $\alpha \in [0,1]$. In this way we generalize the strict pervasiveness of the latent factors, corresponding to the case $\alpha=1$, necessary to ensure the effectiveness of the recovery of the number of factors in the methods by \cite{hallin2007determining} and \cite{onatski2009testing}. We notice that our results would hold even for $r=O(\log(p))$, however, as common in the literature, and as proved in \citet{forni2001generalized} for the GDFM, we consider $r$ as fixed and independent of $p$ to avoid unnecessary complications.

In order to give a rigorous definition of sparsity, we need to introduce further notation. The tangent spaces to $\mathcal{L}(r)$ and $\mathcal{S}(q)$ in two generic matrices $L \in \mathcal L(r)$ and $S \in\mathcal S(q)$, defined in \eqref{var:L} and \eqref{var:S}, are respectively defined as:
\begin{align}
T(L)&=\{M \in \C^{p \times p} \mid M=U Y_1^\dag+Y_2 U^\dag \mid Y_1,Y_2 \in \C^{p \times
r}, U^\dag LU \in \C^{r \times r} \mbox{diagonal}, L\in \mathcal L(r)\},\\
\Omega(S)&=\{N \in \C^{p \times p} \mid \text{ supp}(N)\subseteq
\text{supp}(S), S\in\mathcal S(q)\}.
\end{align}
The following uncertainty principle holds \citep{chandrasekaran2011rank}:
if $L(\theta)$ is nearly sparse, $S(\theta)$ cannot be recovered, and
if $S(\theta)$ is nearly low rank, $L(\theta)$ cannot be recovered.
Therefore, in order to achieve consistency we need to control for the
spikiness of the eigenvalues of $S(\theta)$ and the sparsity pattern
of $L(\theta)$. To this end, we make use of the following rank-sparsity incoherence measures between $T(L)$ and $\Omega(S)$, introduced in \cite{chandrasekaran2011rank}:
$\xi(T(L))=\max_{\substack{N \in T(L)\\ \Vert N\Vert_{2} \leq 1}} {\Vert N \Vert_{\infty}}$, $\mu(\Omega(S))=\max_{\substack{N \in \Omega(S)\\ \Vert N \Vert_{\infty} \leq 1}}
{\Vert N\Vert_{2}}$.
In order to identify $L(\theta)$ and $S(\theta)$, we need to control these rank-sparsity incoherence measures, which in turn implies defining the admissible sparsity patterns for $S(\theta)$. Indeed, as proved in Section \ref{sec:cons}, a necessary condition
to ensure both parametric and algebraic consistency is
\beq\label{incoerenza}
\xi(T({L}(\theta)))\ \mu(\Omega({S}(\theta))\leq \frac{1}{54}, \quad \theta \in [-\pi,\pi],
\eeq
which guarantees that ${L}(\theta)$ is far from sparsity and $S(\theta)$ is far from rank-deficiency. Indeed, the smaller is the product between the dual norms $\xi(T({L}(\theta)))$ and $\mu(\Omega({S}(\theta))$, the closer the two spaces $\mathcal{L}(r)$ and $\mathcal{S}(s)$ are to orthogonality, thus making easier to perform low rank plus sparse decomposition effectively.

Under Assumption \ref{basic}, $T(L(\theta))=T(L^*)$ and $\Omega(S(\theta))=\Omega(S^*)$, for all $\theta\in[-\pi,\pi]$. We then make the following assumption.

\begin{Ass}\label{ass2}
For all $p\in\mathbb N$, there exist:
\begin{inparaenum}
\item [(i)] $\delta \in [0,\frac{1}{2}]$, with $\delta < \alpha$, and $\delta_2>0$, such that $\Vert S^{*}\Vert_{0,v}=\max_{1\le i \leq p} \sum_{j =1}^p \mathbbm{1}(S^*_{ij}=0)\leq \delta_2 p^{\delta}$;
\item [(ii)] $\kappa_L,\kappa_S>0$ with $\frac{\kappa_S}{\kappa_L}\leq \frac{1}{54}$ and $\kappa_S\leq\delta_2$, such that $\xi(T({L^{*}}))=\frac{\sqrt r}{\kappa_L p^{\delta}}$ and $\mu(\Omega({S^{*}}))=\kappa_S p^{\delta}$;
\item[(iii)] $\underline{\delta},\delta_2'$ with $\underline{\delta}\leq \delta+\frac{1}{2}$ such that $\Vert S^{*} \Vert_{1,v} \leq \delta_2' p^{\underline{\delta}}$.
\end{inparaenum}
\end{Ass}

Let us now consider the generalized linear filter setting of Assumption \ref{gen}. In this case, while, as before, $\xi(T({L(\theta)}))=\xi(T({L^{*}}))$, since the low rank variety $\mathcal{L}(r)$ does not change, $\mu(\Omega(S(\theta)))$ is different from $\mu(\Omega(S^*))$ because the underlying sparsity pattern
now depends on the sparsity pattern of
the matrix $S_{\infty}=\sum_{k=-\infty}^\infty \Gamma_\epsilon(k)$. We therefore make the following assumption.

\begin{Ass}\label{ass2bis}
For all $p\in\mathbb N$, there exist:
\begin{inparaenum}
\item [(i)] $\delta \in [0,\frac{1}{2}]$, with $\delta < \alpha$, and $\delta_2>0$, such that $\Vert S_\infty\Vert_{0,v}=\max_{1\le i \leq p} \sum_{j =1}^p \mathbbm{1}(S_{\infty,ij}=0)\leq \delta_2 p^{\delta}$;
\item [(ii)] $\kappa_L,\kappa_S>0$ with $\frac{\kappa_S}{\kappa_L}\leq \frac{1}{54}$ and $\kappa_S\leq\delta_2$, such that $\xi(T({L^{*}}))=\frac{\sqrt r}{\kappa_L p^{\delta}}$ and $\mu(\Omega({S_\infty}))=\kappa_S p^{\delta}$;
\item[(iii)] $\underline{\delta},\delta_2'$ with $\underline{\delta}\leq \delta+\frac{1}{2}$ such that $\Vert S_{\infty} \Vert_{1,v} \leq \delta_2' p^{\underline{\delta}}$.
\end{inparaenum}
\end{Ass}
In both Assumptions \ref{ass2} and \ref{ass2bis}, part (ii) guarantees that \eqref{incoerenza} is satisfied, and the maximum number of nonzeros per row in $S^*$ or in $S_\infty$, which is controlled in part (i), is crucial to this end, as it is a natural upper bound for $\mu(\Omega({S}(\theta))$ \citep{chandrasekaran2011rank}. The condition $\delta < \alpha$ is instead needed to preserve the identification of the two components of the spectral density matrix. Part (iii) controls the $l_{1,v}$ norms of $S^*$ or $S_\infty$, notice that this condition is compatible with part (i).

Hereafter, depending on which filter setting we are considering, we define $q'=\Vert S^*\Vert_{0,v}$ or  $q'=\Vert S_\infty\Vert_{0,v}$. Notice that by means of Assumptions \ref{ass1} or \ref{ass1bis} we are constraining the number of nonzeros in $S(\theta)$ to be at most $q$, while by means of Assumption \ref{ass2}(i) or \ref{ass2bis}(i) we are further constraining the maximum number of nonzeros in each row to be at most $q'$. Obviously $q'\le q$ and moreover, because of our assumptions $q=O(p^{1+\delta})$, hence it can grow with $p$ at a maximum rate $p^{3/2}$.

Finally, by means of the next assumption and similarly to \citet{wu2018asymptotic}, we control the amount of physical dependence of our stochastic process (\ref{mod_base}) across time.

\begin{Ass}\label{ass3}
There exist $\delta_3>0$, $\rho\in [0,1)$, and $\delta'\in [0,\alpha)$,
such that, for all $p\in\mathbb N$:\linebreak
\begin{inparaenum}
\item [(i)] $\Vert \sum_{s=0}^\infty B_{s}\Vert_{1,v} \le \delta_3 p^{\alpha/2} \frac{1}{1-\rho}$;
\item [(ii)] $\Vert \sum_{s=0}^\infty B_{s}\Vert_{\infty,v} \leq  \delta_3 \frac{1}{1-\rho}$;
\item [(iii)] $\Vert \sum_{s=0}^\infty C_{s}\Vert_{1} \leq \delta_3 p^{\delta'/2} \frac{1}{1-\rho}$;\linebreak
\item [(iv)] $\Vert \sum_{s=0}^\infty C_{s}\Vert_{\infty,v} \leq  \delta_3 \frac{1}{1-\rho}$.
\end{inparaenum}
\end{Ass}

For any fixed $p$, Assumption \ref{ass3} implies geometrically decaying auto-covariances as $\rho^{\vert s\vert}$. Similar assumptions are made by \cite{wu2018asymptotic}, but here we generalize them to allow for a low rank plus sparse structure. Three comments are necessary. First, notice that we need to control the entire $\Vert .\Vert_{1}$ norm of the coefficients $C_s$ of the idiosyncratic filters.
Obviously, we could do the same for the coefficients $B_s$ of the common filters, since if we bound their  $l_{1}$ norm their $l_{1,v}$ norm is bounded as well, because $r$ is finite.
Second, part (iv) bounds the column sums of the coefficients of the idiosyncratic filters. This, together with part (iii), implies that the diverging behavior of those coefficients is due to the row sums. This means that for each given $\epsilon_{it}$ the $p$ idiosyncratic innovations have a finite impact.
This is just a useful way of parametrizing the model and we notice that, equivalently, we could also assume the viceversa or let both row and column sums diverge (compatibly with part (iii)). It is also straightforward to see that parts (iii) and (iv) are compatible with Assumption \ref{ass2}(i)
 or \ref{ass2bis}(i). Third, the assumption $\delta'< \alpha$ ensures that
the low rank component dominates the sparse component, in agreement with the idea of common factors.

\section{Consistency}\label{sec:cons}

Suppose now that we observe a sample of $p$-dimensional data vectors with size $T$, i.e., we observe $\{ X_{it}, i=1,\ldots, p,\, t =1,\ldots, T\}$,
and we compute the estimator $\wh{\Sigma}(\theta)=\wh{L}(\theta)+\wh{S}(\theta)$ such that:
\begin{equation}
\l(\wh{L}(\theta),\wh{S}(\theta)\r)=\arg\! \!\!\! \!\!\! \min_{\substack{\underline L(\theta)\in\mathcal{L}(r) \\ \underline S(\theta))\in \mathcal{S}(q)}} \frac{1}{2}\Vert \wt{\Sigma}(\theta)-(\underline L(\theta)+\underline S(\theta))\Vert_{F}^2
+\psi\Vert\underline L(\theta)\Vert_{*}+\rho\Vert\underline S(\theta)\Vert_{1},  \label{func:ob_spec_again}
\end{equation}
where $\wt{\Sigma}(\theta)$ is the smoothed periodogram defined in \eqref{kernel}.

In this section, we prove the algebraic and parametric consistency of the
pair of estimates $(\wh{L}(\theta),\wh{S}(\theta))$, and, in order to do this, we introduce two definitions, taken from \cite{chandrasekaran2012}. First, we say that $(\wh{S}(\theta),\wh{L}(\theta))$
is algebraically consistent if the following conditions hold, for any given $\theta\in[-\pi,\pi]$:
\begin{inparaenum}
\item $\text{\upshape{rk}}(\wh{L}(\theta))=\text{rk}({L}(\theta))$;
\item $\text{\upshape sgn}(\mathrm{Re}(\wh{S}(\theta)_{ij}))=\text{\upshape sgn}(\mathrm{Re}({S}(\theta)_{ij}))$,
for all $i,j=1,\ldots, p$ (by convention we let $\text{\upshape sgn} (0)=0$);
\item $\wh{L}(\theta)+\wh{S}(\theta)$ and $\wh{S}(\theta)$ are positive definite and $\wh{L}(\theta)$ is positive semidefinite.
\end{inparaenum}
Notice that condition 2 is often referred to also as sparsistency \citep{chandrasekaran2012}.

Second, we say that (parametric) consistency holds if the estimates $(\wh{S}(\theta),\wh{L}(\theta))$ are close to $({S}(\theta),{L}(\theta))$, for any given $\theta\in[-\pi,\pi]$, with high probability, according to the following norm:
\begin{equation}
g_\gamma(\wh{L}(\theta),\wh{S}(\theta))=\max\left(\frac{\Vert\wh{S}(\theta)-{S}(\theta)\Vert_{\infty}}{\gamma},\Vert\wh{L}(\theta)-{L}(\theta)\Vert_2\right)\label{ggamma}
\end{equation}
where $\gamma=\frac {\rho}{\psi}$ is the ratio of the thresholds in \eqref{func:ob_spec_again}.

By properly adapting the results of \cite{wu2018asymptotic} to the intermediate spikiness-sparsity regimes described in Section \ref{sec:model}, we prove  uniform consistency over frequencies of the smoothed periodogram pre-estimator  \eqref{kernel}.
\begin{Lemma}\label{cons}
Suppose that the kernel function $K(.)$ is even, bounded, with support $[-1,1]$, and bandwidth $M_T$, such that:
\begin{inparaenum}
\item[(i)] for some $\kappa>0$, $\vert K(s)-1\vert =O(s^{\kappa})$, as $s \rightarrow 0$;
\item[(ii)] $\int_{-\infty}^{\infty} K^2(s) ds< \infty$;
\item[(iii)] $\sum_{s'\in\mathbb Z} \sup_{\vert s'-s''\vert \leq 1} \vert K(s'\vartheta)-K(s''\vartheta)\vert =O(1)$, as $\vartheta \rightarrow 0$;
\item[(iv)]  $c_1 T^{\underline{\zeta}} \leq M_T \leq c_2 T^{\zeta}$, for some $c_1,c_2>0$ and $\zeta,\underline\zeta>0$, with $0<\underline{\zeta}<\zeta<1<\underline \zeta(2\kappa+1)$.
\end{inparaenum}\\
Then, under Assumption \ref{ass3}, there exists some positive real $G$, independent of $p$ and $T$, such that,
as $p,T\to\infty$, for $\theta_h=\frac{h\pi}{M_T}$:
$\mc{P}\left(\max_{\vert h\vert \leq M_T} \frac{1}{p^{\alpha}}\Vert \wt{\Sigma}(\theta_h)-\Sigma(\theta_h)\Vert_{2}\leq G\sqrt{\frac{ M_T \log(M_T)}T}
\right) \rightarrow 1$.
\end{Lemma}

\begin{Rem}[Bandwidth choice]\upshape{
Notice that the bias term, which is of order $\frac 1{M_T^\kappa}$, is not included in the above result, since for all $M_T$ satisfying condition (iv) this term is always dominated by the variance term. Indeed, while the optimal choice balancing variance and squared bias is $M_T=O(T^{1/(2\kappa+1)})$, in condition (iv) we are instead assuming $M_T=O(T^\zeta)$ with $\zeta>\underline\zeta>\frac 1{2\kappa+1}$. In other words, similarly to \citet{wu2018asymptotic}, with this choice of $M_T$ the mean squared error  of the smoothed periodogram is dominated by the variance, while the squared bias becomes negligible, as $T\to\infty$.  Typical values of $\kappa$ are 1 if we choose the Bartlett kernel, or 2 if we choose the Parzen kernel. All following theoretical results are unaffected if we relaxed our bandwidth choice and we picked smaller values of $\zeta$, provided that, when needed, we also account for the bias in the bound in Lemma \ref{cons}.}
\end{Rem}

We are now ready to show parametric and latent rank consistency
of $(\wh{L}(\theta),\wh{S}(\theta))$ under the basic filter setting.

\begin{thm}\label{thmMineUNALSE}
Let $\Omega=\Omega({S}^{*})$ and
$\mathcal T=T({L}^{*})$. Suppose that the assumptions of Lemma \ref{cons} hold,
with Assumptions \ref{assGDFM}, \ref{basic}, \ref{ass1}, and \ref{ass2}.
Set
$\psi=\frac{p^{\alpha}}{\xi(\mathcal T)}\sqrt{\frac{M_T\log M_T}{T}}$ and $\rho=\gamma \psi$, where $\gamma \in [9\xi(\mathcal T),\frac 1{6\mu(\Omega)}]$. In addition, suppose that $\underline{\delta}_T p^{2(\alpha-\underline{\delta})} < T < \overline{\delta}_T p^{6\delta}$ for some $\underline{\delta}_T,\overline{\delta}_T$ such that $0<\underline{\delta}_T<\overline{\delta}_T$, and
the minimum eigenvalue of $L^{*}$ is such that $\lambda_r(L^{*})>G_2 \frac{\psi}{\xi^2(\mathcal T)}$.  Then, there exists a positive real $G$ independent of $p$ and $T$
such that, as $p,T\to\infty$, for $\theta_h=\tfrac{h\pi}{M_T}$:
\begin{compactenum}
\item
 $\mc{P}\l(\max_{\vert h\vert \leq M_T}\frac{1}{p^{\alpha}}\Vert \wh{L}(\theta_h)-{L}(\theta_h)\Vert_{2}\leq G \frac{1}{\xi(\mathcal T)}\sqrt{\frac{M_T\log M_T}{T}}\r) \rightarrow 1$;
\item $\mc{P}\l(\max_{\vert h\vert \leq M_T}\frac{1}{p^{\alpha}}\Vert
\wh{S}(\theta_h)-{S}(\theta_h)\Vert_{\infty} \leq G  \frac{\gamma}{\xi(\mathcal T)}\sqrt{\frac{M_T\log M_T}{T}}\r)
\rightarrow 1$ and\\
$\mc{P}\l(\max_{\vert h\vert \leq M_T}\frac{1}{p^{\alpha}}\Vert \wh{S}(\theta_h)-{S}(\theta_h)\Vert_{2} \leq G q' \sqrt{\frac{M_T\log M_T}{T}}\r)
\rightarrow 1$;
\item $\mc{P}\l(\max_{\vert h\vert \leq M_T}\frac{1}{p^{\alpha}}\Vert \wh{\Sigma}(\theta_h)-\Sigma(\theta_h)\Vert_{2} \leq G \l[\frac{1}{\xi(\mathcal T)}+q'\r]\sqrt{\frac{M_T\log M_T}{T}}\r) \rightarrow 1$.
\end{compactenum}
Moreover,  if $\psi<C$ for some positive real $C$, then, for $\theta_h=\frac{h\pi}{M_T}$:
\begin{inparaenum}
\item [4.] $\mc{P}\l(\text{\upshape{rk}}(\wh{L}(\theta_h))=r\r) \rightarrow 1$;
\end{inparaenum}
and, if the minimum absolute value of the nonzero off-diagonal entries of $S^{*}$ is such that $\Vert S^*\Vert_{min,off}>G_3\frac{\psi}{\mu(\Omega)}$ then, for $\theta_h=\frac{h\pi}{M_T}$:
\begin{inparaenum}
\item [5.]  $\mc{P}\l(\text{\upshape sgn}(\mathrm{Re}(\wh{S}(\theta_h)_{ij})=\text{\upshape {sgn}}(\mathrm{Re}({S}(\theta_h)_{ij}))\r) \rightarrow 1$, for all $i,j=1,\ldots, p$.
\end{inparaenum}
\end{thm}

Similarly, under the generalized linear filter setting we have the following.

\begin{thm}\label{thmMineUNALSEgen}
Let $\Omega=\Omega({S}_{\infty})$ and
$\mc T=T({L}^{*})$. Suppose that the assumptions of Lemma \ref{cons} hold,
with Assumptions \ref{assGDFM}, \ref{gen}, \ref{ass1bis}, and \ref{ass2bis}. Set
$\psi=\frac{p^{\alpha}}{\xi(\mathcal T)}\sqrt{\frac{M_T\log M_T}{T}}$ and $\rho=\gamma \psi$, where $\gamma \in [9\xi(\mathcal T),\frac 1{6\mu(\Omega)}]$. In addition, suppose that $\underline{\delta}_T p^{2(\alpha-\underline{\delta})} < T < \overline{\delta}_T p^{6\delta}$
for some $\underline{\delta}_T,\overline{\delta}_T$ such that $0<\underline{\delta}_T<\overline{\delta}_T$, and that, for all $\theta\in[-\pi,\pi]$,
the minimum eigenvalue of $L(\theta)$ is such that $\lambda_r(L(\theta))>G_2 \frac{\psi}{\xi^2(\mathcal T)}$.  Then, there exists a positive real $G$ independent of $p$ and $T$
such that, as $p,T\to\infty$, for $\theta_h=\tfrac{h\pi}{M_T}$:
\begin{compactenum}
\item
 $\mc{P}\l(\max_{\vert h\vert \leq M_T}\frac{1}{p^{\alpha}}\Vert \wh{L}(\theta_h)-{L}(\theta_h)\Vert_{2}\leq G \frac{1}{\xi(\mathcal T)}\sqrt{\frac{M_T\log M_T}{T}}\r) \rightarrow 1$;
\item $\mc{P}\l(\max_{\vert h\vert \leq M_T}\frac{1}{p^{\alpha}}\Vert
\wh{S}(\theta_h)-{S}(\theta_h)\Vert_{\infty} \leq G  \frac{\gamma}{\xi(\mathcal T)}\sqrt{\frac{M_T\log M_T}{T}}\r)
\rightarrow 1$ and\\
$\mc{P}\l(\max_{\vert h\vert \leq M_T}\frac{1}{p^{\alpha}}\Vert \wh{S}(\theta_h)-{S}(\theta_h)\Vert_{2} \leq G q' \sqrt{\frac{M_T\log M_T}{T}}\r)
\rightarrow 1$;
\item $\mc{P}\l(\max_{\vert h\vert \leq M_T}\frac{1}{p^{\alpha}}\Vert \wh{\Sigma}(\theta_h)-\Sigma(\theta_h)\Vert_{2} \leq G \l[\frac{1}{\xi(\mathcal T)}+q'\r]\sqrt{\frac{M_T\log M_T}{T}}\r) \rightarrow 1$.
\end{compactenum}
Moreover,  if $\psi<C$ for some positive real $C$, then, for $\theta_h=\frac{h\pi}{M_T}$:
\begin{inparaenum}
\item [4.] $\mc{P}\l(\text{\upshape{rk}}(\wh{L}(\theta_h))=r\r) \rightarrow 1$;
\end{inparaenum}
and if, for all $\theta\in[-\pi,\pi]$, the minimum absolute value of the nonzero off-diagonal entries of $S(\theta)$ is such that $\Vert S(\theta)\Vert_{min,off}>G_3\frac{\psi}{\mu(\Omega)}$, then, for $\theta_h=\frac{h\pi}{M_T}$:
\begin{inparaenum}
\item [5.]  $\mc{P}\l(\text{\upshape sgn}(\mathrm{Re}(\wh{S}(\theta_h)_{ij})=
\text{\upshape{sgn}}(\mathrm{Re}({S}(\theta_h)_{ij}))\r) \rightarrow 1$, for all $i,j=1,\ldots, p$.
\end{inparaenum}
\end{thm}

Some important remarks follow.

\begin{Rem}\label{rem:T1}\upshape{
The upper bound  $T<\overline{\delta}_T p^{6\delta}$ is a non-asymptotic condition necessary to to ensure that
the conditions of the Theorem and Assumptions \ref{ass1}-\ref{ass1bis}, requiring the eigenvalue $\lambda_r({L}(\theta))$ to diverge as $p^\alpha$, are satisfied. Indeed, from the conditions of the Theorem and by Assumptions \ref{ass1}-\ref{ass1bis}, we must have:
\[
\lambda_r({L}(\theta))>G_2\frac {\psi}{\xi^2(\mathcal T)}= G_2 \frac{p^\alpha}{\xi^3(\mathcal T)} \sqrt{\frac{M_T\log M_T}{T}} = K p^\alpha,
\]
for some positive real $K$. Therefore, since by Assumption \ref{ass2}-\ref{ass2bis} $\xi(\mathcal T)=\frac{\sqrt r}{\kappa_L p^
\delta}$, and from the conditions of Lemma \ref{cons} $M_T\ge c_1 T^{\underline \zeta}$ with $\underline \zeta>0$, then it must hold that
$p^{3\delta}>M T^{1/2}$, i.e. $T<\overline {\delta}_T p^{6\delta}$ with $\overline {\delta}_T=\frac 1{K^2}$. Let us stress that this is a non-asymptotic condition, that is, it must hold for all $p$ and $T$. Notice that equivalently this requires $p>T^{1/(6\delta)}$ which shows that the less sparse is $S(\theta)$ (higher $\delta$) the larger must be $p$ in order to ensure the eigen-gap to be large enough for identification of the latent rank and the sparsity pattern (parts 4 and 5 of the theorems) to be possible even for finite $p$ and $T$. Notice also that if $\delta =0$ we can still have parametric consistency (parts 1, 2, and 3 of the theorems) as long as $p,T\to\infty$ but $\frac {{T}}{p} \to 0$}.
\end{Rem}

\begin{Rem}\upshape{
The lower bound $T>\underline{\delta}_T p^{2(\alpha-\underline \delta)}$ is also a non-asymptotic condition necessary to to ensure that
the conditions of the Theorem and Assumptions \ref{ass2}-\ref{ass2bis}, requiring that the maximum number of non-zeros per row in $S(\theta)$ to diverge as $p^\delta$, are satisfied. Indeed, from the conditions of the Theorem and Assumptions \ref{ass2}-\ref{ass2bis}
\[
\delta_2'p^{\underline \delta}\ge \Vert S(\theta)\Vert_{1,v}\ge q'\Vert S(\theta)\Vert_{min,off}>G_3 \frac{q'\psi}{\mu(\Omega)}\ge G_3\frac{C'p^\alpha}{\xi(\mathcal T)\mu(\Omega)} \sqrt{\frac{M_T\log M_T}{T}},
\]
since $q'\ge C'$ for some positive real $C'$. Therefore, since $\xi(\mathcal T)\mu(\Omega)\le \frac 1{54}$ by Assumptions \ref{ass2}-\ref{ass2bis} and from the conditions of Lemma \ref{cons} $M_T\ge c_1 T^{\underline \zeta}$ with $0<\underline \zeta<1$, we must have $p^{\underline \delta}> C p^\alpha T^{-1/2}$ for some positive real $C$ or equivalently $T> \underline{\delta}_Tp^{2(\alpha-\underline\delta)}$ with $\underline{\delta}_T=C^{2}$. Notice that
$\underline\delta<\delta'<\alpha$ by Assumption \ref{ass3} and since $\delta'$ is a bound on the $l_1$ norm while $\underline\delta$ is a bound on the $l_{1,v}$ norm. Moreover, $\underline\delta\le\delta+\frac 12$ so when $S(\theta)$ is the least sparse possible ($\delta=\frac 12$) we do not need a very large $T$, the lower bound being $T>\underline{\delta}_T p^{2(\alpha-1)}$, and if the latent eigenvalues are very spiked ($\alpha=1$) the sparsity pattern can be identified without imposing constraints between $T$ and $p$. In the most sparse case ($\delta=0$) we need a $T>\underline{\delta}_T p^{2(\alpha-1/2)}$, which in the spiked case ($\alpha =1$) implies $T$ at least comparable to $p$, while, for the upper bound in Remark \ref{rem:T1} to still hold, we must have both $p$ and $T$ bounded by a constant (recall that these are non-asymptotic conditions). If this lower bound for $T$ is not satisfied we cannot identify the sparsity pattern, but we can still recover the latent rank, and the parametric consistency of both $\wh L(\theta_h)$ and $\wh S(\theta_h)$ still holds.}
\end{Rem}

\begin{Rem}\label{rem:chandra}\upshape{
Parts 1, 2, and 3 provide bounds for the estimation error of the spectral density matrices, which are uniform over all frequencies. Parts 4 and 5 guarantee rank consistency for $\wh L(\theta_h)$ and sparsistency for $\wh S(\theta_h)$. In order for these conditions to be verified we need $\psi=\frac{p^{\alpha}}{\xi(\mathcal T)}\sqrt{\frac{M_T\log M_T}{T}}$ to be finite for all $p$ and $T$, however it is not required for $\psi$ to decrease as $p$ and $T$ increase. Furthermore, when $\psi$ is finite for all $p$ and $T$, then there exists a constant $\varphi$, depending on $\mu(\Omega)$ and $\xi(\mathcal T)$, such that if $g_\gamma(\wh{L}(\theta),\wh{S}(\theta))\le \varphi$ then parts 4 and 5 hold with probability 1 for all $p$ and $T$ (see
\citealp[Propositions 5.2 and 5.3]{chandrasekaran2012}, for details).}
\end{Rem}

\begin{Rem}\upshape{
The error bound in spectral norm for the sparse component in part 2 is larger than the bound $O\l(\sqrt{\frac{\log{p}}{T}}\r)$ reported in \cite{bickel2008covariance}, as we are allowing $q'$ to grow as $p^{\delta}$.}
\end{Rem}

\begin{Rem}\upshape{
From parts 1 and 2, it immediately follows that $\mc P\l(g_\gamma(\wh{L}(\theta),\wh{S}(\theta)) \le G \frac{p^{\alpha}}{\xi(\mathcal T)}\sqrt{\frac{M_T\log M_T}{T}} \r)\to 1$. Therefore, since $\xi(\mathcal T)=O(p^\delta)$ by Assumption \ref{ass2bis}, then, if $p^{\alpha+\delta}\sqrt{\frac{M_T\log M_T}{T}}\to 0$, as  $p,T\to\infty$, then $g_\gamma(\wh{L}(\theta),\wh{S}(\theta))\to 0$, with probability tending to 1. Notice that, to achieve parametric consistency, we must have $T$ growing faster than $p$, and in the worst case, i.e., $\alpha=1$ and $\delta=\frac 12$, this means that we need $p^{3}{\frac{M_T\log M_T}{T}}\to 0$. However, notice also that if, as common in the literature on the estimation of large matrices, we were to consider the norm relative to the dimension $p$, then the error bound would be $p^{\alpha+\delta-1}\sqrt{\frac{M_T\log M_T}{T}}$, which in the worst case requires $ p{\frac{M_T\log M_T}{T}}\to 0$. On the other hand in the GDFM case, i.e., when $\alpha=1$ and $\delta=0$, we would have the same bound, $\sqrt{\frac{M_T\log M_T}{T}}$,
which was stated in Theorem \ref{th:introTh}. This bound is also the same one derived for the classical smoothed periodogram estimator in Lemma \ref{cons}. }
\end{Rem}

Finally, we have a useful Corollary about the inverses
of $\wh{S}(\theta_h)$ and $\wh{\Sigma}(\theta_h)$.
\begin{coroll}\label{pos_def}
Under the assumptions of Theorem \ref{thmMineUNALSE} or Theorem \ref{thmMineUNALSEgen}, there exists a positive real $G$ independent of $p$ and $T$
such that, as $p,T\to\infty$,
for $\theta_h=\frac{h\pi}{M_T}$:
\begin{compactenum}
\item $\wh{\Sigma}(\theta_h)$ is positive definite if $\frac{\lambda_{p}(\Sigma(\theta_h))}{p^{\alpha}}>G \l(\frac{1}{\xi(\mc T)\sqrt{T}}+\frac{q'}{\sqrt{T}}\r)\sqrt{\frac{M_T\log M_T}{T}}$; \item $\wh{S}(\theta_h)$ is positive definite if $\frac{\lambda_{p}(S(\theta_h))}{p^{\alpha}}>G \frac{q'}{\sqrt{T}}\sqrt{\frac{M_T\log M_T}{T}}$;
\item $\wh{\Sigma}^{-1}(\theta_h)$ is positive definite if $\frac{\lambda_{p}(\Sigma(\theta_h))}{p^{\alpha}}\geq 2G \l(\frac{1}{\xi(\mc T)\sqrt{T}}+\frac{q'}{\sqrt{T}}\r) \sqrt{\frac{M_T\log M_T}{T}}$;
\item $\wh{S}^{-1}(\theta_h)$ is positive definite if $\frac{\lambda_{p}(S(\theta_h))}{p^{\alpha}}\geq 2G \frac{q'}{\sqrt{T}}\sqrt{\frac{M_T\log M_T}{T}}$.
\end{compactenum}
In addition, it also holds for $\theta_h=\tfrac{h\pi}{M_T}$:
\begin{compactenum}
\item [5.] $\mc{P}\left(\max_{\vert h\vert \leq M_T}\frac{1}{p^{\alpha}}\Vert\wh{\Sigma}(\theta_h)^{-1}-\Sigma(\theta_h)^{-1}\Vert_{2} \leq G \l(\frac{1}{\xi(\mc T)\sqrt{T}}+\frac{q'}{\sqrt{T}}\r)\sqrt{\frac{M_T\log M_T}{T}}\right) \rightarrow 1$;
\item [6.] $\mc{P}\left(\max_{\vert h\vert \leq M_T}\frac{1}{p^{\alpha}}\Vert\wh{S}(\theta_h)^{-1}- S(\theta_h)^{-1}\Vert_{2} \leq G \frac{q'}{\sqrt{T}}\sqrt{\frac{M_T\log M_T}{T}} \right)\rightarrow 1$.
\end{compactenum}
\end{coroll}

\begin{Rem}[Unshrinking]\upshape{
We stress that above defined estimates may suffer from systematic sub-optimality for what concerns estimated eigenvalues. In particular, if $p$ is large and the latent eigenvalues are spiked, the singular value thresholding procedure may lead to the over-shrinkage of latent eigenvalues. For this reason, following \cite{farne2020large}, we perform the un-shrinkage of the estimated latent eigenvalues, i.e., we give back the threshold to $\wh L(\theta_h)$. The new idiosyncratic estimate is then obtained keeping fixed the off-diagonal sparsity pattern recovered, and deriving its diagonal by difference from the diagonal of $\wh{\Sigma}(\theta_h)$. The resulting matrix estimators are our  UNshrunk ALgebraic Spectral Estimators (UNALSE). More specifically, for any $\theta_h=\frac{h\pi}{M_T}$, consider the spectral decomposition $\wh L(\theta_h)=\wh{W}(\theta_h)\wh D(\theta_h)\wh{W}^\dag(\theta_h)$,
 then we define:
\begin{align}
&\wh{L}_{\text{\tiny{UNALSE}}}(\theta_h)=\wh{W}(\theta_h)\l(\wh{D}(\theta_h)+{\psi} {I}_r\r)\wh{W}^\dag(\theta_h),\nn\\
&\text{diag}(\wh{S}_{\text{\tiny{UNALSE}}}(\theta_h))=\text{diag}(\wh{\Sigma}(\theta_h))-\text{diag}(\wh{L}(\theta_h)),\quad \text{off-diag}(\wh{S}_{\text{\tiny{UNALSE}}}(\theta_h))=\text{off-diag}(\wh{S}(\theta_h)),\nn\\
&\wh{\Sigma}_{\text{\tiny{UNALSE}}}(\theta_h)=\wh{L}_{\text{\tiny{UNALSE}}}(\theta_h)+\wh{S}_{\text{\tiny{UNALSE}}}(\theta_h),\nn
\end{align}
where ${\psi}>0$ is the same as in Theorems \ref{thmMineUNALSE} or \ref{thmMineUNALSEgen}. The above defined UNALSE estimates have two relevant optimality properties.  First, they have the smallest possible Frobenius loss from the targets into the recovered matrix varieties.
Second, they have the maximally concentrated eigenvalues into the
class of algebraically consistent estimators, for any given sample size $T$. Under our assumptions,
the above optimality properties of the un-shrinkage procedure of \cite{farne2020large},
to which we refer for the details, hold straightforwardly.}
\end{Rem}

\section{Threshold selection}\label{sec:thresholds}

In solving problem (\ref{func:ob_spec}), the choice of the eigenvalue
threshold $\psi$ and the sparsity threshold $\rho$ is a nontrivial issue.
Differently from the covariance matrix context, in fact, the magnitude of the eigenvalues
can vary a lot across frequencies, which may cause the optimization of (\ref{func:ob_spec}) to be strongly sensitive to the magnitude of both thresholds.

Let us suppose that $\wh{L}_{\psi,\rho}(\theta_h)$, $\wh{S}_{\psi,\rho}(\theta_h)$,
$\wh{\Sigma}_{\psi,\rho}(\theta_h)=\wh{L}_{\psi,\rho}(\theta_h)+\wh{S}_{\psi,\rho}(\theta_h)$
are the solutions of (\ref{func:ob_spec}) with thresholds $\psi$ and $\rho$, under the assumptions of Theorem \ref{thmMineUNALSE}.
The dual norm of the composite loss (\ref{ggamma}) is considered, from which we define at each $\theta_h$ the following criterion:
\begin{equation}MC_h(\psi,\rho)=\max
\left\{\frac{\wh{r}\Vert\wh{L}_{\psi,\rho}(\theta_h)\Vert_2}{\wh{\beta}_{\psi,\rho}(\theta_h)},
\frac{\frac{\psi}{\rho}\Vert\wh{S}_{\psi,\rho}(\theta_h)\Vert_{1,v}}{(1-\wh{\beta}_{\psi,\rho}(\theta_h))}\right\},
\label{MC}
\end{equation}
where $\Vert\wh{S}_{\psi,\rho}(\theta_h)\Vert_{1,v}=\max_{i=1,\ldots,p}\sum_{j=1}^p |\wh{S}_{ij,\psi,\rho}(\theta_h)|$ and $\wh{\beta}_{\psi,\rho}(\theta_h)=
\text{tr}(\wh{L}_{\psi,\rho}(\theta_h))/\text{tr}(\wh{\Sigma}_{\psi,\rho}(\theta_h))$ is the estimated proportion of latent variance.
The optimal threshold pair $(\breve{\psi}_h,\breve{\rho}_h)$ is thus selected as the mini-max $(\breve{\psi}_h,\breve{\rho}_h)=\arg \min_{\psi,\rho}{MC_h(\psi,\rho)}$,
where $\psi$ and $\rho$ vary across pre-specified grids.
This threshold selection method penalizes solution pairs with too dispersed latent eigenvalues and too many residual nonzeros in single rows, by comparing two appropriately re-scaled versions of the spectral norm of the low rank solution and the row-wise maximum norm of the residual solution.

In order to ensure the effectiveness of the above criterion,
the threshold grids need to be properly set up at each $\theta_h$,
according to the unknown underlying algebraic structure.
We thus recall from Theorem \ref{thmMineUNALSEgen} that $\psi=\sqrt{\frac{p M_T\log M_T}{T}}\frac{1}{\xi(T)}$, 
assuming the intermediate value $\alpha=1/2$ and recalling from \cite{chandrasekaran2011rank} that
$\text{\upshape inc}(L(\theta_h))\leq \xi(T)\leq 2 \text{\upshape inc}(L(\theta_h))$, where $\text{\upshape inc}(L(\theta_h))$ is the incoherence of $L(\theta_h)$, defined as
$\text{\upshape inc}(L(\theta_h))=\max_{i=1,\ldots,p}\Vert \mathcal{P} e_i \Vert$, with $e_i$ the canonical basis vector ($i$th column of the $p$ dimensional identity matrix), and the operator $\mathcal{P}$ projecting each $e_i$ onto the row/column space of $L(\theta_h)$.
The extreme incoherence values are $\text{\upshape inc}(L(\theta_h))=1$, when any vector of the standard basis belongs to the row/column space of $L(\theta_h)$, and $\text{\upshape inc}(L(\theta_h))=\sqrt{\frac{r}{p}}$, when $L(\theta_h)$ is a Hadamard matrix.
In light of this, and since $\sqrt{\frac{M\log{M_T}}{T}}>\sqrt{\frac{1}{T}}$,
we initialize the grid for the eigenvalue threshold ${\breve{\psi}}$ as the sequence
of $n_{thr}$ equi-spaced real numbers
from $\sqrt{\frac{p}{T}}\frac{1}{2\wt{inc}}$ to $\sqrt{\frac{p}{T}}\frac{1}{\wt{inc}}$, where, for any given value of $r_{thr}$,
we set $\wt{inc}=\sqrt[4]{\frac{r_{thr}}{p}}$, which is the geometric mean of the minimum and maximum incoherence values, i.e. $\sqrt{\frac{r_{thr}}{p}}$ and $1$ respectively (see \cite{chandrasekaran2011rank}).

We start by setting $r_{thr}=1$. Then, we apply the solution algorithm of (\ref{func:ob_spec})
(see Section \ref{sec:res}) with the grid for $\breve{\psi}$ as defined above.
Note that we iteratively adapt the grid for $\breve{\psi}$ during the optimization process by dividing, at each step $k$ of the iteration, the grid components  by the Gini index of the eigenvalues of $\mathcal E_{Y,k}(\theta_h)=Y_{k-1}(\theta_h)- \frac{1}{2}(Y_{k-1}(\theta_h)+Z_{k-1}(\theta_h)-\wt{\Sigma}(\theta_h))$ (see point 2b of the solution algorithm). This is done in order to adapt the eigenvalue thresholds
to the underlying degree of spikiness of latent eigenvalues.

Then, if the ranks of $\wh{L}_{\breve{\psi},\breve{\rho}}(\theta_h)$ across thresholds vary too much or the eigenvalue threshold selected by criterion (\ref{MC}) lies in the grid extremes,
the value of $r_{thr}$ must be changed. In particular,
it must be decreased if the recovered rank is very large uniformly across thresholds,
and increased if it is very small or zero.
When the eigenvalue threshold selected by the MC criterion is far away from the boundaries,
and the recovered rank is constant and stable across thresholds,
we stop and select $r_{thr}$. 

Concerning the sparsity threshold, we decide to set $\breve{\gamma}$ as the sequence of $n_{thr}$ equi-spaced real numbers from $s_{thr} \times p^{-1/2}$ to $s_{thr} \times p^{-1/4}$.
The two functions of $p$, $p^{-1/2}$ and $p^{-1/4}$, represent two plausible extremes for residual nonzero proportions, while $s_{thr}$ is a magnitude parameter. We first set $s_{thr}$ to $1$, and we run the solution algorithm of (\ref{func:ob_spec}) setting $\breve{\rho}=\breve{\gamma}\sqrt{\frac{p}{T}}\frac{1}{\wt{inc}}$.

Similarly, if the sparsity threshold selected by criterion (\ref{MC}) lies in the grid extremes,
we decrease or increase it, in order to obtain a non-diagonal solution
with a reasonable and stable number of nonzeros.
When the value of $\breve{\gamma}$ selected by the MC criterion is far from the grid extremes and the number of nonzeros is approximately stable across thresholds, we stop and select $s_{thr}$.
The described process allows to map the problem of selecting thresholds $\psi$ and $\rho$
into the more intuitive problem of selecting $r_{thr}$ and $s_{thr}$.
Apart from it, the two thresholds can also be chosen manually, as long as
the recovered rank and sparsity pattern are constant in a neighborhood
of the chosen threshold pairs.
\section{Simulation study}\label{sec:sim}


In order to test the performance of UNALSE under an exhaustive range of situations, we consider three different simulated scenarios:
\begin{inparaenum}
  \item[A]: simulated with basic filters;
  \item [B]: simulated with general filters and a very sparse, almost negligible, residual pattern;
  \item [C]: simulated with general filters and a less sparse, more relevant, residual pattern.
\end{inparaenum}
For each Scenario, we setup five Settings, with different dimensions, sample sizes, spectral magnitudes, latent ranks, condition numbers, and sparsity degrees.
In particular, the Settings from $1$ to $3$ present $p=100$ and $T=1000$. Setting 4 presents $p=T=150$, Setting 5 presents $p=200$ and $T=100$.

Our spectra follow a reverse S-shape:
the minimum latent eigenvalue decreases across frequencies for all settings,
the minimum off-diagonal nonzero residual entry (in absolute value)
varies with a similar trend across frequencies.
The proportion of latent variance differs across settings.
In absolute terms, the largest one is for Setting 4, followed by Setting 5, and Settings 1, 2 and 3.
At the same time, considering the proportion of residual covariance,
Settings 1 and 4 are the most sparse, followed by Settings 2 and 5,
while Setting 3 is the least sparse.
In Scenario C we have a situation with a small minimum latent eigenvalue and minimum residual non-zero off-diagonal entry (in absolute terms) across frequencies, and the capability of UNALSE to recover the sparsity pattern can be meaningfully tested.

We fix the frequency grid as $\theta_h=\frac{\pi h}{12}$, $h=0,\ldots,5$. We simulate $N=100$ replications of the data $X_t$ having a spectral density with a low rank plus sparse structure designed according to a given Scenario and Setting (see Section D in the supplementary appendix for details on the simulation mechanism), computing the pre-estimator of the spectral density matrix (\ref{kernel}) on each simulated dataset. To this end, we adopt the classical choice of a Bartlett kernel \citep{forni2000generalized,forni2017dynamic} and we set $M_T=\lfloor \sqrt T\rfloor$.

For each replication, we apply the solution algorithm (see Section \ref{sec:res}) with the threshold selection procedure described in Section \ref{sec:thresholds}, thus getting $100$ optimal spectral density matrix estimates $\wh{L}^{(b)}_{\text{\tiny{UNALSE}}}(\theta_h)$, $\wh{S}^{(b)}_{\text{\tiny{UNALSE}}}(\theta_h)$, $\wh{\Sigma}^{(b)}_{\text{\tiny{UNALSE}}}(\theta_h)$, for $b=1,\ldots,100$, $h=0,\ldots,5$. In the literature, the only existing competitor is for the low rank component and relies on the dynamic principal components of \cite{brillinger2001time} (see also \cite{forni2000generalized}).
We call it $\wh{L}_{DYN}(\theta_h)$. Note however that this estimator requires a pre-specified rank, which, in a high-dimensional setting, can be determined via the information criterion by \citet{hallin2007determining} or the test by \citet{onatski2009testing}.

Results are reported using the frequencies ${f}_h=\frac{\theta_h}{\pi}$, $h=0,\ldots,5$. For each $f_h$ and for each replication, we calculate some relevant statistics and some relative metrics to evaluate the quality of the rank and sparsity pattern recovery.
\begin{compactenum}[(i)]
\item the latent variance proportion $\wh{\beta}(f_h)=\frac 1{100}\sum_{b=1}^{100}\frac{tr(\wh{L}^{(b)}(f_h))}{tr(\wh{\Sigma}^{(b)}(f_h))}$;
\item the binary indicator sum for the correct estimated rank
$\wh{R}=\frac 1 6\sum_{b=1}^{100}\sum_{h=0}^5\mathbb I(\hat{r}^{(b)}(f_h)=r)$.
\item the nonzero predictive value: $nzpv(f_h) =\frac 1{100}\sum_{b=1}^{100}\frac{\sum_{i=1}^p \sum_{j=i+1}^p \mathbb I \{\wh{S}^{(b)}_{ij}(f_h) \ne 0 \bigcup {S}_{ij}(f_h) \ne 0\}}{\sum_{i=1}^p \sum_{j=i+1}^p {\mathbb I\{\wh{S}_{ij}^{(b)}(f_h) \ne 0\}}}$;
\item the positive predictive value: $ppv(f_h) =\frac 1{100}\sum_{b=1}^{100}\frac{\sum_{i=1}^p \sum_{j=i+1}^p \mathbb I \{\wh{S}^{(b)}_{ij}(f_h) > 0 \bigcup {S}_{ij}(f_h) > 0\}}{\sum_{i=1}^p \sum_{j=i+1}^p {\mathbb I\{{S}_{ij}^{(b)}(f_h) > 0\}}}$;
\item the negative predictive value: $npv(f_h) =\frac 1{100}\sum_{b=1}^{100}\frac{\sum_{i=1}^p \sum_{j=i+1}^p \mathbb I \{\wh{S}^{(b)}_{ij}(f_h) < 0 \bigcup {S}_{ij}(f_h) < 0\}}{\sum_{i=1}^p \sum_{j=i+1}^p {\mathbb I\{{S}_{ij}^{(b)}(f_h) < 0\}}}$;
\item the maximum of the sum of the binary indicator of nonzero recovered residual entries of each row: $mnz_i(f_h)=\max_{j=1,\ldots,p, j \ne i}\sum_{b=1}^{100}\mathbb I\{\wh{S}^{(b)}_{ij}(f_h)\ne 0\}$.
\end{compactenum}

In order to evaluate the properties of our estimates, we calculate also the Frobenius loss of each estimate from the relative target, rescaled by the dimension:
\begin{compactenum}[(i)]
\item the low rank component Frobenius loss
$err_{\wh{L}}(f_h)=\frac 1{100}\sum_{b=1}^{100}\Vert \wh{L}^{(b)}(f_h)-{L}(f_h) \Vert_F/p$;
\item the ratio between the overall UNALSE and the input Frobenius loss with respect to the target:
$err_{ratio}(f_h)=\frac 1{100}\sum_{b=1}^{100} \frac{\Vert \wh{\Sigma}^{(b)}(f_h)-{\Sigma}(f_h) \Vert_F}{\Vert \wt{\Sigma}^{(b)}(f_h)-{\Sigma}(f_h) \Vert_F}$.

\end{compactenum}

For all above quantities, we calculate also the standard deviation across the $100$ trials. We present here results only for Scenarios A and C and Settings 3 and 4, while all other results are available in Section D of the supplementary appendix.

%

First of all, we compare our estimates of the latent rank $r$ with those obtained with the test by \cite{onatski2009testing}, based on sample dynamic eigenvalues of the smoothed periodogram estimator (see Table \ref{tab:addlabel}).  We observe that the test presents some empirical level issues, particularly when the eigenvalues are not so spiked. For Setting 1, for instance, the observed proportion of correct decisions is as low as $38\%$. For the other settings of all scenarios, we note that the same proportion increases considerably, consistently with the increased spikiness of latent eigenvalues. At the same time, the observed outcome is always less than our $100\%$.

\begin{table}[t!]
  \centering
  \caption{Estimation of latent rank $r$.}
  \footnotesize{
    \begin{tabular}{c | cc | cc | cc}
    \hline
    \hline
    &	\multicolumn{6}{c}{\textbf{Scenario}}\\
    \hline
     & \multicolumn{2}{c|}{A} & \multicolumn{2}{c|}{B} & \multicolumn{2}{c}{C} \\
     \cline{2-7}
{\textbf{Setting}}     & Onatski & UNALSE & Onatski & UNALSE & Onatski & UNALSE\\
     \hline
     	1 & 38 & 100 & 90 & 100 & 82 & 100\\
	2 & 95 & 100 & 93 & 100 & 86 & 100\\
	3 & 93 & 100 & 94 & 100 & 90 & 100\\
	4 & 90 & 100 & 94 & 100 & 92 & 100\\
	5 & 96 & 100 & 98 & 100 & 95 & 100\\
\hline
\hline
    \end{tabular}
    }
  \label{tab:addlabel}
\end{table}

          \begin{figure}[t!]
          \caption{Estimated latent variance proportion $\wh{\beta}(f_h)$ - Scenario A.}\label{fig:beta_A.1.U}
          \centering
               \begin{tabular}{cc}
                             {\footnotesize Setting 3}&{\footnotesize Setting 4}\\
              \includegraphics[width=.2\textwidth]{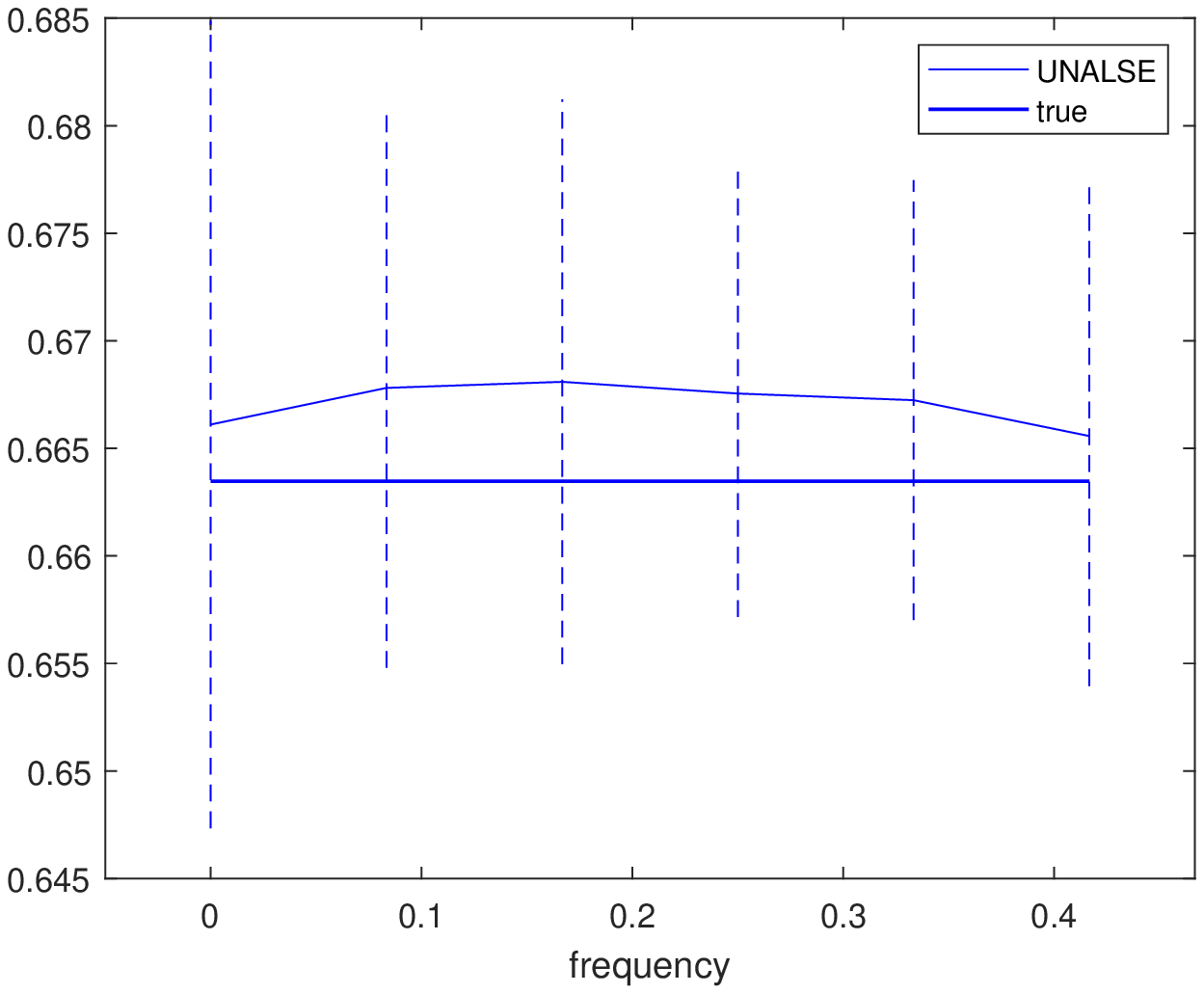}&
              \includegraphics[width=.2\textwidth]{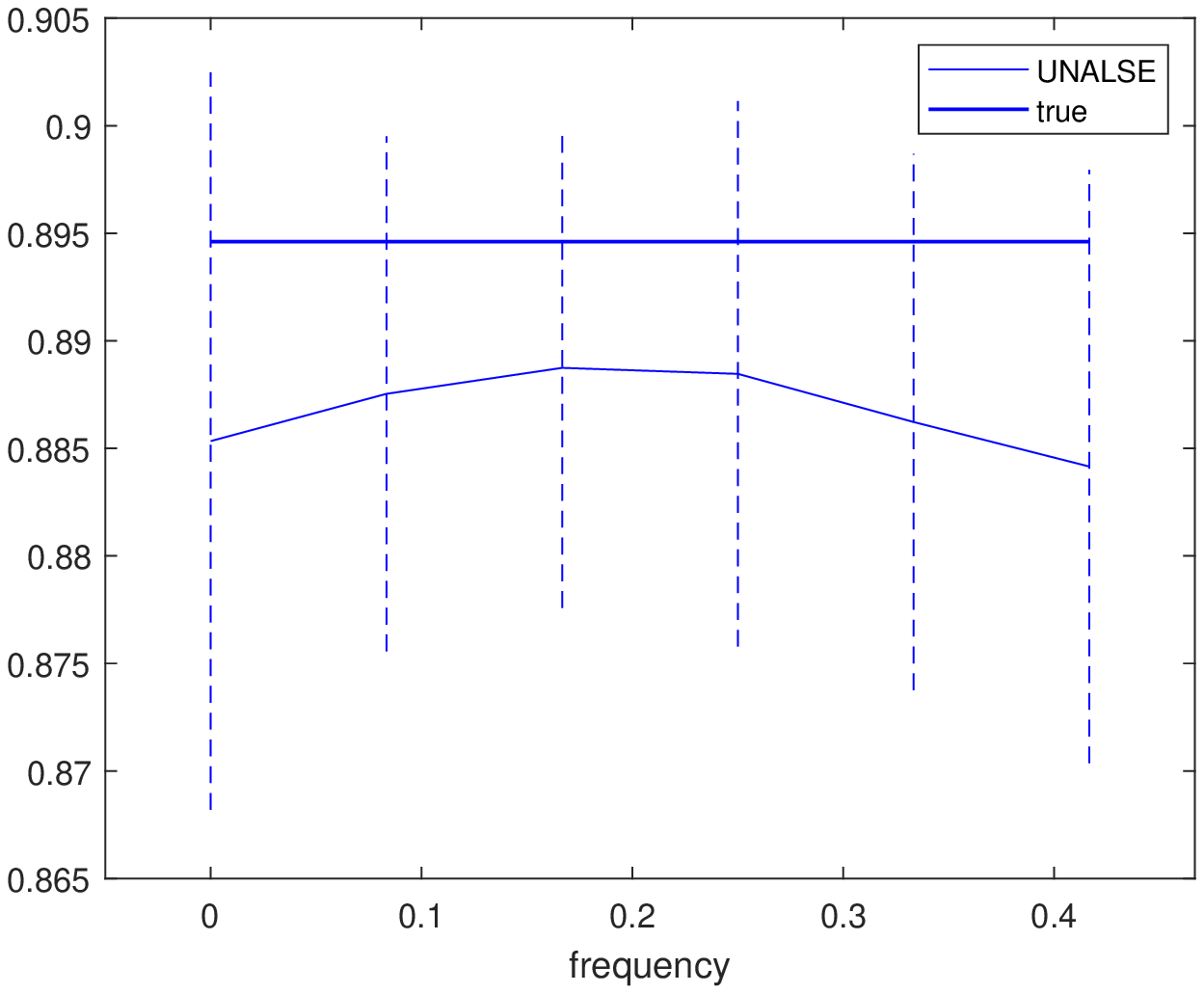}\\
              \end{tabular}
          \end{figure}

In Figure \ref{fig:beta_A.1.U} we show the estimated latent variance proportion $\wh{\beta}(f_h)$. We
notice that UNALSE estimates systematically better than DYN the true $\beta(f_h)$ across frequencies (results not reported),
in particularly for Setting 3, which has $p/T=0.1$. This is due to the bad properties of sample eigenvalues with respect to the presence of non-spiked latent eigenvalue structures.
In addition, UNALSE can recover the residual sign pattern at each frequency (see Figure \ref{fig:pred_pos_A.1.U}),
and the positive and negative predictive values are pretty similar.
When the overall magnitude is larger, the true predictive rate of nonzeros across frequencies grows considerably: indeed, in Setting 3, it overcomes $80\%$ at all frequencies. This happens because there are many nonzeros of sufficient magnitude, even compared to the low rank component.

 \begin{figure}[t!]
          \caption{Positive and negative predictive values $ppv(f_h)$ and $npv(f_h)$ - Scenario A.}\label{fig:pred_pos_A.1.U}
          \centering
               \begin{tabular}{cccc}
	      {\footnotesize $ppv(f_h)$ Setting 3}&{\footnotesize $ppv(f_h)$ Setting 4}&
	      {\footnotesize $npv(f_h)$ Setting 3}&{\footnotesize $npv(f_h)$ Setting 4}\\
              \includegraphics[width=.2\textwidth]{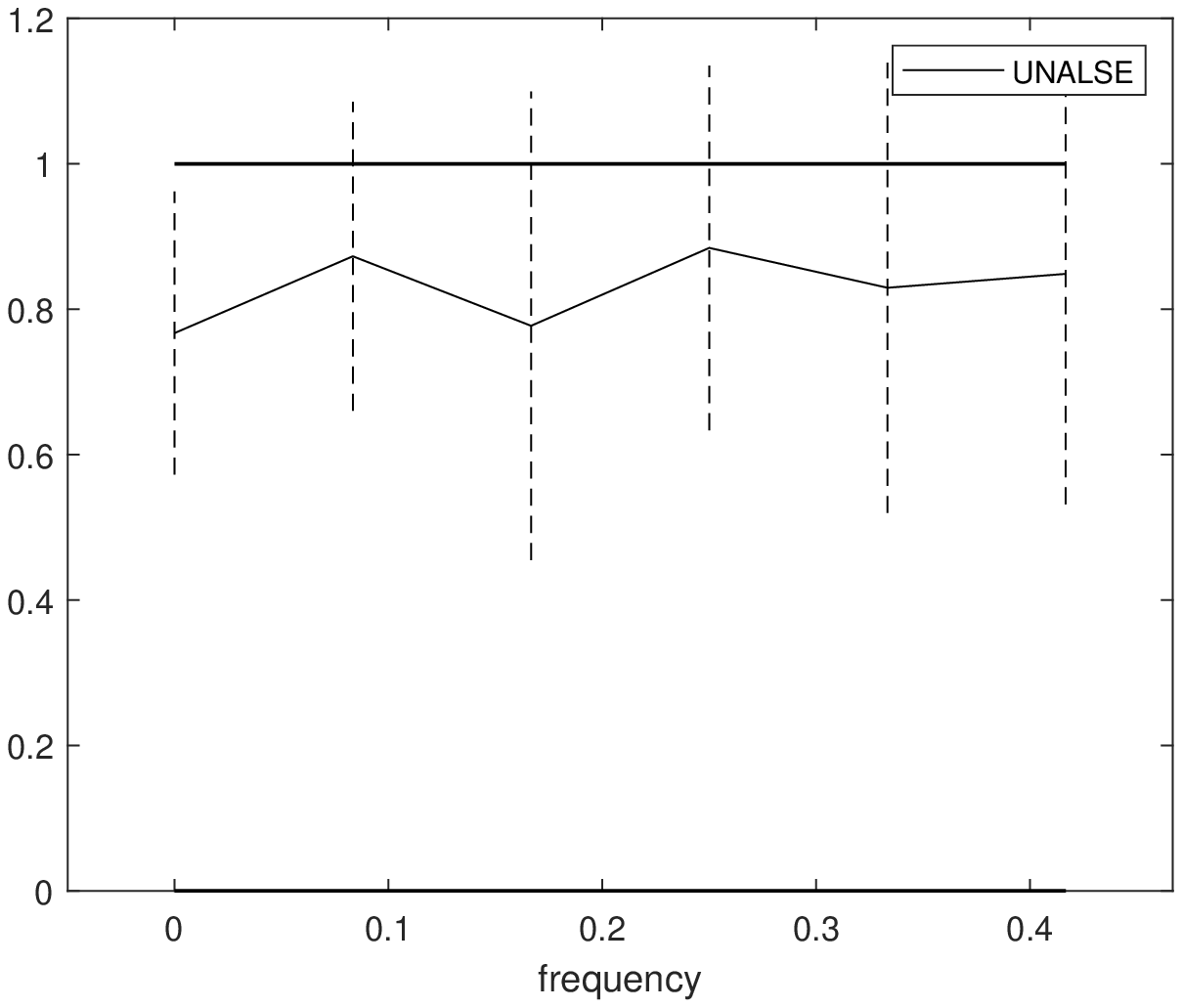}&
              \includegraphics[width=.2\textwidth]{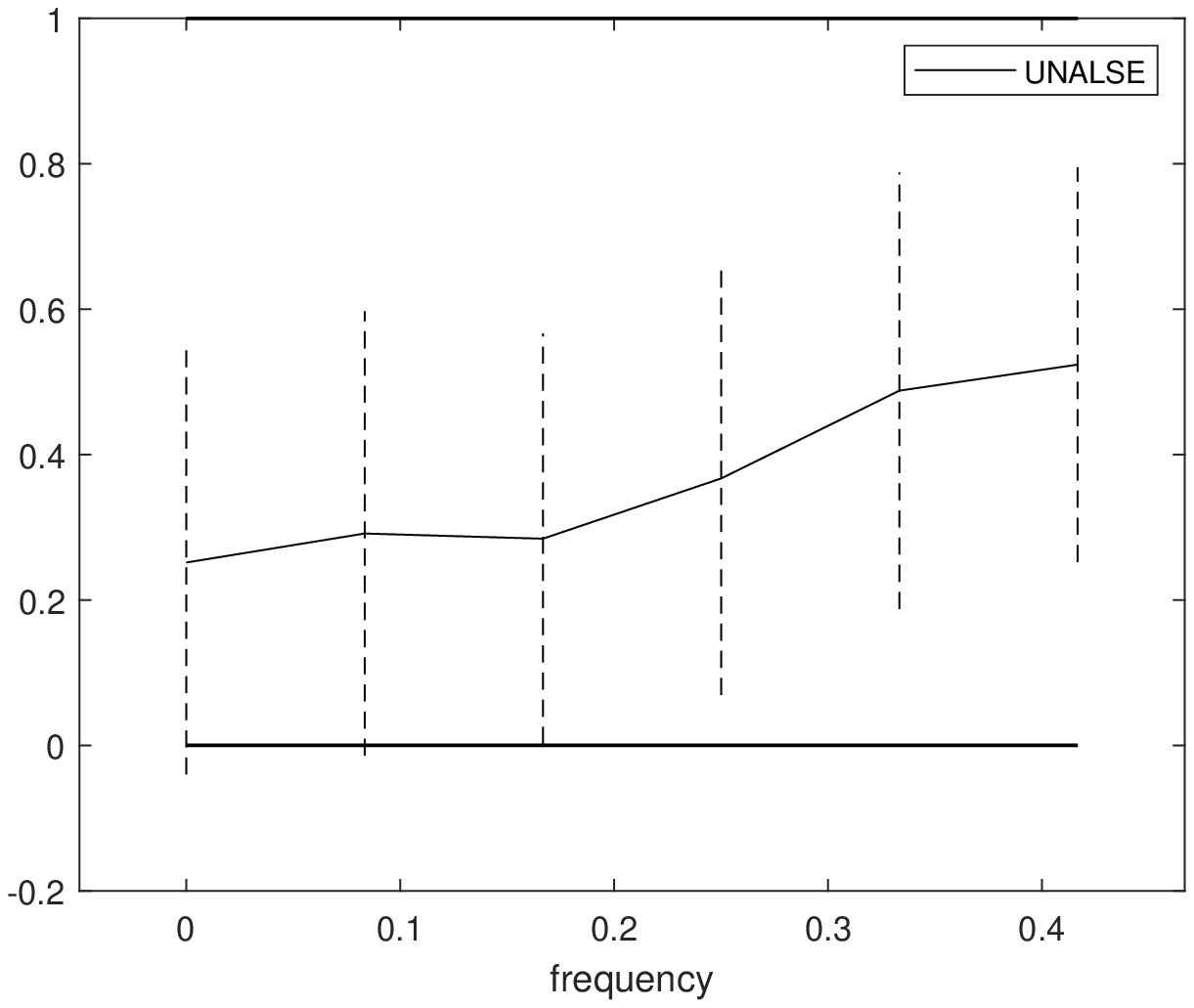}&
               \includegraphics[width=.2\textwidth]{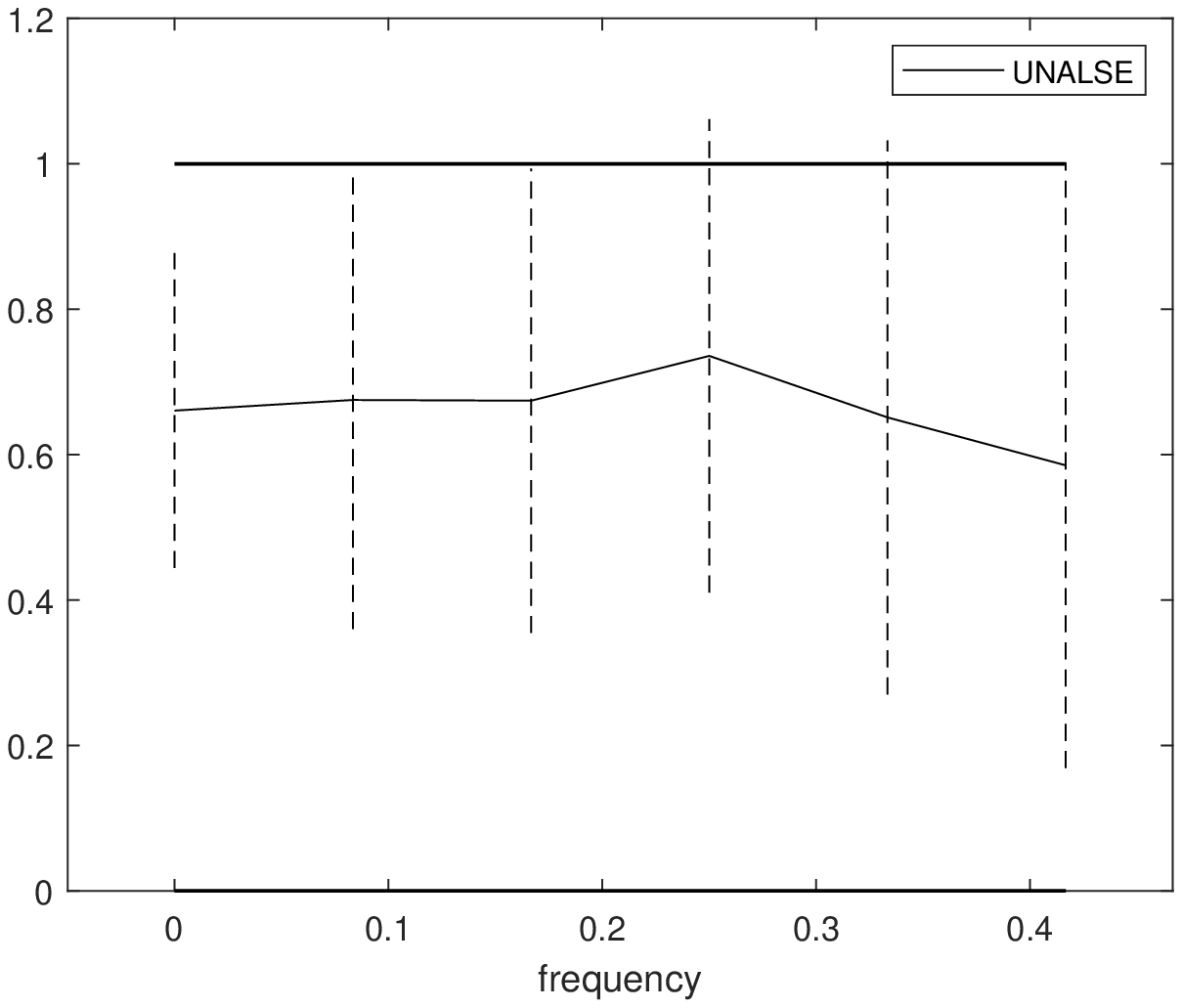}&
              \includegraphics[width=.2\textwidth]{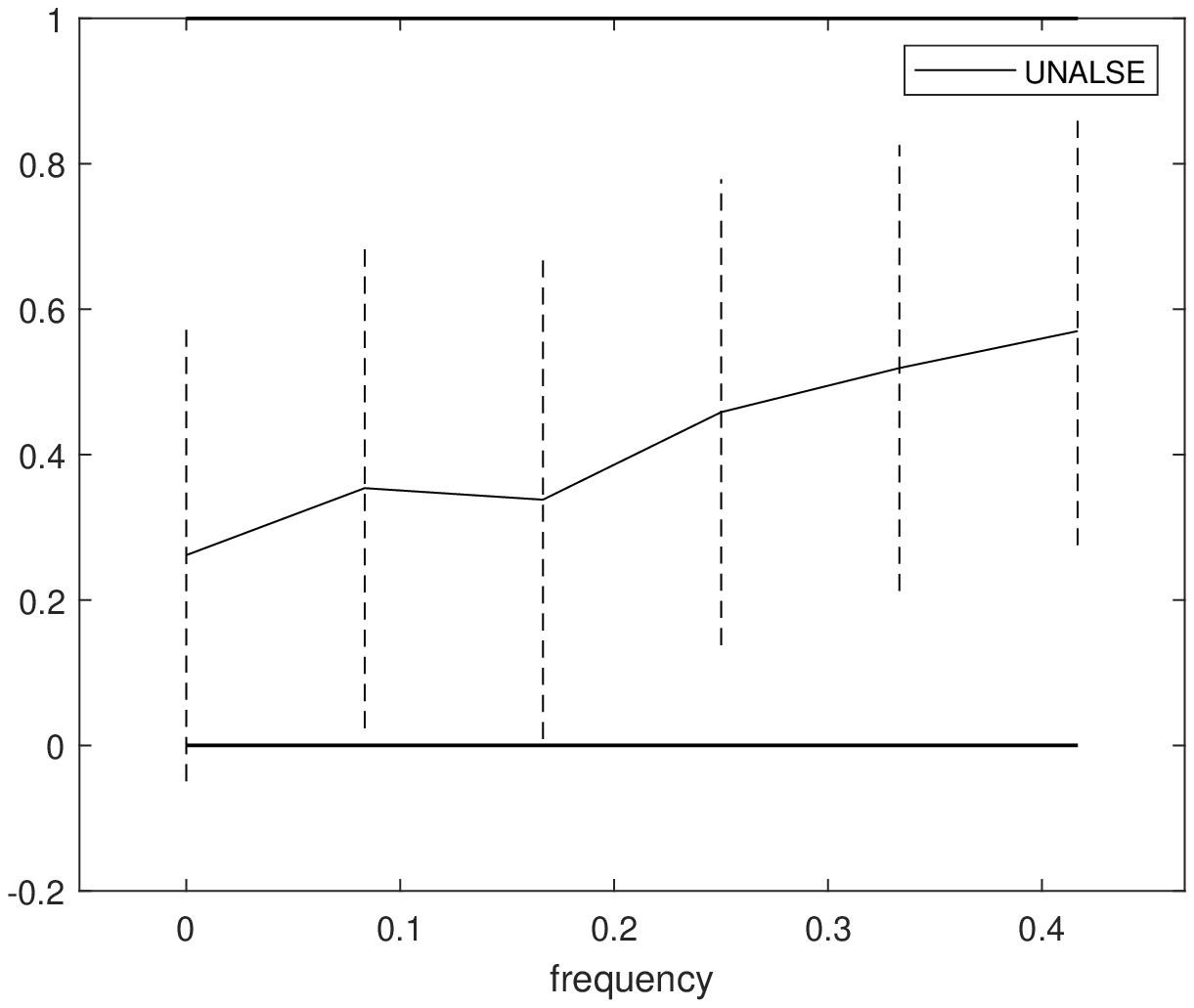}\\
              \end{tabular}
          \end{figure}

Concerning the error metrics (see Figure \ref{fig:err_L_A.1.U}), $err_{\wh{L}}(f_h)$ is slightly worse for UNALSE than for DYN at low frequencies, and quite better at high frequencies. Again, this is consistent with the bad properties of the dynamic principal components under weak factors in not so large dimensions. The same pattern is visible for $err_{ratio}(f_h)$.

\begin{figure}[t!]
          \caption{$err_{\wh{L}}(f_h)$ and $err_{ratio}(f_h)$ - Scenario A.}\label{fig:err_L_A.1.U}
          \centering
               \begin{tabular}{cccc}
	      {\footnotesize $err_{\wh{L}}(f_h)$ Setting 3}&{\footnotesize $err_{\wh{L}}(f_h)$ Setting 4}&
	      {\footnotesize $err_{ratio}(f_h)$ Setting 3}&{ \footnotesize $err_{ratio}(f_h)$ Setting 4}\\
              \includegraphics[width=.2\textwidth]{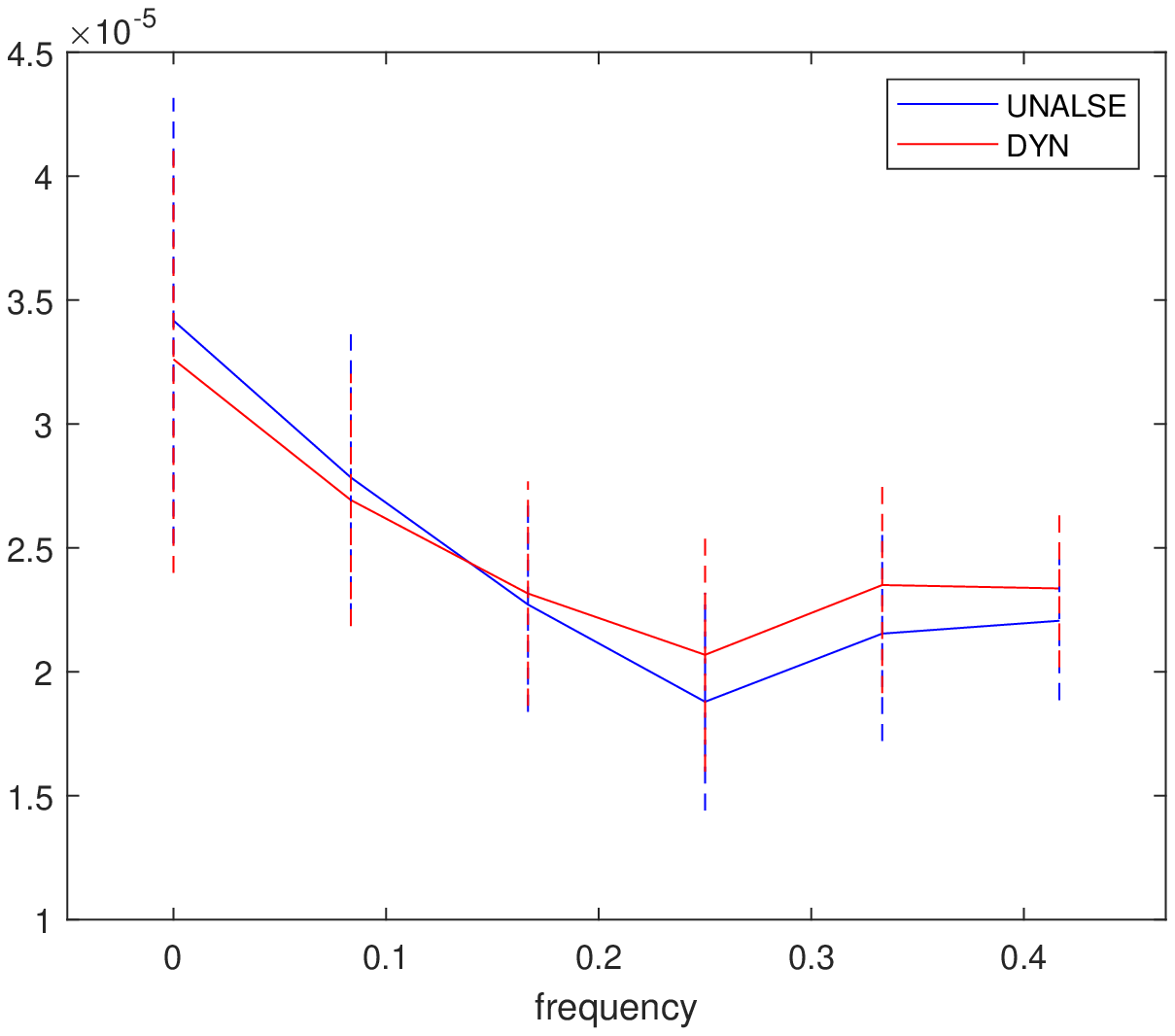}&
              \includegraphics[width=.2\textwidth]{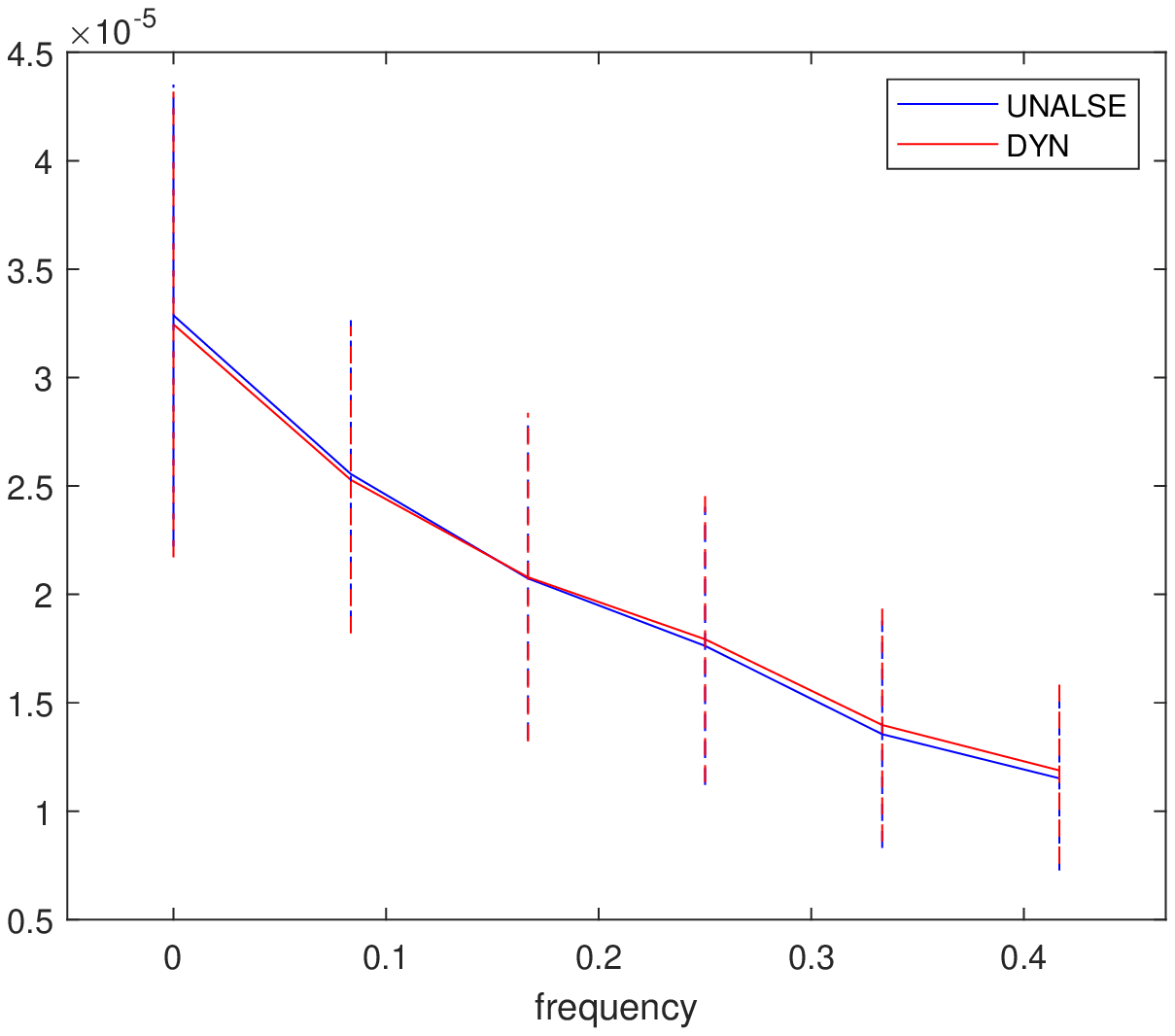}&	
              \includegraphics[width=.2\textwidth]{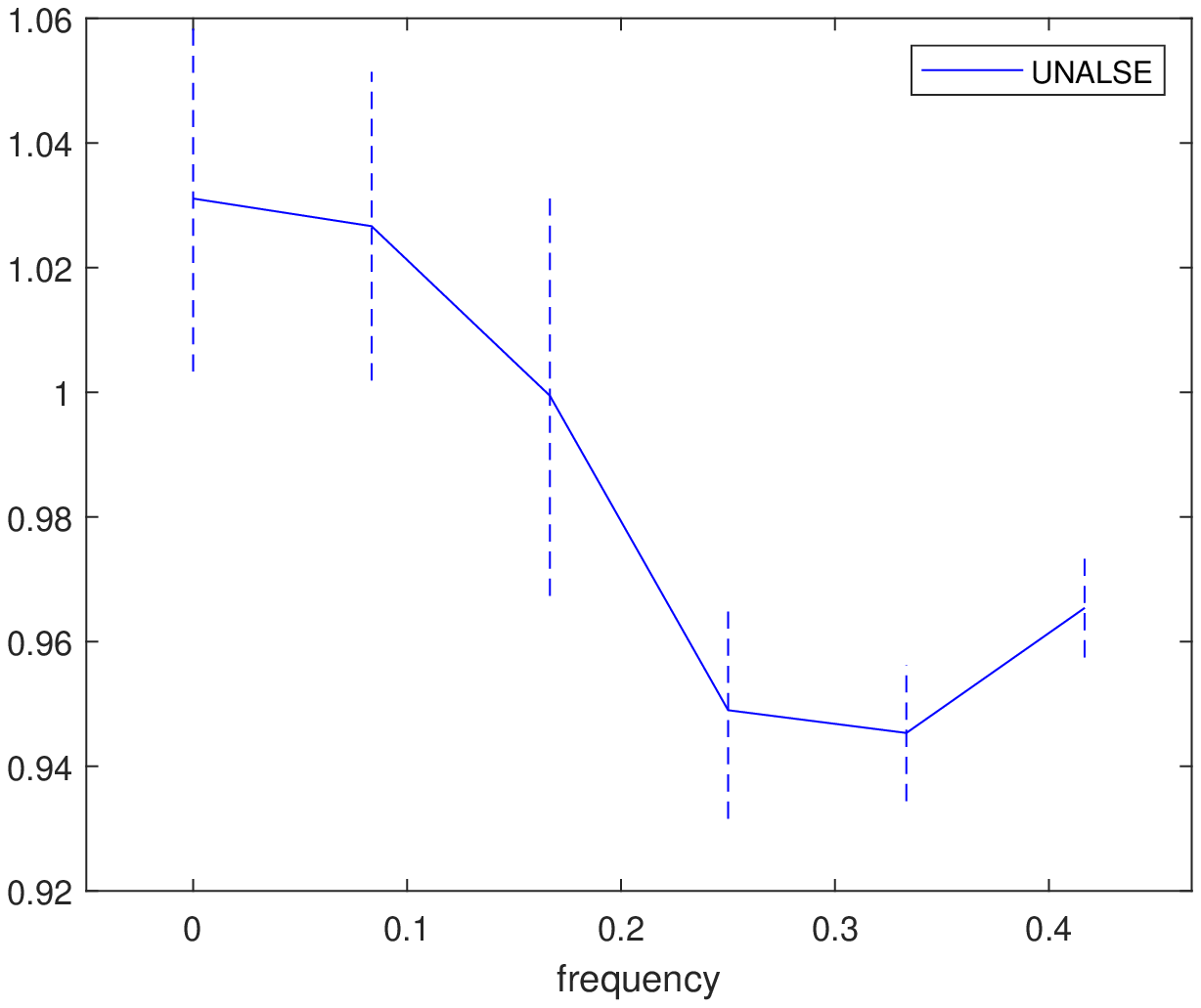}&
              \includegraphics[width=.2\textwidth]{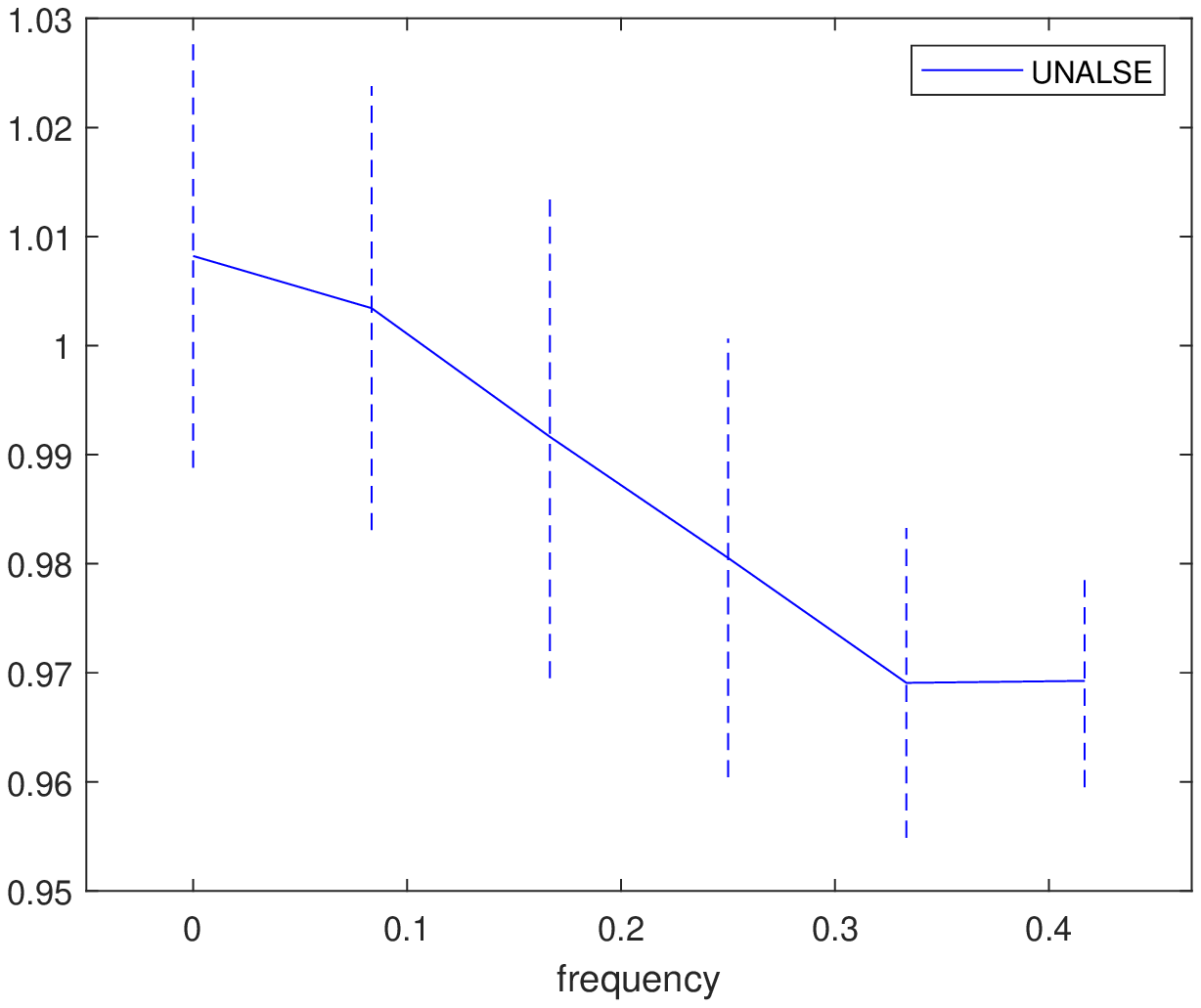}\\
              \end{tabular}
\end{figure}


Scenario C also shows that the proportion of latent variance is estimated very well by UNALSE, which is very close, for Setting 3, even to a target lower than $40\%$ (see Figure \ref{fig:beta_C.1.U}).
Figure \ref{fig:err_L_C.1.U} shows that UNALSE is still very good regarding $err_{\wh{L}}(f_h)$ and $err_{ratio}(f_h)$, particularly for Setting 3 at high frequencies.

\begin{figure}[t!]
          \caption{Estimated latent variance proportion $\wh{\beta}(f_h)$ - Scenario C.}\label{fig:beta_C.1.U}
          \centering
               \begin{tabular}{cc}
                             {\footnotesize Setting 3}&{\footnotesize Setting 4}\\
              \includegraphics[width=.2\textwidth]{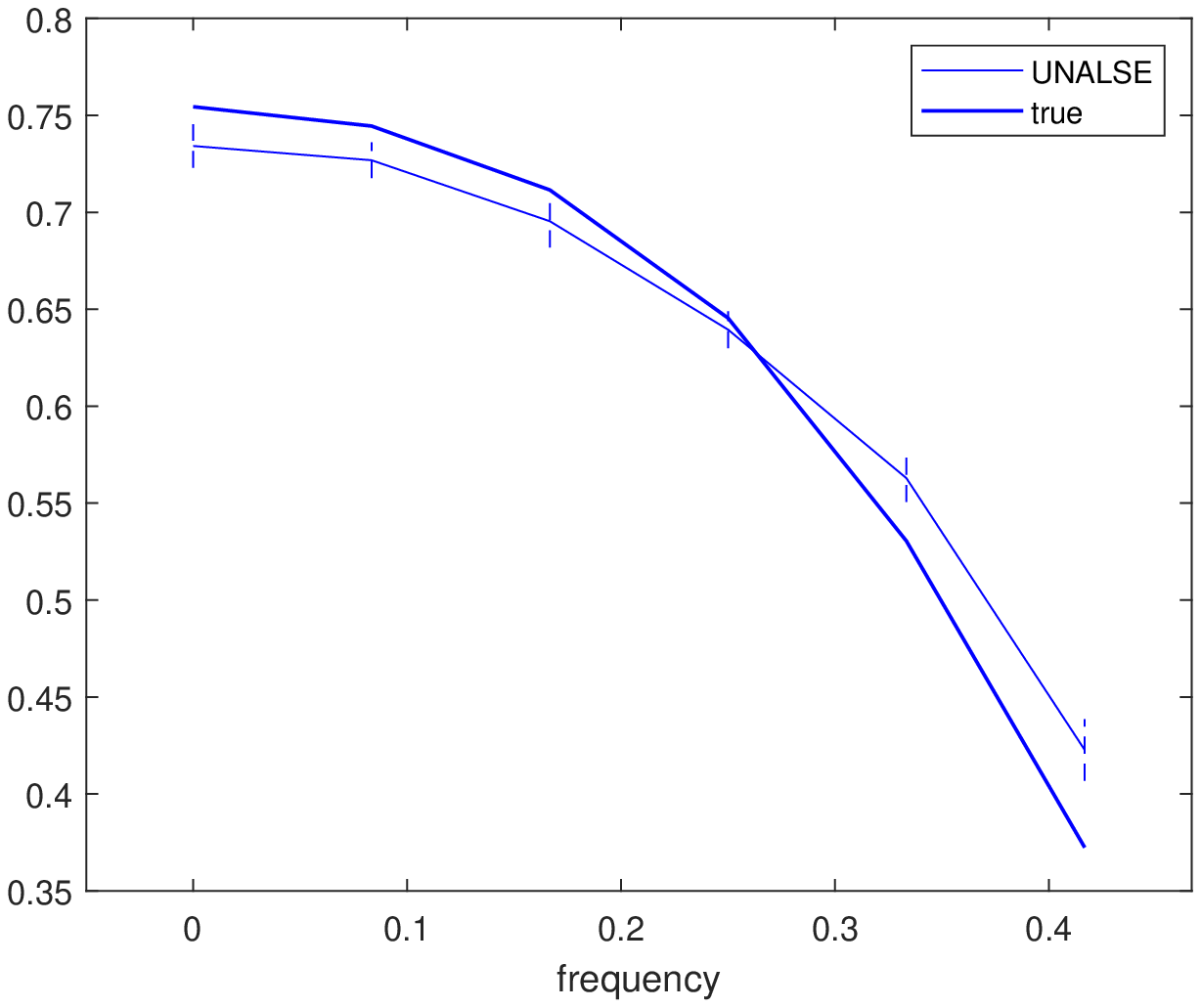}&
              \includegraphics[width=.2\textwidth]{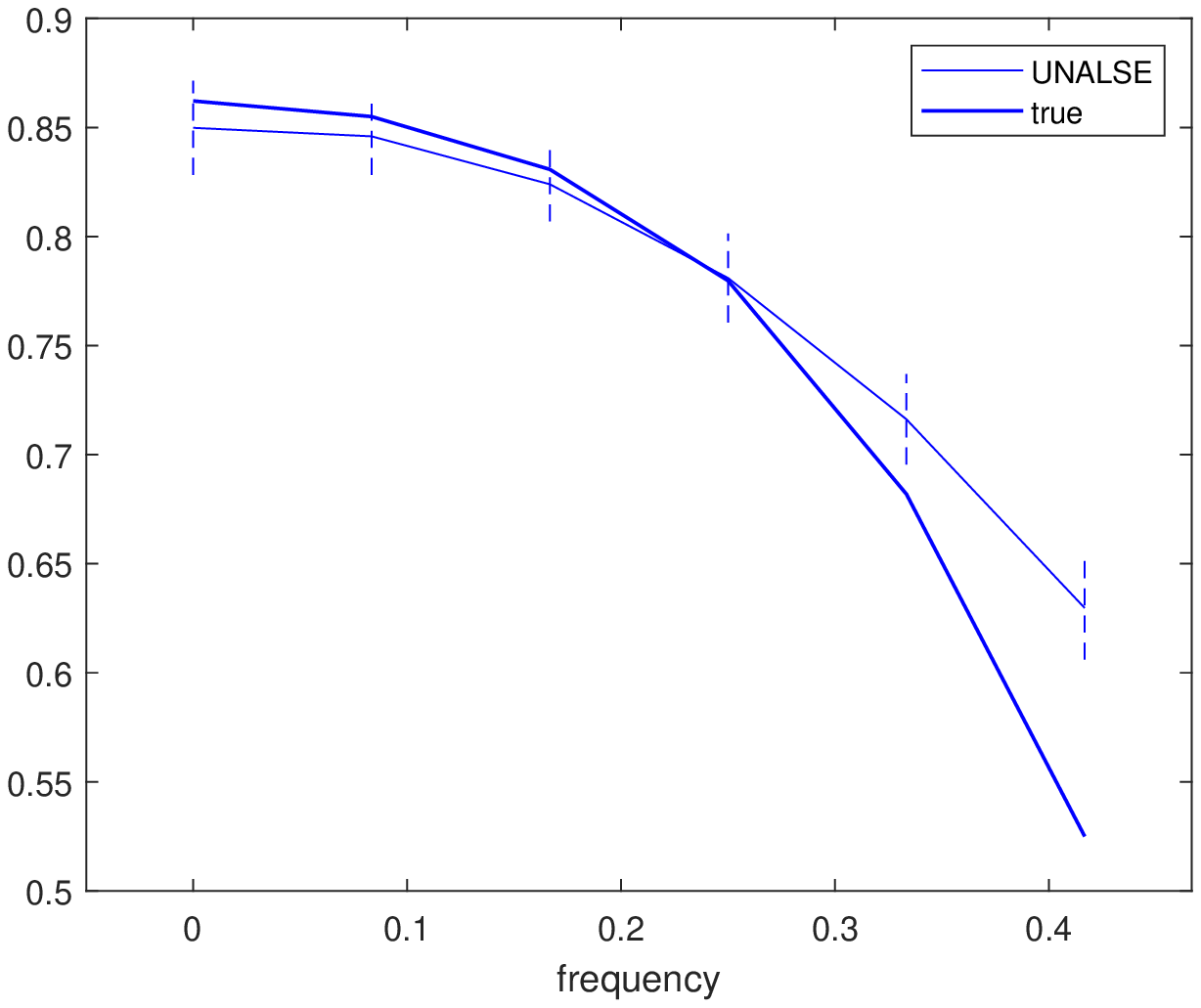}\\
              \end{tabular}
\end{figure}

\begin{figure}[t!]
          \caption{$err_{\wh{L}}(f_h)$ and $err_{ratio}(f_h)$- Scenario C.}\label{fig:err_L_C.1.U}
          \centering
          \begin{tabular}{cccc}
          {\footnotesize $err_{\wh{L}}(f_h)$ Setting 3}&{\footnotesize $err_{\wh{L}}(f_h)$ Setting 4}&{\footnotesize $err_{ratio}(f_h)$ Setting 3}&{\footnotesize 	$err_{ratio}(f_h)$ Setting 4}\\
              \includegraphics[width=.2\textwidth]{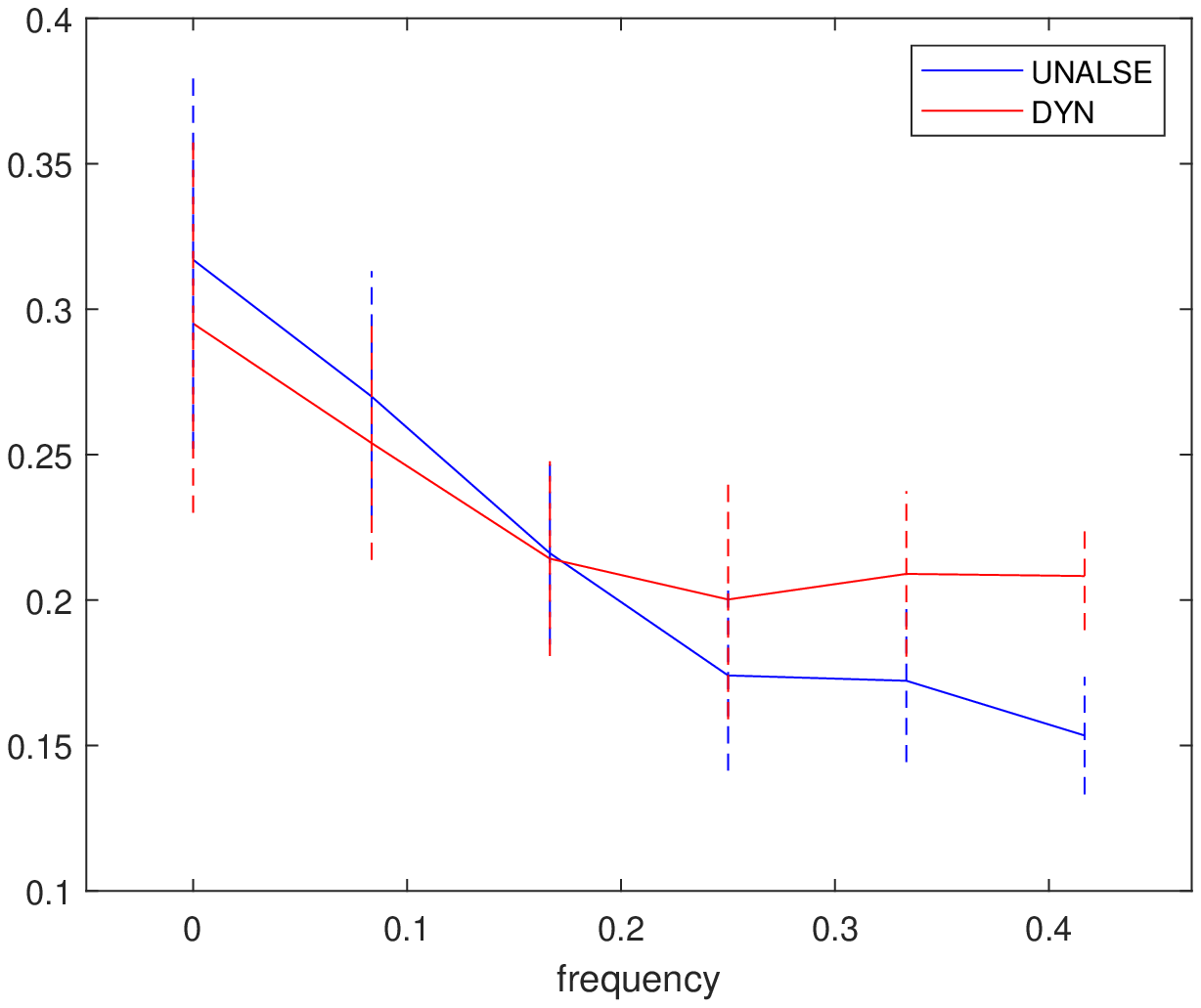}&
              \includegraphics[width=.2\textwidth]{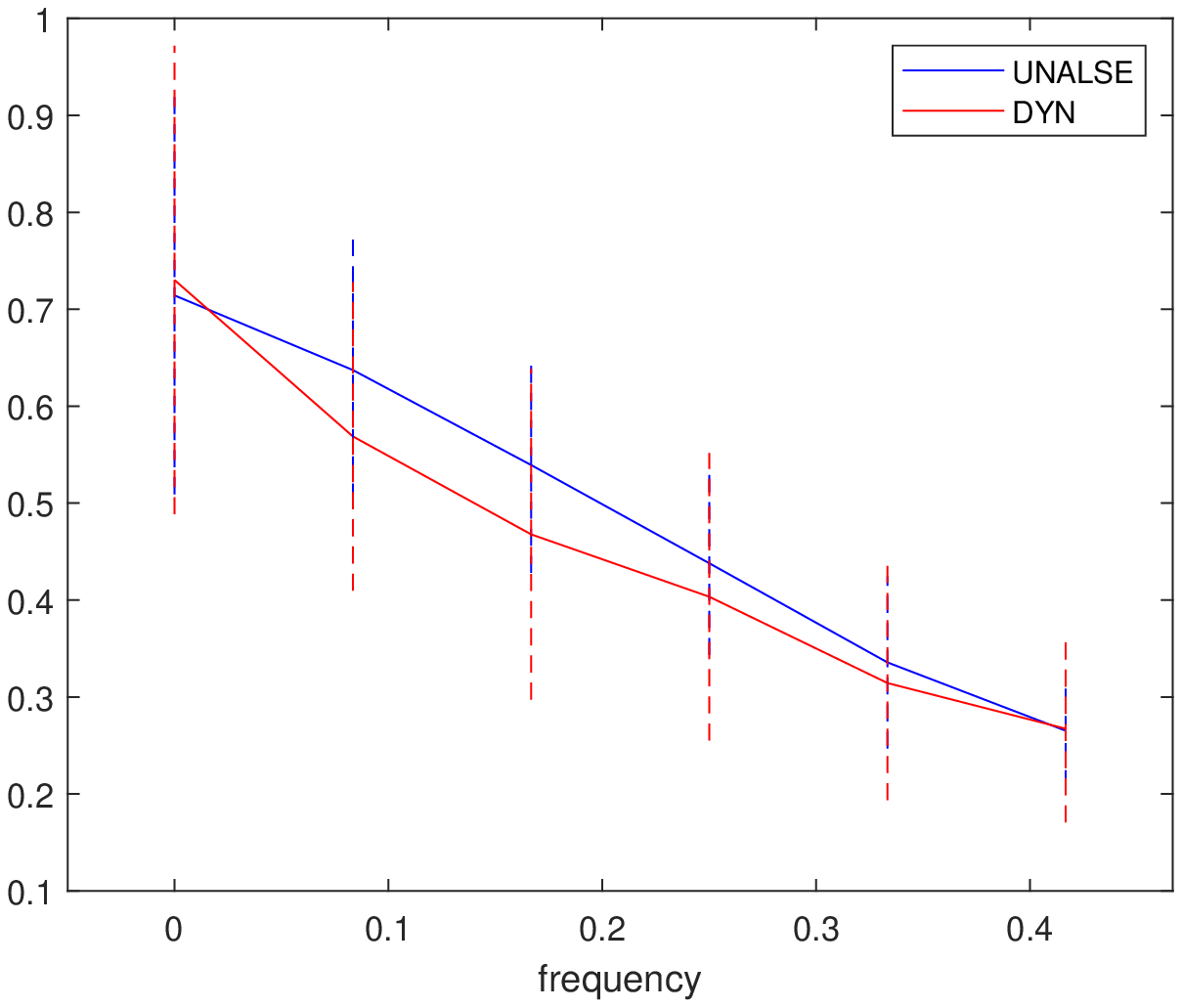}&
              \includegraphics[width=.2\textwidth]{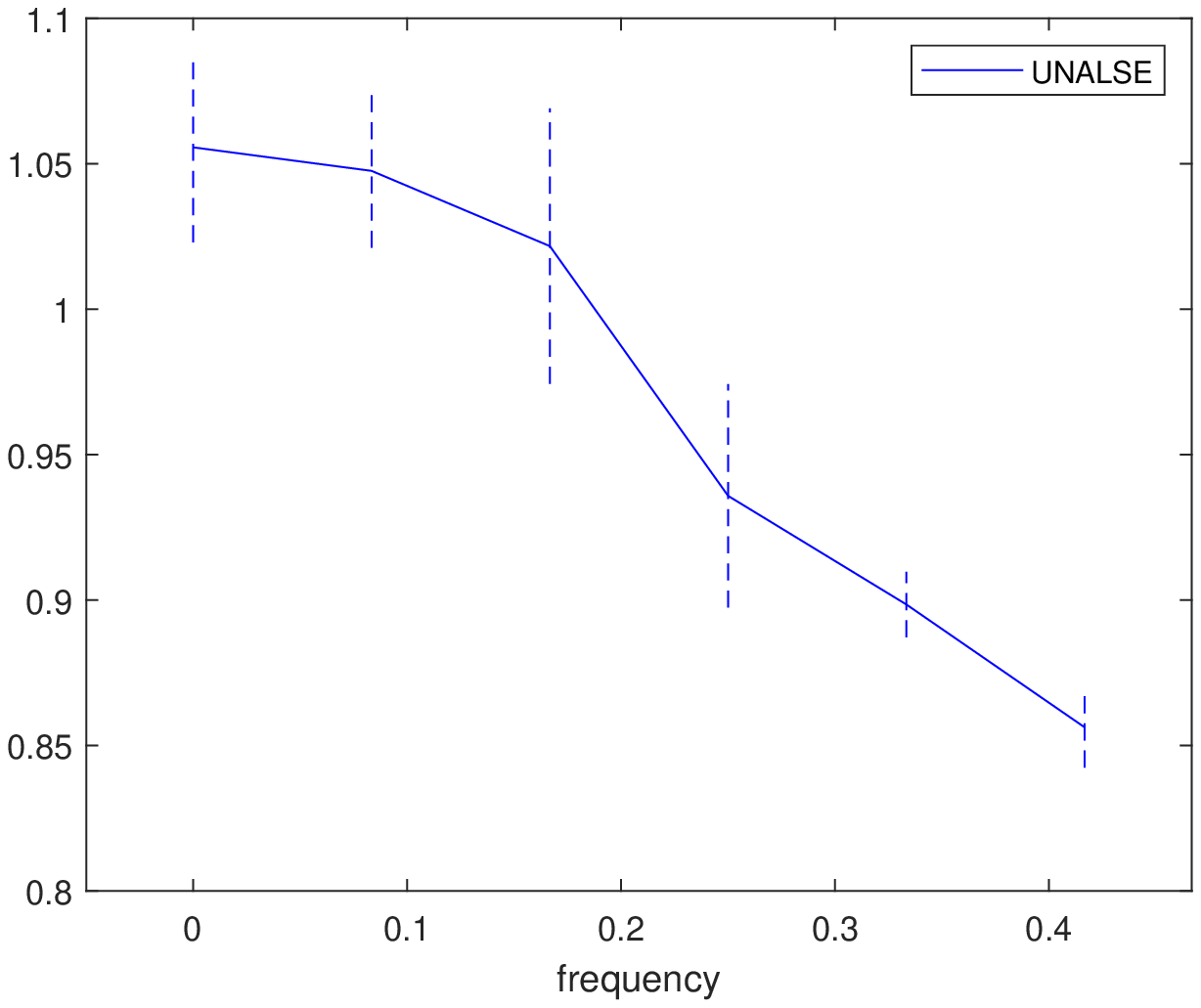}&
              \includegraphics[width=.2\textwidth]{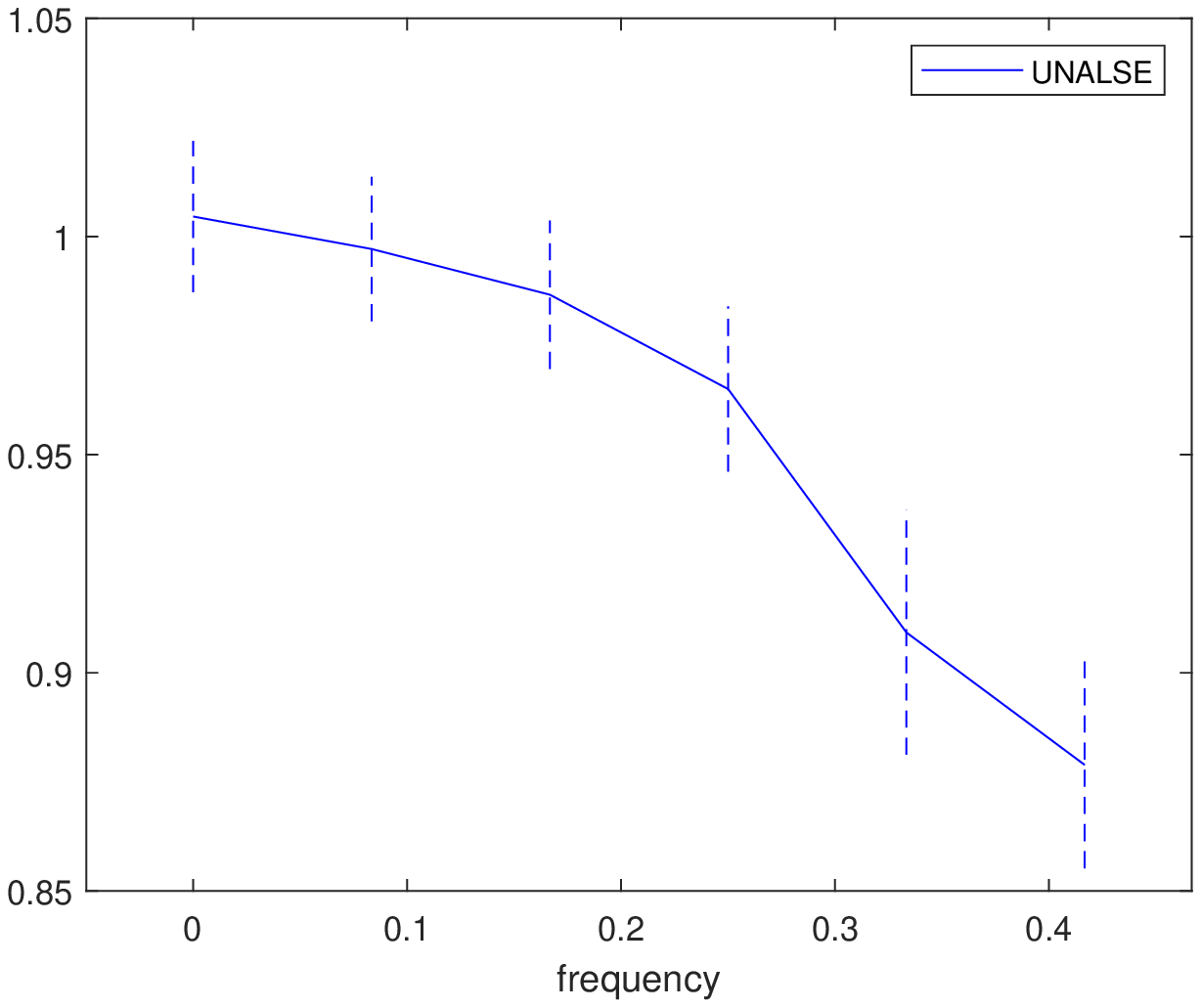}\\
              \end{tabular}
\end{figure}

\begin{figure}[t!]
          \caption{Nonzero predictive values $nzpv(f_h)$ - Scenario C.}\label{fig:pred_pp_C.3.U}
          \centering
          \begin{tabular}{cc}
          {\footnotesize Setting 3}&{\footnotesize Setting 4}\\
              \includegraphics[width=.2\textwidth]{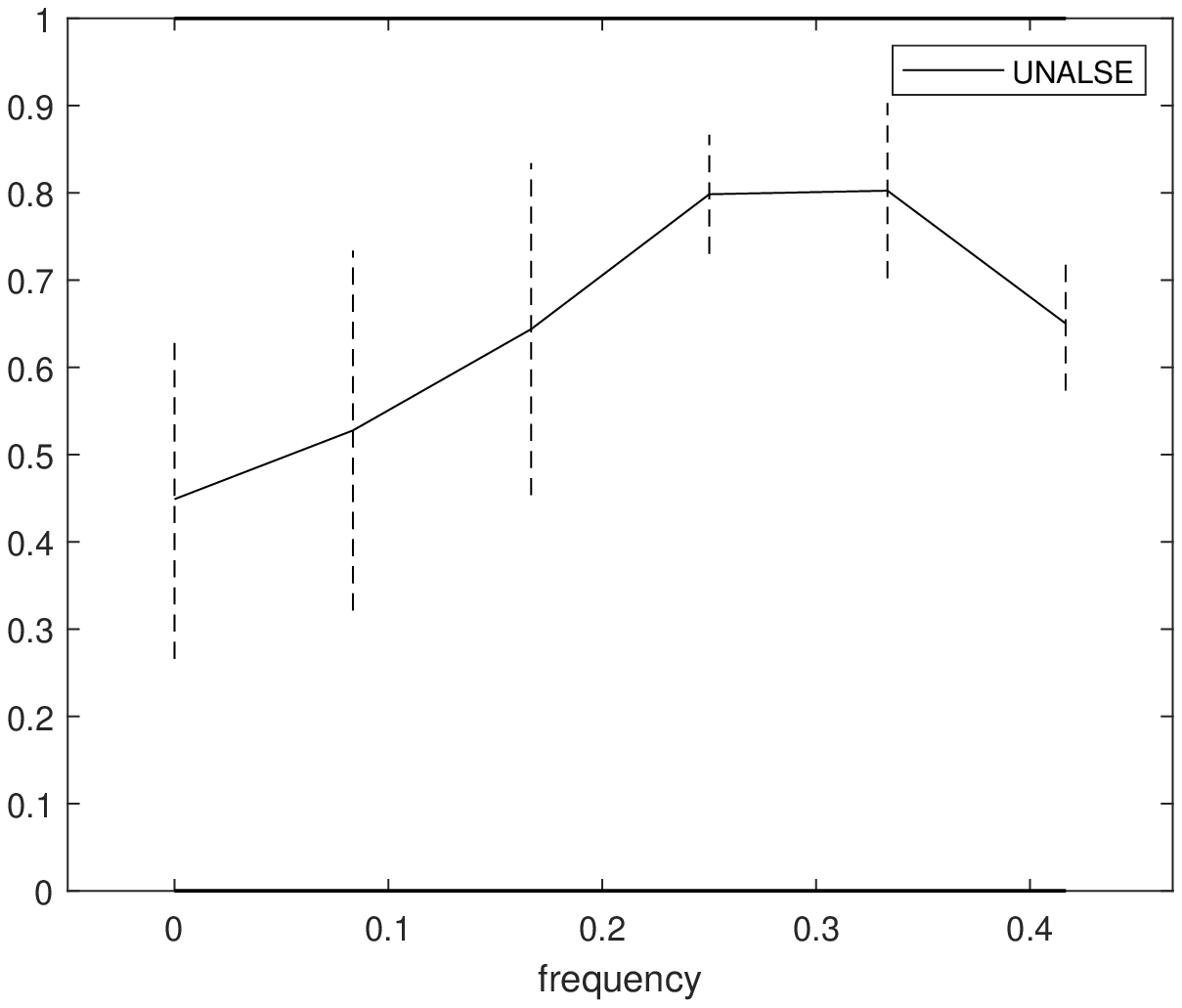}&
              \includegraphics[width=.2\textwidth]{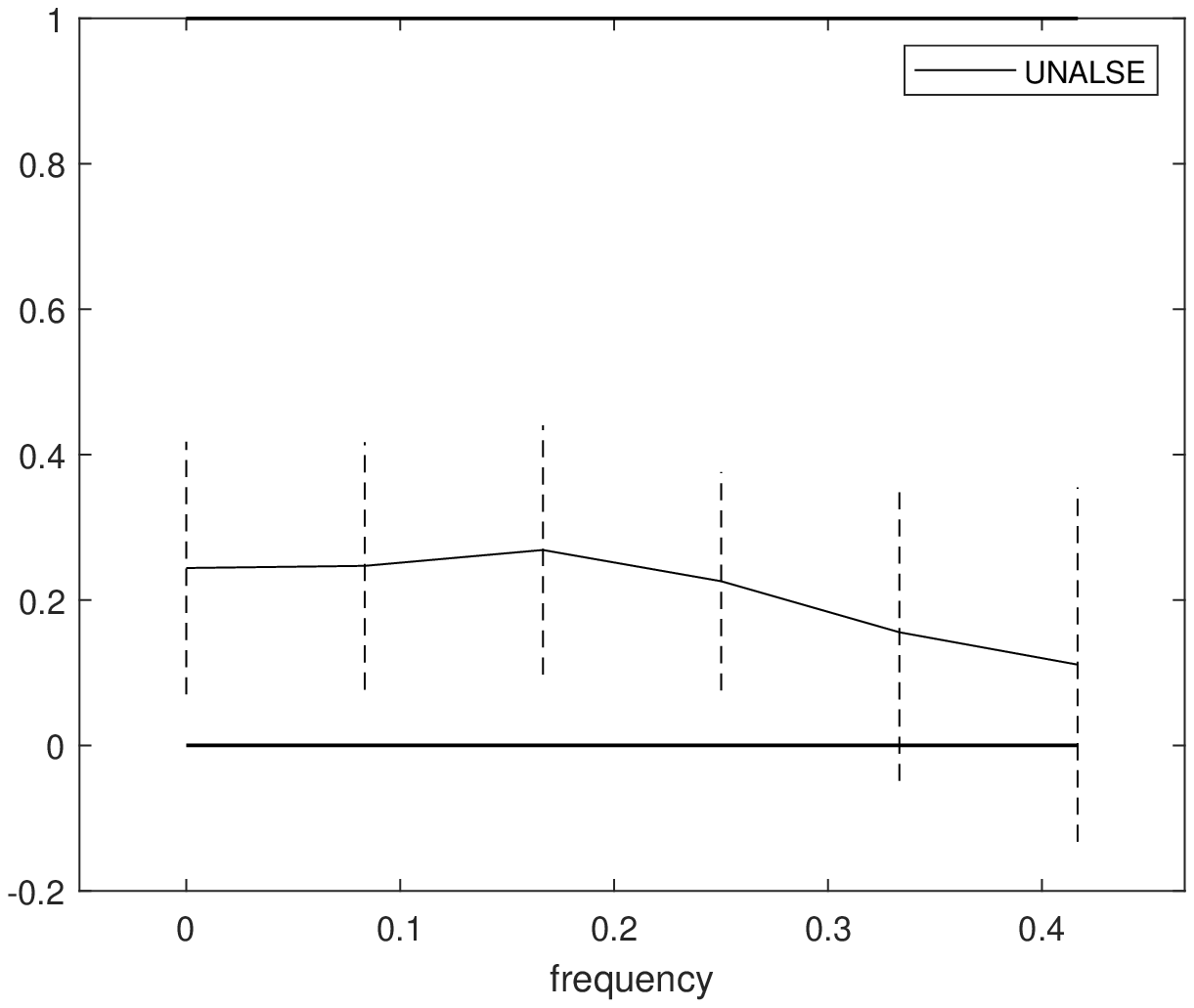}\\
              \end{tabular}
 \end{figure}

\begin{figure}[t!]
          \caption{$mnz_i(f_h)$, $h=0,1,2$ - Scenario C - Setting 3.}               \label{fig:cumulate_C.3.U}
          \centering
          \begin{tabular}{ccc}
          {\footnotesize $f_1=0$}&{\footnotesize $f_2=0.0833$}&{\footnotesize $f_3=0.1667$}\\[-2pt]
              \includegraphics[width=.2\textwidth]{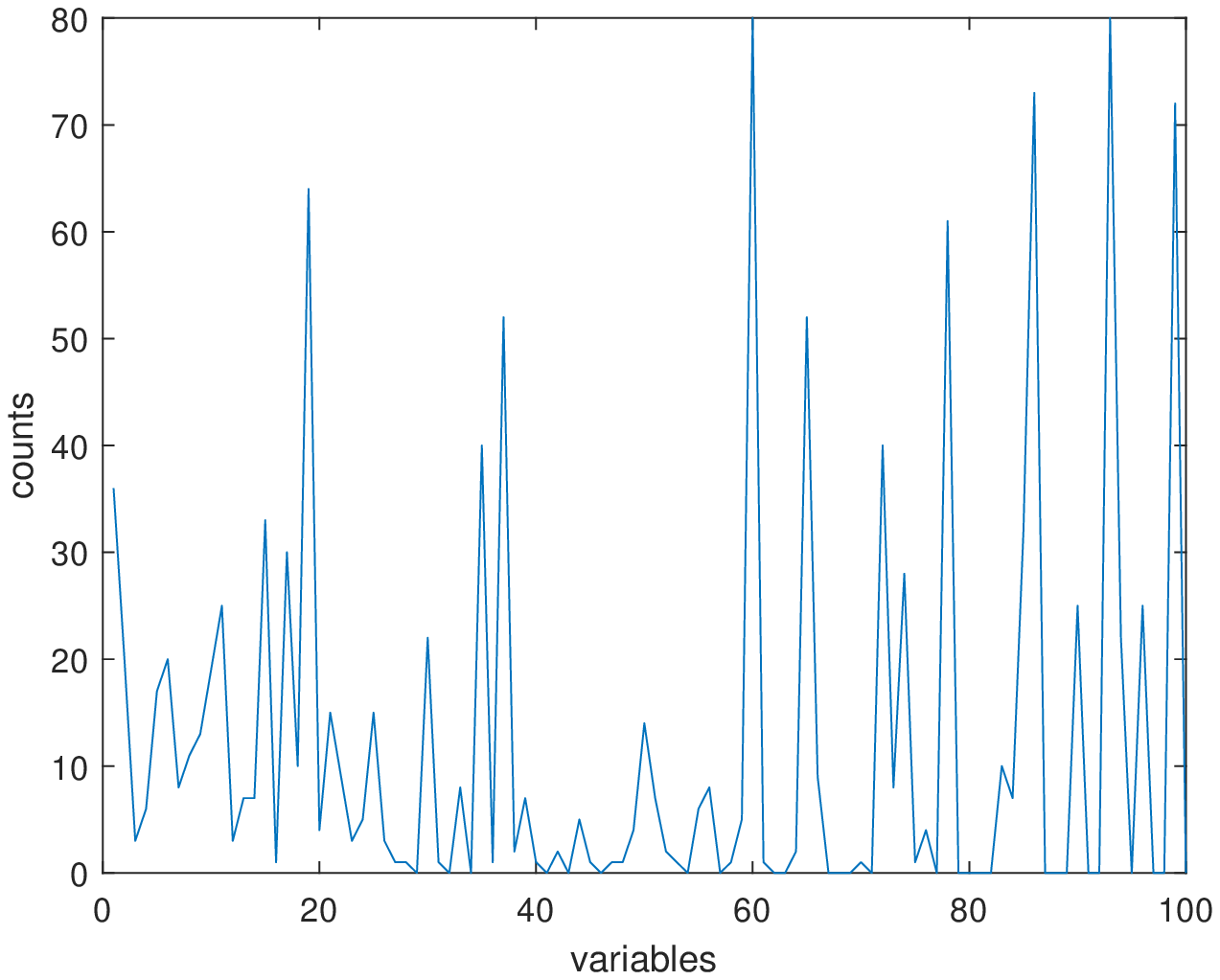}&
              \includegraphics[width=.2\textwidth]{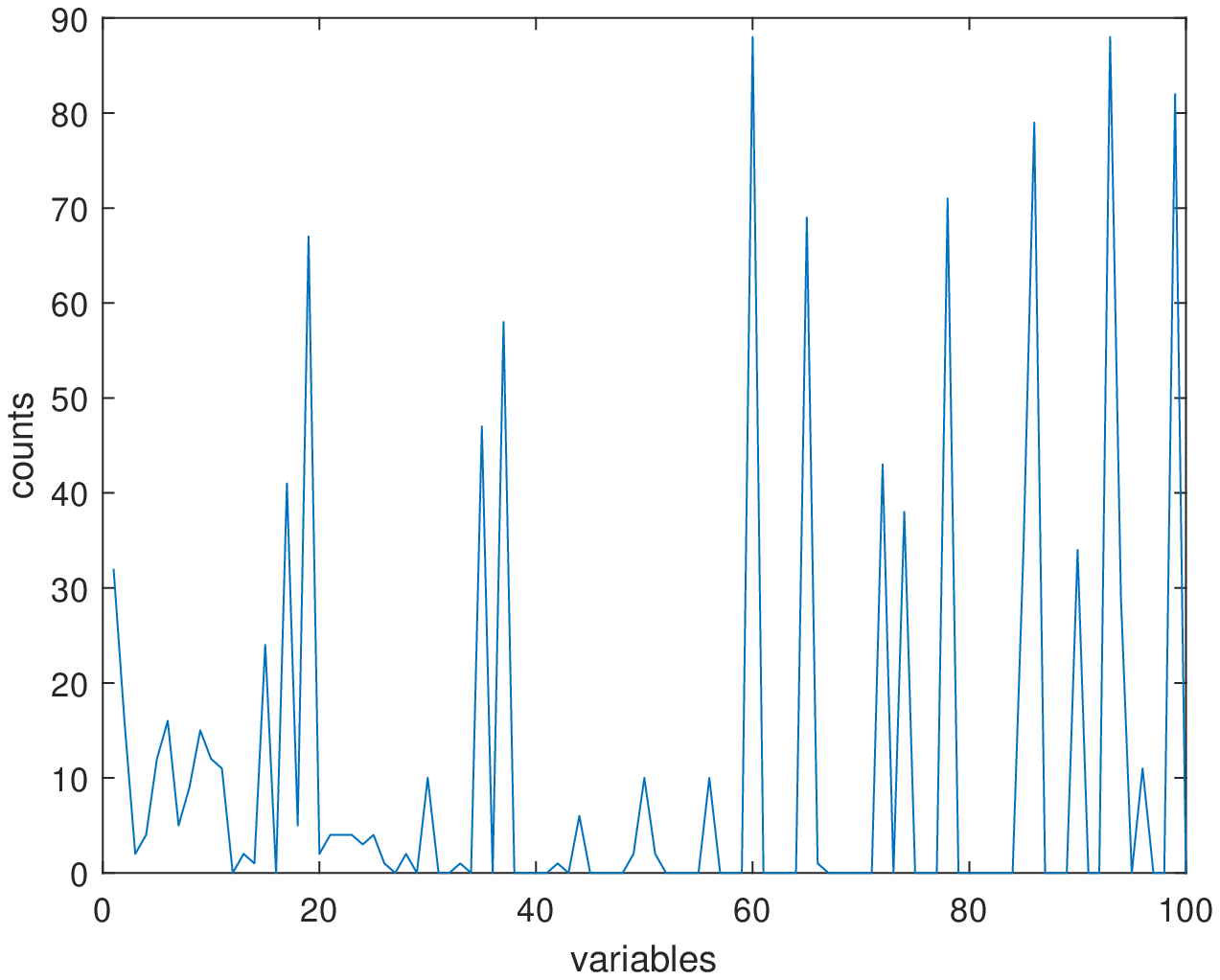}&
              \includegraphics[width=.2\textwidth]{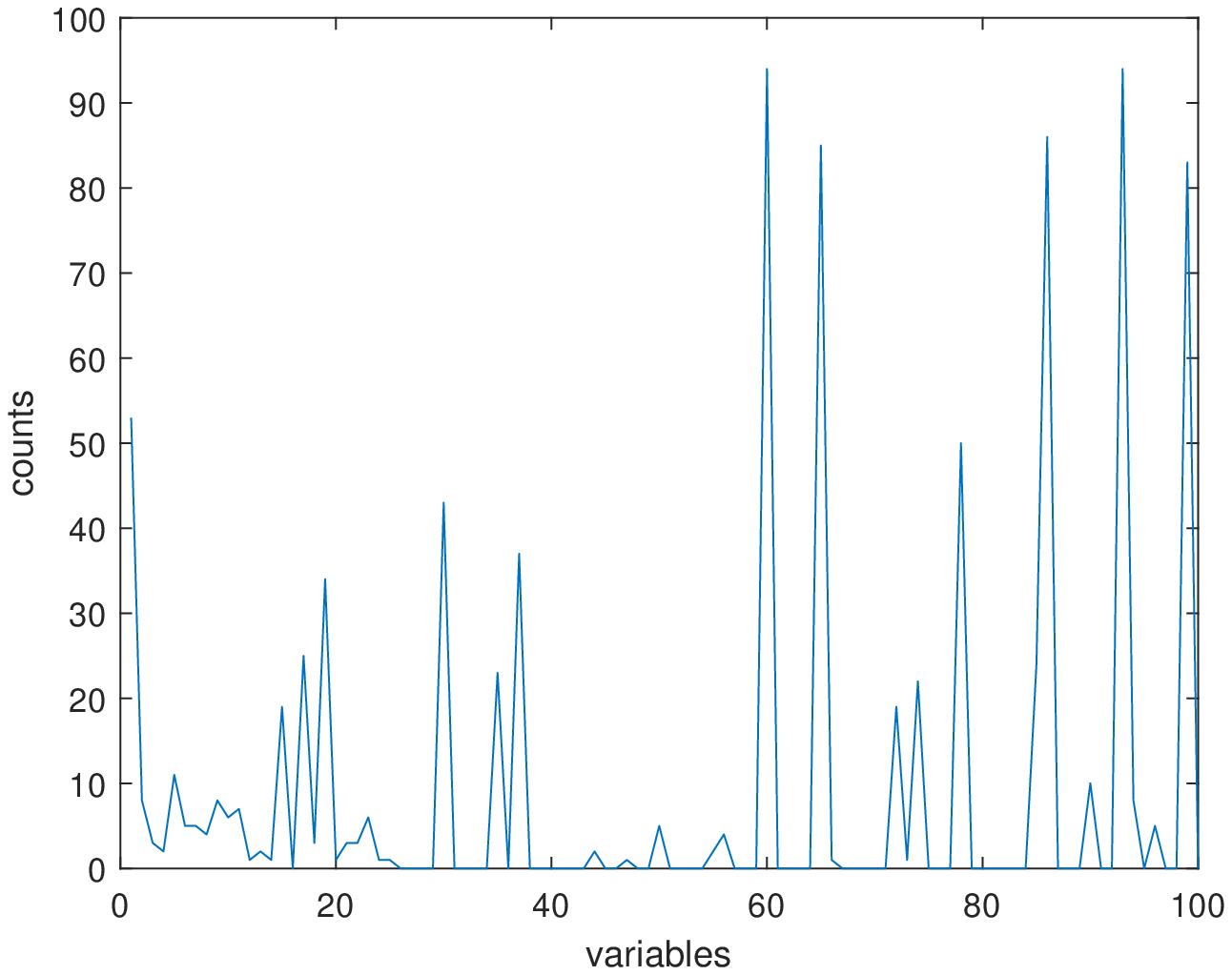}\\
              \end{tabular}
\end{figure}

Concerning the sparsity pattern, the nonzeros recovered by UNALSE are stable across frequencies and the nonzero predictive value overcomes $70\%$ at central frequencies for Setting 3 and it is still acceptable for Setting 4 (see Figure \ref{fig:pred_pp_C.3.U}).
Finally, Figure \ref{fig:cumulate_C.3.U} shows how the indicator $mnz_i(f_h)$, $h=0,1,2$,
presents a very similar pattern across variables for the first three frequencies, proving that the
sparsity pattern is consistent over frequencies. Note that this property is common to all Scenarios and Settings, even when the predictive value is not good.

\section{Real data analysis}\label{sec:real}
We consider a  dataset of $p=101$ quarterly macroeconomic indicators regarding the US economy observed over $T=210$ time points spanning the period 1960:Q2-2012:Q3 \citep[see, e.g.,][]{mccracken2020fred}. Throughout, we compute the smoothed periodogram by setting $M_T=\lfloor \sqrt{T}\rfloor=14$ and using the Bartlett kernel. The analysis that follows shows that the underlying spectral structure seems to be quite relevant at certain frequencies and heterogeneous.

\begin{figure}[t!]
          \caption{US macroeconomic data - Co-movements and sparsity.}\label{fig:macro1}
                  \centering
          \begin{tabular}{cccc}
              \includegraphics[width=.2\textwidth]{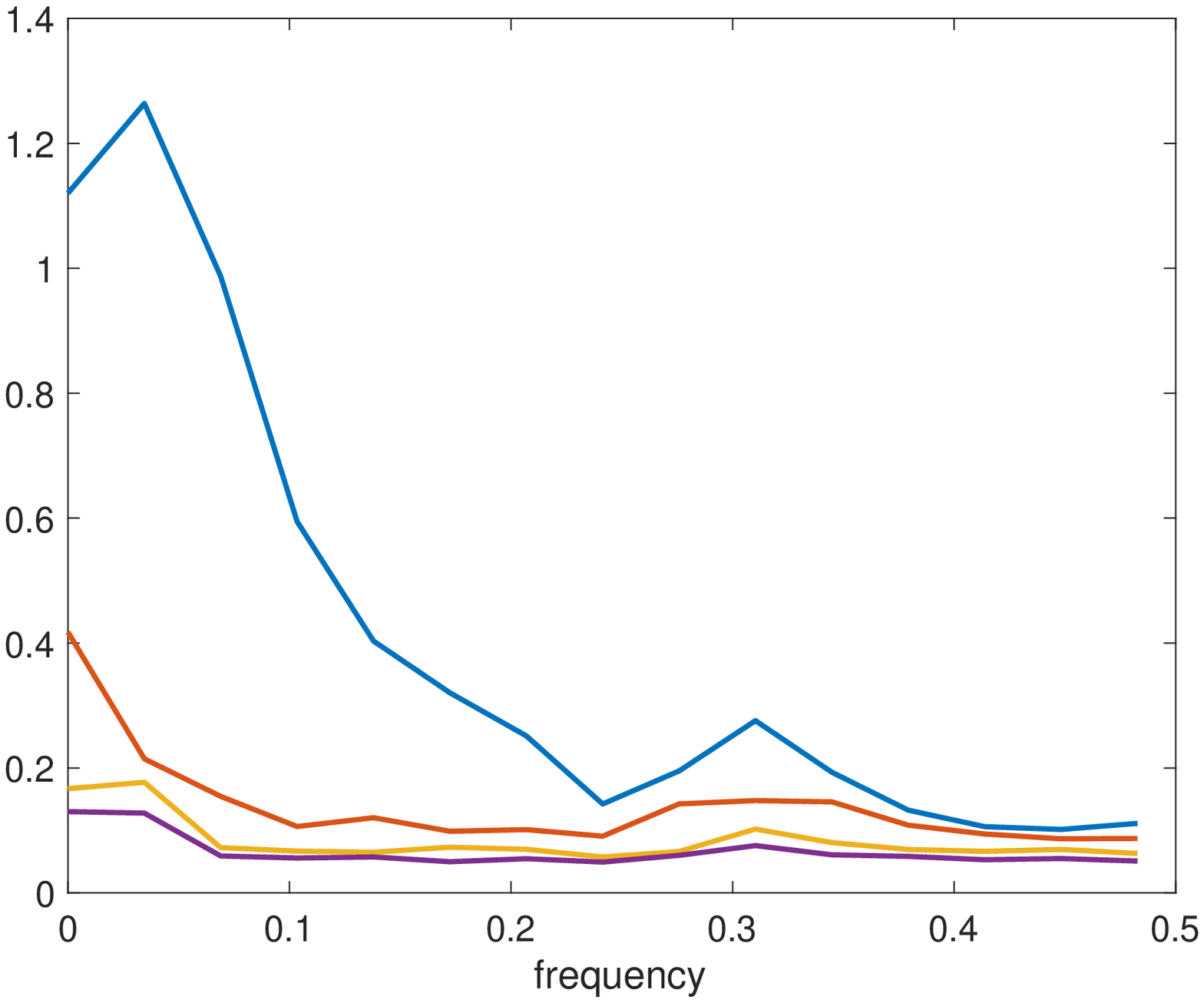}&
              \includegraphics[width=.2\textwidth]{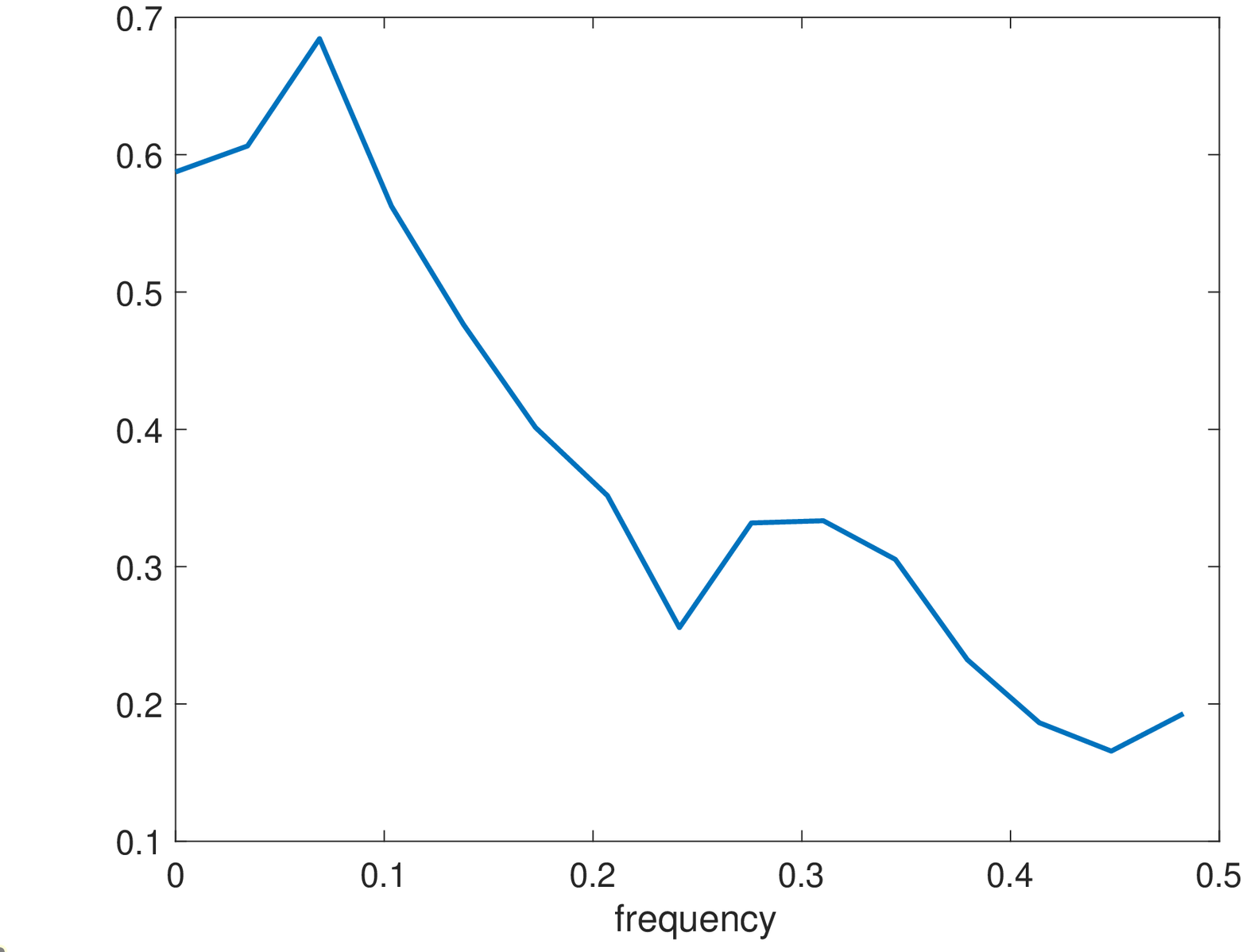}&
                            \includegraphics[width=.2\textwidth]{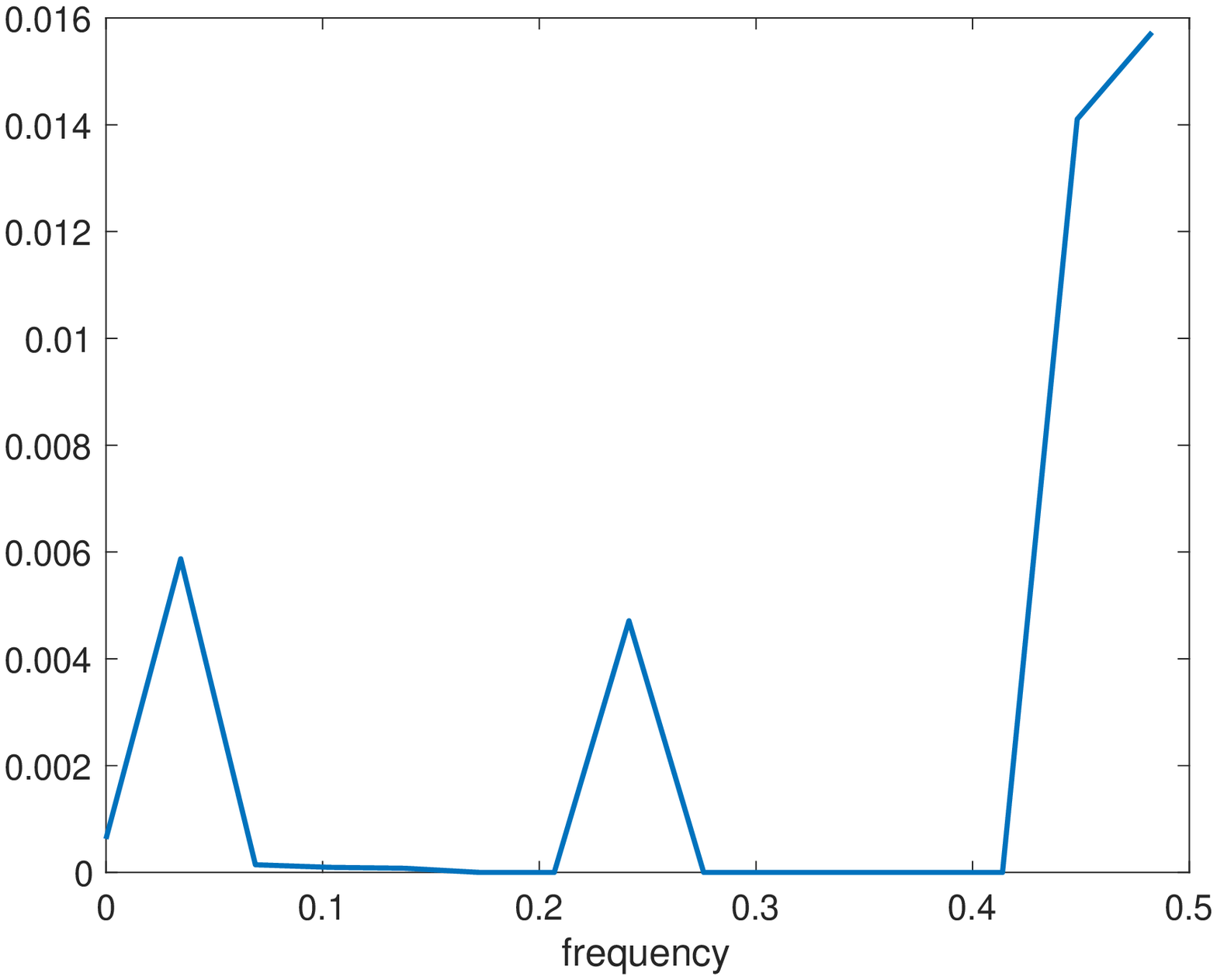}&
              \includegraphics[width=.2\textwidth]{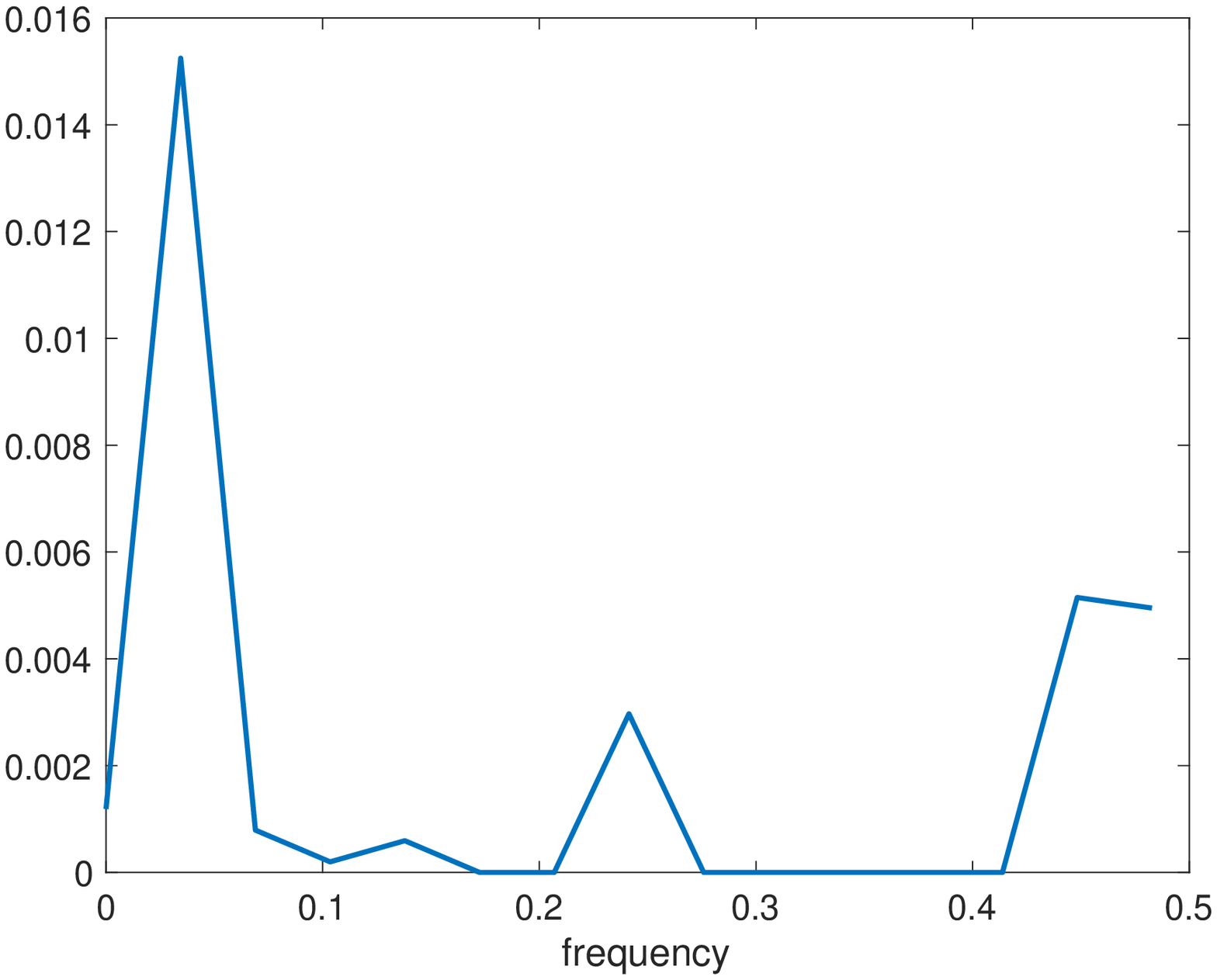}\\
              {\footnotesize $\frac{\wh{\lambda}_j(\wt{\Sigma}(f_h))}p$, $j=1,2,3,4$.}&{\footnotesize $\wh{\beta}(f_h)$.}&
              	      {\footnotesize $\wh{\zeta}(f_h)$.}&{\footnotesize fraction of nonzeros in $\wh S(f_h)$.}\\
              \end{tabular}
\end{figure}

In the top left panel of Figure \ref{fig:macro1} we show the four largest eigenvalues, rescaled by $p$, of the smoothed periodogram estimator.
The top eigenvalue shows a decreasing shape from $f=0$ and two auxiliary peaks at $f=0.07$ and $f=0.3$, corresponding to periods of about 3.5 years and 9 months, respectively. Note that 3.5 years is around the typical period of a business cycle. The estimated rank by UNALSE is $\wh r=2$ at all frequencies. The top right panel of Figure \ref{fig:macro1} shows the proportion of latent variance $\wh{\beta}(f_h)=\frac{\text{rk}(\wh{L}(f_h))}{\text{rk}(\wh{\Sigma}(f_h))}$ so the contribution of $\wh{L}(f_h)$, which follows the pattern of the leading eigenvalues of $\wt{\Sigma}(f_h)$, hence it captures the business cycle frequency. The bottom left panel of Figure \ref{fig:macro1} reports the proportion of residual covariance $\wh{\zeta}(f_h)=\frac{\sum_{i=1}^p\sum_{j=i+1}^p\vert \wh{S}_{ij}(f_h)\vert}{\sum_{i=1}^p\sum_{j=i+1}^p\vert \wh{\Sigma}_{ij}(f_h)\vert}$, summarizing the contribution of $\wh{S}(f_h)$, which has the main contribution at a higher frequency $f=0.45$, corresponding to a period of 6 months. Secondary maxima are at $f=0.03$, i.e., a period of 7 years, and $f=0.25$, corresponding to a period of 1 year. Finally, the fraction of nonzeros has a similar pattern (see bottom right panel of Figure \ref{fig:macro1}).

\begin{figure}[t!]
          \caption{US macroeconomic data - $\wh{L}(f_h)$.}\label{fig:heatL}
                  \centering
               \begin{tabular}{ccc}
               \includegraphics[width=.2\textwidth]{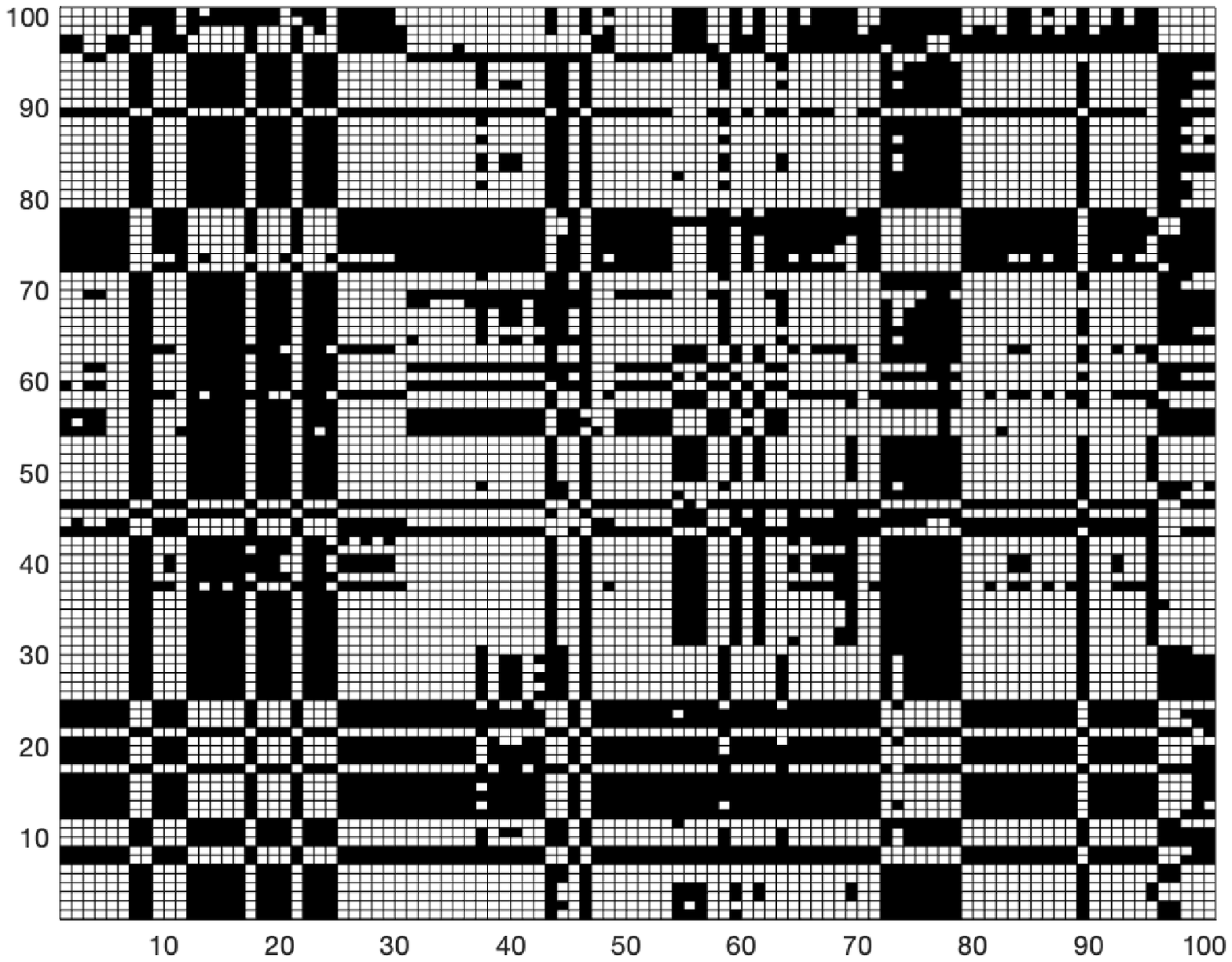}&
              \includegraphics[width=.2\textwidth]{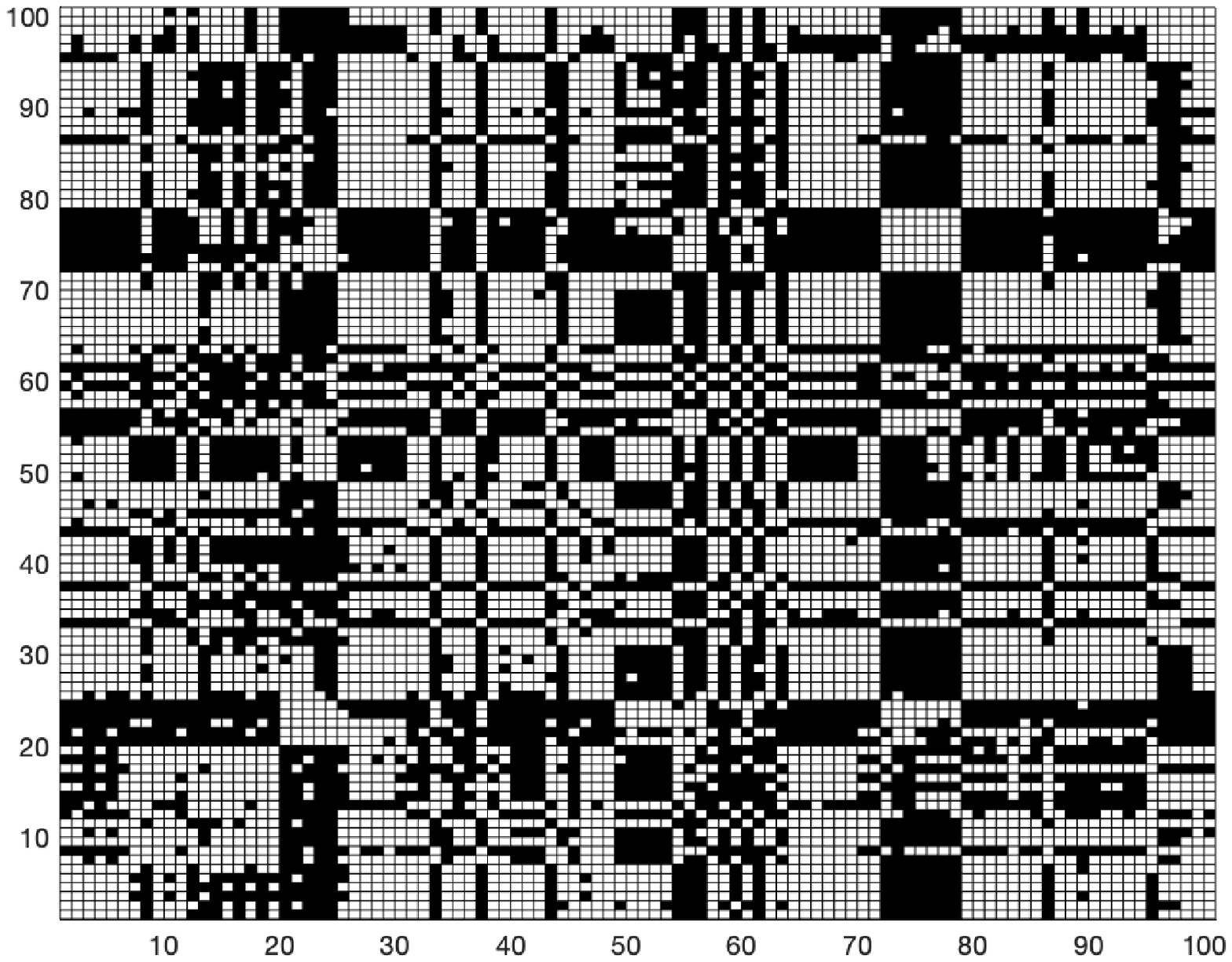}&
              \includegraphics[width=.2\textwidth]{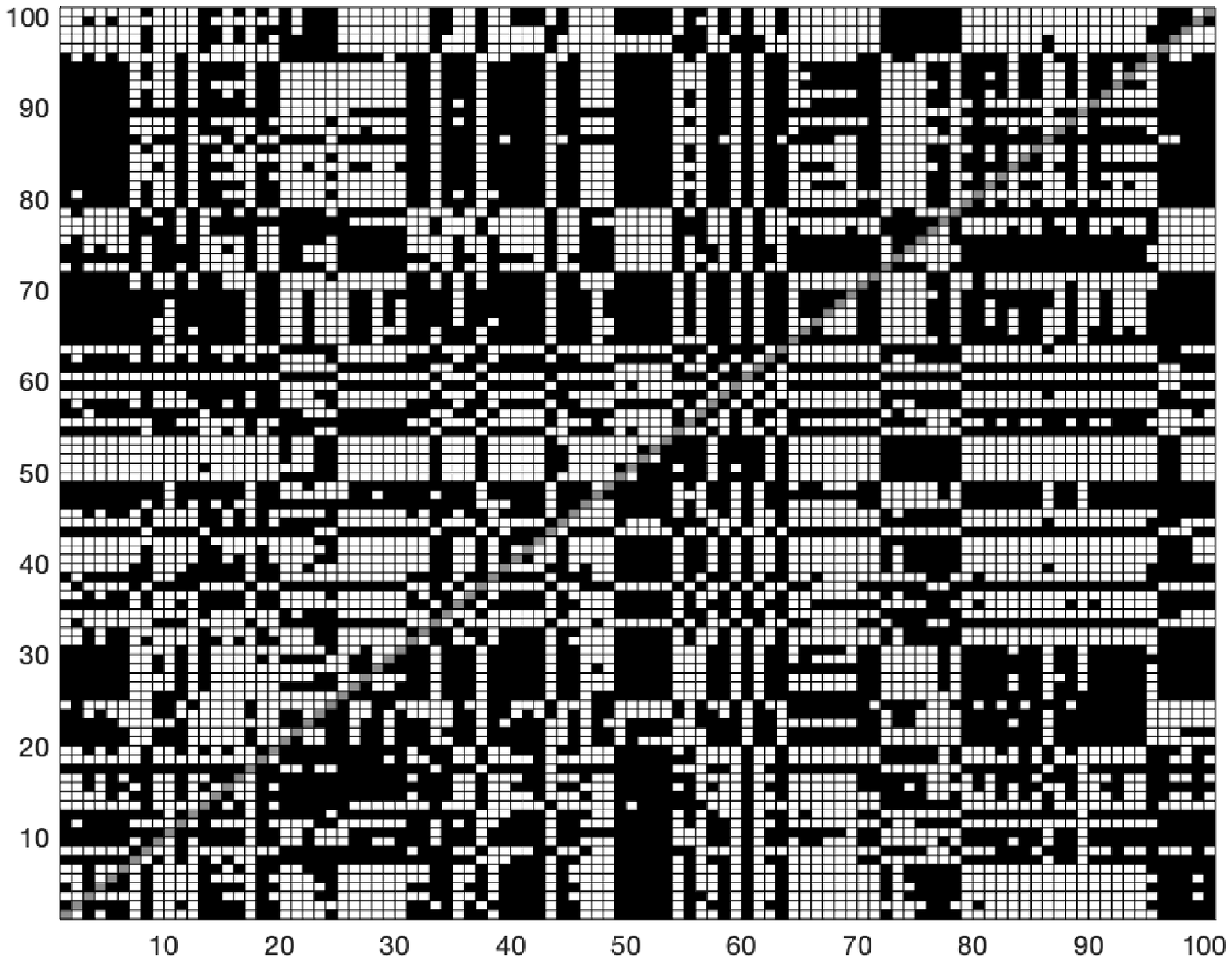}\\[-10pt]
              {\footnotesize $\wh{L}(f_h=0)$.}& {\footnotesize $\text {Re}(\wh{L}(f_h=0.07))$.}&{\footnotesize $\text {Im}(\wh{L}(f_h=0.07))$.}\\
              \end{tabular}
\end{figure}

In Figure \ref{fig:heatL} we show heat-maps of $\wh{L}(f_h)$ at frequencies $0$ and $0.07$.
The variables having the largest (in term of explained variance) common component are all related to the labor market:
\begin{inparaenum}[(i)]
\item Civilian Employment;
\item All Employees in Service-Providing Industries, Total Private Industries, and in Trade, Transportation \& Utilities.
\end{inparaenum}
%

\begin{figure}[t!]
          \caption{US macroeconomic data - $\wh{S}(f_h)$.}\label{fig:heatS}
                  \centering
          \begin{tabular}{cccc}
              \includegraphics[width=.2\textwidth]{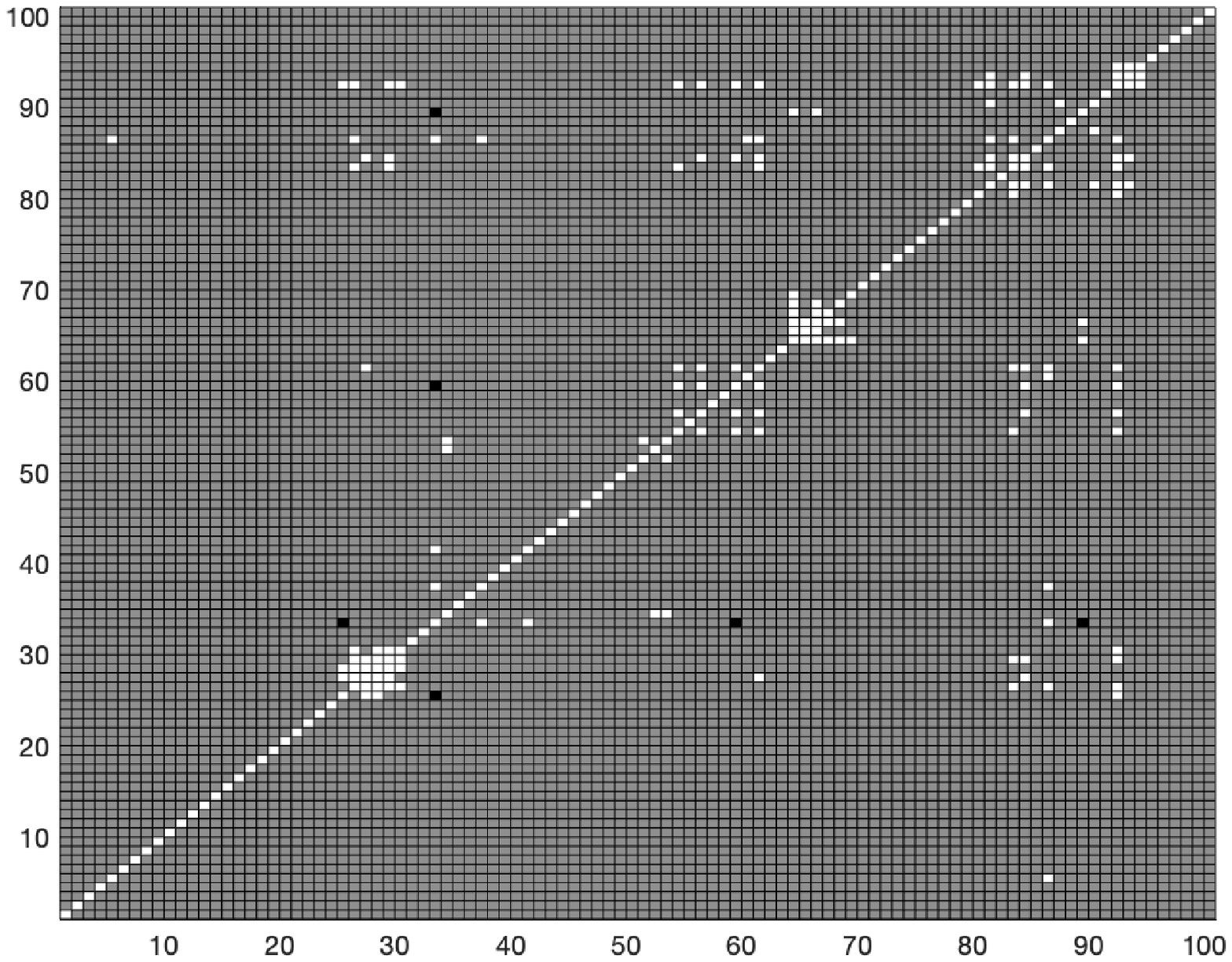}&
              \includegraphics[width=.2\textwidth]{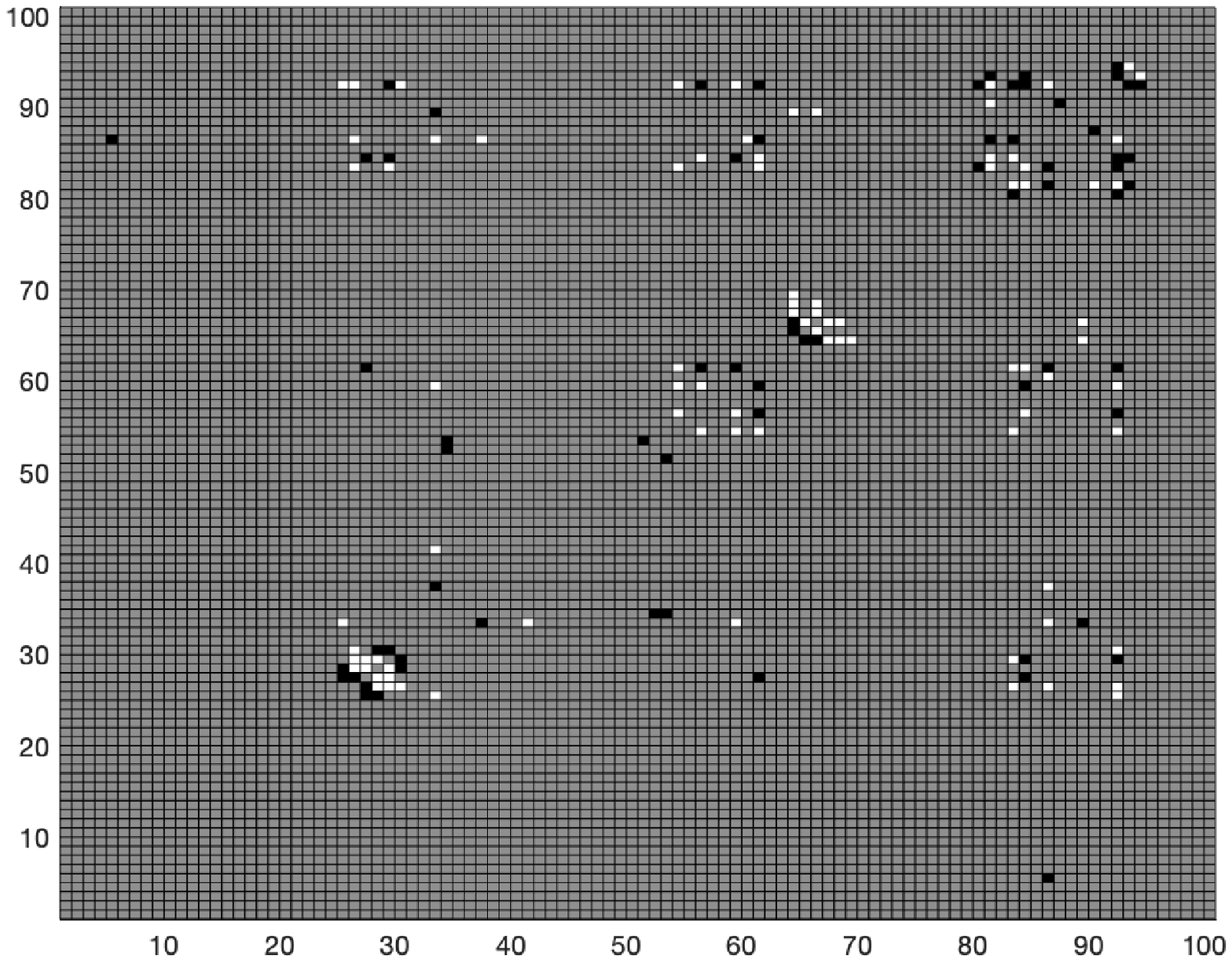}&
              \includegraphics[width=.2\textwidth]{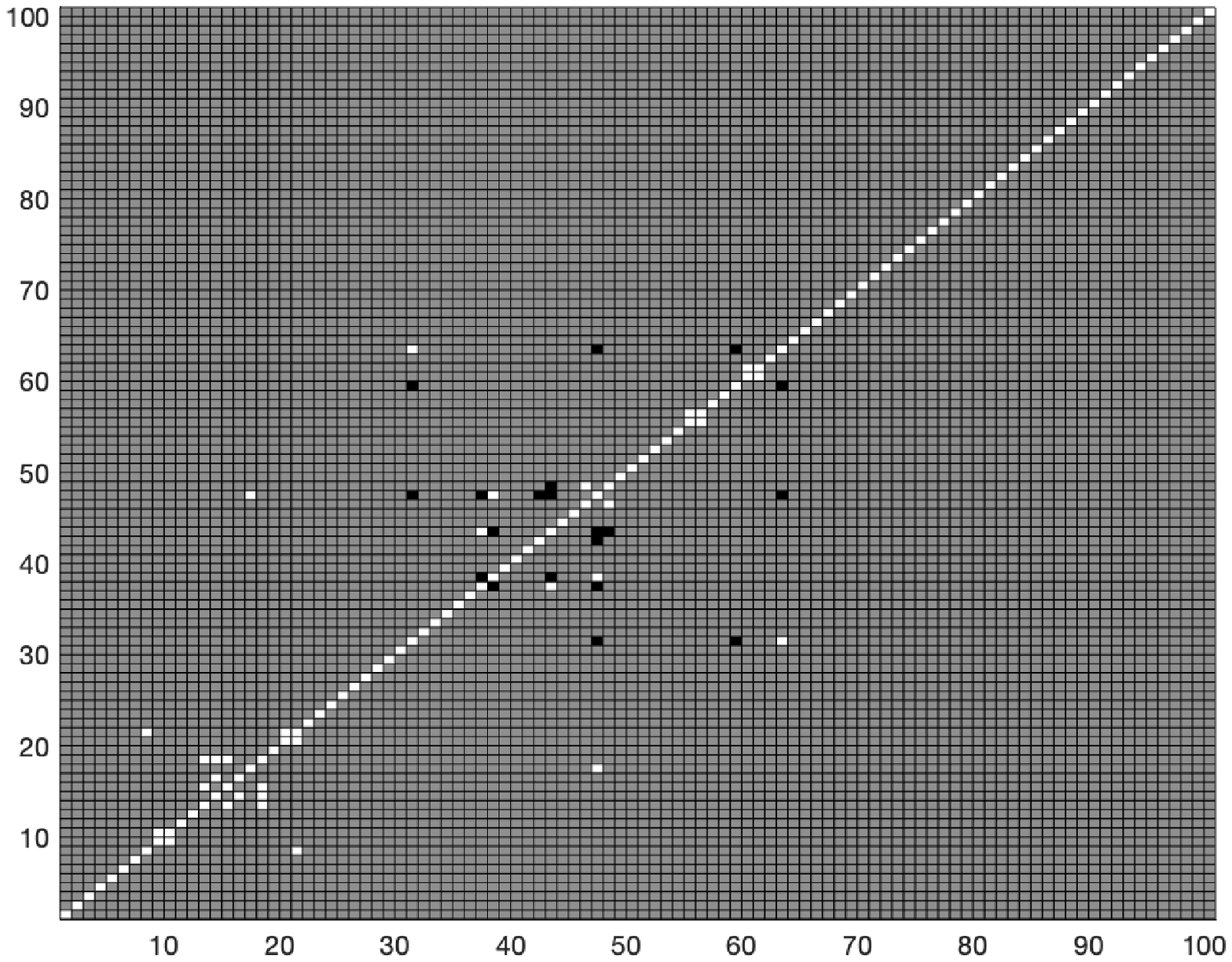}&
              \includegraphics[width=.2\textwidth]{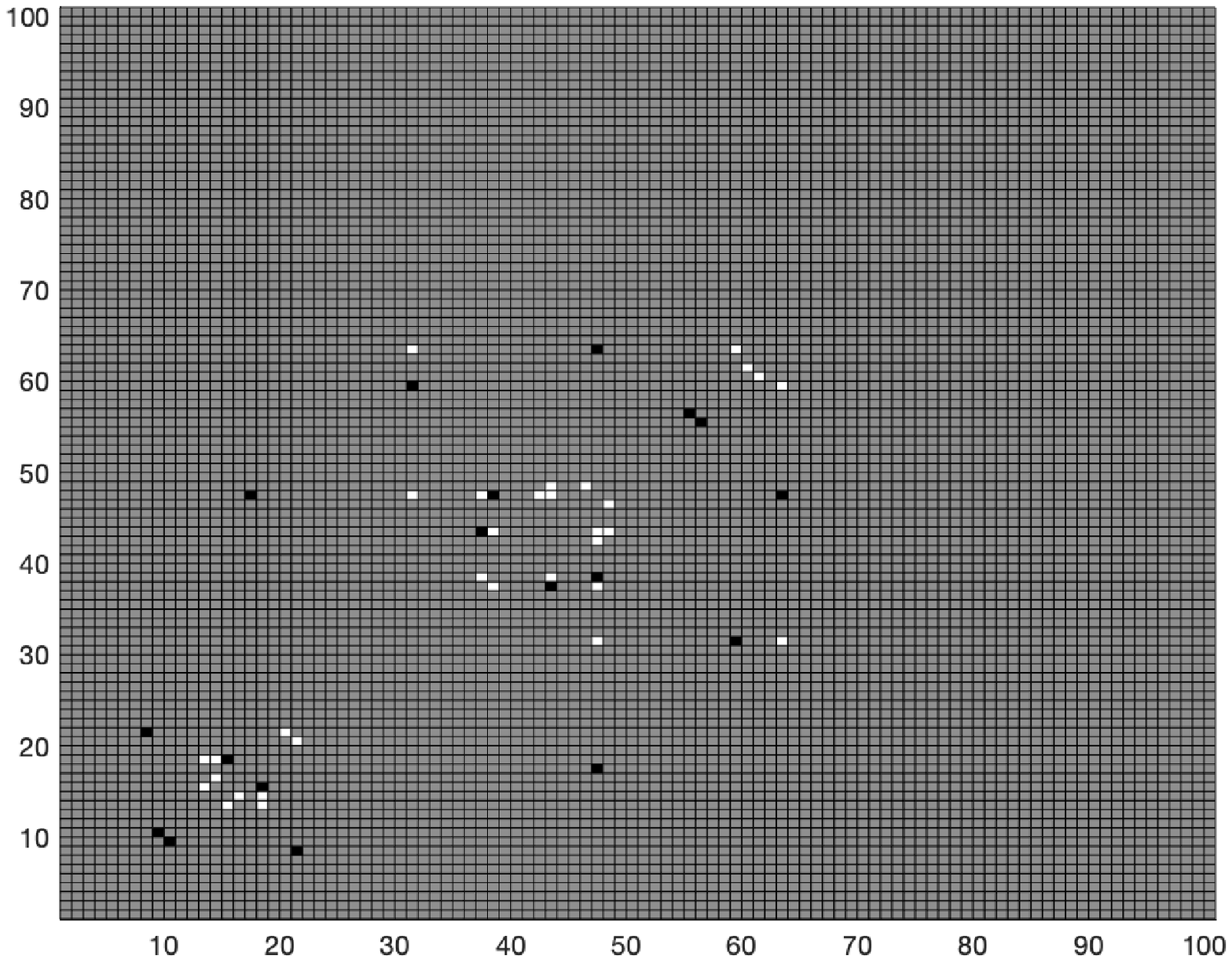}\\[-10pt]
              {\footnotesize $\text{Re}(\wh{S}(f_h=0.03))$.}&{\footnotesize $\text {Im}(\wh{S}(f_h=0.03))$.}&
               {\footnotesize $\text{Re}(\wh{S}(f_h=0.45))$.}&{\footnotesize $\text {Im}(\wh{S}(f_h=0.45))$.}\\
               \end{tabular}
\end{figure}

In Figure \ref{fig:heatS} we show heat-maps of $\wh{S}(f_h)$ at frequencies $0.03$ and $0.45$.
It is worth mentioning some prominent co-spectral relationship at selected frequencies.
At frequency $f=0.03$, the following pairs display strong co-dependence:
\begin{inparaenum}[(i)]
\item Consumer Loans at All Commercial Banks
and Total Consumer Credit (Owned and Securitized);
\item 3-Year Treasury Constant Maturity Rate and 10-Year Treasury Constant Maturity Rate;
\item Compensation Per Hour in the Manufacturing Sector and in the Business Sector.
\end{inparaenum}


At frequency $f=0.45$, we observe a strong relationship for the following variable pairs:
\begin{inparaenum}[(i)]
\item Consumer Price Index for All Urban Consumers of All Items Less Energy and
of All Items Less Food \& Energy;
\item Real Imports and Exports of Goods \& Services Per-Capita;
\item Real Exports of Goods \& Services and Gross Private Domestic Investment Price Index;
\item Unit Labor Cost in the Business Sector and Output Per Hour of All Persons in the Business Sector.
\end{inparaenum}

\section{Conclusions}\label{concl}

In this paper, we consistently estimate the spectral density matrix under the assumption of a dynamic factor model with a sparse idiosyncratic multivariate spectrum for the data, where the relative pervasiveness of the common and the idiosyncratic components are allowed to vary across frequencies.
We prove that the nuclear norm plus $l_1$ norm heuristics consistently recovers across frequencies the spectral components and their sum, as well as the dynamic rank and the residual sparsity pattern.
We call the resulting estimators UNALSE (UNshrunk ALgebraic Spectral Estimator).

The empirical implications of the UNALSE approach are discussed on a US macroeconomic dataset, showing that UNALSE is able to catch the driving variables of the latent dynamics as well as the particular strength of specific relationships at all frequencies. This opens up the way to enhanced dynamic factor scores estimation and temporal network analysis.

\appendix
\small

\newpage
\setcounter{figure}{0}
\section{Proofs}
\paragraph{Proof of Lemma \ref{cons}}\label{proof_lemma}

First, due to the uncorrelation between $u_t$ and $e_t$, we can decompose the expected value
of the $(ij)-$ entry of the raw periodogram $\wh{\Sigma}^{X}_{Raw,ij}(\theta_h)$ at
each frequency $\theta_h=\pi \frac{h}{M_T}$, $h\leq |M_T|$, as follows:
$$
\E(\wh{\Sigma}^{X}_{Raw,ij}(\theta_h))=\E(d_i^{X}(\theta_h) \overline{d_j^{X}(\theta_h)})=\E(d_i^{\chi}(\theta_h)
\overline{d_j^{\chi}(\theta_h)})+\E(d_i^{\epsilon}(\theta_h)\overline{d_j^{\epsilon}(\theta_h)}),
$$
where $d_i^{X}(\theta_h)=\frac{1}{\sqrt{T}}\sum_{t=1}^T X_{i,t} \mathrm{e}^{-\mathrm{i}\theta_h t}$,
$d_i^{\chi}(\theta_h)=\frac{1}{\sqrt{T}}\sum_{t=1}^T \chi_{i,t} \mathrm{e}^{-\mathrm{i}\theta_h t}$,
$d_i^{\epsilon}(\theta_h)=\frac{1}{\sqrt{T}}\sum_{t=1}^T \epsilon_{i,t} \mathrm{e}^{-\mathrm{i}\theta_h t}$ are the Fourier transforms at frequency $\theta_h$ of $X_t$, $\chi_t$ and $\epsilon_t$ respectively.

At each frequency $\theta_h=\pi \frac{h}{M_T}$, $h\leq |M_T|$,
we define the common component of $\wh{\Sigma}^{X}_{ij,Raw}(\theta_h)$ as
$\wh{\Sigma}^{\chi}_{ij,Raw}(\theta_h)= d_i^{\chi}(\theta_h)\overline{d_j^{\chi}(\theta_h)}$,
and the idiosyncratic component of $\wh{\Sigma}^{X}_{ij,Raw}(\theta_h)$ as
$\wh{\Sigma}^{\epsilon}_{ij,Raw}(\theta_h)= d_i^{\epsilon}(\theta_h)\overline{d^{\epsilon}_j(\theta_h)}$.

Then, we derive the product $T \wh{\Sigma}^{\chi}_{ij,Raw}(\theta_h)$ as
\begin{align}
T d_i^{\chi}(\theta_h)\overline{d_j^{\chi}(\theta_h)}&=\l(\sum_{t=1}^T \sum_{s=0}^{\infty} \sum_{k=1}^r b_{ik,s} u_{k,t-s}
\mathrm{e}^{-\mathrm{i}\theta_h t}\r)\l(\sum_{t=1}^T \sum_{s=0}^{\infty} \sum_{k=1}^r b_{jk,s} u_{k,t-s} \mathrm{e}^{\mathrm{i}\theta_h t}\r)\nn\\
&=\sum_{s=0}^{\infty} \sum_{k=1}^r b_{ik,s} \sum_{s'=0}^{\infty} \sum_{k'=1}^r b_{jk',s'}
\sum_{t=1}^T \sum_{t'=1}^T u_{k,t-s} u_{k',t'-s'} \mathrm{e}^{-\mathrm{i}\theta_h (t-t')},\nn
\end{align}
and the product $T \wh{\Sigma}^{\epsilon}_{ij,Raw}(\theta_h)$ as
\begin{align}
T d_i^{\epsilon}(\theta_h)\overline{d_j^{\epsilon}(\theta_h)}&=\l(\sum_{t=1}^T \sum_{s=0}^{\infty} \sum_{k=1}^r c_{ik,s} e_{k,t-s}
\mathrm{e}^{-\mathrm{i}\theta_h t}\r)\l(\sum_{t=1}^T \sum_{s=0}^{\infty} \sum_{k=1}^r c_{jk,s} e_{k,t-s} \mathrm{e}^{\mathrm{i}\theta_h t}\r)\nn\\
&=\sum_{s=0}^{\infty} \sum_{k=1}^r c_{ik,s} \sum_{s'=0}^{\infty} \sum_{k'=1}^r c_{jk',s'}
\sum_{t=1}^T \sum_{t'=1}^T e_{k,t-s} e_{k',t'-s'} \mathrm{e}^{-\mathrm{i}\theta_h (t-t')}.\nn
\end{align}

Similarly, at each $\theta_h$, $h\leq |M_T|$,
we define the common component of the input spectral density matrix
$\wt{\Sigma}(\theta_h)$ (see (\ref{kernel})) as
\begin{align}
\wt{\Sigma}_{ij}(\theta_h)^{\chi}&=\frac{1}{2\pi} \sum_{h=-M_T}^{M_T} K^{*}(\theta_h) \wh{\Sigma}_{ij,Raw}^{\chi}(\theta_h)\nn\\
&=\frac{1}{2\pi} \sum_{h=-M_T}^{M_T} K^{*}(\theta_h) \frac{1}{T} \sum_{s=0}^{\infty} \sum_{k=1}^r b_{ik,s} \sum_{s'=0}^{\infty} \sum_{k'=1}^r b_{jk',s'} \sum_{t=1}^T \sum_{t'=1}^T  u_{i,t-s} u_{j,t'-s'} \mathrm{e}^{-\mathrm{i}\theta_h (t-t')},\nn
\end{align}
and the idiosyncratic component of the input $\wt{\Sigma}(\theta_h)$ as
\begin{align}
\wt{\Sigma}_{ij}(\theta_h)^{\epsilon}&=\frac{1}{2\pi}\sum_{h=-M_T}^{M_T} K^{*}(\theta_h) \wh{\Sigma}_{ij,Raw}^{\epsilon}(\theta_h)\nn\\
&=\frac{1}{2\pi}\sum_{h=-M_T}^{M_T} K^{*}(\theta_h) \frac{1}{T} \sum_{s=0}^{\infty} \sum_{k=1}^p c_{ik,s} \sum_{s'=0}^{\infty} \sum_{k'=1}^p c_{jk',s'} \sum_{t=1}^T \sum_{t'=1}^T  e_{i,t-s} e_{j,t'-s'} \mathrm{e}^{-\mathrm{i}\theta_h (t-t')},\nn
\end{align}
where $K^{*}(\theta_h)$ is the spectral window associated to the lag window $K\left(\frac{k}{M_T}\right)$ satisfying the kernel assumption (see e.g. \cite{priestley1981spectral}).

We can now apply the framework of \cite{wu2018asymptotic},
disentangling the common and the idiosyncratic components of the input estimator.
Let $P_{tt',ij}^{\chi}=\sum_{s=0}^{\infty} \sum_{k=1}^r b_{ik,s} u_{i,t-s} \sum_{s'=0}^{\infty} \sum_{k'=1}^r b_{jk',s'} u_{j,t'-s'}$ and $P_{tt',ij}^{\epsilon}=\sum_{s=0}^{\infty} \sum_{k=1}^p c_{ik,s} e_{i,t-s} \sum_{s'=0}^{\infty} \sum_{k'=1}^p c_{jk',s'} e_{j,t'-s'}$.

Then
\begin{align}
&\max_{1 \leq i,j \leq p} \E \max_{h\leq |M_T|} \vert \wt{\Sigma}_{ij}(\theta_h)-{\Sigma}_{ij}(\theta_h)\vert ^2 \label{sum_both}\\
&=\max_{1 \leq i,j \leq p} \E \max_{h\leq |M_T|}\frac{1}{2\pi T}\Bigg\vert \sum_{h=-M_T}^{M_T} K^{*}(\theta_h) \sum_{t=1}^T \sum_{t'=1}^T (P_{tt',ij}^{\chi}+P_{tt',ij}^{\epsilon}) \mathrm{e}^{-\mathrm{i}\theta_h (t-t')}-({\sigma}_{ij}(\theta_h)^{\chi}+{\sigma}_{ij}(\theta_h)^{\epsilon})\Bigg \vert^2 \nn\\
&=\max_{1 \leq i,j \leq p} \E \max_{h\leq |M_T|}\frac{1}{2\pi T} \Bigg \vert\left(\sum_{h=-M_T}^{M_T} K^{*}(\theta_h) \sum_{t=1}^T \sum_{t'=1}^T P_{tt',ij}^{\chi}\mathrm{e}^{-\mathrm{i}\theta_h (t-t')}-{\sigma}_{ij}(\theta_h)^{\chi}\right)\nn\\
&+\max_{1 \leq i,j \leq p} \E \max_{h\leq |M_T|}\left(\sum_{h=-M_T}^{M_T} K^{*}(\theta_h) \sum_{t=1}^T \sum_{t'=1}^T P_{tt',ij}^{\epsilon}\mathrm{e}^{-\mathrm{i}\theta_h (t-t')}-{\sigma}_{ij}(\theta_h)^{\epsilon}\right)\Bigg \vert^2\nn\\
&\leq \max_{1 \leq i,j \leq p} \E \max_{h\leq |M_T|} \frac{1}{2\pi T}
\Bigg \vert \sum_{h=-M_T}^{M_T} K^{*}(\theta_h) \sum_{t=1}^T \sum_{t'=1}^T P_{tt',ij}^{\chi}\mathrm{e}^{-\mathrm{i}\theta_h (t-t')}-{\sigma}_{ij}(\theta_h)^{\chi}\Bigg \vert^2\nn\\
&+\Bigg \vert \sum_{h=-M_T}^{M_T} K^{*}(\theta_h) \sum_{t=1}^T \sum_{t'=1}^T P_{tt',ij}^{\epsilon}\mathrm{e}^{-\mathrm{i}\theta_h (t-t')}-{\sigma}_{ij}(\theta_h)^{\epsilon}
\Bigg \vert^2\nn
\end{align}
which in turn is equal to

\begin{align}
&=\max_{1 \leq i,j \leq p} \E \max_{h\leq |M_T|} \Bigg \vert \sum_{s=0}^{\infty} \sum_{k=1}^r b_{ik,s} \sum_{s'=0}^{\infty} \sum_{k'=1}^r b_{jk',s'} \nn \\
&\left[\frac{1}{2\pi T} \sum_{h=-M_T}^{M_T} K^{*}(\theta_h)
\sum_{t=1}^T \sum_{t'=1}^T u_{k,t-s} u_{k',t'-s'}
\mathrm{e}^{-\mathrm{i}\theta_h (t-t')}-\mathrm{e}^{-\mathrm{i}\theta_h (s-s')}\right]\Bigg \vert^2 \label{sum_com}\\
&+\max_{1 \leq i,j \leq p} \E \max_{h\leq |M_T|} \Bigg \vert \sum_{s=0}^{\infty} \sum_{k=1}^p c_{ik,s} \sum_{s'=0}^{\infty} \sum_{k'=1}^p c_{jk',s'} \nn \\
&\left[\frac{1}{2\pi T} \sum_{h=-M_T}^{M_T} K^{*}(\theta_h)
\sum_{t=1}^T \sum_{t'=1}^T e_{k,t-s} e_{k',t'-s'}
\mathrm{e}^{-\mathrm{i}\theta_h (t-t')}-\mathrm{e}^{-\mathrm{i}\theta_h(s-s')}\right]\Bigg \vert^2 \label{sum_idio}
\end{align}

Under Assumption \ref{ass3} and the kernel assumption,
due to Theorem 1 in \cite[Section 4.2]{wu2018asymptotic},
which states that
$$
\max_{1 \leq i,j \leq p} \E \max_{h\leq |M_T|}
\Bigg \vert \frac{1}{2\pi T}\sum_{h=-M_T}^{M_T} K^{*}(\theta_h)
\sum_{t=1}^T \sum_{t'=1}^T u_{k,t-s} u_{k',t'-s'}
\mathrm{e}^{-\mathrm{i}\theta_h (t-t')}-\mathrm{e}^{-\mathrm{i}\theta_h (s-s')} \Bigg \vert^2 =O\left(\frac{M_T\log(M_T)}{T}\right),
$$
for the common component of (\ref{sum_both}),
which is (\ref{sum_com}), it holds
\begin{align}
&\max_{1 \leq i,j \leq p} \E \max_{h\leq |M_T|}(\Bigg \vert\sum_{s=0}^{\infty} \sum_{k=1}^r b_{ik,s} \sum_{s'=0}^{\infty} \sum_{k'=1}^r b_{jk',s'}\nn\\
&\left[\frac{1}{2\pi T} \sum_{h=-M_T}^{M_T} K^{*}(\theta_h)
\sum_{t=1}^T \sum_{t'=1}^T u_{k,t-s} u_{k',t'-s'}
\mathrm{e}^{-\mathrm{i}\theta_h (t-t')}-\mathrm{e}^{-\mathrm{i}\theta_h (s-s')}\right]\Bigg\vert^2\nn\\
&\leq \Bigg\Vert\sum_{s} B_{s}\Bigg\Vert_{\infty,v}^2\; \Bigg\Vert\sum_{s} B_{s}\Bigg\Vert_{\infty,v}^2 \frac{M_T\log(M_T)}{T}=O\l(\frac{M_T\log(M_T)}{T}\r).
\end{align}

Similarly, due to the same reasons, for the idiosyncratic component of (\ref{sum_both}),
which is (\ref{sum_idio}), it holds
\begin{align}
&\max_{1 \leq i,j \leq p} \E \max_{h\leq |M_T|}
\Bigg \vert\sum_{s=0}^{\infty} \sum_{k=1}^r c_{ik,s} \sum_{s'=0}^{\infty} \sum_{k'=1}^r c_{jk',s'} \nn \\
&\left[\frac{1}{2\pi T} \sum_{h=-M_T}^{M_T} K^{*}(\theta_h)
\sum_{t=1}^T \sum_{t'=1}^T e_{k,t-s} e_{k',t'-s'}
\mathrm{e}^{-\mathrm{i}\theta_h (t-t')}-\mathrm{e}^{-\mathrm{i}\theta_h (s-s')}\right]\Bigg\vert^2\nn\\
&\leq \Bigg\Vert\sum_{s} C_{s}\Bigg\Vert_{\infty,v}^2\; \Bigg\Vert\sum_{s} C_{s}\Bigg\Vert_{\infty,v}^2 \frac{M_T\log(M_T)}{T}=O\l(\frac{M_T\log(M_T)}{T}\r).
\end{align}

Let us now consider the expected maximum overall Frobenius loss of the input across frequencies:
\begin{align}\label{sum_spec}
&\E \max_{h\leq |M_T|} \Vert \wh{\Sigma}(\theta_h)-\Sigma(\theta_h)\Vert_{F}=\sqrt{\sum_{i=1}^p \sum_{j=1}^p (\wh{\Sigma}_{ij}(\theta_h)-\Sigma_{ij}(\theta_h))^2}\nn\\
&\leq \E \max_{h\leq |M_T|} \sqrt{\sum_{i=1}^p \sum_{j=1}^p (\wh{\Sigma}_{ij}^{\chi}(\theta_h)-\Sigma_{ij}^{\chi}(\theta_h))^2}+
\E \max_{h\leq |M_T|} \sqrt{\sum_{i=1}^p \sum_{j=1}^p (\wh{\Sigma}_{ij}^{\epsilon}(\theta_h)-\Sigma_{ij}^{\epsilon})^2}.
\end{align}
For the first term in the rhs of inequality \eqref{sum_spec} it holds
\begin{align}
&\E \max_{h \leq |M_T|} \Bigg \vert \sum_{i=1}^p \sum_{s=0}^{\infty} \sum_{k=1}^r b_{ik,s} \sum_{j=1}^p \sum_{s'=0}^{\infty}\sum_{k'=1}^r b_{jk',s'}\nn\\
&\left[\frac{1}{2\pi T}\sum_{h=-M_T}^{M_T} K^{*}(\theta_h) \sum_{t=1}^T \sum_{t'=1}^T u_{k,t-s} u_{k',t'-s'} \mathrm{e}^{-\mathrm{i}\theta_h (t-t')}-\mathrm{e}^{-\mathrm{i}\theta_h (s-s')}\right]\Bigg \vert^2\nn\\
&\le \Bigg \vert \sum_{i=1}^p \sum_{s=0}^{\infty} \sum_{k=1}^r  b_{ik,s} \sum_{j=1}^p \sum_{s'=0}^{\infty} \sum_{k'=1}^r b_{jk',s'}\Bigg \vert^2  G\left(\frac{M_T\log(M_T)}{T}\right)\nn\\
&\leq \Bigg \vert \sum_{i=1}^p \sum_{s=0}^{\infty} \sum_{k=1}^r b_{ik,s} \Bigg \vert ^2\; \Bigg \vert \sum_{j=1}^p \sum_{s'=0}^{\infty} \sum_{k'=1}^r b_{jk',s'}\Bigg \vert^2 G\left(\frac{M_T\log(M_T)}{T}\right)\nn\\
&\leq  {r^2} \Bigg \vert \Bigg \vert \sum_{s} B_{s} \Bigg \vert \Bigg \vert_{1,v}^2 \; \Bigg \vert \Bigg \vert \sum_{s} B_{s} \Bigg \vert \Bigg \vert_{1,v}^2 G\left(\frac{M_T\log(M_T)}{T}\right)\nn\\
& \leq G p^{2\alpha} \left(\frac{M_T\log(M_T)}{T}\right)=O\left(\frac{p^{2\alpha} M_T \log(M_T)}{T}\right)
\end{align}
because Assumptions \ref{ass3} and the kernel assumption hold.

Similarly, due to the same reasons, for the second term on the rhs of equation \eqref{sum_spec} it holds:
\begin{align}
&\E \max_{h \leq |M_T|} \Bigg \vert \sum_{i=1}^p \sum_{s=0}^{\infty} \sum_{k=1}^p c_{ik,s} \sum_{j=1}^p  \sum_{s'=0}^{\infty} \sum_{k'=1}^p c_{jk',s'}\nn\\
&\left[\frac{1}{2\pi T}\sum_{h=-M_T}^{M_T} K^{*}(\theta_h) \sum_{t=1}^T \sum_{t'=1}^T e_{k,t-s} e_{k',t'-s'}
\mathrm{e}^{-\mathrm{i}\theta_h (t-t')}-\mathrm{e}^{-\mathrm{i}\theta_h (s-s')}\right]\Bigg \vert^2\nn\\
&\le \Bigg \vert \sum_{i=1}^p \sum_{s=0}^{\infty} \sum_{k=1}^p c_{ik,s} \sum_{j=1}^p \sum_{s'=0}^{\infty} \sum_{k'=1}^p c_{jk',s'}\Bigg \vert^2  G\left(\frac{M_T\log(M_T)}{T}\right)\nn\\
&\leq \Bigg \vert \sum_{i=1}^p \sum_{s=0}^{\infty} \sum_{k=1}^p c_{ik,s} \Bigg \vert ^2\; \Bigg \vert \sum_{j=1}^p \sum_{s'=0}^{\infty} \sum_{k'=1}^p c_{jk',s'}\Bigg \vert^2 G\left(\frac{M_T\log(M_T)}{T}\right)\nn\\
&\leq \Bigg \vert \Bigg \vert \sum_{s} C_{s} \Bigg \vert \Bigg \vert_{1,v}^2 \; \Bigg \vert \Bigg \vert \sum_{s} C_{s} \Bigg \vert \Bigg \vert_{1,v}^2 G\left(\frac{M_T\log(M_T)}{T}\right)\nn\\
& \leq G_1 p^{2\delta'} \left(\frac{M_T\log(M_T)}{T}\right)=O\left(\frac{p^{2\delta'} M_T \log(M_T)}{T}\right).
\end{align}

%

Since $\delta' < \alpha$, from (\ref{sum_spec}) we can derive that
$\E \max_{h \leq |M_T|} \Vert \wh{\Sigma}(\theta_h)-\Sigma(\theta_h)\Vert_{F} \leq O\left(p^{\alpha}
\sqrt{\frac{M_T \log(M_T)}{T}}\right)$,
from which the thesis
$\E \max_{h \leq |M_T|} \Vert \wh{\Sigma}(\theta_h)-\Sigma(\theta_h)\Vert_{2} \leq O\left(p^{\alpha}
\sqrt{\frac{M_T \log(M_T)}{T}}\right)$ follows.

\paragraph{Proof of Theorem \ref{thmMineUNALSE}}\label{proof_cons_basic}

Following \cite{luo2011high},
we note that under Assumptions \ref{ass1} and \ref{ass2},
setting
$\psi=\frac{p^{\alpha}}{\xi(T)}\frac{1}{\sqrt{T}}$
with $\rho=\gamma \psi$ (where $\gamma \in [9\xi(T),1/(6\mu(\Omega))]$),
and further assuming that
$\underline{\delta}_T p^{2(\alpha-\underline{\delta})} < T < \overline{\delta}_T p^{6\delta}$ for some $\underline{\delta}_T,\overline{\delta}_T$ such that $0<\underline{\delta}_T<\overline{\delta}_T$
and that the minimum eigenvalue of $L^{*}$, $\lambda_r(L^{*})$, is larger than $G_2 \frac{\psi}{\xi^2(T)}$, 
Propositions 12, 13, and 14 in \cite{luo2011high} can be directly applied to our setting, proving that,
for each $\theta_h=\frac{h\pi}{M_T}$, $h \in -[M_T],\ldots,[M_T]$, the pair $(\wh{L}(\theta_h),\wh{S}(\theta_h))$ minimizing (\ref{func:ob_spec})
satisfies the following theses:
\ben
\item[i)] $g_\gamma(\wh{S}(\theta_h)-{S}(\theta_h),\wh{L}(\theta_h)-{S}(\theta_h))$ is upper bounded;
\item[ii)] $\wh{L}(\theta_h)$ is rank-consistent: $\text{\upshape rk}(\wh{L}(\theta_h))=\text{\upshape rk}({L}(\theta_h))$;
\een
with probability depending on the random loss $\Vert \wh{\Sigma}(\theta_h)-\Sigma(\theta_h) \Vert_2$.
In particular, the upper bound on $g_\gamma(\wh{S}(\theta_h)-{S}(\theta_h),\wh{L}(\theta_h)-{S}(\theta_h))$
depends on $\Vert \wh{\Sigma}(\theta_h)-\Sigma(\theta_h) \Vert_2$.
Considering the frequency grid $\theta_h=\pi \frac{h}{M_T}$, $h\leq |M_T|$,
Lemma \ref{cons} states that under Assumptions \ref{ass3} and the kernel assumption it holds
$\max_{h} \Vert \wh{\Sigma}(\theta_h)-\Sigma(\theta_h)\Vert_{2}\leq O\left(p^{\alpha} \sqrt{\frac{M_T \log(M_T)}{T}}\right)$, which leads, setting
$\psi=\frac{p^{\alpha}}{\xi(\mathcal T)}\sqrt{\frac{M_T\log M_T}{T}}$,
to
\begin{equation}
g_\gamma(\wh{S}(\theta_h)-{S}(\theta_h),\wh{L}(\theta_h)-{L}(\theta_h))\leq C \psi. \label{ggamma2}
\end{equation}
From (\ref{ggamma2}),
all the claims of Theorem \ref{thmMineUNALSE} follow, because
\bea
||\wh{L}(\theta_h)-{L}(\theta_h)||_2 &\leq& C \psi; \nn\\
||\wh{S}(\theta_h)-{S}(\theta_h)||_{\infty} &\leq&
C \gamma \psi \leq \xi(T) \psi; \nn\\
||\wh{S}(\theta_h)-{S}(\theta_h)||_{2} &\leq& C q'
\gamma \psi \leq q' \xi(T) \psi; \nn\\
||\wh{\Sigma}(\theta_h)-{\Sigma}(\theta_h)||_{2}
&\leq& C \psi +  q' \xi(T) \psi. \nn
\eea




Unlike Assumptions \ref{ass1}-\ref{ass2} and the lower bound on $\lambda_r(L^{*})$,
if the lower bound on $S_{\text{\tiny{min,off}}}$ does not hold, there is no consequence on the identification of the two underlying algebraic varieties, i.e. on parametric consistency
and rank recovery. 
The only consequence lies in the fact that some nonzero elements
of $\wh{S}(\theta_h)$ are not recovered.
This fact can be appreciated by looking at the proofs of Propositions 5.2 and 5.3
in \cite{chandrasekaran2012}, directly exploited by \cite{luo2011high}.

If, instead, the condition $\Vert S \Vert_{\text{\tiny{min,off}}}>G_3\frac{\psi}{\mu(\Omega)}$ holds in addition to all the assumptions and conditions of Theorem \ref{thmMineUNALSE},
the same Propositions in \cite{chandrasekaran2012}
allow to conclude that the recovered sparsity pattern is also consistent: $\text{\upshape sgn}(\mathrm{Re}(\wh{S}(\theta_h)))=\text{\upshape sgn}(\mathrm{Re}({S}(\theta_h)))$.

\paragraph{Proof of Theorem \ref{thmMineUNALSEgen}}\label{proof_cons_gen}
Under the assumptions of Theorem \ref{thmMineUNALSEgen}, the assumptions and conditions of Theorem \ref{thmMineUNALSE} are automatically satisfied. While the assumptions on kernel window and temporal dependence are exactly the same as in the basic filter setting of Assumption \ref{basic}, the assumptions on latent eigenvalues and residual sparsity pattern are reshaped to cope with the general filters prescribed by Assumption \ref{gen}. The assumption needed to ensure the identifiability of underlying algebraic varieties is also reshaped accordingly. Therefore, all the claims of Theorem \ref{thmMineUNALSE} are still valid, in a much more general context where the latent eigenvalues and the residual sparsity pattern are intermediately spiked, and the latent coefficient matrices may have different condition numbers across frequencies while the residual coefficient matrices may have different sparsity patterns. The GDFM setting is a special case into this context, where the latent eigenvalues are spiked with $p$ and the maximum row-wise number of residual nonzeros is bounded by a constant. Therefore, Theorem \ref{th:introTh} is proved as a special case of Theorem \ref{thmMineUNALSEgen} with $\alpha=1$.

\paragraph{Proof of Corollary \ref{pos_def}}\label{proof_inv}

Let us define $\phi=G (\frac{p^{\alpha}}{\xi(T)}\frac{1}{\sqrt{T}}+q'\frac{p^{\alpha}}{\sqrt{T}})\sqrt{\frac{M_T\log M_T}{T}}$
and $\phi_S=G (q'\frac{p^{\alpha}}{\sqrt{T}})\sqrt{\frac{M_T\log M_T}{T}}$.
Suppose that
$\wh{L}=\wh{L}_{\text{\tiny{UNALSE}}}(\theta_h)$, $\wh{S}=\wh{S}_{\text{\tiny{UNALSE}}}(\theta_h)$, $\wh{\Sigma}=\wh{\Sigma}_{\text{\tiny{UNALSE}}}(\theta_h)$, for each $\theta_h=\pi\frac{h}{M_T}$, $h\leq |M_T|$ and $\Sigma=\Sigma(\theta_h)$ for each $\theta \in [-\pi,\pi]$.
Weyl's Theorem prescribes that, for any matrix $\Sigma$, we have
$$\vert \wh{\lambda}_i-\lambda\vert  \leq \Vert \wh{\Sigma}-\Sigma\Vert_{2}, \;\forall i=1,\ldots,p,$$
where $\wh{\lambda}_i$, $i=1,\ldots,p$, are the sample eigenvalues.
This result relates the rate of sample eigenvalues to the matrix spectral loss rate.
The triangular inequality gives
\begin{equation}
\vert \lambda_p(\wh{L}+\wh{S})-\lambda_p(\Sigma)\vert
\leq \vert \lambda_p(\wh{L}+\wh{S})\vert  + \vert -\lambda_p(\Sigma)\vert
=\vert \lambda_p(\wh{L}+\wh{S})\vert  +\lambda_p(\Sigma),
\end{equation}
because $\Sigma$ is positive definite.
Thus, $$\vert \lambda_p(\wh{L}+\wh{S})\vert  \geq  \vert \lambda_p(\wh{L}+\wh{S})-\lambda_p(\Sigma)\vert  - \lambda_p(\Sigma).$$
Since for Weyl's theorem $\vert \lambda_p(\wh{L}+\wh{S})-\lambda_p(\Sigma)\vert  \leq \phi$
we have 
\begin{equation}\lambda_p(\wh{L}+\wh{S})>0 \Longleftrightarrow \lambda_p(\Sigma)> \phi.\end{equation}
This proves the first part of the claim.

In order to achieve the same rate $\phi$ for the inverse spectral rate $\Vert (\wh{L}+~\wh{S})^{-1}-~\Sigma^{-1}\Vert_{2}$,
it is necessary that $\lambda_p(\Sigma) \geq 2\phi$.
In fact, the triangular inequality gives
\begin{equation}\Vert (\wh{L}+\wh{S})^{-1}-\Sigma^{-1}\Vert_{2} \leq \Vert (\wh{L}+\wh{S})^{-1}\Vert_{2} + \lambda_p(\Sigma)^{-1}\label{base}\end{equation}

By summing and subtracting $\Sigma$ and using the triangular inequality
\begin{equation}\Vert (\wh{L}+\wh{S})^{-1}\Vert_{2}= \Vert (\wh{L}+\wh{S}-\Sigma+\Sigma)^{-1}\Vert_{2} \leq
\Vert (\wh{L}+\wh{S}-\Sigma)^{-1}\Vert_{2}+\Vert \Sigma^{-1}\Vert_{2}
\leq\Vert (\wh{L}+\wh{S}-\Sigma)^{-1}\Vert_{2}+\lambda_p(\Sigma)^{-1}.
\end{equation}


For the triangular inequality, we have
\begin{equation}
\vert \lambda_p((\wh{L}+\wh{S})^{-1})-\lambda_p(\Sigma)^{-1}\vert \leq \vert \lambda_p((\wh{L}+\wh{S})^{-1})\vert  + \vert -\lambda_p(\Sigma)^{-1}\vert \leq \vert \lambda_p((\wh{L}+\wh{S})^{-1})\vert  + \lambda_p(\Sigma)^{-1}
\end{equation}
since $\Sigma$ is positive definite.

At the same time, we want that 
$$\Vert (\wh{L}+\wh{S})^{-1}-\Sigma^{-1}\Vert_{2}\leq \phi.$$




Hence, inequality (\ref{base}) becomes
$$\phi^{-1} \leq \vert \lambda_p((\wh{L}+\wh{S})^{-1})\vert +2\lambda_p(\Sigma)^{-1}.$$
We can write $$\vert \lambda_p((\wh{L}+\wh{S})^{-1})\vert  \geq \phi^{-1}-2\lambda_p(\Sigma)^{-1},$$
which allows to conclude that \begin{equation}\Vert \wh{\Sigma}^{-1}-\Sigma^{-1}\Vert_{2}\leq \phi \Longleftrightarrow \phi ^{-1} \geq 2 \lambda_p(\Sigma)^{-1}.\label{inv}\end{equation}

Using (\ref{inv}),
it is possible to derive the rate for $(\wh{L}+\wh{S})^{-1}$.
By the property (see \cite{luo2011high}, pp. 31-32):
\begin{equation}
\Vert (M+N)^{-1}-M^{-1}\Vert_{2} \leq \Vert M^{-1}\Vert_{2} \cdot\Vert N\Vert_{2}\cdot\Vert (M+N)^{-1}\Vert_{2}
\end{equation}
we obtain
\begin{align}
& \Vert (\wh{L}+\wh{S})^{-1}-\Sigma^{-1}\Vert_{2}= \Vert (\wh{L}+\wh{S})^{-1}[\wh{L}+\wh{S}-\Sigma]\Sigma^{-1}\Vert_{2} \leq \nn \\
& \leq \Vert (\wh{L}+\wh{S})^{-1}\Vert_{2}\cdot \Vert [\wh{L}+\wh{S}-\Sigma]\Vert_{2}\cdot \Vert \Sigma^{-1}\Vert_{2}\leq \frac{2}{\lambda_p(\Sigma)^2}\Vert [\wh{L}+\wh{S}-\Sigma]\Vert_{2}. \nn
\end{align}

Hence, we have \begin{equation}\Vert \wh{\Sigma}^{-1}-\Sigma^{-1}\Vert_{2} \leq G(q'\xi(T)+1)\psi=\phi\end{equation}

The same reasoning can be carried out for $\wh{S}$ by simply replacing $\phi$ by $\phi_S$, thus obtaining
\begin{equation}\Vert \wh{S}^{-1}-S^{-1}\Vert_{2} \leq G(q'\xi(T))\psi=\phi_S.\end{equation}

\section{Admissible sparsity regimes}

\begin{Rem}[Admissible sparsity regimes]\upshape{
Finally, we analyze in detail the admissible sparsity regimes for the residual spectral density matrix $S(\theta)$ at a given frequency. The parameters involved in this analysis are the latent eigenvalues rate $\alpha$, the maximum number of nonzero elements per row, $q'$, the minimum absolute nonzero off-diagonal element $\Vert S(\theta)\Vert_{min,off}$, the row-wise maximum $l_1$ norm $\Vert S(\theta)\Vert_{1,v}$, and the sample size $T$.}

First of all, we note that the following inequality holds:
\beq q' \Vert S(\theta)\Vert_{min,off} \leq \Vert S(\theta)\Vert_{1,v} \label{sp1}.\eeq
At the same time, from Theorems 4.1 and 4.2 we know that sparsistency requires that
$\Vert S(\theta)\Vert_{min,off} >\frac{\psi}{\mu(\Omega)}$, from which
we can write
$\Vert S(\theta)\Vert_{min,off} > \frac{p^{\alpha}}{\sqrt{T}}\frac{1}{\xi(\mathcal{T})\mu(\Omega)}$,
that becomes \beq \Vert S(\theta)\Vert_{min,off} \sqrt{T} \gtrsim p^{\alpha} \label{sp2}\eeq
because $\xi(\mathcal{T})\mu(\Omega)=O(1)$.

From the assumptions of Theorems 4.1 and 4.2, we know that $q'\lesssim p^{\delta}$,
$\delta < \alpha$,
and that
$\underline{\delta}_T p^{2(\alpha-\underline{\delta})} < T < \overline{\delta}_T p^{6\delta}$ for some $\underline{\delta}_T,\overline{\delta}_T$ such that $0<\underline{\delta}_T<\overline{\delta}_T$.
For some $\delta_2>0$, we require $\Vert S(\theta)\Vert_{1,v} \leq \delta_2 p^{\underline{\delta}}$, with $\underline{\delta} \leq \delta+0.5$, $\underline{\delta}\leq\delta'$, $\underline{\delta} <\alpha$.
Letting $\Vert S(\theta)\Vert_{min,off}=O(p^{\vartheta})$ and $\sqrt{T}=O(p^{\iota})$,
by combining \eqref{sp1} and \eqref{sp2} we obtain
\begin{equation}p^{\alpha-\iota} < p^{\vartheta} \leq p^{\underline{\delta}-\delta}.
\label{sparsity_regimes}
\end{equation}


As $\iota \leq \delta$, sparsistency is not possible, because \eqref{sparsity_regimes} leads to $\underline{\delta} \geq \alpha$.
As $\iota > \delta$, i.e., as $T$ grows, sparsistency becomes possible.
If $\iota=\frac{3}{2}\delta$, sparsistency requires for instance that $\alpha < \underline{\delta}+\frac{1}{2} \delta$. This condition is compatible, among others, with the setting $\alpha=\frac{3}{4}$, $\underline{\delta}=\frac{2}{3}$, $\delta=\frac{1}{3}$, because by (\ref{sparsity_regimes}) $\iota=\frac{3}{2}\delta$ leads to $\alpha<\frac{5}{6}$.

In order to understand the admissible relative scalings of $\vartheta$, $\delta$, $\underline{\delta}$,
we need to study the sign of the quantity $\alpha-\iota$.
If for instance $\iota=\frac{3}{2}\delta$ and $\delta=\frac{1}{3}$, we can notice that $\alpha-\iota =\alpha-\frac{3}{2} \delta = \alpha - \frac{3}{2} \frac{1}{3} < 0$ if $\alpha < \frac{1}{2}$.
Therefore, from (\ref{sparsity_regimes}) we know that the condition $\underline{\delta}-\delta<0$ is admissible in that case. As an example, setting $\delta=\frac{1}{3}$ and $\alpha=\frac{3}{8}$, we obtain from (\ref{sparsity_regimes}) that $\frac{3}{8}-\frac{1}{2}<\underline{\delta}-\frac{1}{3}$, which means $\underline{\delta}>-\frac{1}{8}+\frac{1}{3}=\frac{5}{24}$, smaller than $\delta=\frac{1}{3}$ and $\alpha=\frac{3}{8}$.

To sum up, if for some $\varepsilon>0$ $T \gtrsim p^{\delta+\varepsilon}$, i.e., if $\iota>\delta$, sparsistency is possible. If $\alpha-\iota>0$, it is necessary that $\vartheta>0$, with $\underline{\delta}$ large and $\delta$ small to respect the condition $\alpha-\iota < \vartheta < \underline{\delta}-\delta$ with $\underline{\delta}<\alpha$ and $\delta<\alpha$.
This situation corresponds to a sparsity pattern with few large nonzero entries.

As $\iota$ increases, i.e. as the sample size $T$ grows, it is more likely that $\alpha-\iota<0$. In that case, $\vartheta$ may be smaller than $0$, and it may also hold $\underline{\delta}<\delta$. This corresponds to a sparsity pattern with a larger number of smaller nonzero entries (than before).

Note that the former is usually a case with large $\alpha$, while the latter requires a smaller $\alpha$. This implies that spiked latent eigenvalues require large residual entries to achieve sparsistency, while smaller latent eigenvalues require small residual entries to ensure latent rank recovery.

\end{Rem}

\section{Simulation mechanism}\label{subsec:sim}


Our purpose is to obtain data with a low rank plus sparse spectral density matrix.
First, we fix the dimension $p$, the sample size $T$, the latent rank $r$ and the condition number of $L^{*}$, $c$. Concerning basic filters, the first problem we encounter regards the generation of matrices $U_L$ and $\Lambda_u$. We apply the generation algorithm in \cite{farne2016algorithm} to simulate a latent multivariate spectrum with fixed condition number $c$ and trace equal to $\tau \beta p$, where $\beta$ is the latent variance proportion (constant across frequencies) and $\tau$ is a scale parameter. The procedure consists in applying the Gram-Schmidt algorithm to a permutation matrix, randomly drawing $r$ random eigenvectors (to become the orthonormal columns of the matrix $U_L$), and then building $\Lambda_u$ as a diagonal $r\times r$ matrix with ordered equidistant elements such that ${\Lambda_u}_{1,1}/{\Lambda_u}_{r,r}=c$. We thus obtain $L^{*}=U_L \Lambda_u U_L'$.
For given time coefficients $\lambda_s$ 
such that $\sum_{s=0}^{n_{l}}\lambda_s^2=1$
(where $n_{l}$ is the chosen number of lags),
we obtain the coefficient matrices accordingly as $B_{s}= U_L \sqrt{\Lambda_u}\lambda_s$,
which allow us to generate the common component $\chi_t$, $t=1,\ldots,T$, as a vector moving average (VMA) with $n_{l}$ lags.

Following Assumption 3.2, we generate the sparse component in the following way. First, we generate the diagonal of $S^{*}$ from a Dirichlet distribution with parameter $(1-\beta) \tau p$. Then, we re-order these residual variances matching the respective magnitude order of the variances in $L^{*}$. Exploiting Cauchy-Schwartz inequality, we then randomly generate the off-diagonal elements $S_{ij}^{*}$ from a uniform distribution $Unif(0,\delta \sqrt{S_{ii}^{*} S_{jj}^{*}})$, where
$\delta$ is a tuning parameter. The generated residual off-diagonal elements are then ordered, and the survival threshold is set to their maximum times a proportion parameter $\delta_{bis}$. The residual coefficient matrices are then obtained accordingly as $C_{s}= U_S \sqrt{\Lambda_e}\lambda_s$, where $U_S \Lambda_e U_S'$ is the spectral decomposition of $S^{*}$, allowing to generate the idiosyncratic component $\epsilon_t$, $t=1,\ldots,T$, as a VMA as well.

Concerning general filters, we follow Assumption 3.3. We start from the basic filters $B_{s}= U_L \sqrt{\Lambda_u}\lambda_s$. We then replace the scalar $\lambda_s$ by a $r \times r$ diagonal matrix $D_{L,s}$, built as follows. We set a perturbation proportion $\kappa_{pert}=0.1$. We generate the diagonal of $D_{L,s}$ as $\lambda_s\mathbbm{1}_r -\kappa_{pert}\lambda_s\mathbbm{1}_r+2\kappa_{pert}\lambda_s W_r$, where $\mathbbm{1}_r$ is a vector composed by $r$ ones, and $W_r$ is a vector of $r$ random numbers drawn from a uniform between $0$ and $1$. The general filters of the low rank component are then obtained as $B_{s}=U_L D_{L,s} \sqrt{\Lambda_u}$. This procedure generates a varying latent variance proportion $\beta(\theta_h)$ across chosen frequencies.

Concerning the residual filters, our procedure generates the diagonals of $\Gamma_{\epsilon}(s)$, $s=0,1, \ldots, n_{l}$, from a Dirichlet distribution with parameter $(1-\beta) \tau p \lambda_s$. We then order the elements of each of the $n_{l}+1$ diagonals matching the respective magnitude order of the variances in $\Gamma_{\chi}(0)$.
Starting from those diagonals, we repeat the above exposed thresholding procedure $n_{l}+1$ times. For each of the subsequent sparsified $\Gamma_{\epsilon}(s)$, we derive the spectral decomposition $\Gamma_{\epsilon}(s)=U_{S,s}\Lambda^2_{S,s} U_{S,s}'$, and we set the generalized filters as $C_{s}=U_{S,s} \Lambda_{S,s}$. Note that this passage is very delicate, as the matrices of eigenvectors may lead to non-sparse residual spectra if the nonzeros are too many, due to rounding errors.
Anyway, as prescribed in Assumption 3.3, each generated $\Gamma_{\epsilon}(s)$ has in principle
a different sparsity pattern and a different number of nonzeros $q_s$
due to the Cauchy-Schwartz inequality.

The spectral density matrices at each frequency $\theta_h$ are then computed.
We set $|h|\leq M_T$, $\theta_h=\frac{h \pi}{M_T}$,
and we then calculate the low rank and residual transfer function matrices as $B(\theta_h)=\sum_{s=0}^{n_{l}} B_s \mathrm{e}^{(-i s \theta_h)}$ and $C(\theta_h)=\sum_{s=0}^{n_{l}} C_s \mathrm{e}^{(-i s \theta_h)}$.
The latent and residual spectral density matrices are thus derived at each frequency $\theta_h$ as $L(\theta_h)=B(\theta_h)\Lambda_u B(\theta_h)'$ and $S(\theta_h)=C(\theta_h)\Lambda_e C(\theta_h)'$.
The basic filter specifications lead to real spectra. On the contrary, the general filters lead to complex spectra, as the residual spectral components are complex.


Once we have generated the common component $\chi_t$ and the residual component $\epsilon_t$ as in Section \ref{subsec:sim}, we can generate our series $x_t$, $t=1,\ldots,T$, according to equation (5).
A relevant choice to control the spectral shape lies in the vector moving average coefficients $\lambda_s$, $s=0,\ldots,n_{l}$. We observe that the spectral shapes across frequencies are characterized by the following pattern:
\begin{itemize}
\item $n_{l}=1$, a positive $\lambda_0$ and negative $\lambda_1$ lead to a "reverse S-shape";
\item $n_{l}=1$, a negative $\lambda_0$ and positive $\lambda_1$ lead to a " S-shape";
\item $n_{l}=2$, a positive $\lambda_0$, $\lambda_1=0$ and a negative $\lambda_2$ lead to a "U-shape";
\item $n_{l}=2$, a negative $\lambda_0$, $\lambda_1=0$ and a positive $\lambda_2$ lead to a "reverse U-shape";
\end{itemize}
We note that when the non-null coefficients are equal to $0.5$, the spectral shape has the maximum variability possible across frequencies, i.e., the spectral peak attains its maximum and the spectral drop is $0$. Any situation with any of the coefficients equal to $1$ leads instead to a constant spectrum across frequencies.
In the simulation study of Section 6,
we set $\lambda_0=0.8$ and $\lambda_1=0.2$.
We select as our target the reverse S-shape without loss of generality,
because the conditions of Theorems 4.1 and 4.2 must be satisfied frequency-wise,
and the reverse S-shape across frequencies is by far the most common in real time series.

\clearpage
\section{Additional simulation results}\label{subsec:sim}
\subsection{Scenario A}
\begin{figure}[h!]
          \caption{Estimated latent variance proportion $\wh{\beta}(f_h)$ - Scenario A.}\label{fig:beta_A.1.U}
          \centering
          \begin{tabular}{ccc}
          {\footnotesize Setting 1}&{\footnotesize Setting 2}&{\footnotesize Setting 3}\\
              \includegraphics[width=.2\textwidth]{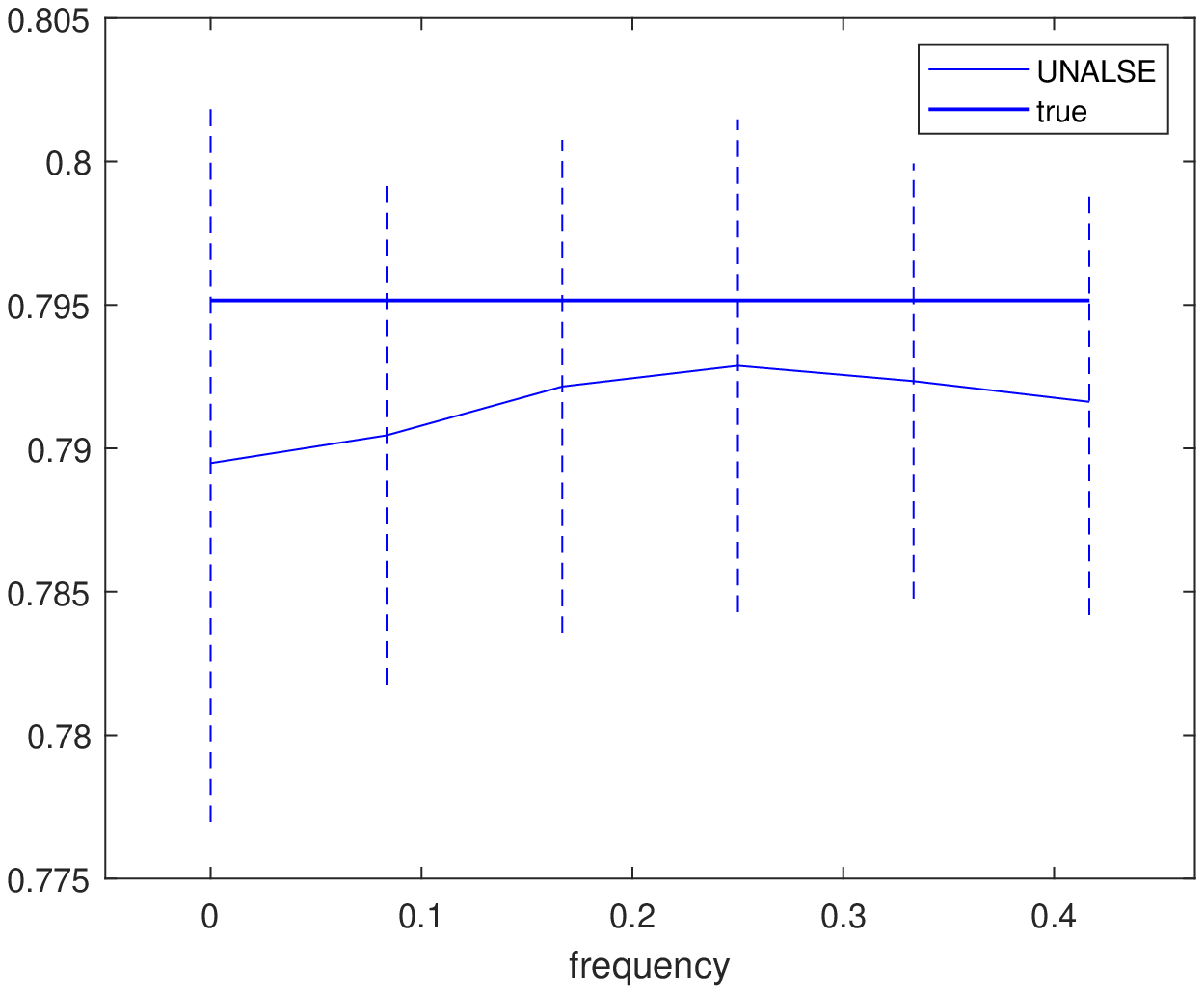}&
              \includegraphics[width=.2\textwidth]{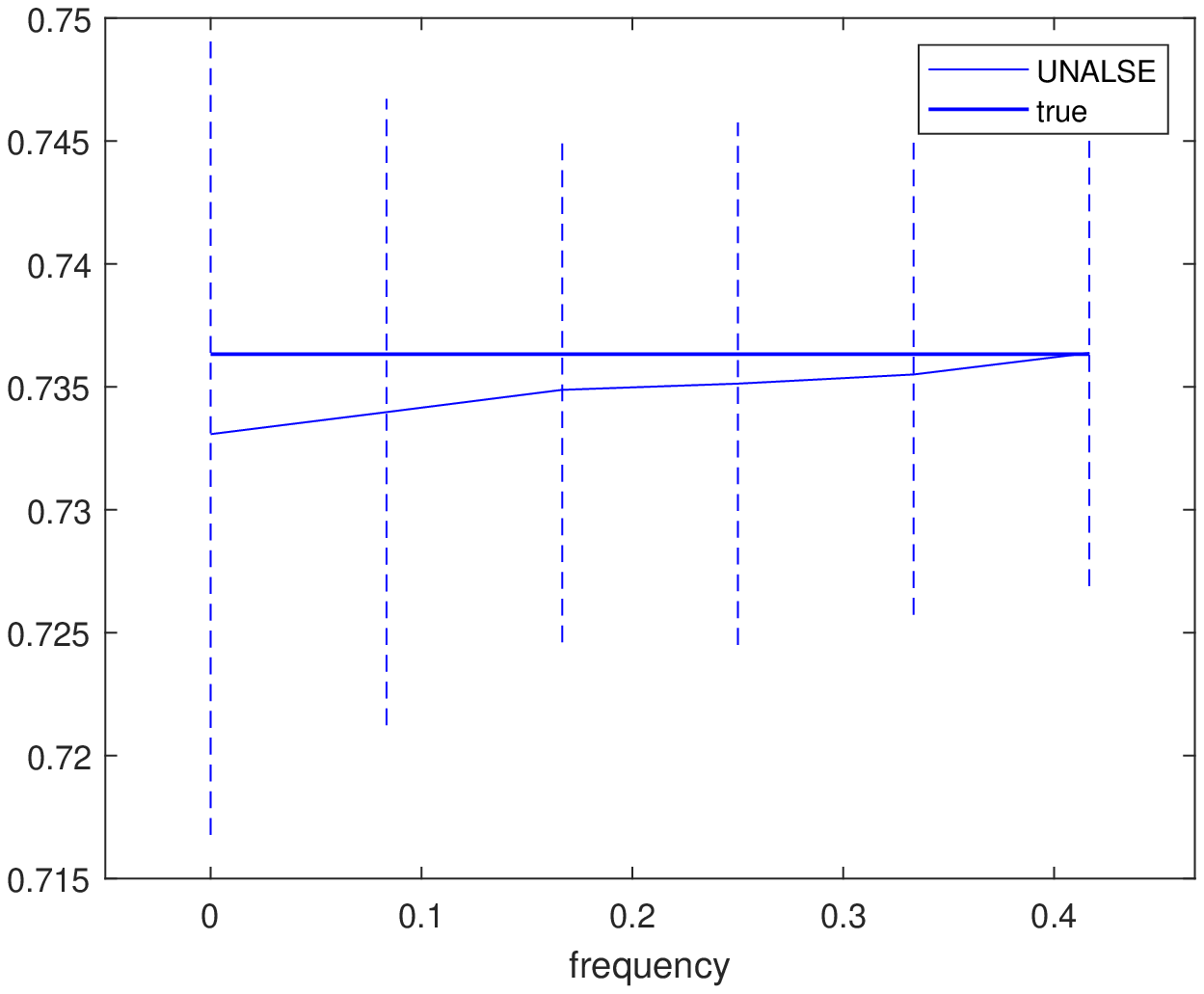}&
              \includegraphics[width=.2\textwidth]{theta_base_UNALSE_3}\\
              \end{tabular}
               \begin{tabular}{cc}
                             {\footnotesize Setting 4}&{\footnotesize Setting 5}\\
              \includegraphics[width=.2\textwidth]{theta_base_UNALSE_4}&
              \includegraphics[width=.2\textwidth]{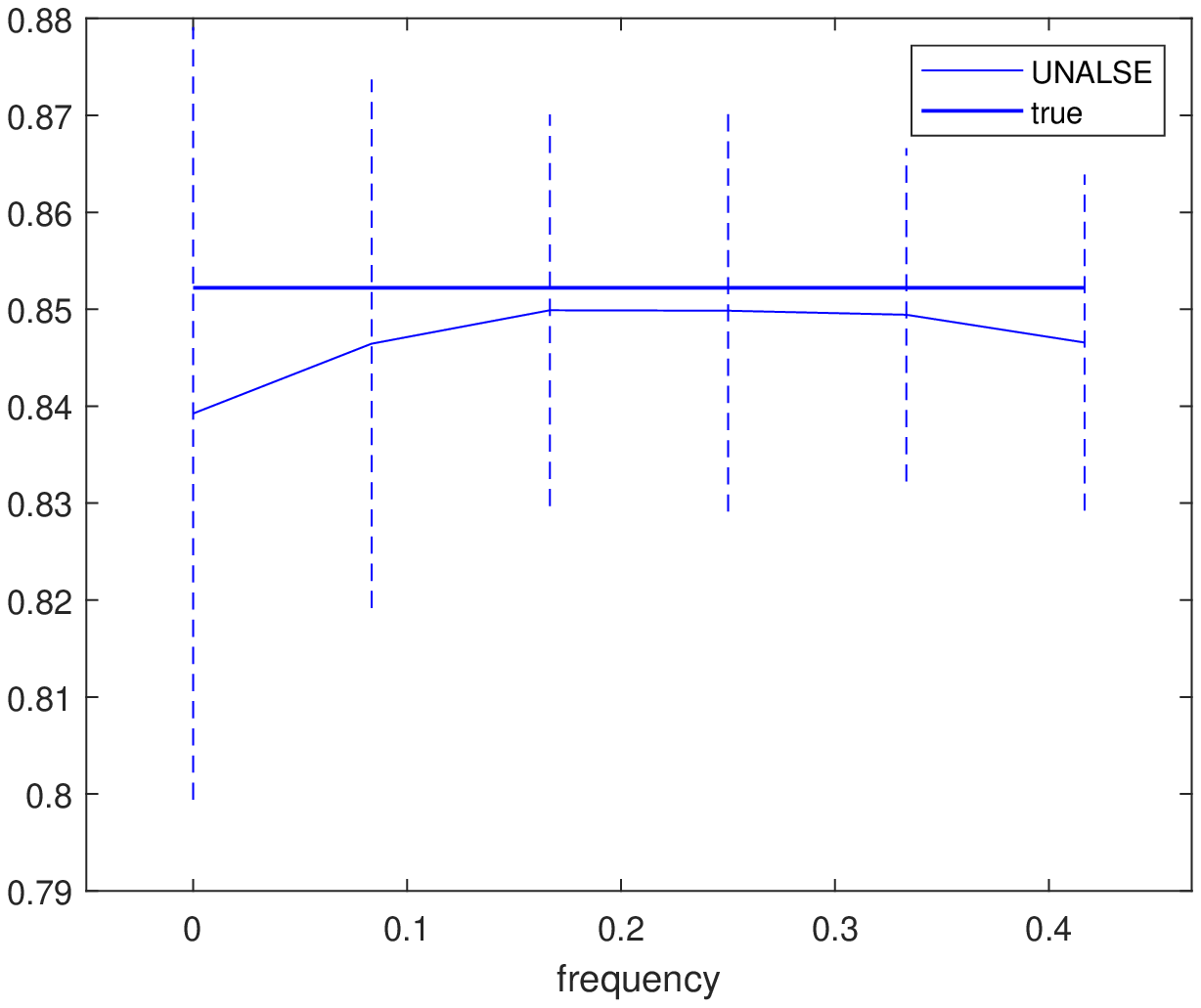}\\
              \end{tabular}
          \end{figure}

          \begin{figure}[h!]
          \caption{Positive predictive value $ppv(f_h)$ - Scenario A.}\label{fig:pred_pos_A.1.U}
          \centering
          \begin{tabular}{ccc}
          {\footnotesize Setting 1}&{\footnotesize Setting 2}&{\footnotesize Setting 3}\\
              \includegraphics[width=.2\textwidth]{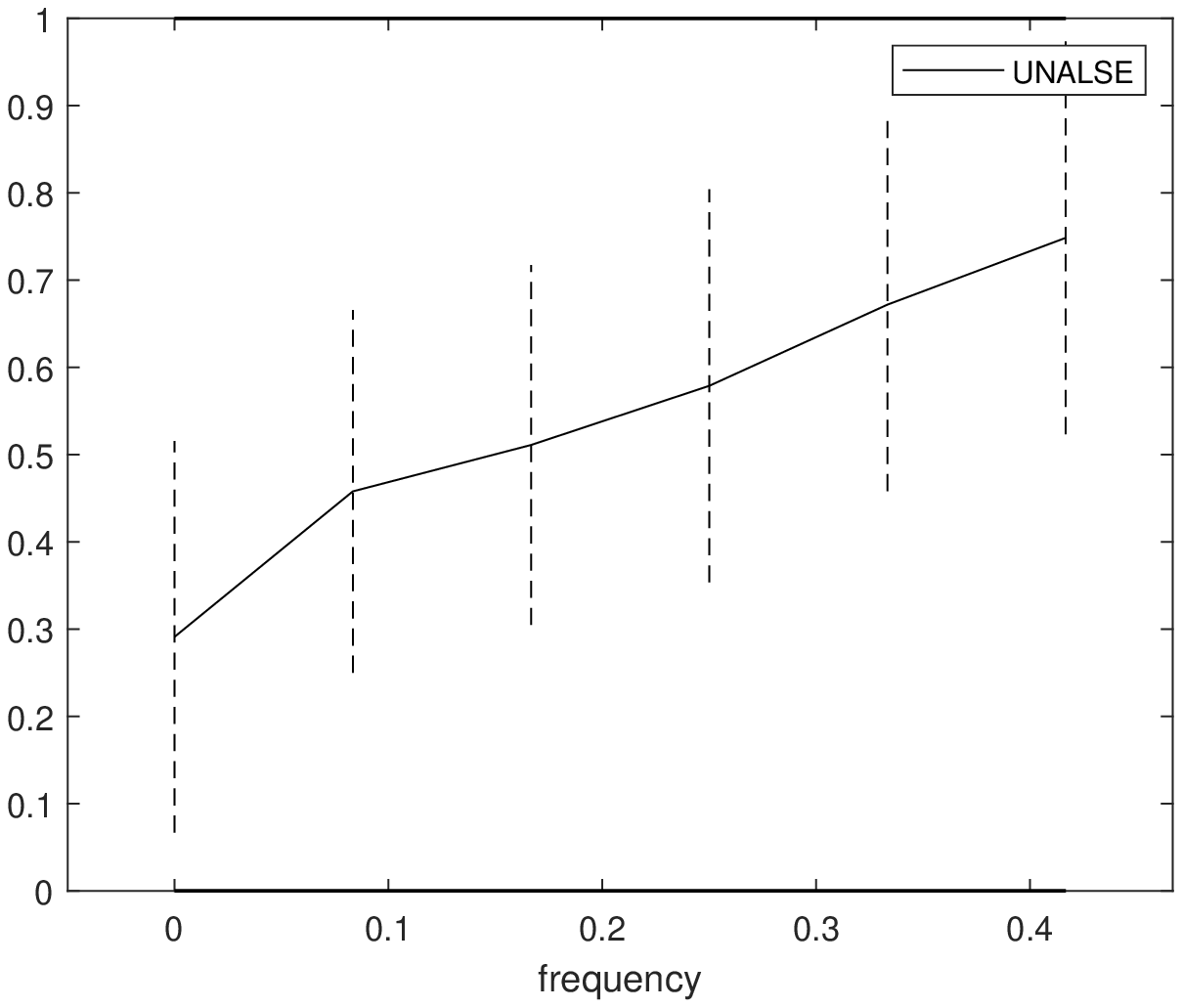}&
              \includegraphics[width=.2\textwidth]{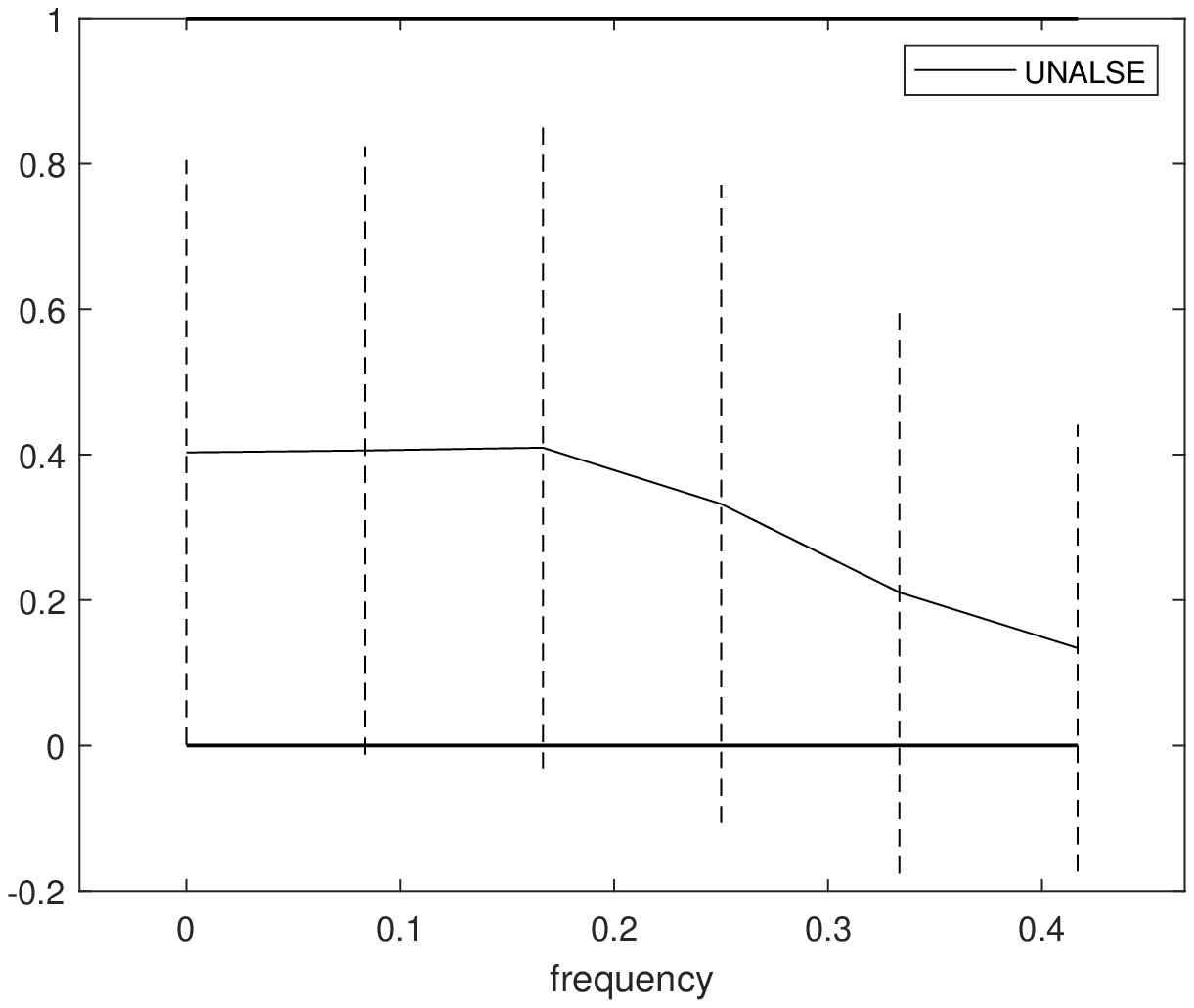}&
              \includegraphics[width=.2\textwidth]{pred_pos_rate_base_UNALSE_3}\\
              \end{tabular}
               \begin{tabular}{cc}
	      {\footnotesize Setting 4}&{\footnotesize Setting 5}\\
              \includegraphics[width=.2\textwidth]{pred_pos_rate_base_UNALSE_4}&
              \includegraphics[width=.2\textwidth]{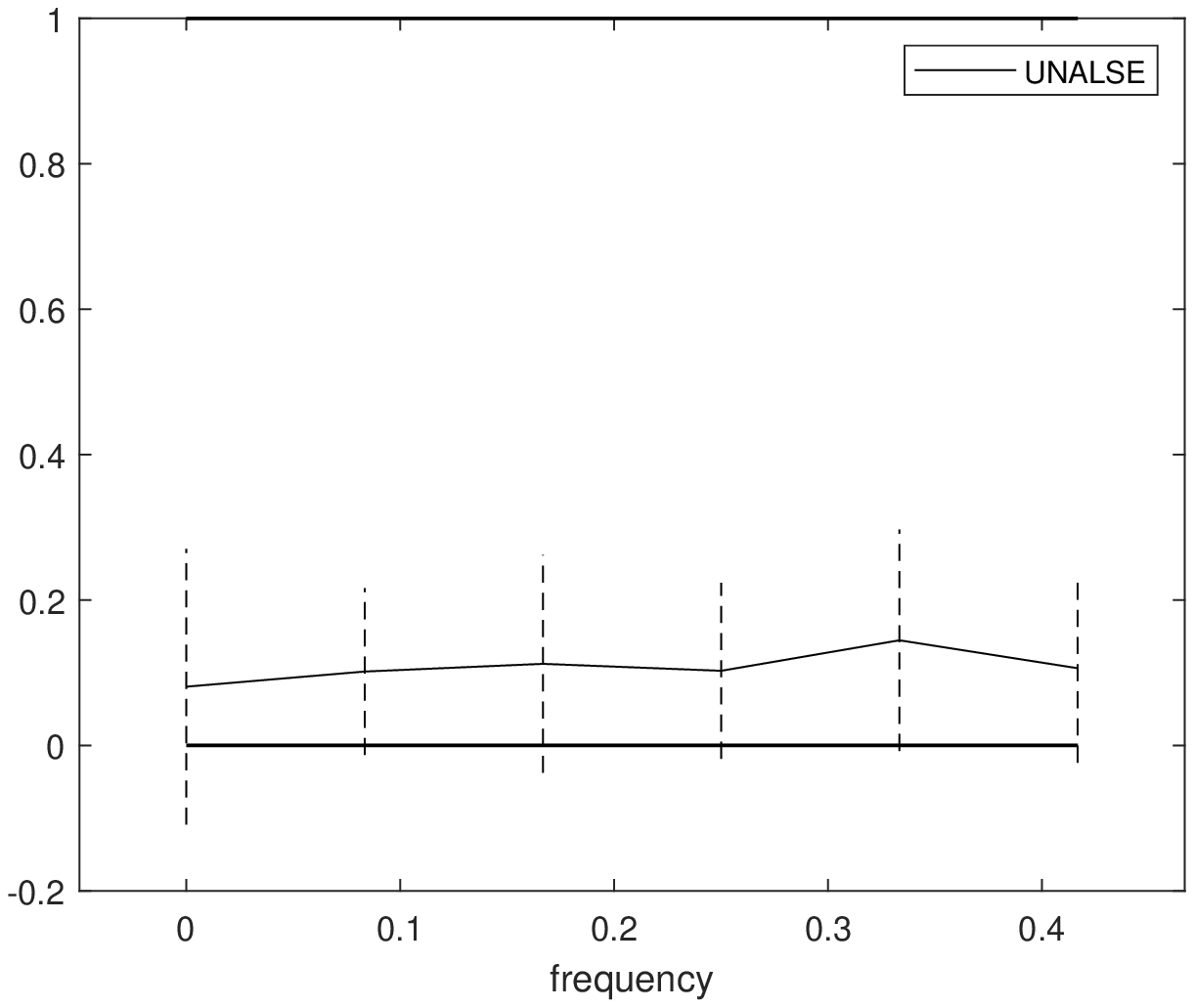}\\
              \end{tabular}
          \end{figure}

 \begin{figure}[h!]
          \caption{Negative predictive value $npv(f_h)$ - Scenario A.}\label{fig:pred_neg_A.1.U}
          \centering
          \begin{tabular}{ccc}
          {\footnotesize Setting 1}&{\footnotesize Setting 2}&{\footnotesize Setting 3}\\
              \includegraphics[width=.2\textwidth]{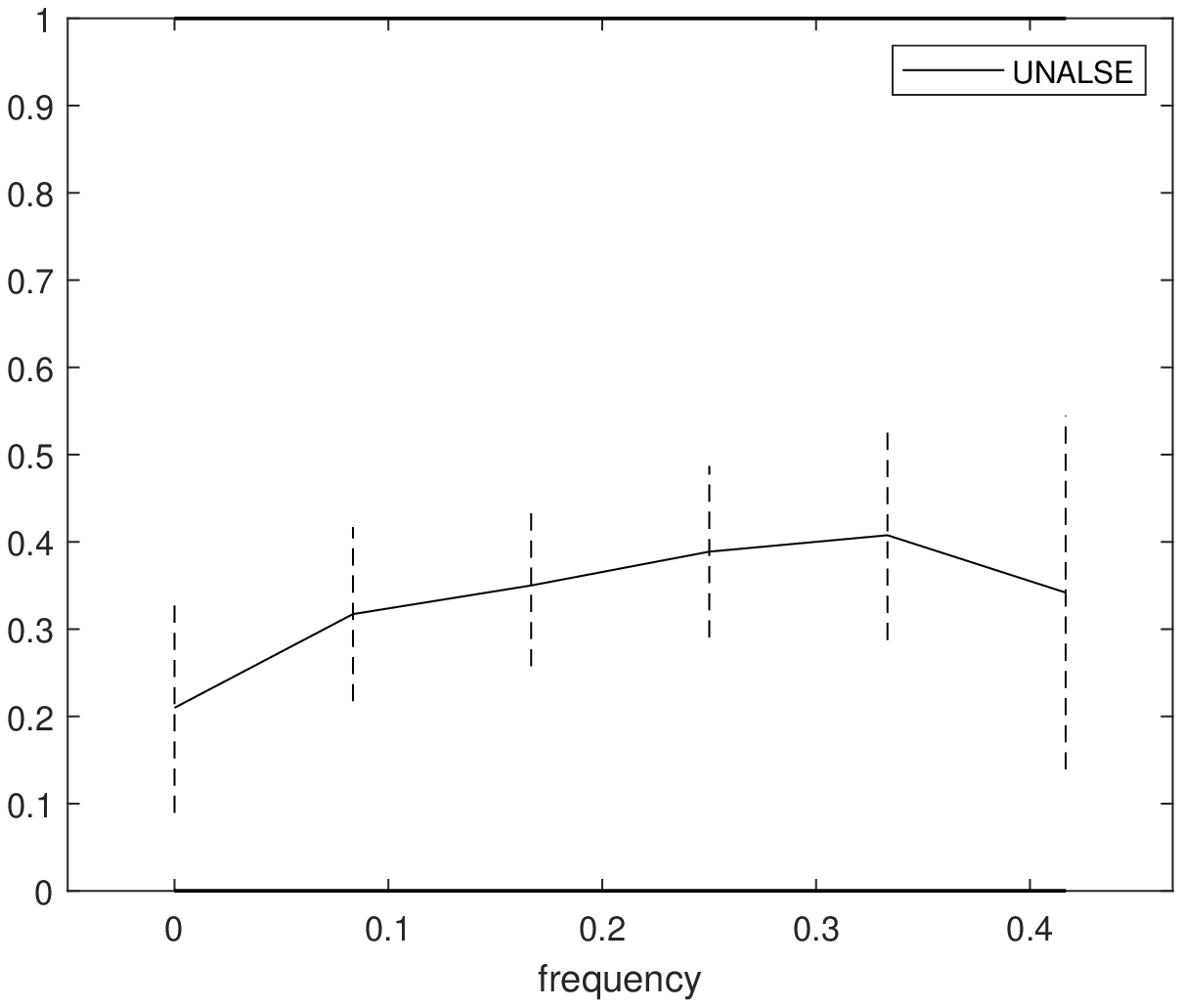}&
              \includegraphics[width=.2\textwidth]{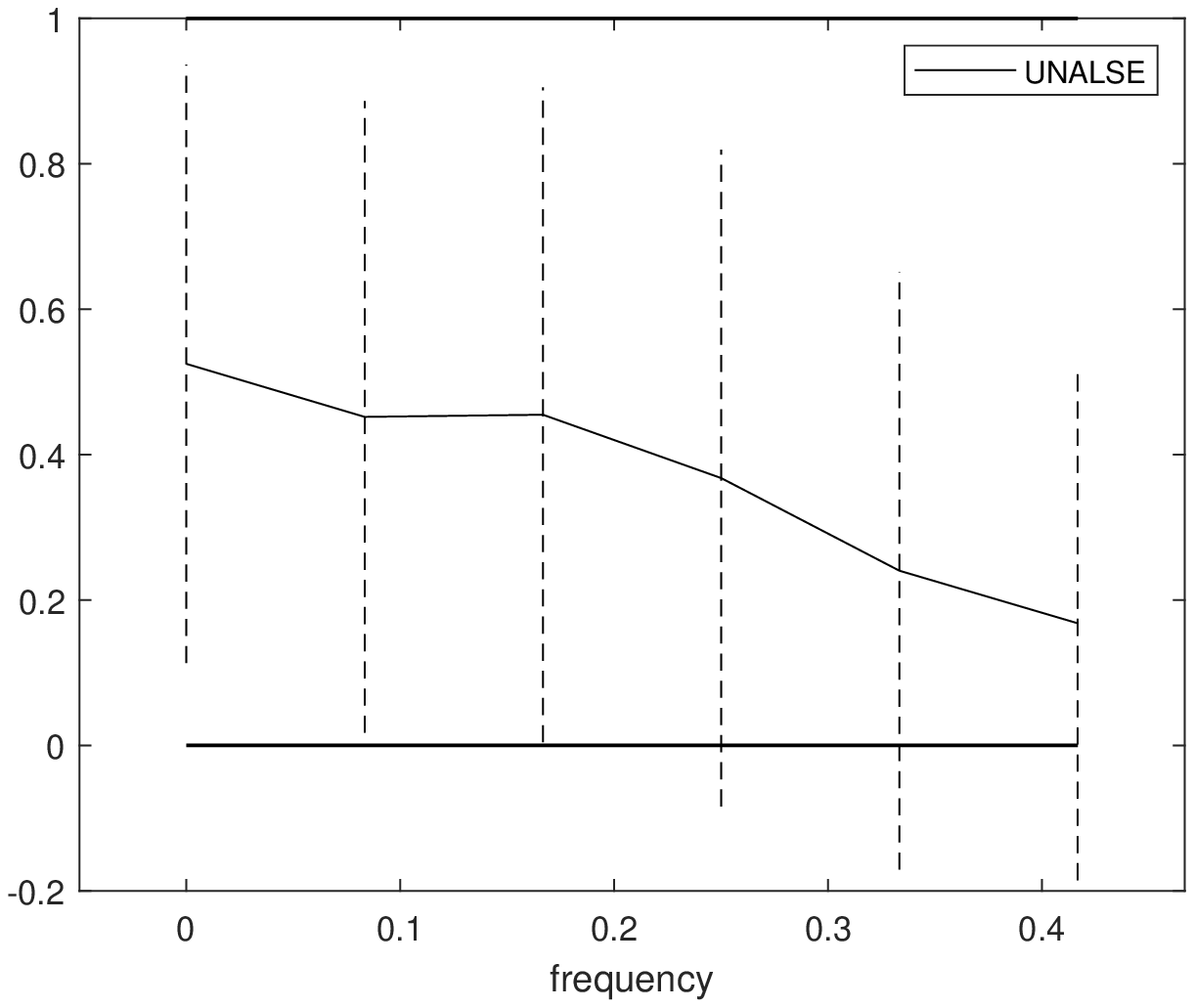}&
              \includegraphics[width=.2\textwidth]{pred_neg_rate_base_UNALSE_3}\\
              \end{tabular}
               \begin{tabular}{cc}
	      {\footnotesize Setting 4}&{\footnotesize Setting 5}\\
              \includegraphics[width=.2\textwidth]{pred_neg_rate_base_UNALSE_4}&
              \includegraphics[width=.2\textwidth]{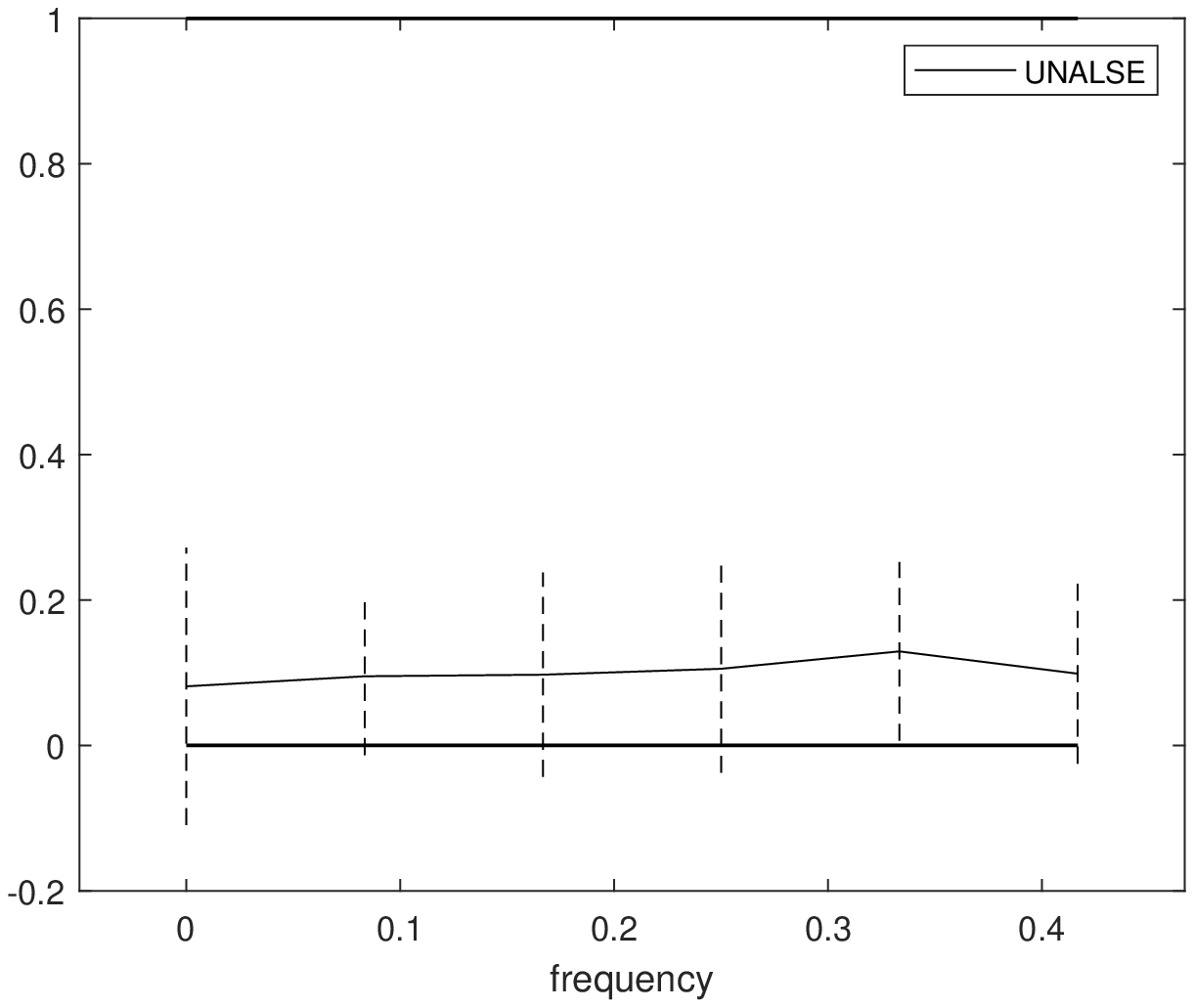}\\
              \end{tabular}
          \end{figure}

          \begin{figure}[h!]
          \caption{$err_{\wh{L}}(f_h)$ - Scenario A.}\label{fig:err_L_A.1.U}
          \centering
          \begin{tabular}{ccc}
          {\footnotesize Setting 1}&{\footnotesize Setting 2}&{\footnotesize Setting 3}\\
              \includegraphics[width=.2\textwidth]{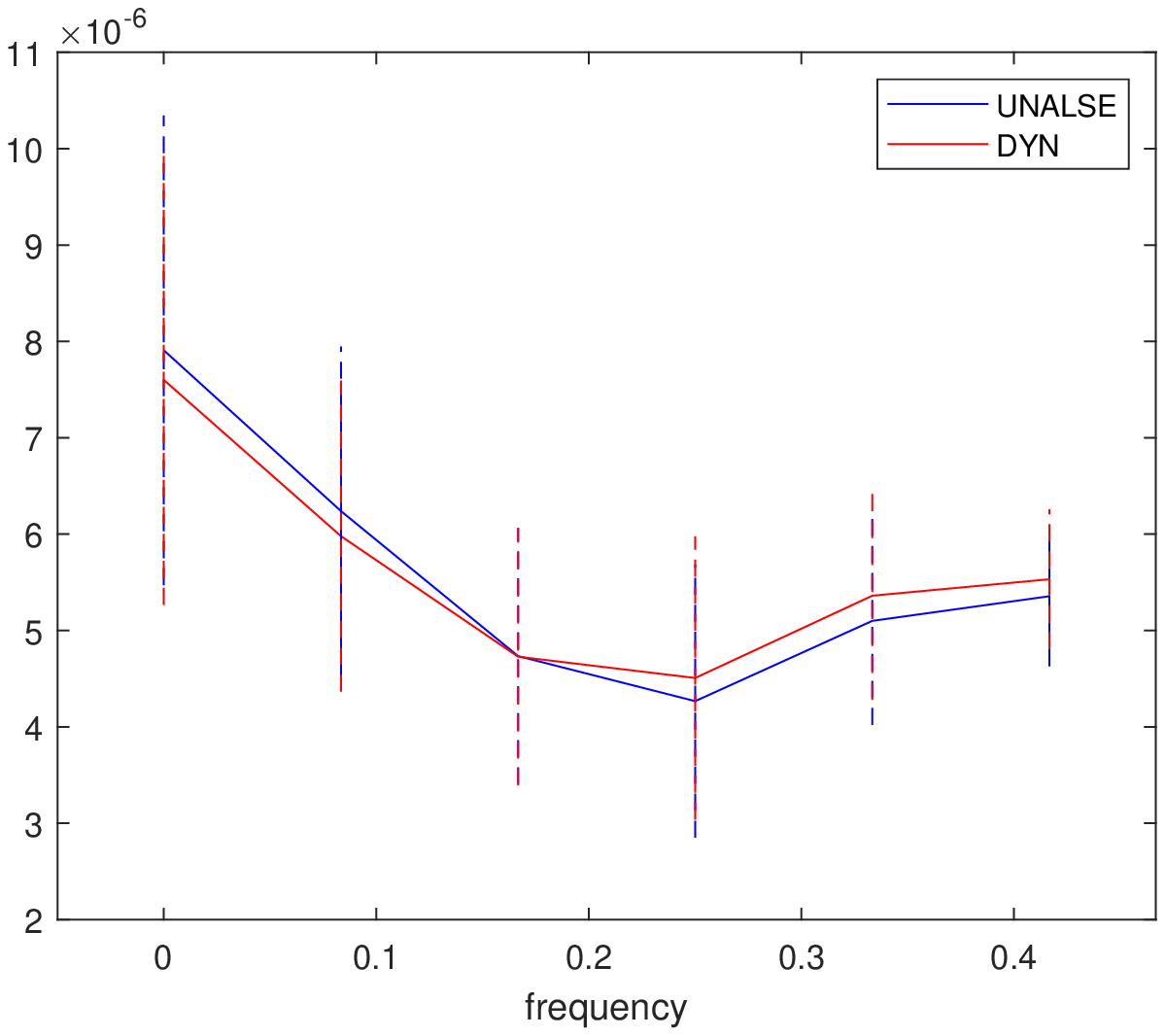}&
              \includegraphics[width=.2\textwidth]{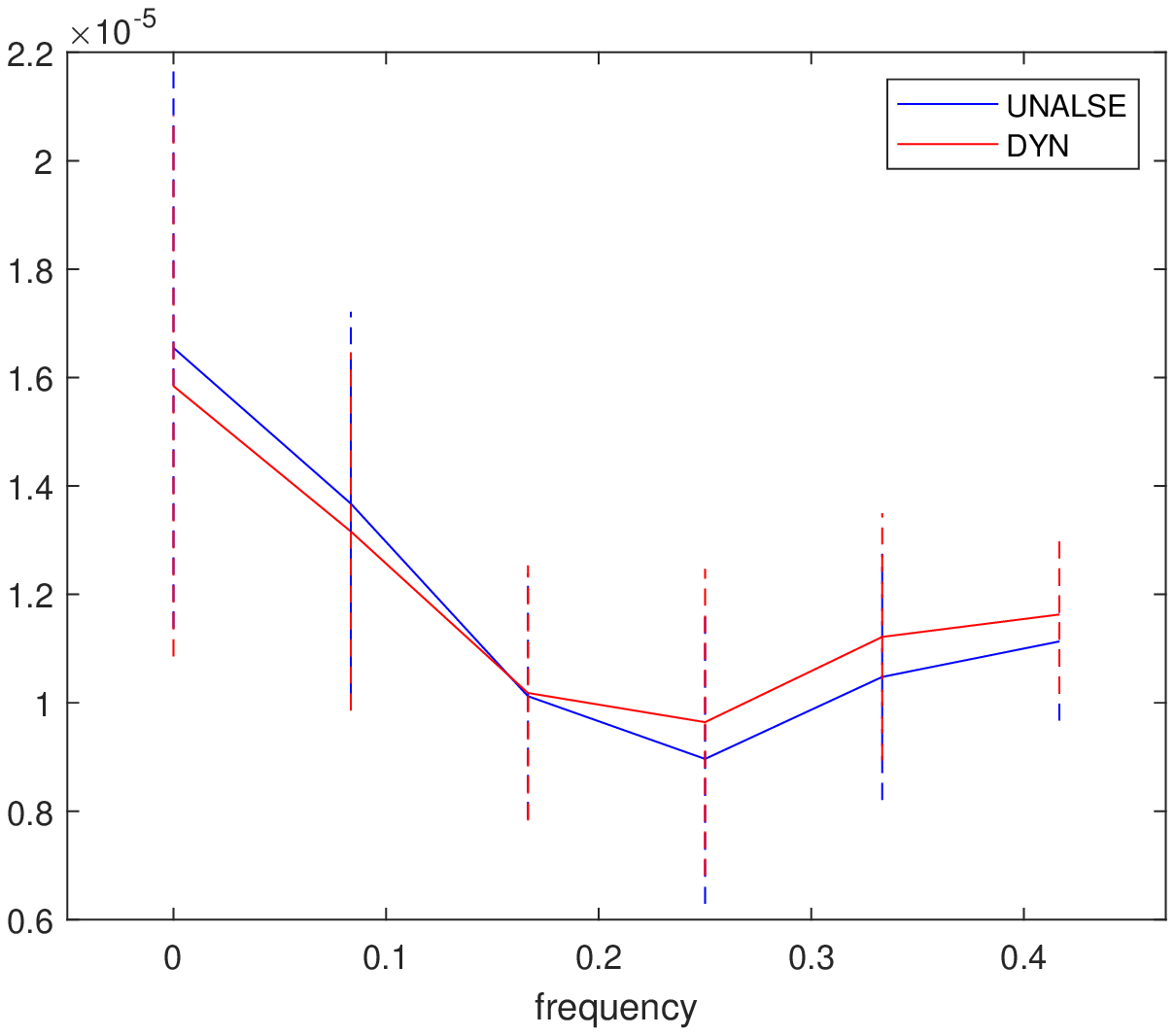}&
              \includegraphics[width=.2\textwidth]{err_L_base_3}\\
              \end{tabular}
               \begin{tabular}{cc}
	      {\footnotesize Setting 4}&{\footnotesize Setting 5}\\
              \includegraphics[width=.2\textwidth]{err_L_base_4}&
              \includegraphics[width=.2\textwidth]{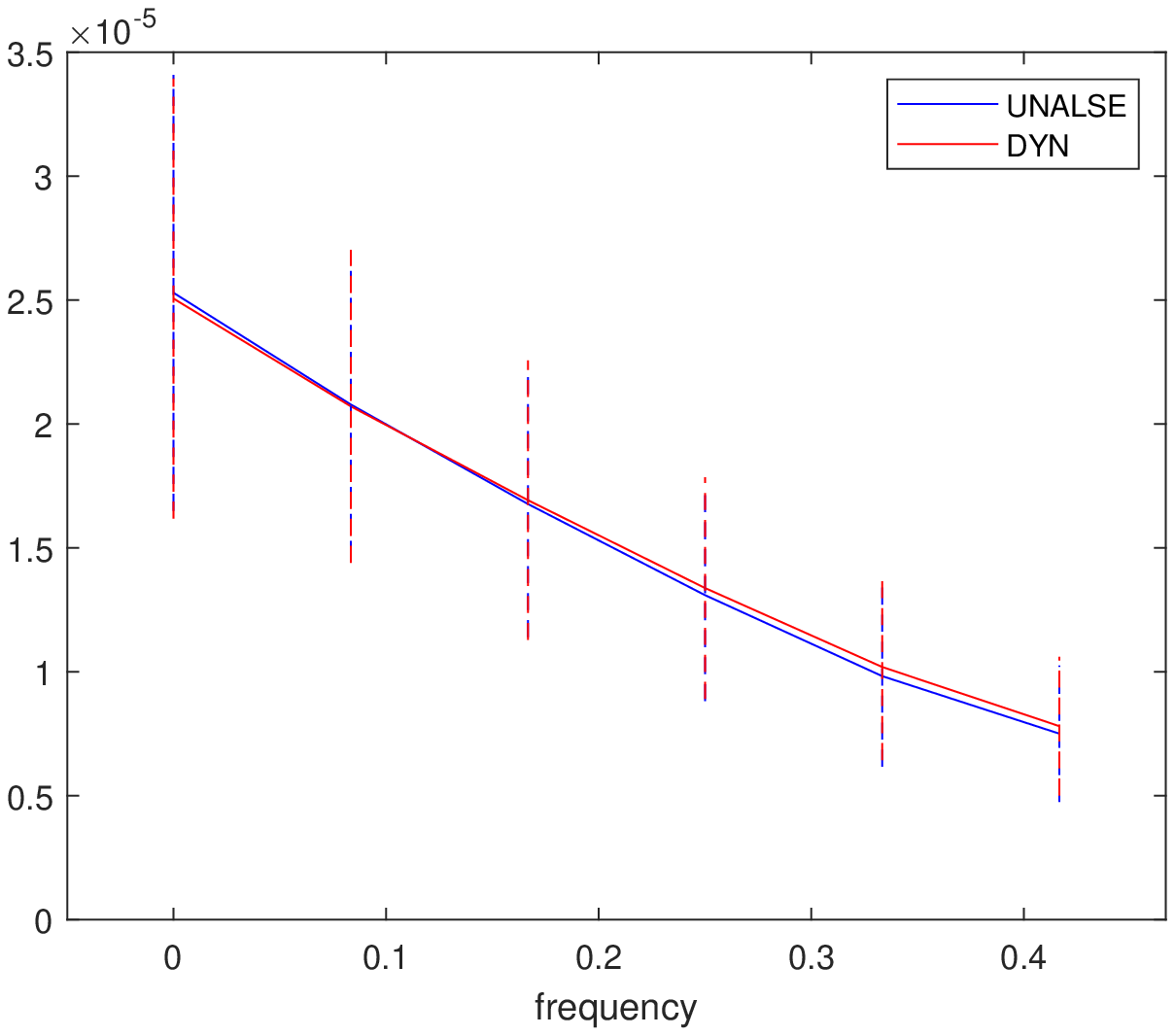}\\
              \end{tabular}
          \end{figure}

\begin{figure}[h!]
          \caption{$err_{ratio}(f_h)$ - Scenario A.}\label{fig:err_ratio_A.1.U}
          \centering
          \begin{tabular}{ccc}
          {\footnotesize Setting 1}&{\footnotesize Setting 2}&{\footnotesize Setting 3}\\
              \includegraphics[width=.2\textwidth]{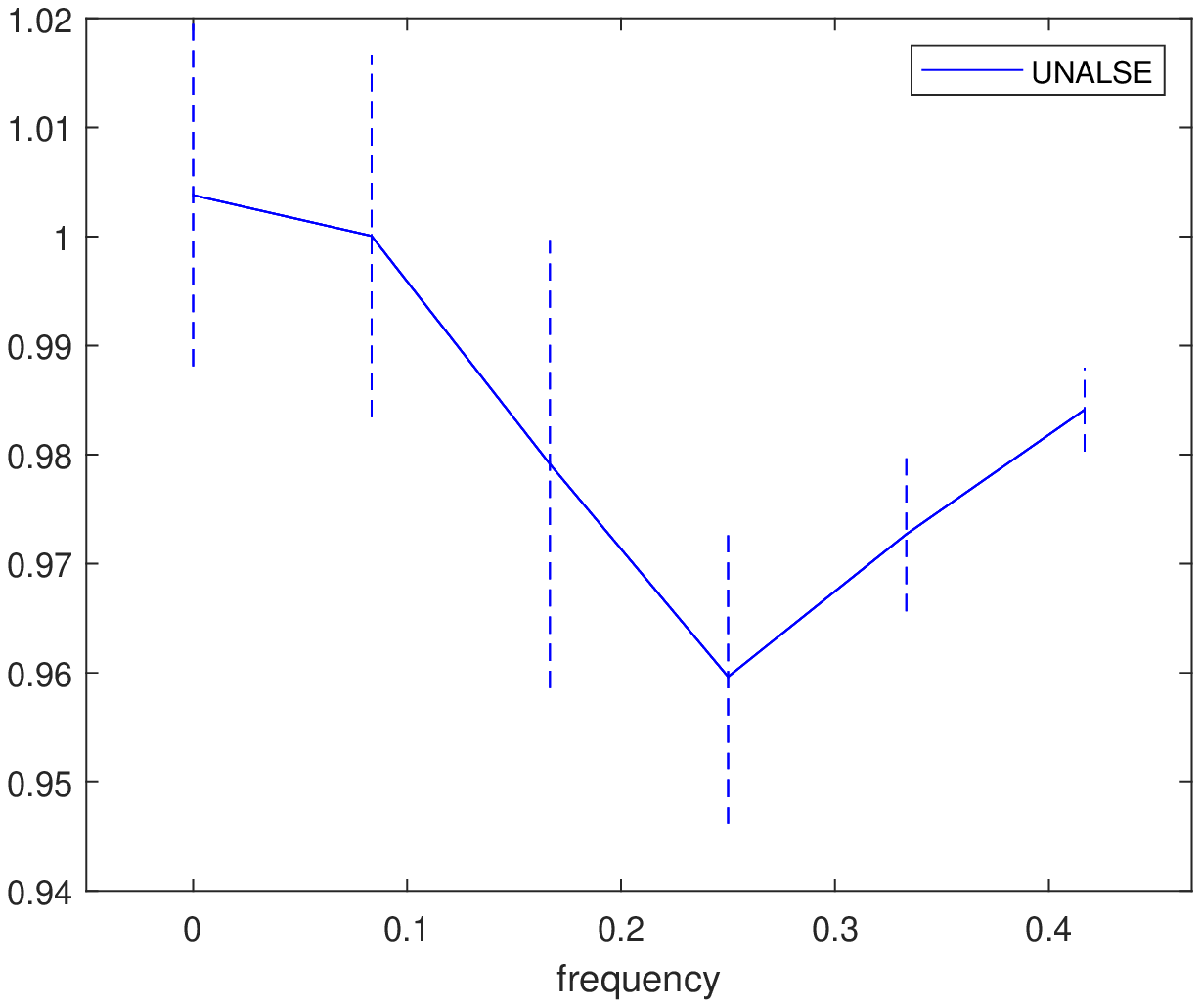}&
              \includegraphics[width=.2\textwidth]{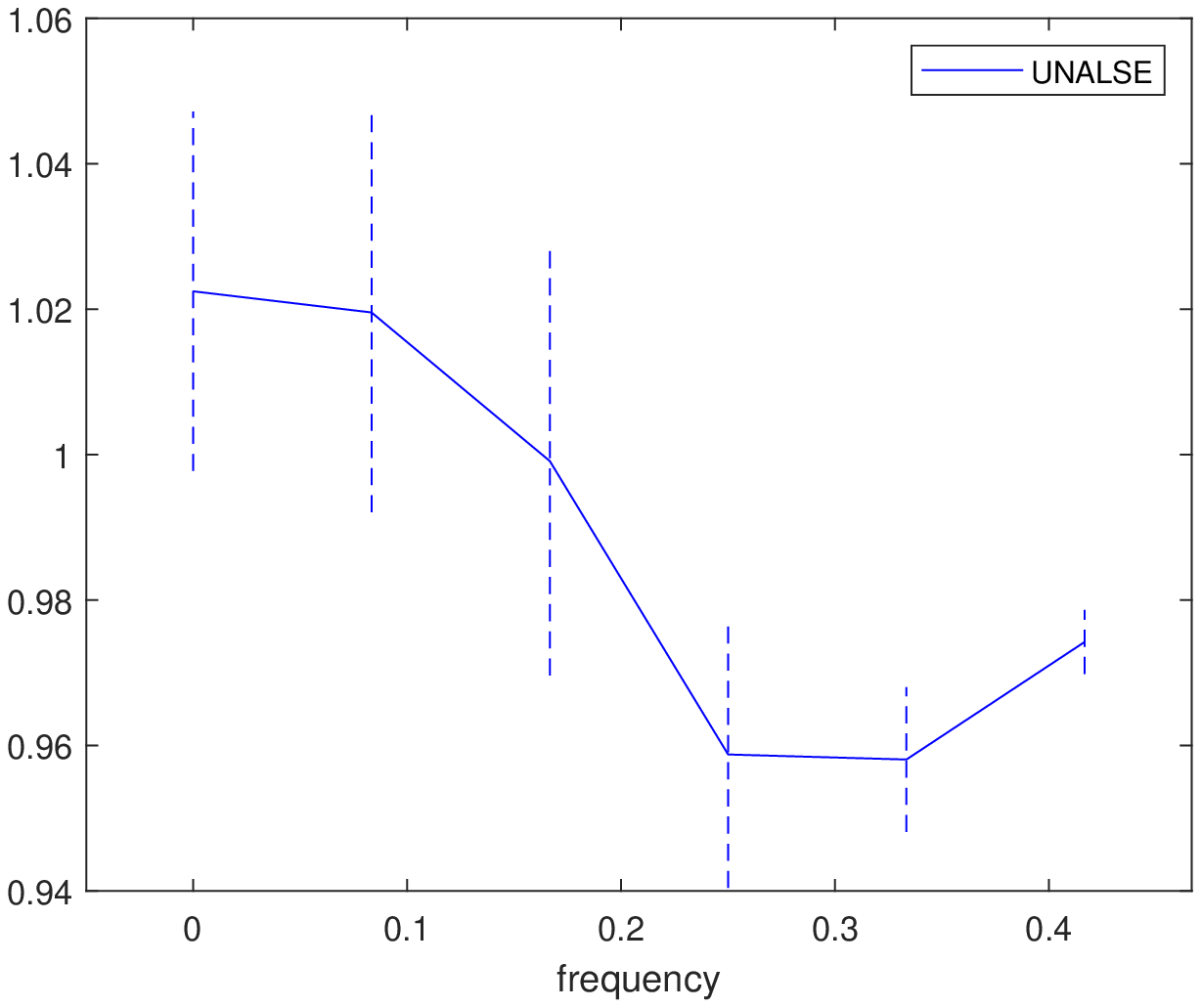}&
              \includegraphics[width=.2\textwidth]{err_ratio_base_3}\\
              \end{tabular}
               \begin{tabular}{cc}
	      {\footnotesize Setting 4}&{\footnotesize Setting 5}\\
              \includegraphics[width=.2\textwidth]{err_ratio_base_4}&
              \includegraphics[width=.2\textwidth]{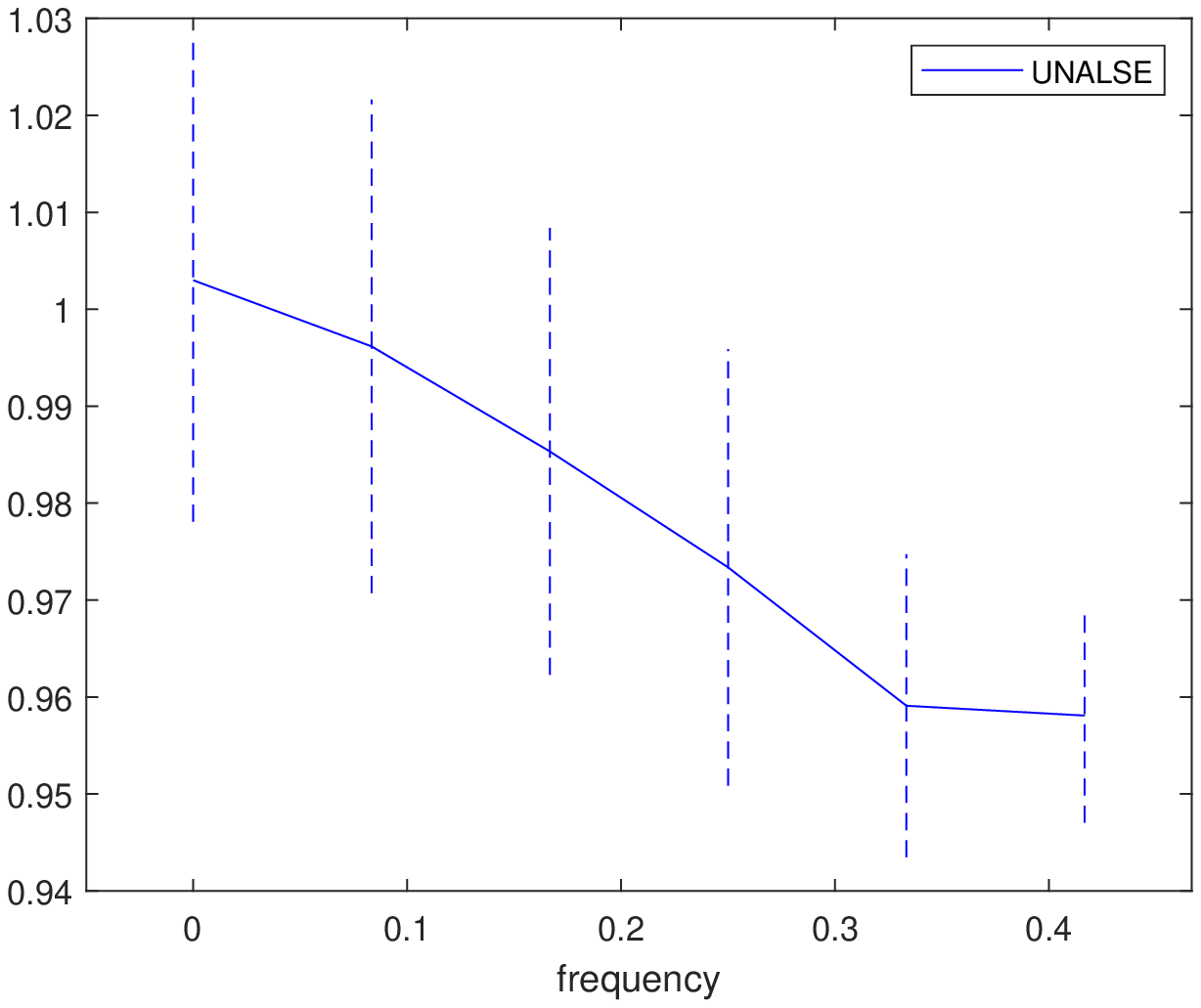}\\
              \end{tabular}
\end{figure}

\clearpage
\subsection{Scenario B}
\begin{figure}[h!]
          \caption{Estimated latent variance proportion $\wh{\beta}(f_h)$ - Scenario B.}\label{fig:beta_B.1.U}
          \centering
          \begin{tabular}{ccc}
          {\footnotesize Setting 1}&{\footnotesize Setting 2}&{\footnotesize Setting 3}\\
              \includegraphics[width=.2\textwidth]{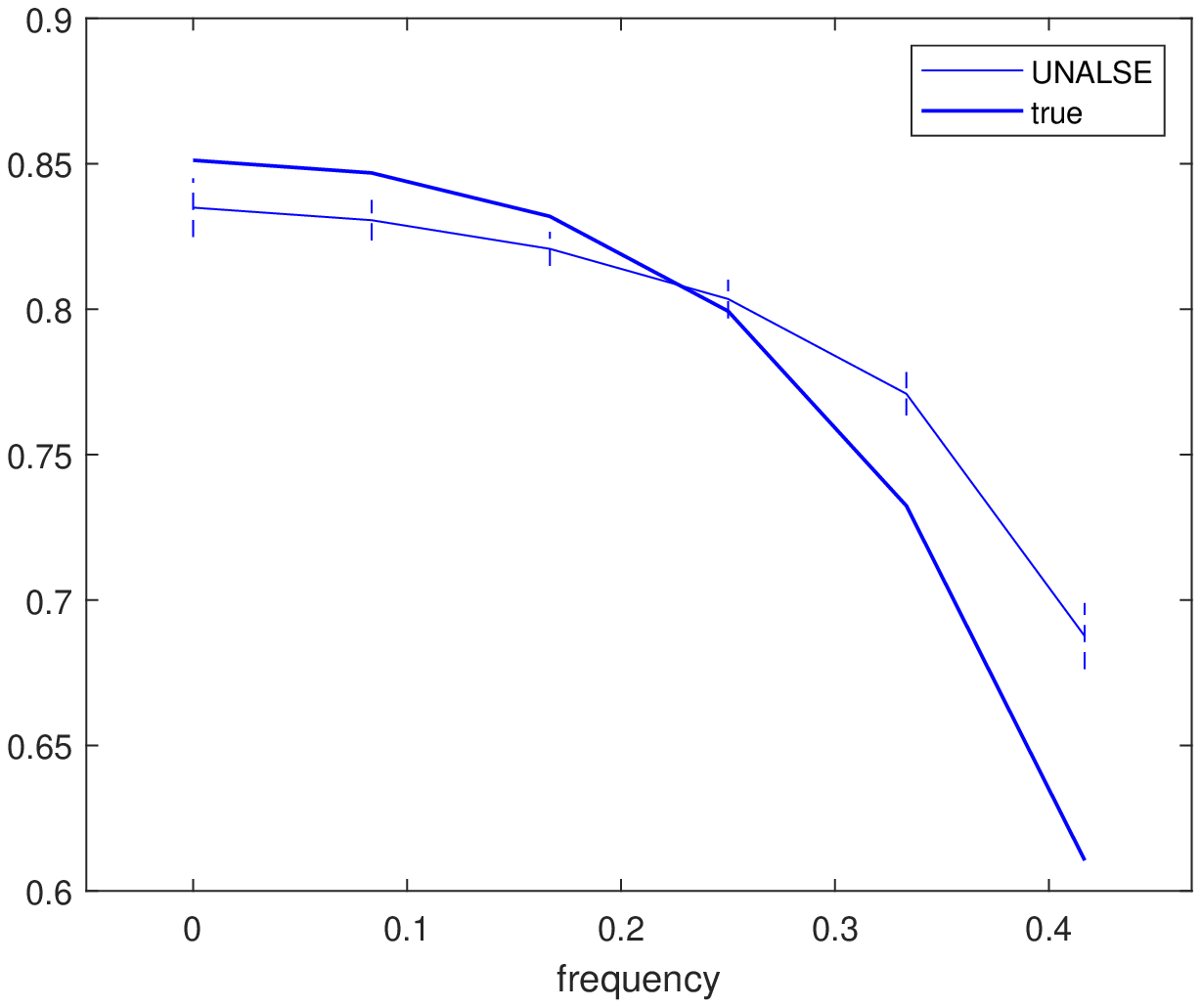}&
              \includegraphics[width=.2\textwidth]{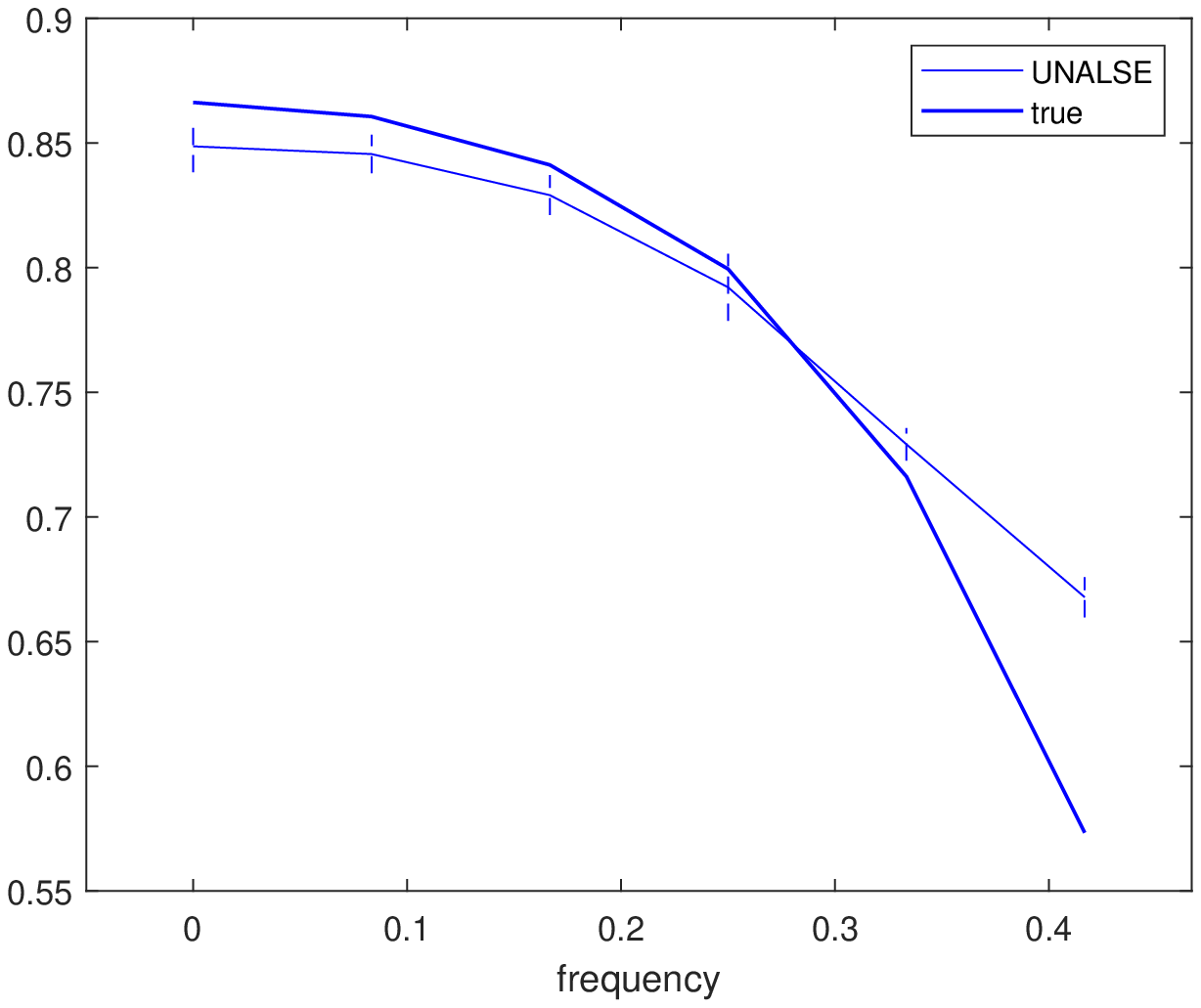}&
              \includegraphics[width=.2\textwidth]{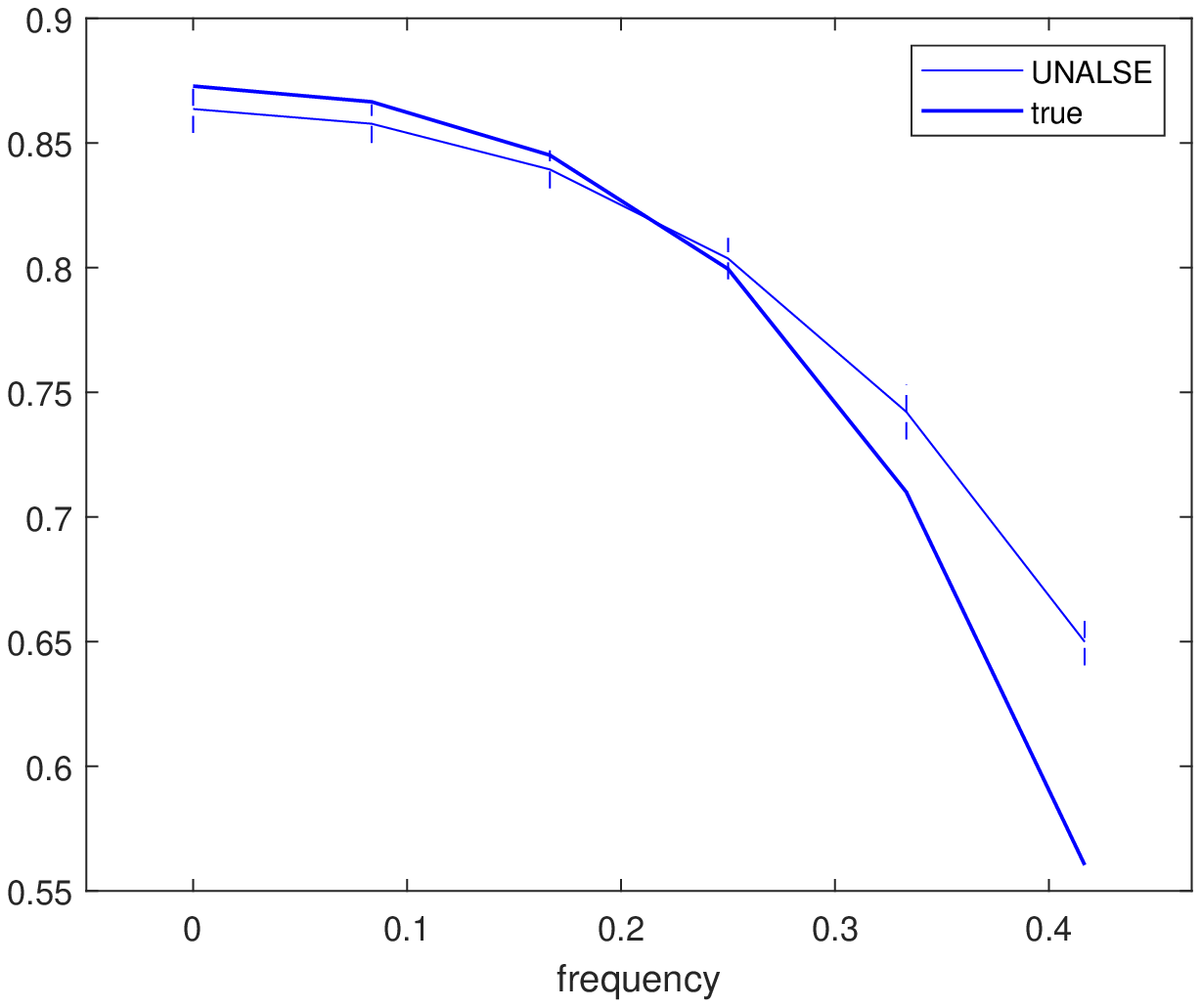}\\
              \end{tabular}
               \begin{tabular}{cc}
                             {\footnotesize Setting 4}&{\footnotesize Setting 5}\\
              \includegraphics[width=.2\textwidth]{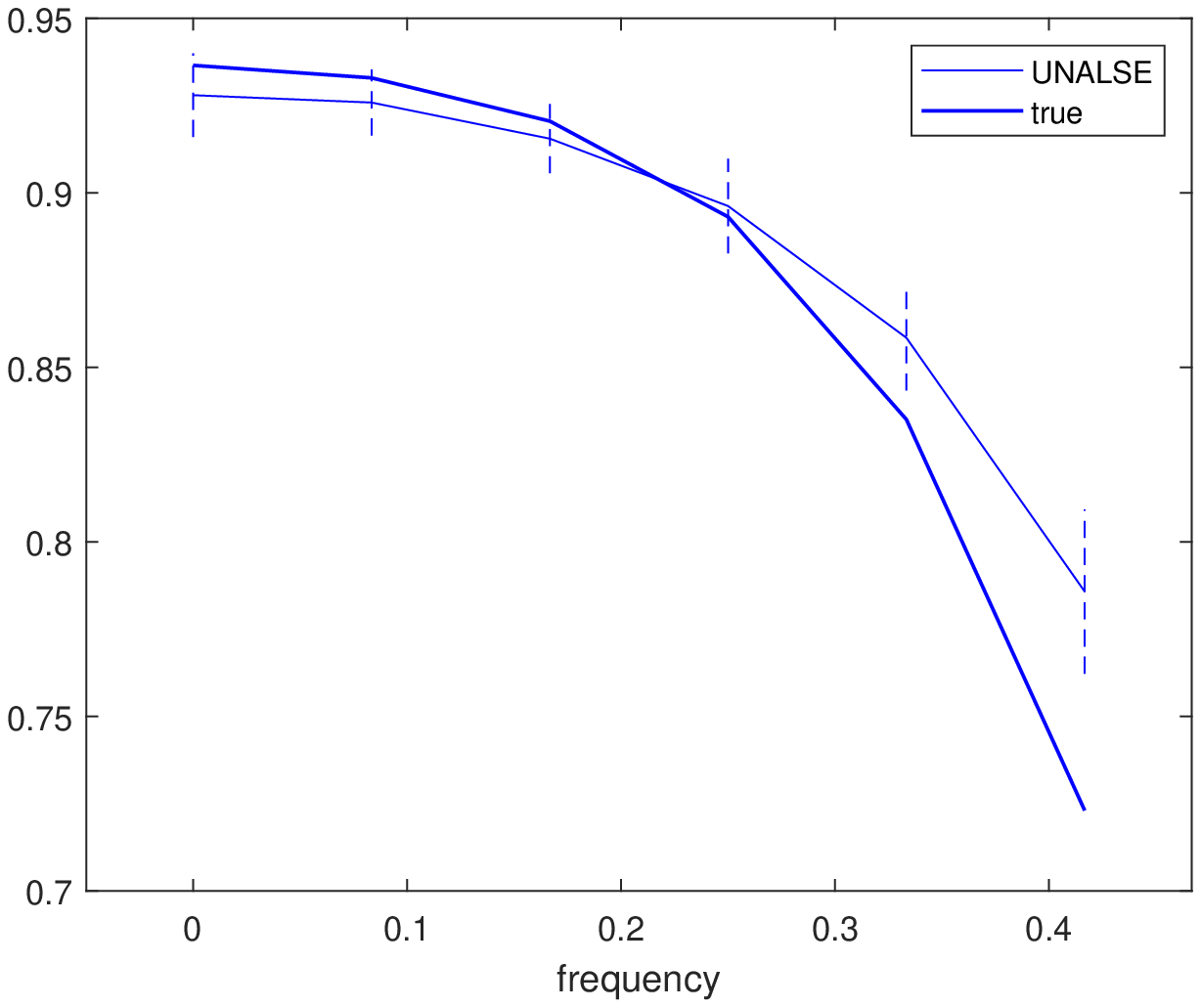}&
              \includegraphics[width=.2\textwidth]{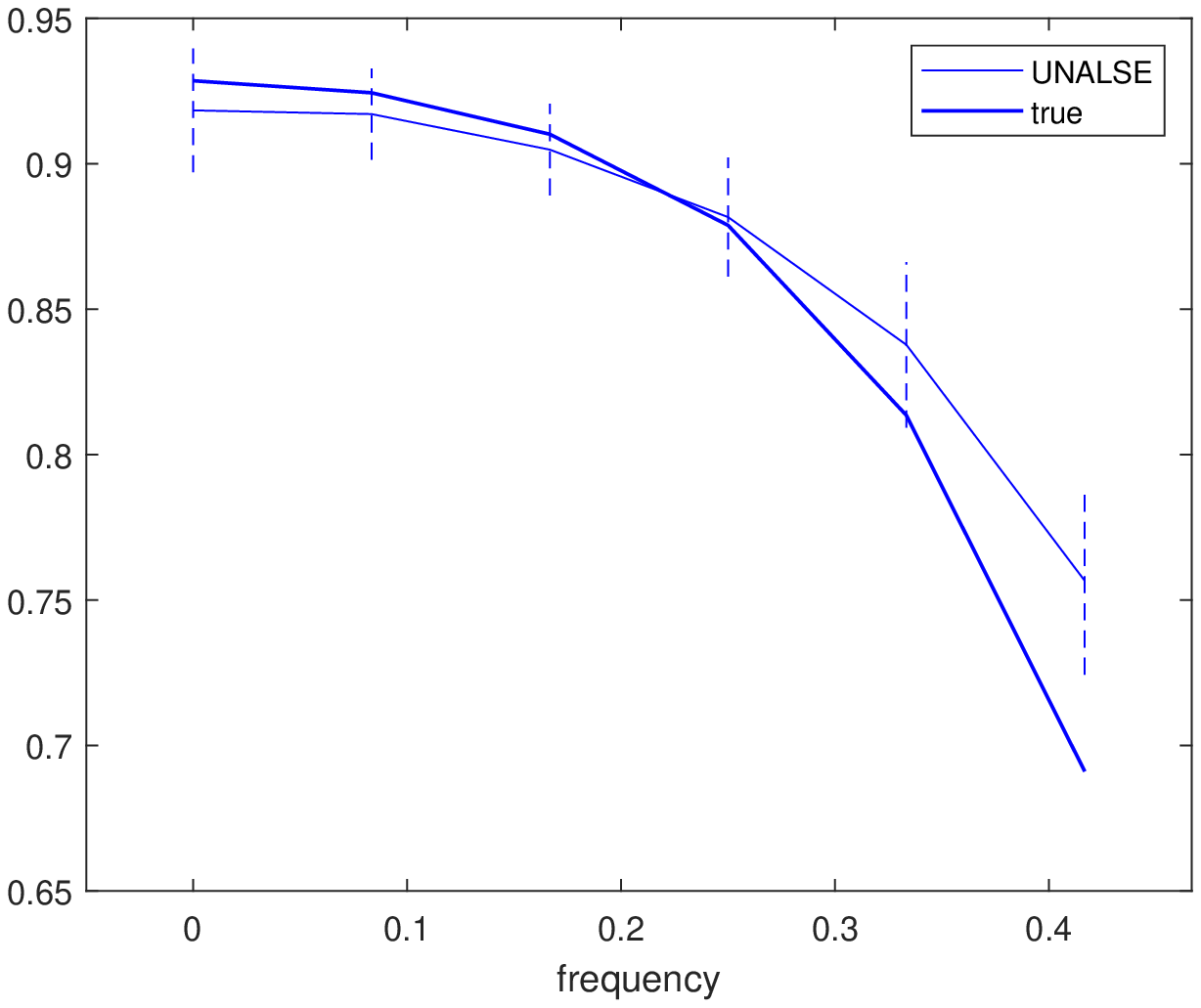}\\
              \end{tabular}
\end{figure}

\begin{figure}[h!]
          \caption{Nonzero predictive values $nzpv(f_h)$ - Scenario B.}
          \label{fig:pred_pos_B.1.U}
          \centering
          \begin{tabular}{ccc}
       {\footnotesize Setting 3}&{\footnotesize Setting 4}&{\footnotesize Setting 5}\\
              \includegraphics[width=.2\textwidth]{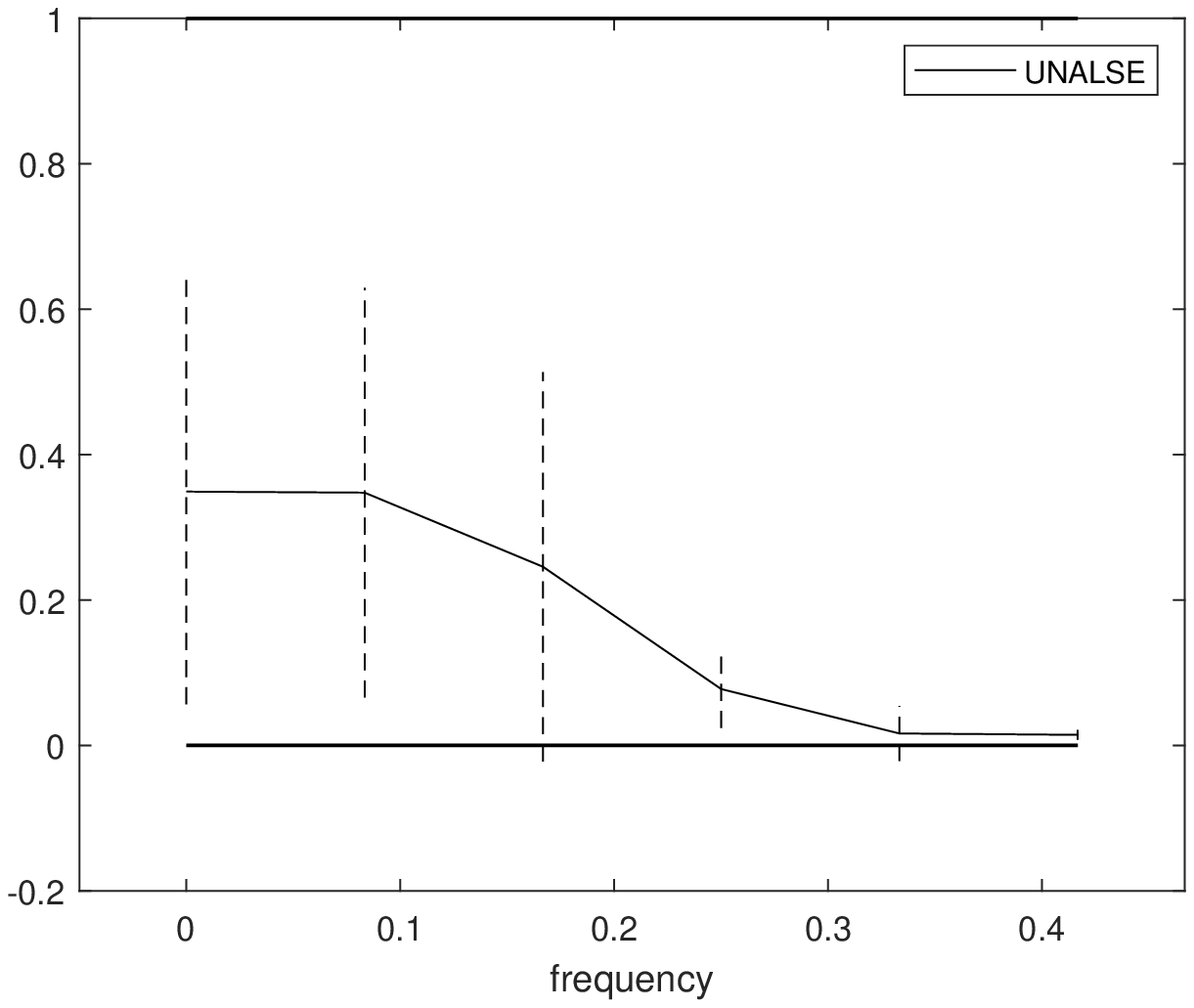}&
              \includegraphics[width=.2\textwidth]{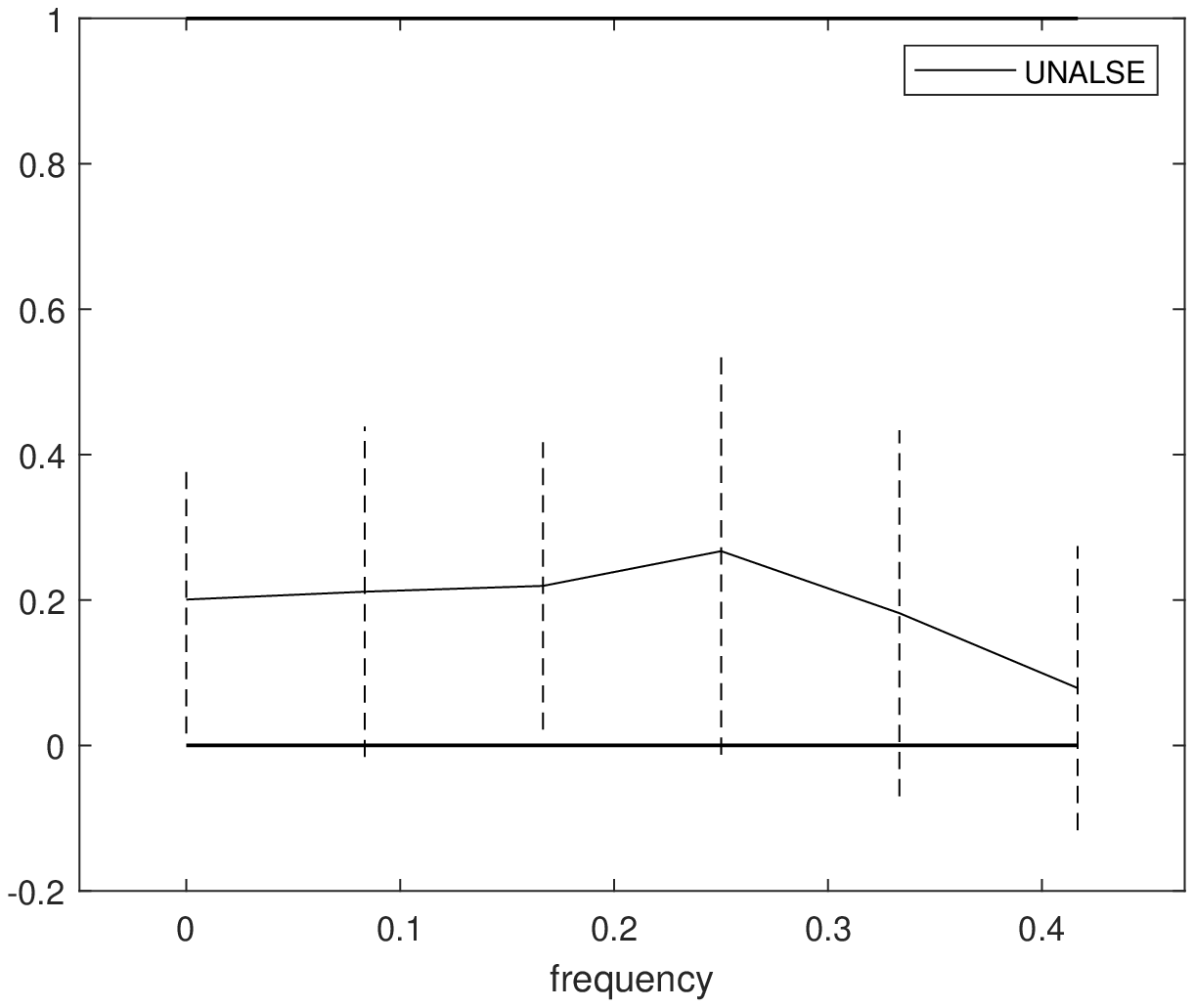}&
              \includegraphics[width=.2\textwidth]{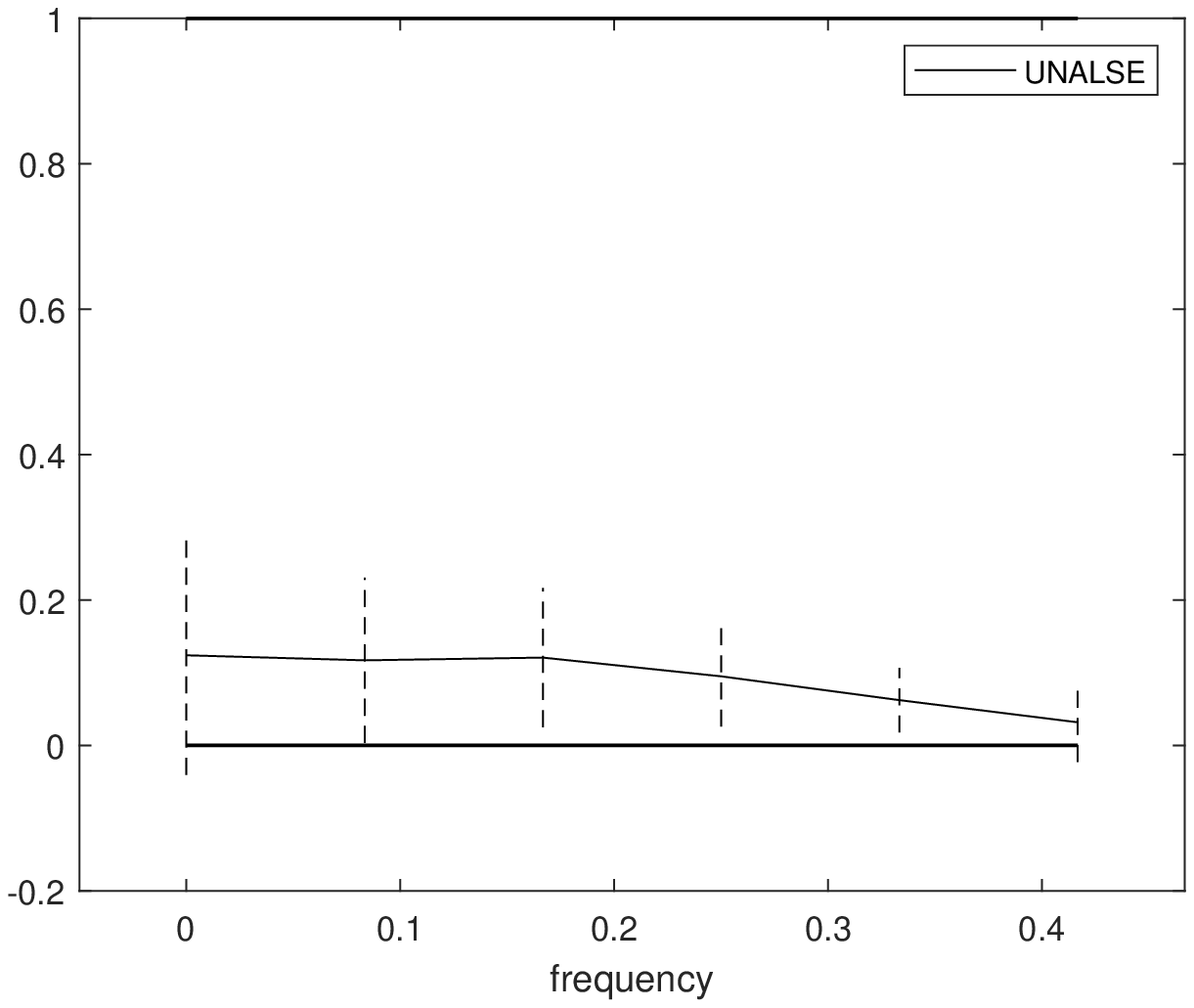}\\
              \end{tabular}
\end{figure}

\begin{figure}[h!]
          \caption{$err_{\wh{L}}(f_h)$ - Scenario B.}\label{fig:err_L_B.1.U}
          \centering
          \begin{tabular}{ccc}
          {\footnotesize Setting 1}&{\footnotesize Setting 2}&{\footnotesize Setting 3}\\
              \includegraphics[width=.2\textwidth]{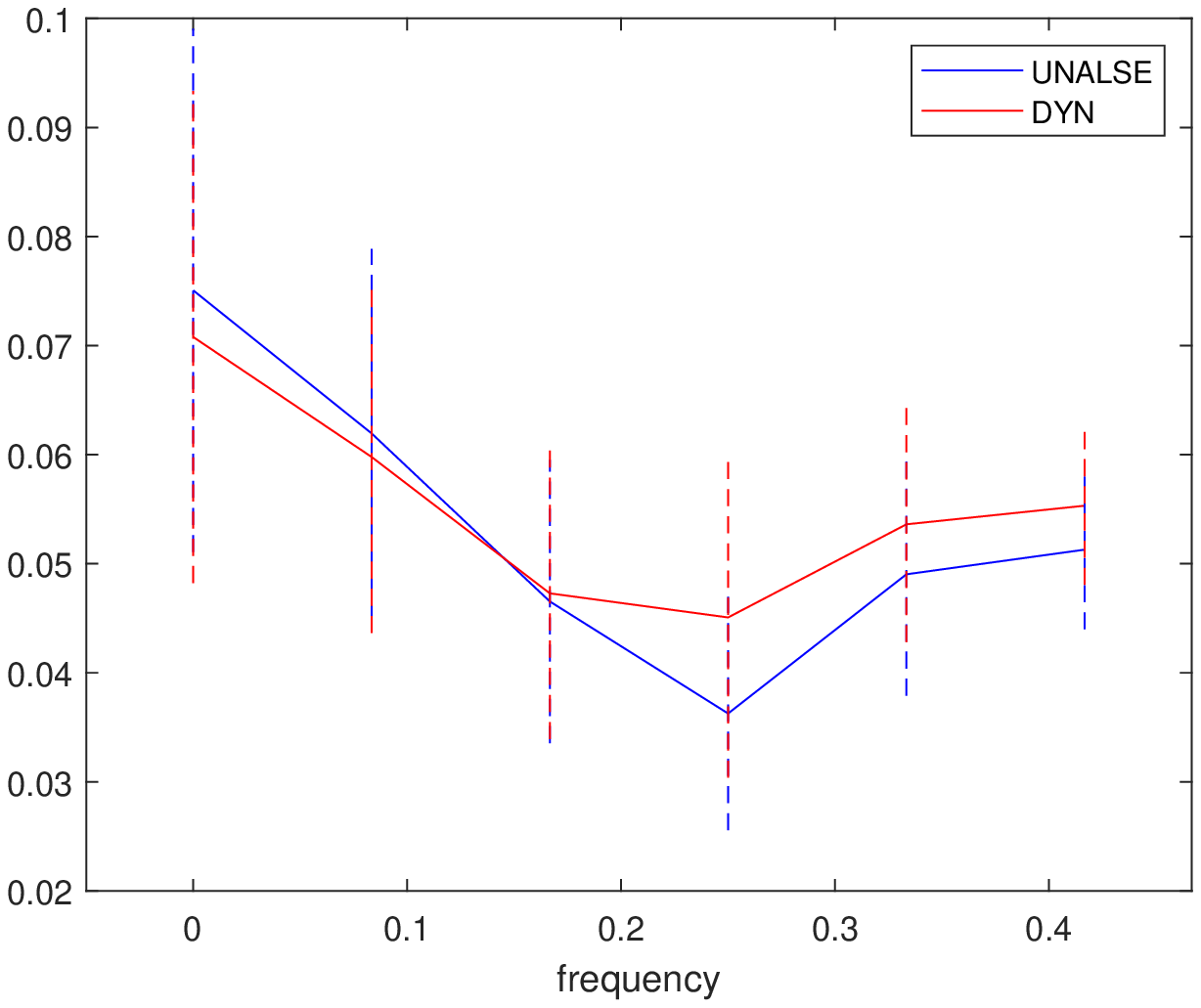}&
              \includegraphics[width=.2\textwidth]{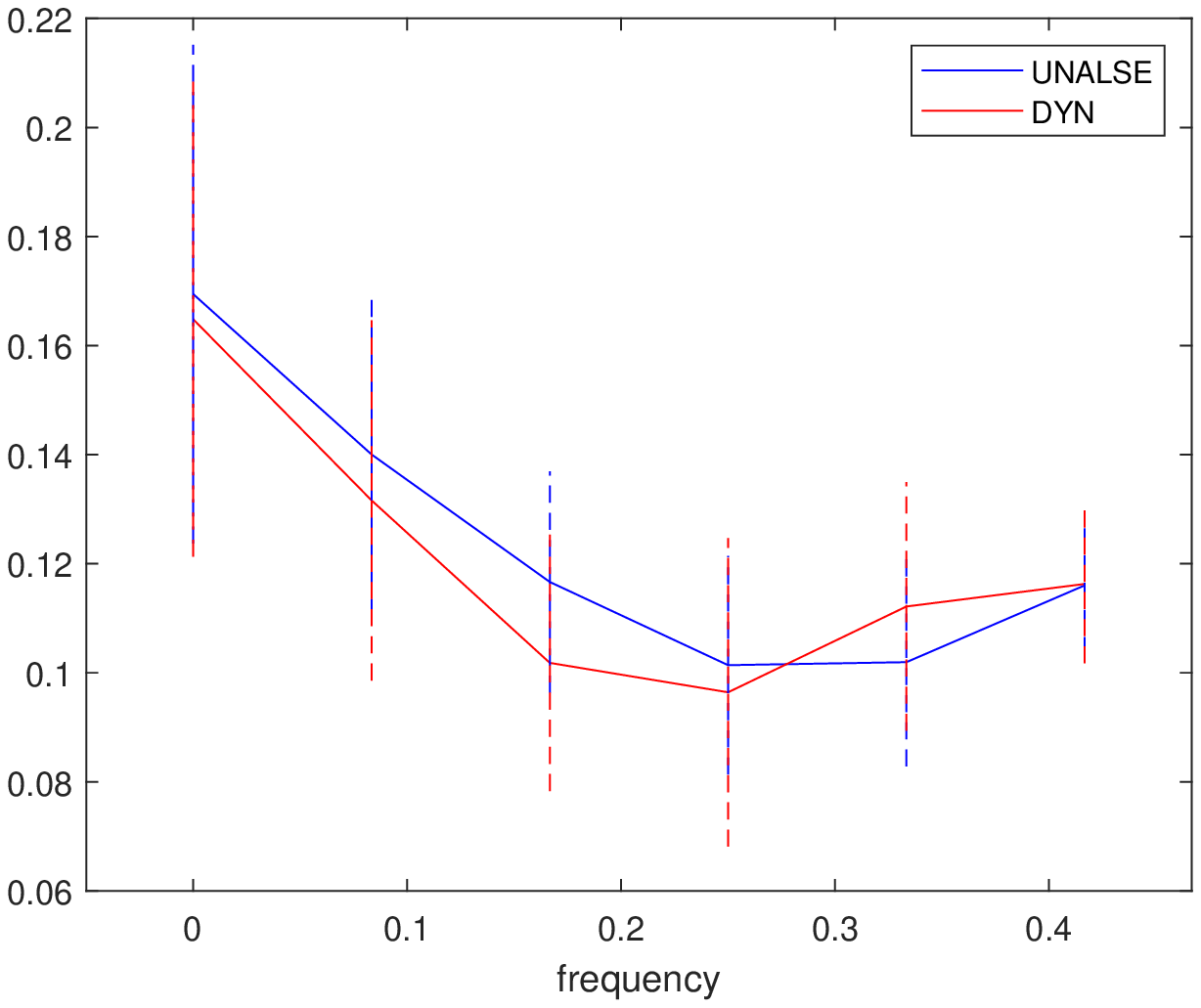}&
              \includegraphics[width=.2\textwidth]{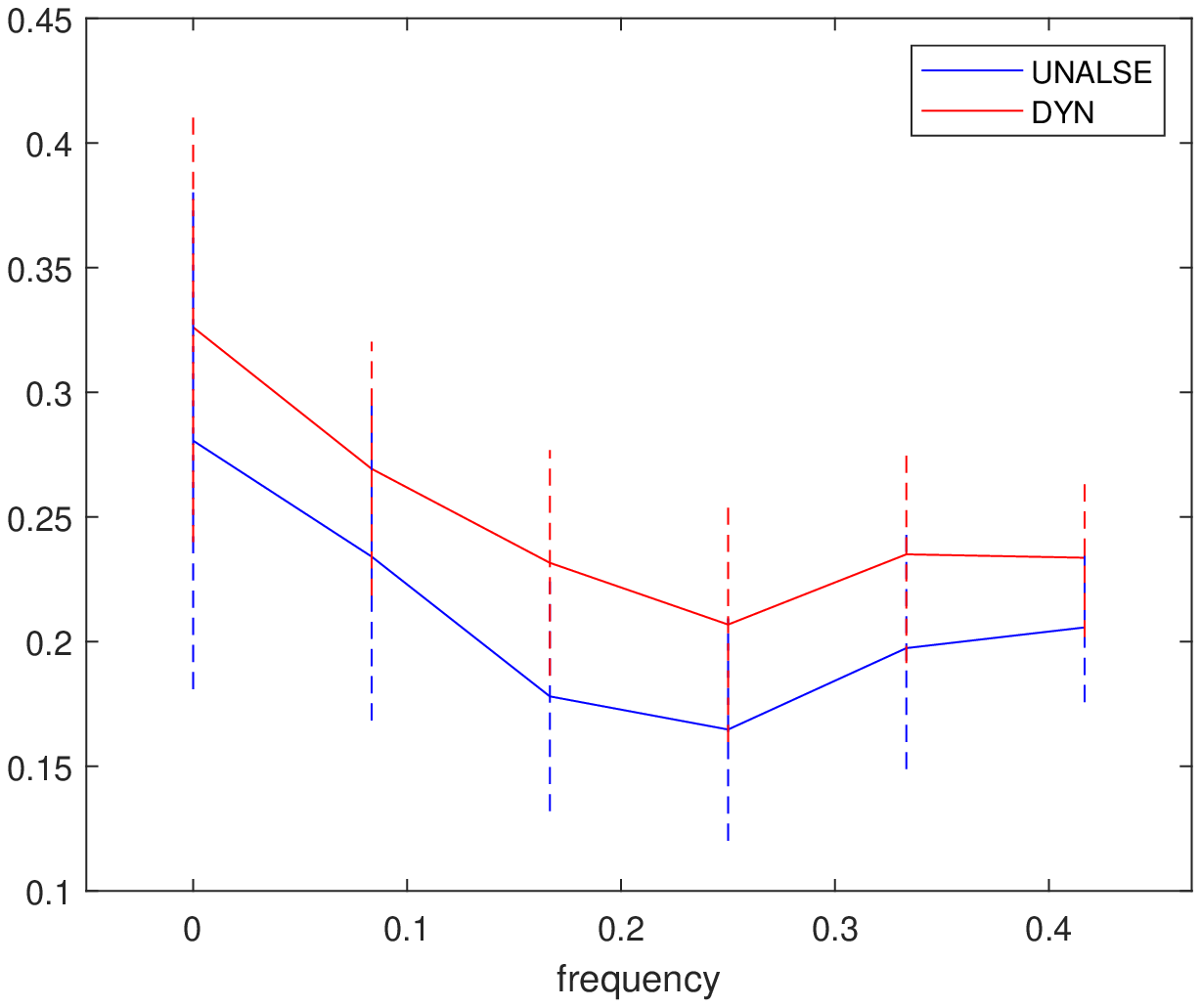}\\
              \end{tabular}
               \begin{tabular}{cc}
	      {\footnotesize Setting 4}&{\footnotesize Setting 5}\\
              \includegraphics[width=.2\textwidth]{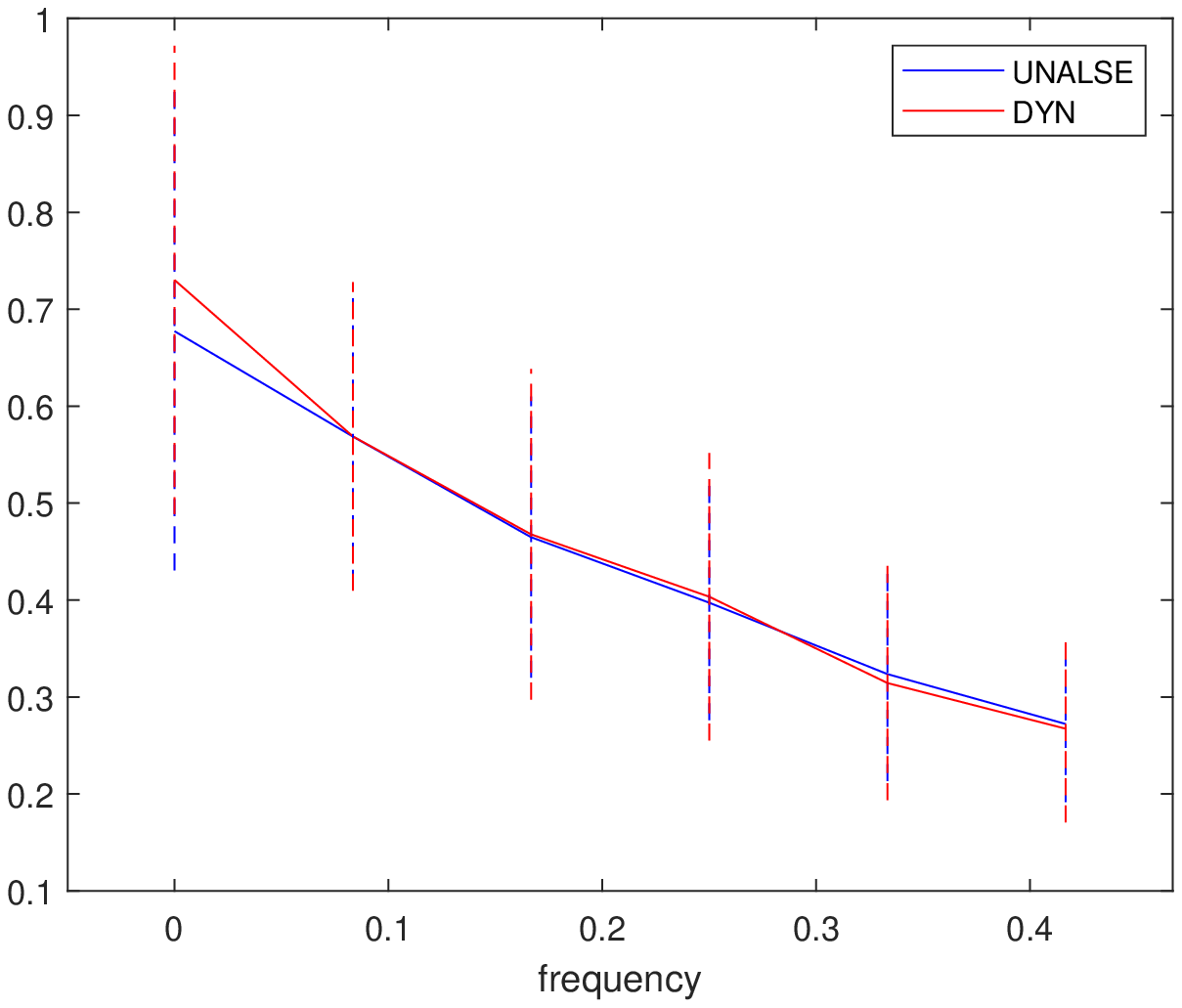}&
              \includegraphics[width=.2\textwidth]{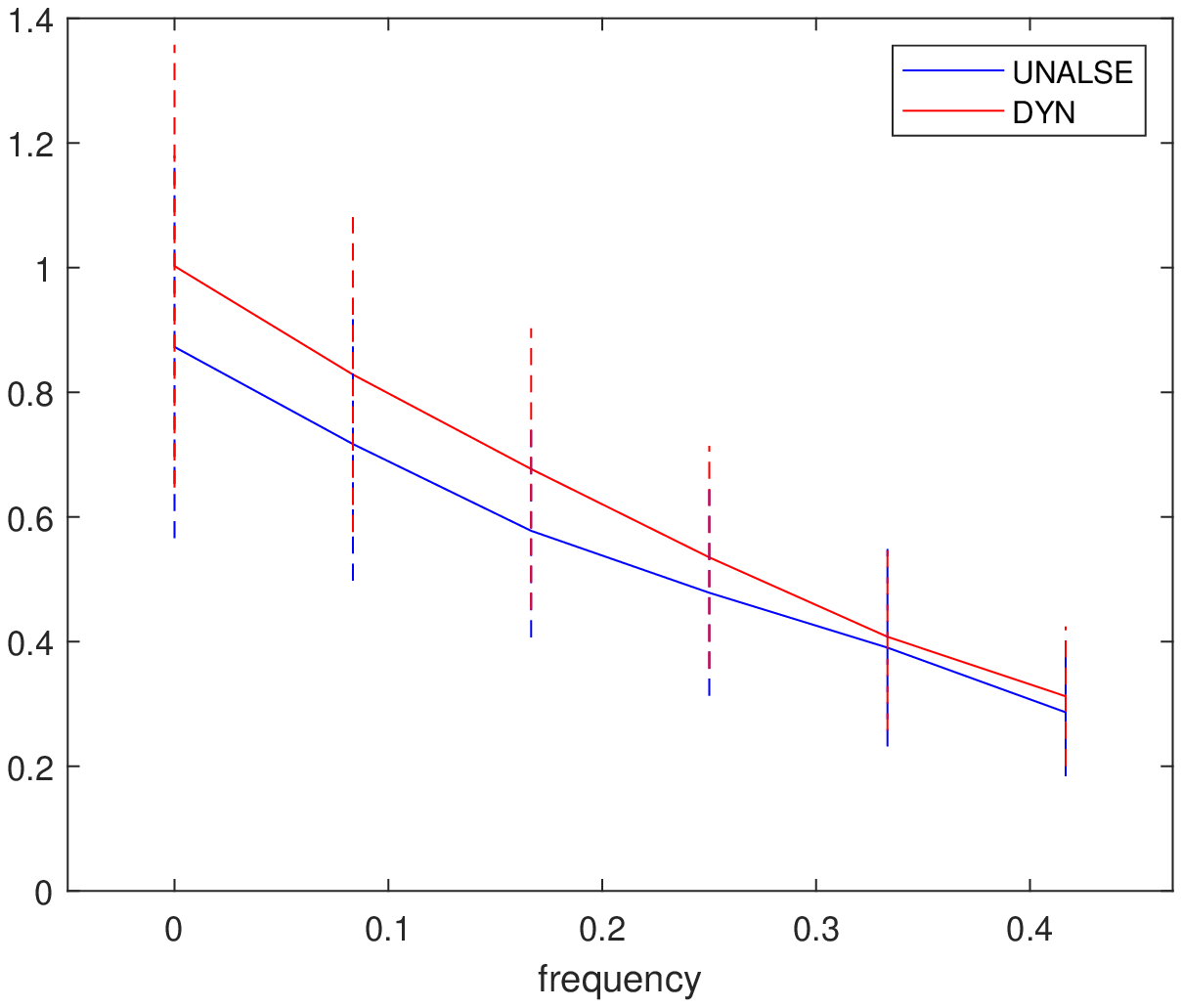}\\
              \end{tabular}
\end{figure}

\begin{figure}[h!]
          \caption{$err_{ratio}(f_h)$ - Scenario B.}\label{fig:err_ratio_B.1.U}
          \centering
          \begin{tabular}{ccc}
          {\footnotesize Setting 1}&{\footnotesize Setting 2}&{\footnotesize Setting 3}\\
              \includegraphics[width=.2\textwidth]{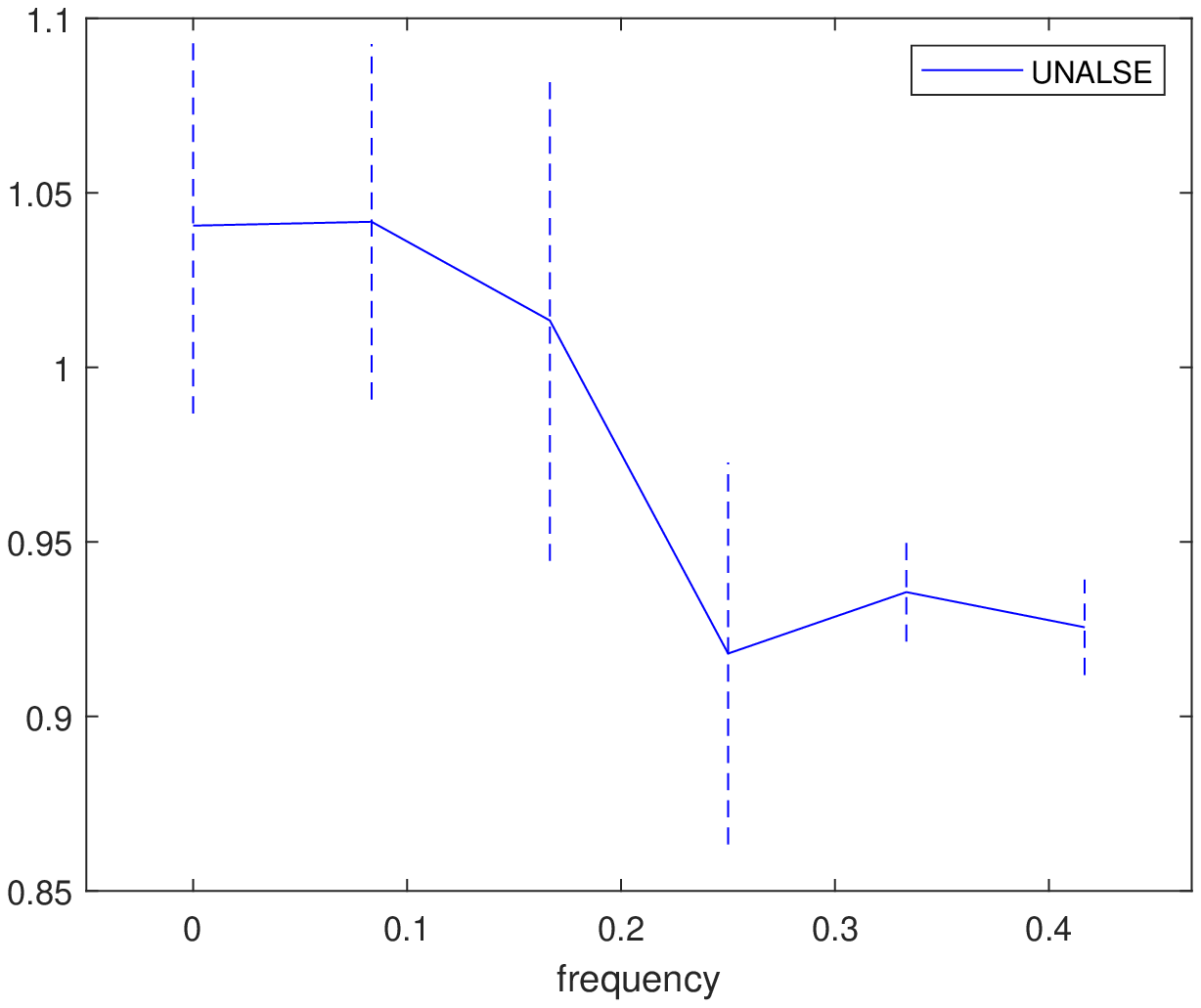}&
              \includegraphics[width=.2\textwidth]{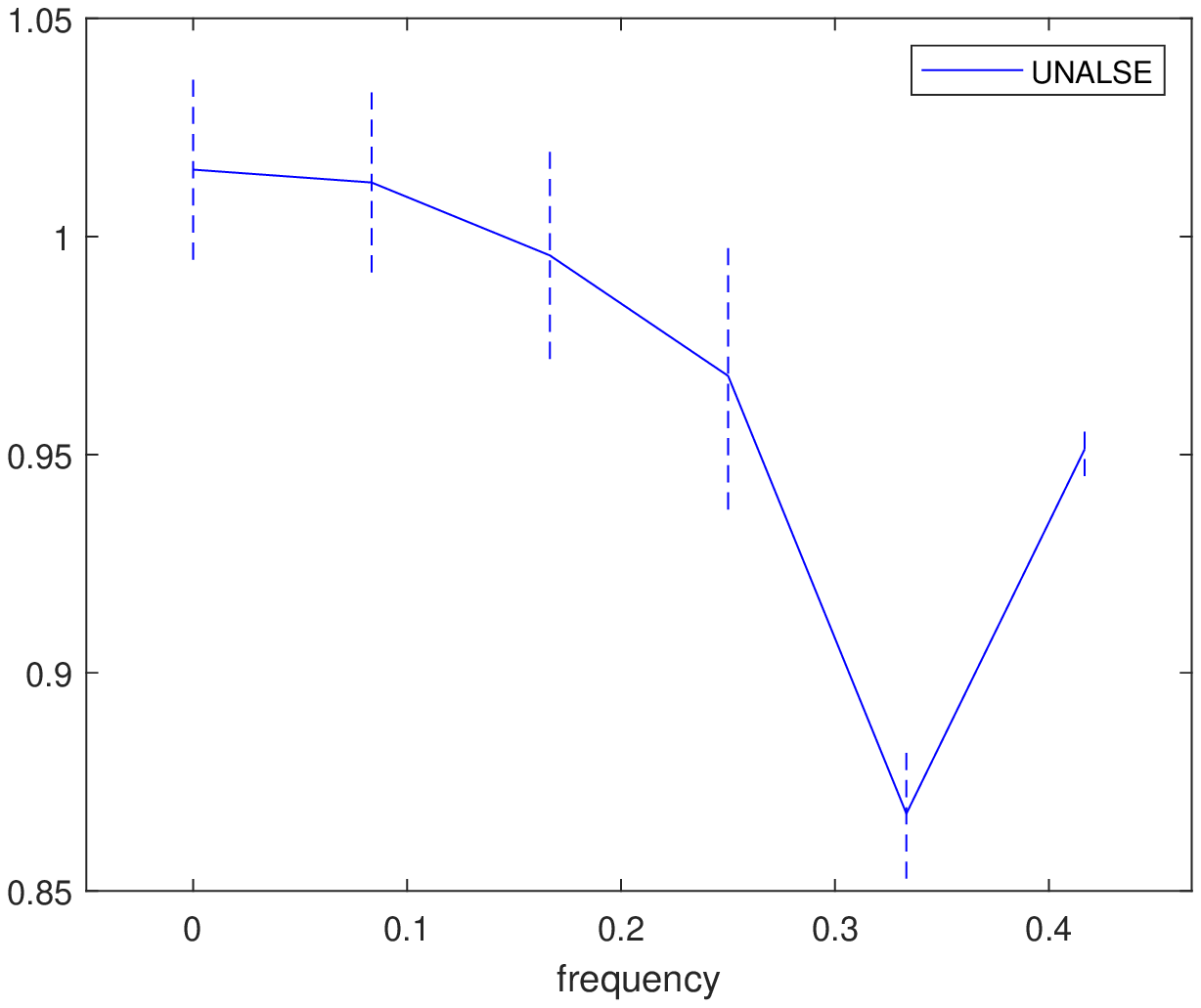}&
              \includegraphics[width=.2\textwidth]{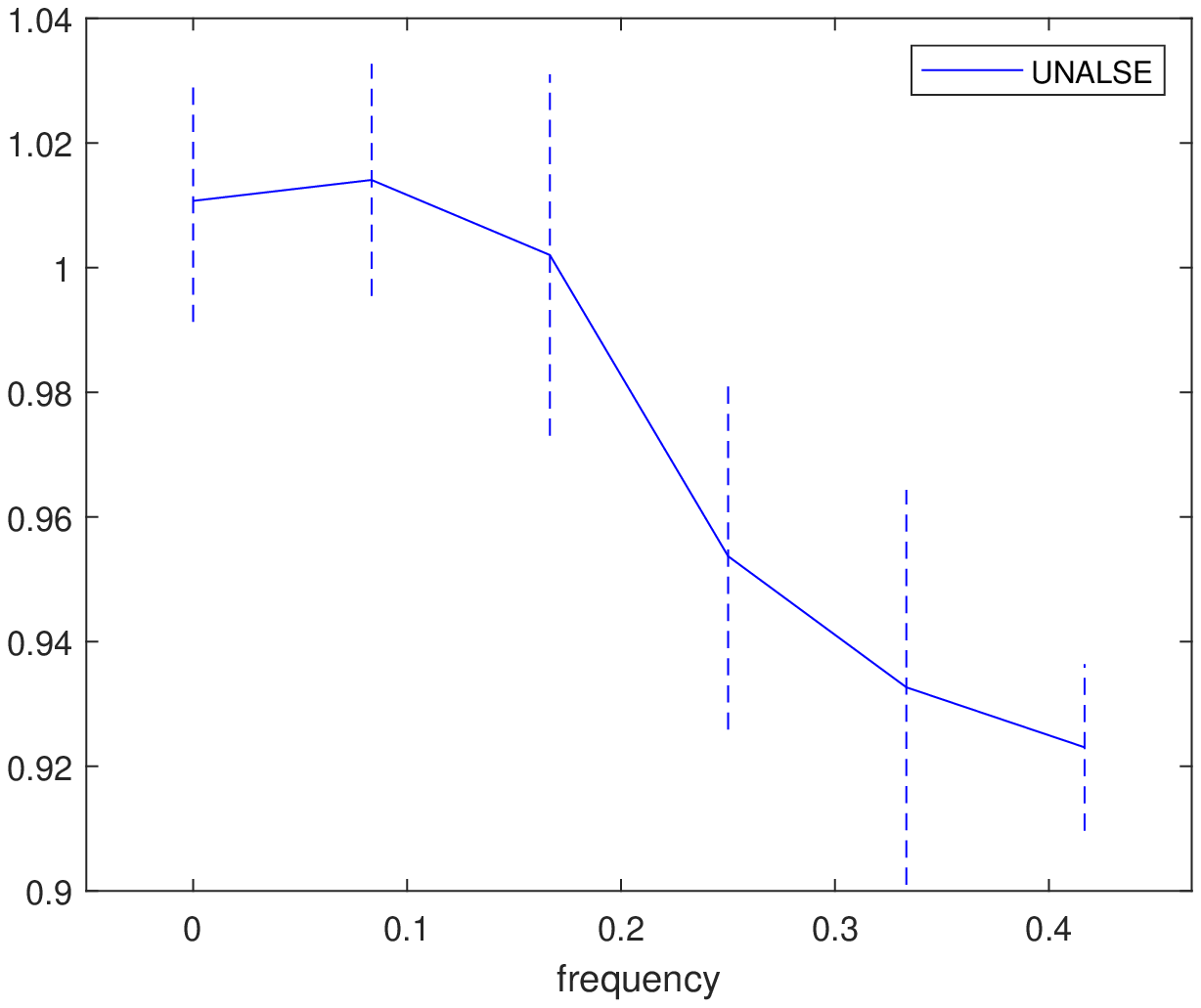}\\
              \end{tabular}
               \begin{tabular}{cc}
	      {\footnotesize Setting 4}&{\footnotesize Setting 5}\\
              \includegraphics[width=.2\textwidth]{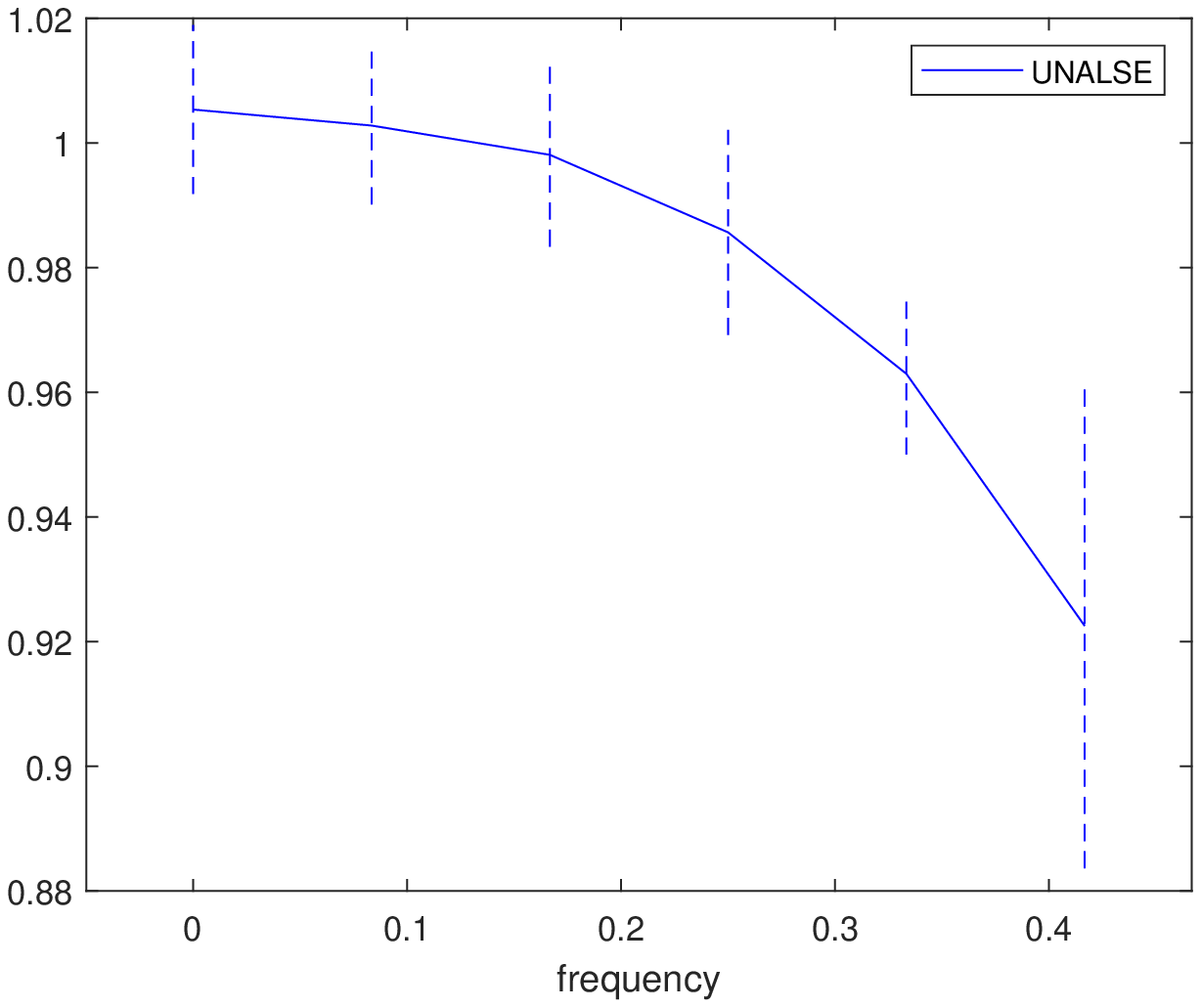}&
              \includegraphics[width=.2\textwidth]{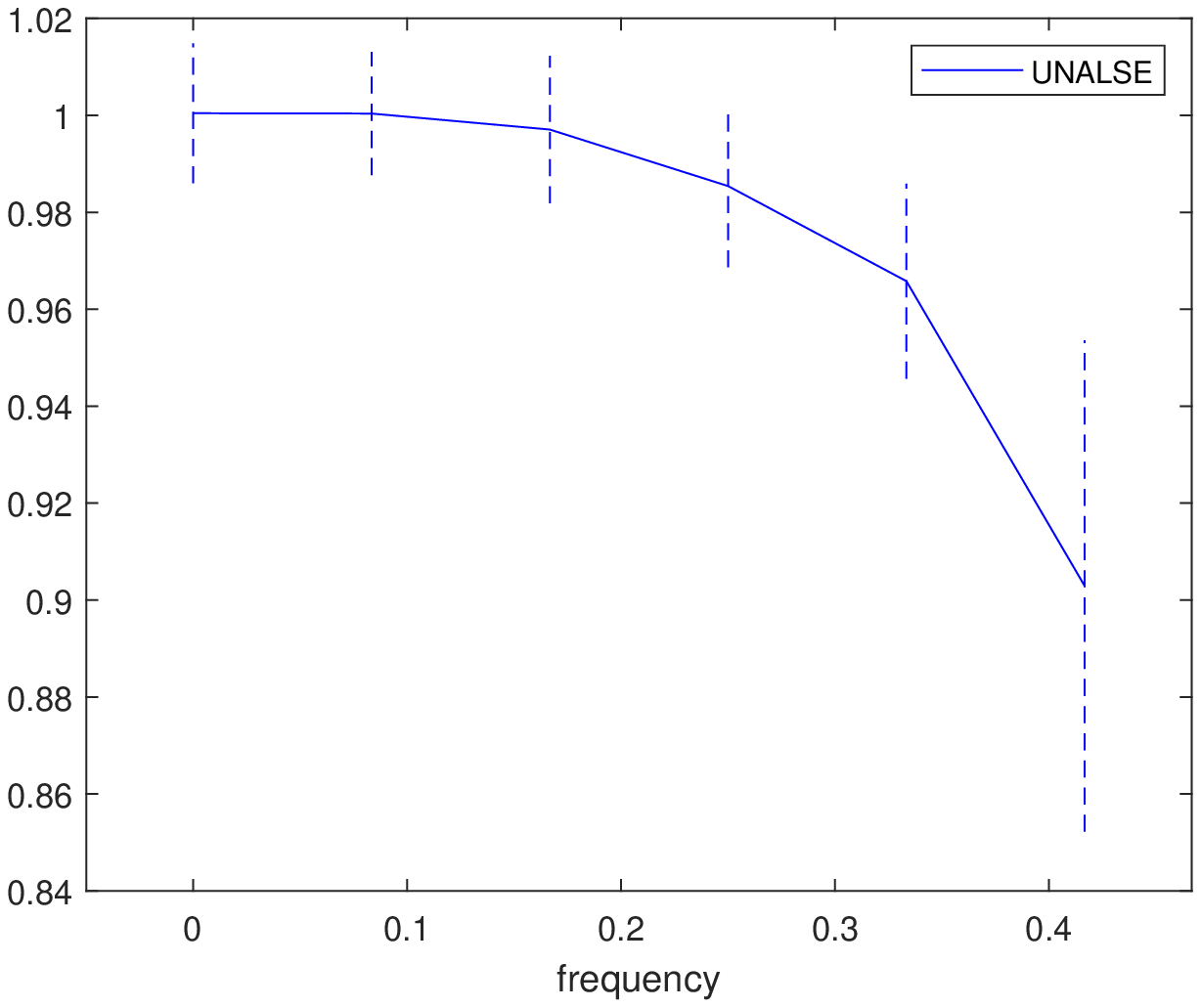}\\
              \end{tabular}
\end{figure}

\clearpage
\subsection{Scenario C}

\begin{figure}[h!]
          \caption{Estimated latent variance proportion $\wh{\beta}(f_h)$ - Scenario C.}\label{fig:beta_C.1.U}
          \centering
          \begin{tabular}{ccc}
          {\footnotesize Setting 1}&{\footnotesize Setting 2}&{\footnotesize Setting 3}\\
              \includegraphics[width=.2\textwidth]{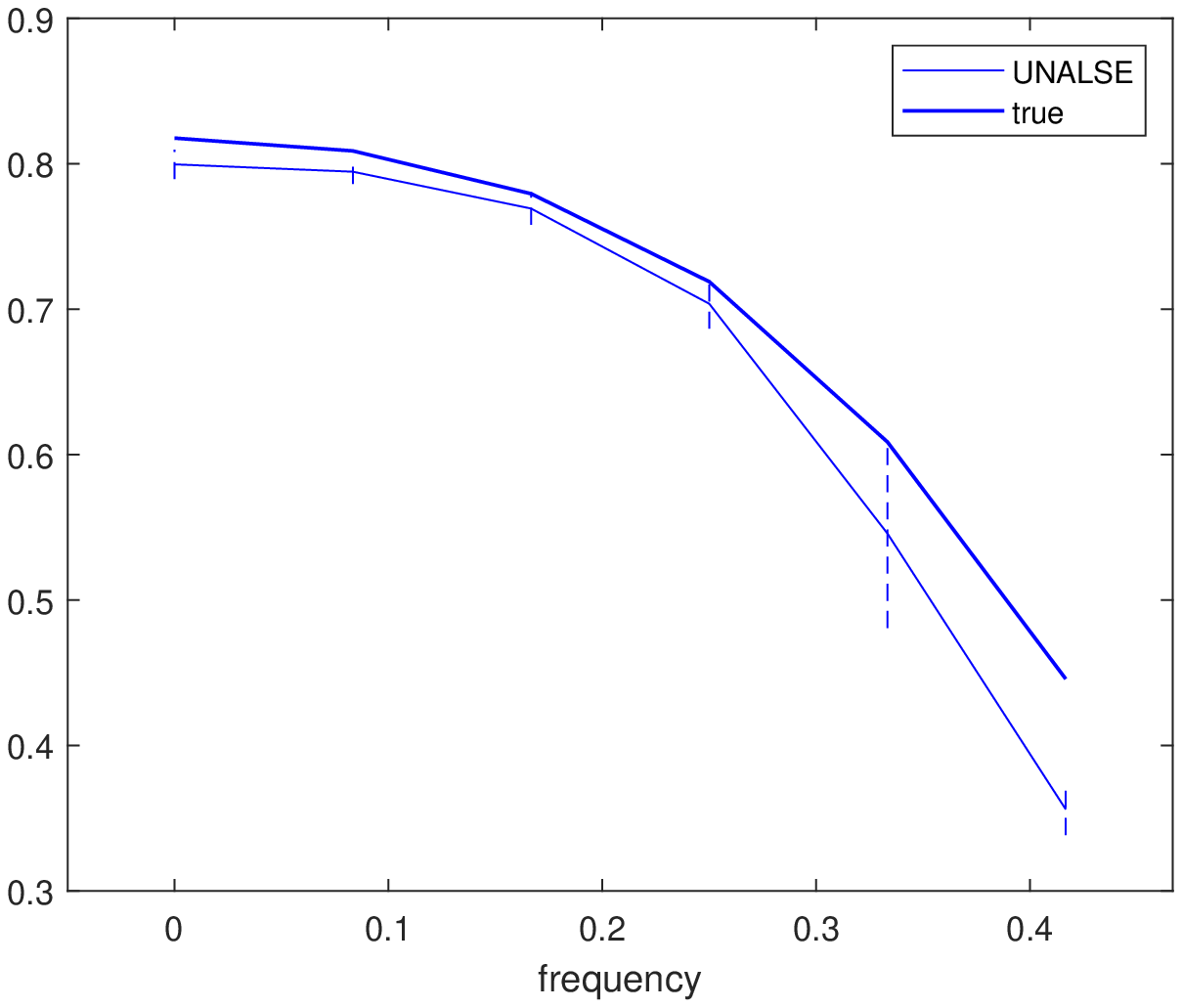}&
              \includegraphics[width=.2\textwidth]{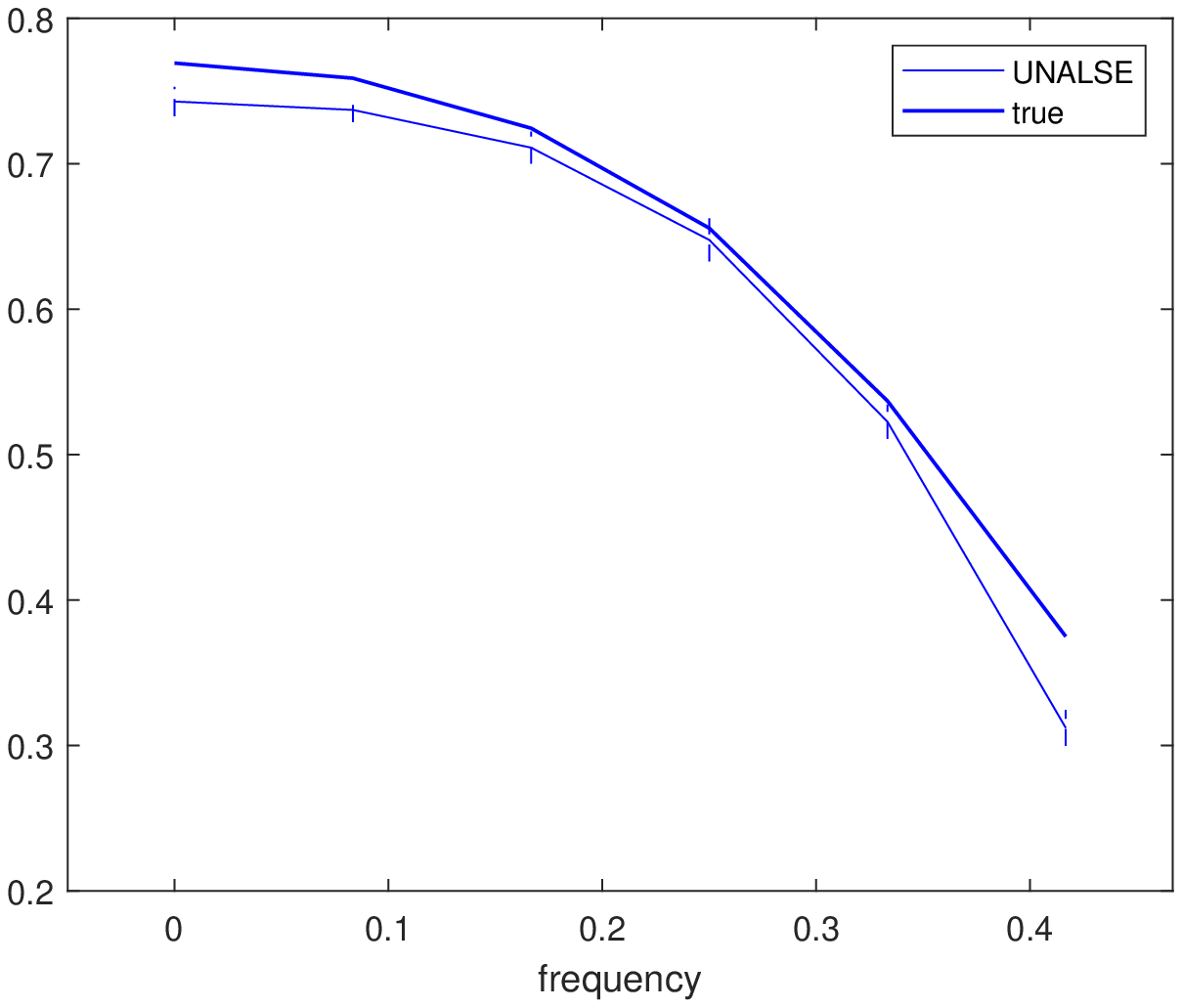}&
              \includegraphics[width=.2\textwidth]{theta_max_UNALSE_3}\\
              \end{tabular}
               \begin{tabular}{cc}
                             {\footnotesize Setting 4}&{\footnotesize Setting 5}\\
              \includegraphics[width=.2\textwidth]{theta_max_UNALSE_4}&
              \includegraphics[width=.2\textwidth]{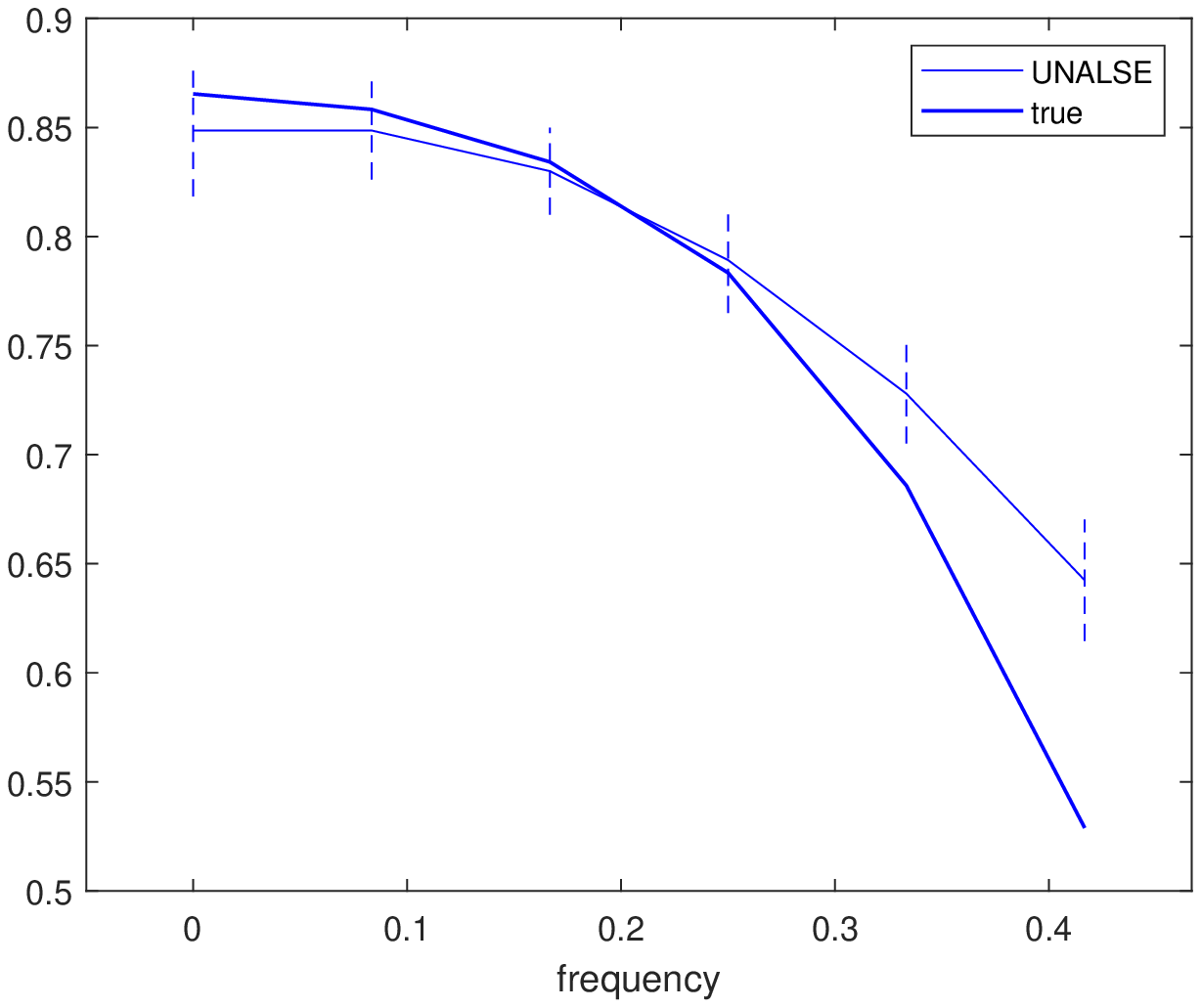}\\
              \end{tabular}
\end{figure}

\begin{figure}[h!]
          \caption{$err_{\wh{L}}(f_h)$ - Scenario C.}\label{fig:err_L_C.1.U}
          \centering
          \begin{tabular}{ccc}
          {\footnotesize Setting 1}&{\footnotesize Setting 2}&{\footnotesize Setting 3}\\
              \includegraphics[width=.2\textwidth]{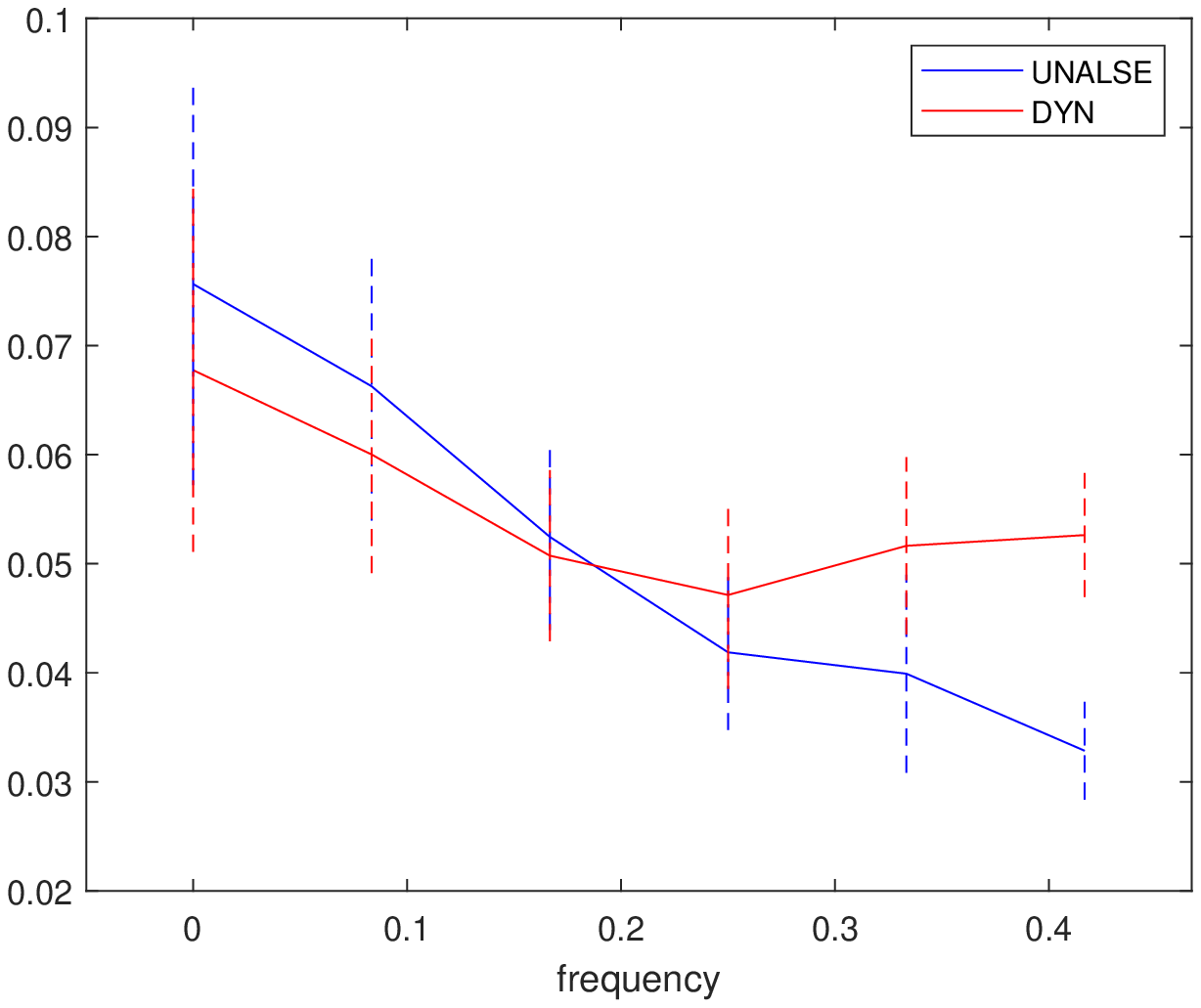}&
              \includegraphics[width=.2\textwidth]{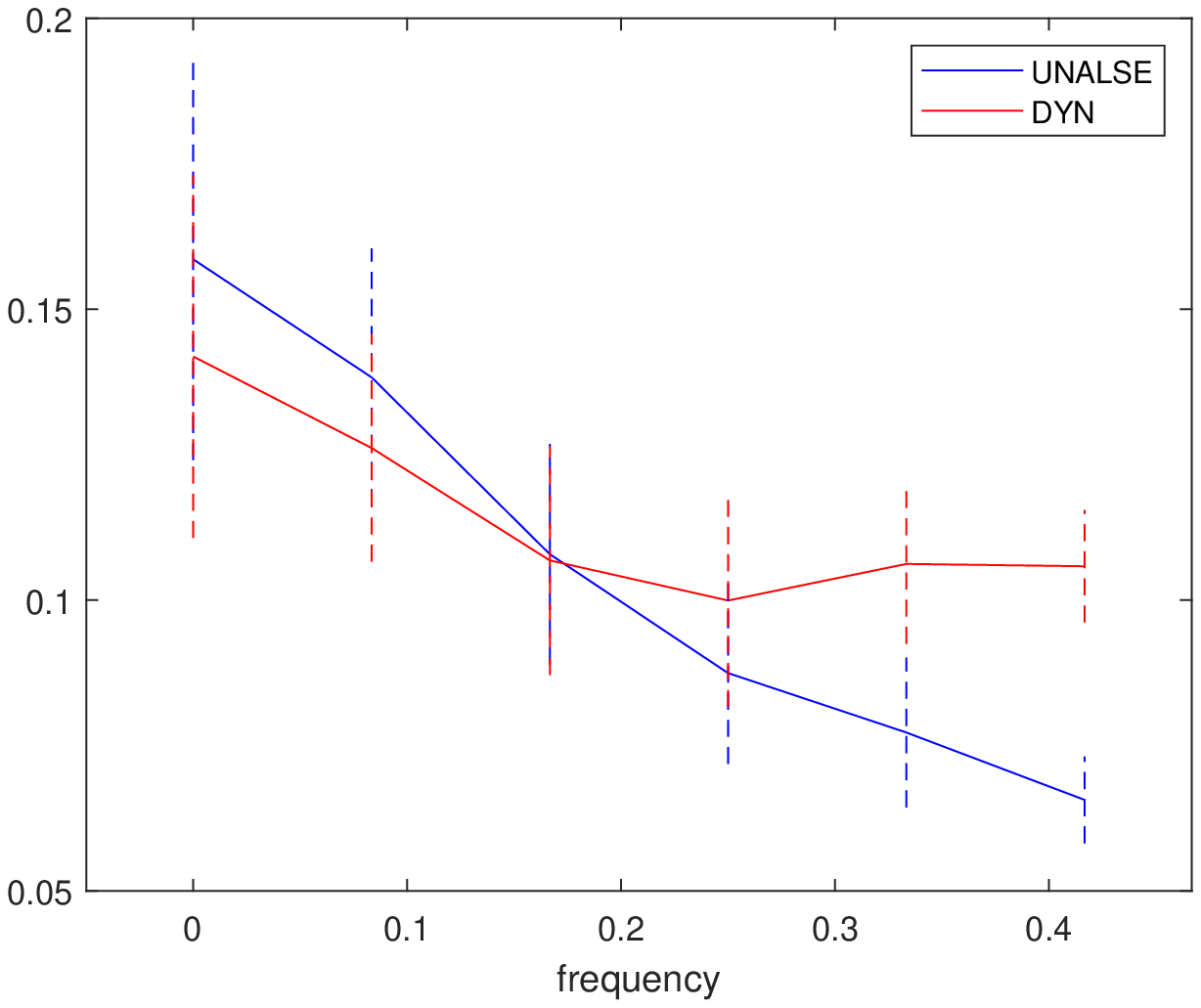}&
              \includegraphics[width=.2\textwidth]{err_L_max_3}\\
              \end{tabular}
               \begin{tabular}{cc}
	      {\footnotesize Setting 4}&{\footnotesize Setting 5}\\
              \includegraphics[width=.2\textwidth]{err_L_max_4}&
              \includegraphics[width=.2\textwidth]{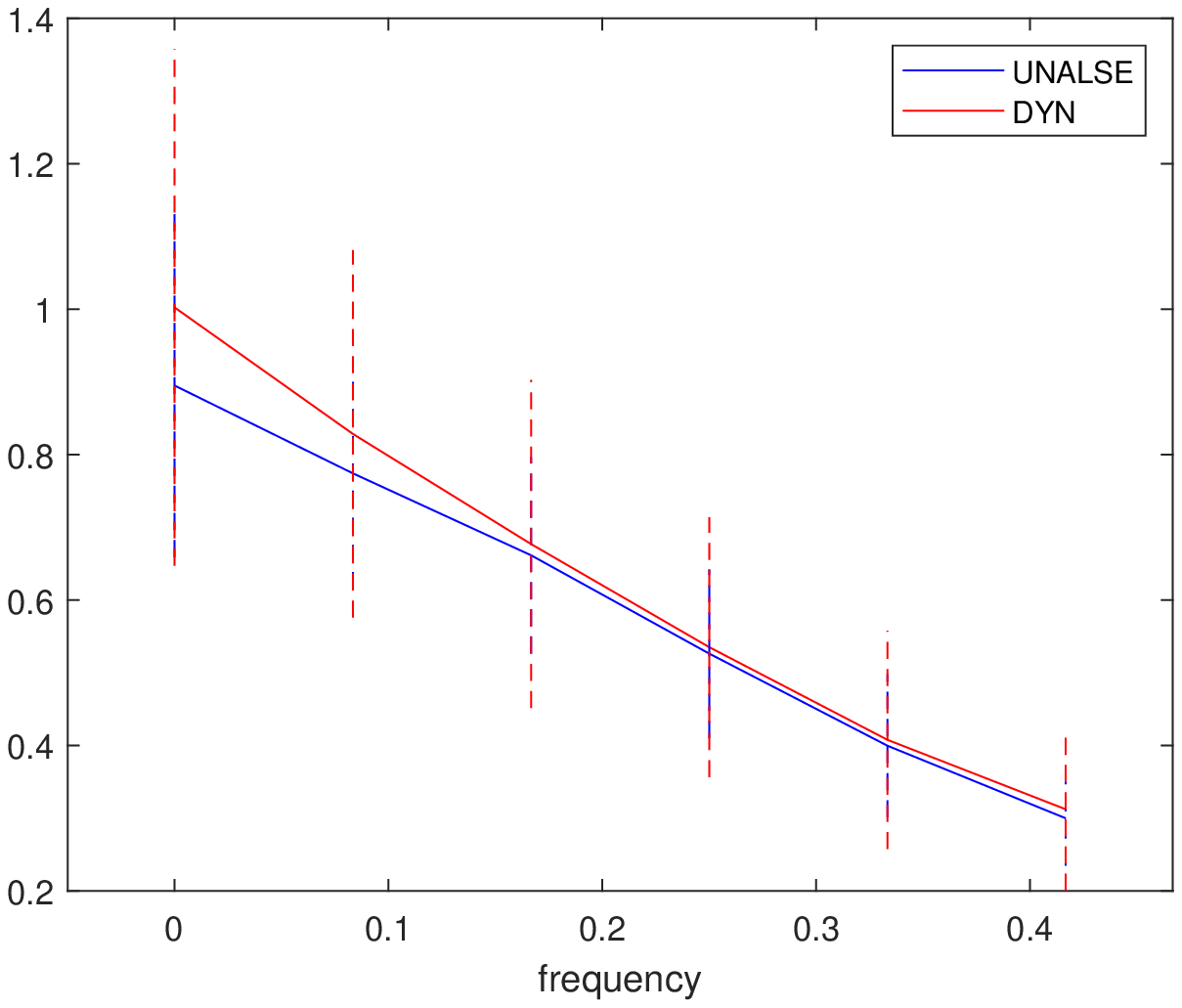}\\
              \end{tabular}
\end{figure}

\begin{figure}[h!]
          \caption{$err_{ratio}(f_h)$ - Scenario C.}\label{fig:err_ratio_C.1.U}
                  \centering
          \begin{tabular}{ccc}
          {\footnotesize Setting 1}&{\footnotesize Setting 2}&{\footnotesize Setting 3}\\
              \includegraphics[width=.2\textwidth]{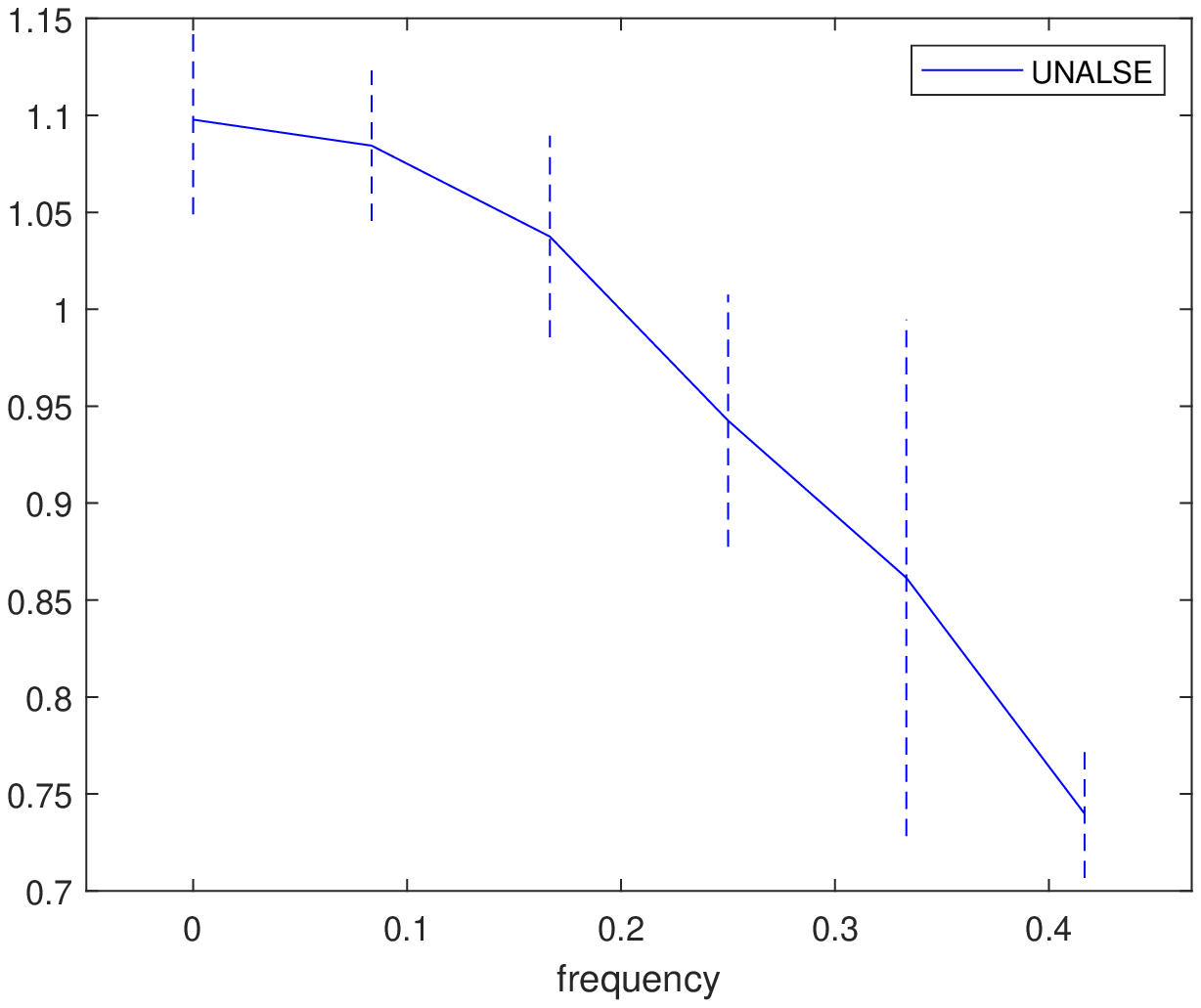}&
              \includegraphics[width=.2\textwidth]{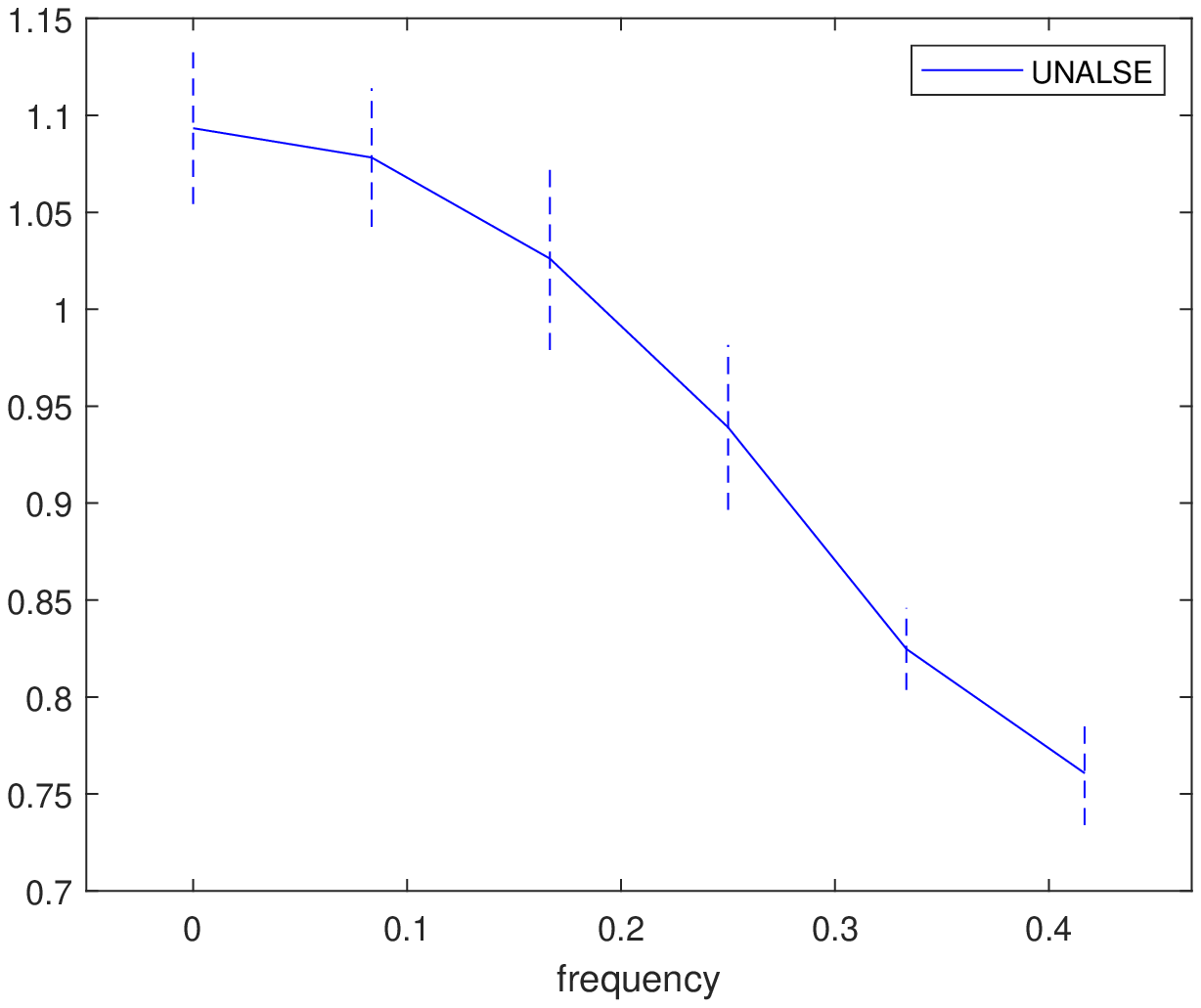}&
              \includegraphics[width=.2\textwidth]{err_ratio_max_3}\\
              \end{tabular}
               \begin{tabular}{cc}
	      {\footnotesize Setting 4}&{\footnotesize Setting 5}\\
              \includegraphics[width=.2\textwidth]{err_ratio_max_4}&
              \includegraphics[width=.2\textwidth]{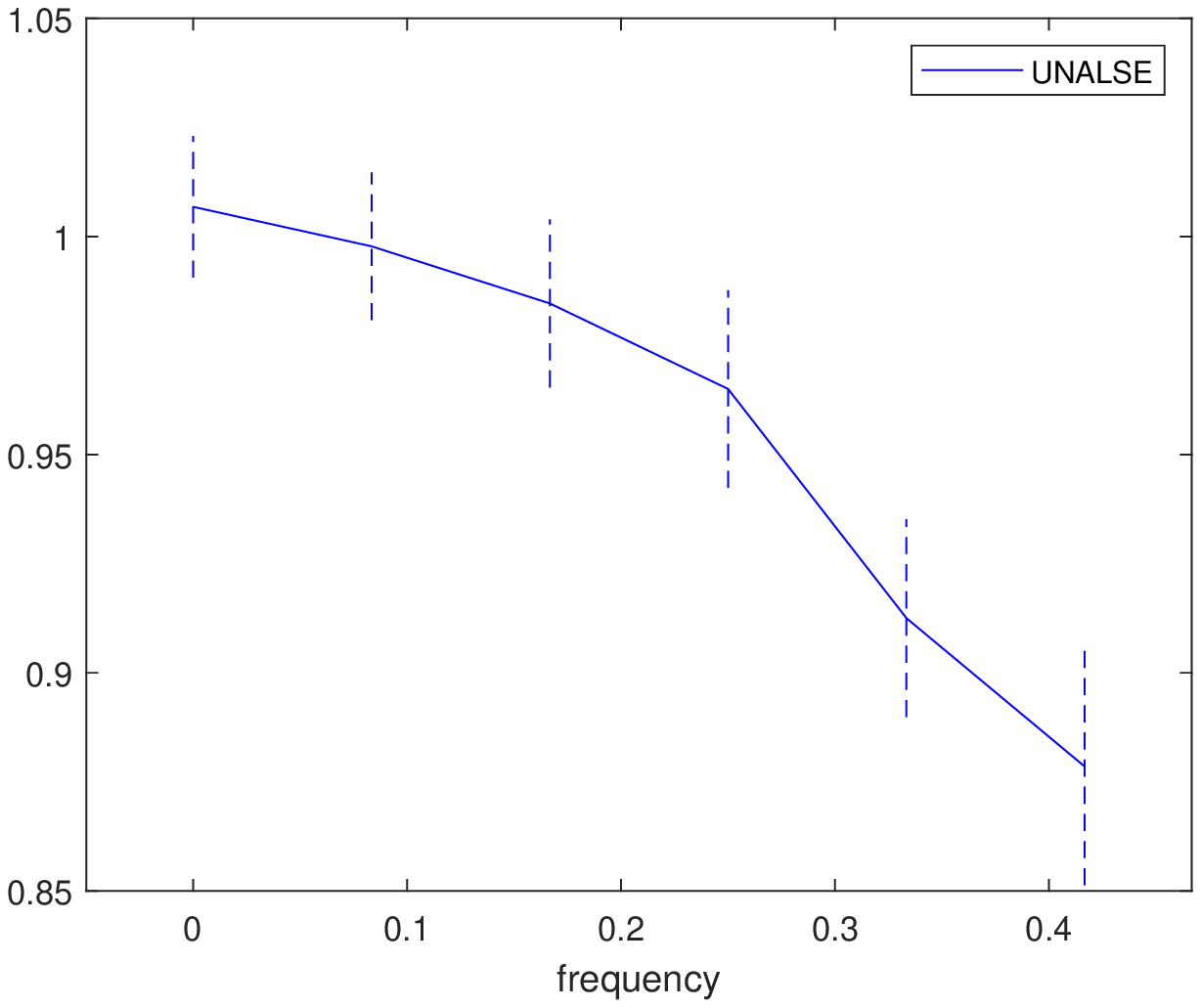}\\
              \end{tabular}
\end{figure}
\clearpage
\singlespace
\bibliographystyle{apalike}

\bibliography{spec_bib}

\begin{thebibliography}{}

\bibitem[Altissimo et~al., 2010]{altissimo2010new}
Altissimo, F., Cristadoro, R., Forni, M., Lippi, M., and Veronese, G. (2010).
\newblock New {Eurocoin}: Tracking economic growth in real time.
\newblock {\em The Review of Economics and Statistics}, 92(4):1024--1034.

\bibitem[Barigozzi and Brownlees, 2019]{barigozzi2019nets}
Barigozzi, M. and Brownlees, C. (2019).
\newblock {NETS}: Network estimation for time series.
\newblock {\em Journal of Applied Econometrics}, 34(3):347--364.

\bibitem[Barigozzi and Hallin, 2017]{barigozzi2017network}
Barigozzi, M. and Hallin, M. (2017).
\newblock A network analysis of the volatility of high dimensional financial
  series.
\newblock {\em Journal of the Royal Statistical Society: Series C (Applied
  Statistics)}, 66(3):581--605.

\bibitem[Barigozzi et~al., 2021]{barigozzi2021large}
Barigozzi, M., Lippi, M., and Luciani, M. (2021).
\newblock Large-dimensional dynamic factor models: Estimation of
  impulse--response functions with {$I (1)$} cointegrated factors.
\newblock {\em Journal of Econometrics}, 221(2):455--482.

\bibitem[Bickel and Levina, 2008]{bickel2008covariance}
Bickel, P.~J. and Levina, E. (2008).
\newblock Covariance regularization by thresholding.
\newblock {\em The Annals of Statistics}, 36(6):2577--2604.

\bibitem[B{\"o}hm and von Sachs, 2008]{bohm2008structural}
B{\"o}hm, H. and von Sachs, R. (2008).
\newblock Structural shrinkage of nonparametric spectral estimators for
  multivariate time series.
\newblock {\em Electronic Journal of Statistics}, 2:696--721.

\bibitem[B{\"o}hm and von Sachs, 2009]{bohm2009shrinkage}
B{\"o}hm, H. and von Sachs, R. (2009).
\newblock Shrinkage estimation in the frequency domain of multivariate time
  series.
\newblock {\em Journal of Multivariate Analysis}, 100(5):913--935.

\bibitem[Breitung and Candelon, 2006]{breitung2006testing}
Breitung, J. and Candelon, B. (2006).
\newblock Testing for short-and long-run causality: A frequency-domain
  approach.
\newblock {\em Journal of Econometrics}, 132(2):363--378.

\bibitem[Brillinger, 2001]{brillinger2001time}
Brillinger, D.~R. (2001).
\newblock {\em {Time Series: Data Analysis and Theory}}.
\newblock SIAM.

\bibitem[Cai et~al., 2010]{cai2010singular}
Cai, J.-F., Cand{\`e}s, E.~J., and Shen, Z. (2010).
\newblock A singular value thresholding algorithm for matrix completion.
\newblock {\em SIAM Journal on Optimization}, 20(4):1956--1982.

\bibitem[Chandrasekaran et~al., 2012]{chandrasekaran2012}
Chandrasekaran, V., Parrilo, P.~A., and Willsky, A.~S. (2012).
\newblock Latent variable graphical model selection via convex optimization.
\newblock {\em The Annals of Statistics}, 40(4):1935--1967.

\bibitem[Chandrasekaran et~al., 2011]{chandrasekaran2011rank}
Chandrasekaran, V., Sanghavi, S., Parrilo, P.~A., and Willsky, A.~S. (2011).
\newblock Rank-sparsity incoherence for matrix decomposition.
\newblock {\em SIAM Journal on Optimization}, 21(2):572--596.

\bibitem[Chaudhuri and Lo, 2015]{chaudhuri2015spectral}
Chaudhuri, S.~E. and Lo, A.~W. (2015).
\newblock Spectral analysis of stock-return volatility, correlation, and beta.
\newblock In {\em 2015 IEEE Signal Processing and Signal Processing Education
  Workshop (SP/SPE)}, pages 232--236. IEEE.

\bibitem[Corbae et~al., 2002]{corbae2002band}
Corbae, D., Ouliaris, S., and Phillips, P.~C. (2002).
\newblock Band spectral regression with trending data.
\newblock {\em Econometrica}, 70(3):1067--1109.

\bibitem[Dahlhaus, 2000a]{dahlhaus2000graphical}
Dahlhaus, R. (2000a).
\newblock Graphical interaction models for multivariate time series.
\newblock {\em Metrika}, 51(2):157--172.

\bibitem[Dahlhaus, 2000b]{dahlhaus2000likelihood}
Dahlhaus, R. (2000b).
\newblock A likelihood approximation for locally stationary processes.
\newblock {\em The Annals of Statistics}, 28(6):1762--1794.

\bibitem[Daubechies et~al., 2004]{daubechies2004iterative}
Daubechies, I., Defrise, M., and De~Mol, C. (2004).
\newblock An iterative thresholding algorithm for linear inverse problems with
  a sparsity constraint.
\newblock {\em Communications on Pure and Applied Mathematics},
  57(11):1413--1457.

\bibitem[Davis et~al., 2016]{davis2016sparse}
Davis, R.~A., Zang, P., and Zheng, T. (2016).
\newblock Sparse vector autoregressive modeling.
\newblock {\em Journal of Computational and Graphical Statistics},
  25(4):1077--1096.

\bibitem[Donoho, 2006]{donoho2006most}
Donoho, D.~L. (2006).
\newblock For most large underdetermined systems of linear equations the
  minimal $l_1$ norm solution is also the sparsest solution.
\newblock {\em Communications on Pure and Applied Mathematics}, 59(6):797--829.

\bibitem[Eichler, 2007]{eichler2007granger}
Eichler, M. (2007).
\newblock Granger causality and path diagrams for multivariate time series.
\newblock {\em Journal of Econometrics}, 137(2):334--353.

\bibitem[Fan et~al., 2013]{fan2013large}
Fan, J., Liao, Y., and Mincheva, M. (2013).
\newblock Large covariance estimation by thresholding principal orthogonal
  complements.
\newblock {\em Journal of the Royal Statistical Society: Series B (Statistical
  Methodology)}, 75(4):603--680.

\bibitem[Farn{\'e}, 2016]{farne2016algorithm}
Farn{\'e}, M. (2016).
\newblock An algorithm to simulate {VMA} processes having a spectrum with fixed
  condition number.
\newblock {\em Communications in Statistics-Simulation and Computation},
  45(5):1664--1675.

\bibitem[Farn{\`e} and Montanari, 2020]{farne2020large}
Farn{\`e}, M. and Montanari, A. (2020).
\newblock A large covariance matrix estimator under intermediate spikiness
  regimes.
\newblock {\em Journal of Multivariate Analysis}, 176:104577.

\bibitem[Farn{\`e} and Montanari, 2021]{farne2018bootstrap}
Farn{\`e}, M. and Montanari, A. (2021).
\newblock A bootstrap method to test {G}ranger-causality in the frequency
  domain.
\newblock {\em Computational Economics}.

\bibitem[Fazel et~al., 2001]{fazel2001rank}
Fazel, M., Hindi, H., and Boyd, S.~P. (2001).
\newblock A rank minimization heuristic with application to minimum order
  system approximation.
\newblock In {\em American Control Conference, 2001. Proceedings of the 2001},
  volume~6, pages 4734--4739. IEEE.

\bibitem[Fiecas et~al., 2019]{fiecas2019spectral}
Fiecas, M., Leng, C., Liu, W., and Yu, Y. (2019).
\newblock Spectral analysis of high-dimensional time series.
\newblock {\em Electronic Journal of Statistics}, 13(2):4079--4101.

\bibitem[Fiecas and Ombao, 2011]{fiecas2011generalized}
Fiecas, M. and Ombao, H. (2011).
\newblock The generalized shrinkage estimator for the analysis of functional
  connectivity of brain signals.
\newblock {\em The Annals of Applied Statistics}, 5(2A):1102--1125.

\bibitem[Fiecas and Ombao, 2016]{fiecas2016modeling}
Fiecas, M. and Ombao, H. (2016).
\newblock Modeling the evolution of dynamic brain processes during an
  associative learning experiment.
\newblock {\em Journal of the American Statistical Association},
  111(516):1440--1453.

\bibitem[Fiecas and von Sachs, 2014]{fiecas2014data}
Fiecas, M. and von Sachs, R. (2014).
\newblock Data-driven shrinkage of the spectral density matrix of a
  high-dimensional time series.
\newblock {\em Electronic Journal of Statistics}, 8(2):2975--3003.

\bibitem[Forni et~al., 2000]{forni2000generalized}
Forni, M., Hallin, M., Lippi, M., and Reichlin, L. (2000).
\newblock The generalized dynamic-factor model: Identification and estimation.
\newblock {\em The Review of Economics and Statistics}, 82(4):540--554.

\bibitem[Forni et~al., 2005]{forni2005generalized}
Forni, M., Hallin, M., Lippi, M., and Reichlin, L. (2005).
\newblock The generalized dynamic factor model: one-sided estimation and
  forecasting.
\newblock {\em Journal of the American Statistical Association},
  100(471):830--840.

\bibitem[Forni et~al., 2017]{forni2017dynamic}
Forni, M., Hallin, M., Lippi, M., and Zaffaroni, P. (2017).
\newblock Dynamic factor models with infinite-dimensional factor space:
  asymptotic analysis.
\newblock {\em Journal of Econometrics}, 199(1):74--92.

\bibitem[Forni and Lippi, 2001]{forni2001generalized}
Forni, M. and Lippi, M. (2001).
\newblock The generalized dynamic factor model: representation theory.
\newblock {\em Econometric theory}, 17(6):1113--1141.

\bibitem[Giannone et~al., 2017]{giannone2017economic}
Giannone, D., Lenza, M., and Primiceri, G.~E. (2017).
\newblock Economic predictions with big data: The illusion of sparsity.
\newblock CEPR discussion paper 12256.

\bibitem[Granger, 1969]{granger1969investigating}
Granger, C.~W. (1969).
\newblock Investigating causal relations by econometric models and
  cross-spectral methods.
\newblock {\em Econometrica}, 37(3):424--438.

\bibitem[Hallin and Lippi, 2013]{hallin2013factor}
Hallin, M. and Lippi, M. (2013).
\newblock Factor models in high-dimensional time series—a time-domain
  approach.
\newblock {\em Stochastic Processes and their Applications}, 123(7):2678--2695.

\bibitem[Hallin and Li{\v{s}}ka, 2007]{hallin2007determining}
Hallin, M. and Li{\v{s}}ka, R. (2007).
\newblock Determining the number of factors in the general dynamic factor
  model.
\newblock {\em Journal of the American Statistical Association},
  102(478):603--617.

\bibitem[Harvey, 1978]{harvey1978linear}
Harvey, A.~C. (1978).
\newblock Linear regression in the frequency domain.
\newblock {\em International Economic Review}, 19(2):507--512.

\bibitem[Joyeux, 1992]{joyeux1992tests}
Joyeux, R. (1992).
\newblock Tests for seasonal cointegration using principal components.
\newblock {\em Journal of Time Series Analysis}, 13(2):109--118.

\bibitem[Luo, 2011a]{luo2011high}
Luo, X. (2011a).
\newblock High dimensional low rank and sparse covariance matrix estimation via
  convex minimization.
\newblock {\em arXiv:1111.1133}.

\bibitem[Luo, 2011b]{luo2011recovering}
Luo, X. (2011b).
\newblock Recovering model structures from large low rank and sparse covariance
  matrix estimation.
\newblock {\em arXiv:1111.1133}.

\bibitem[Mar{\v{c}}enko and Pastur, 1967]{marvcenko1967distribution}
Mar{\v{c}}enko, V.~A. and Pastur, L.~A. (1967).
\newblock Distribution of eigenvalues for some sets of random matrices.
\newblock {\em Mathematics of the USSR-Sbornik}, 1(4):457.

\bibitem[McCracken and Ng, 2020]{mccracken2020fred}
McCracken, M. and Ng, S. (2020).
\newblock {FRED-QD}: {A} quarterly database for macroeconomic research.
\newblock Technical Report 26872, National Bureau of Economic Research.

\bibitem[M{\"u}ller and Watson, 2018]{muller2018long}
M{\"u}ller, U.~K. and Watson, M.~W. (2018).
\newblock Long-run covariability.
\newblock {\em Econometrica}, 86(3):775--804.

\bibitem[Ombao et~al., 2005]{ombao2005slex}
Ombao, H., Von~Sachs, R., and Guo, W. (2005).
\newblock {SLEX} analysis of multivariate nonstationary time series.
\newblock {\em Journal of the American Statistical Association},
  100(470):519--531.

\bibitem[Ombao et~al., 2001]{ombao2001automatic}
Ombao, H.~C., Raz, J.~A., von Sachs, R., and Malow, B.~A. (2001).
\newblock Automatic statistical analysis of bivariate nonstationary time
  series.
\newblock {\em Journal of the American Statistical Association},
  96(454):543--560.

\bibitem[Onatski, 2009]{onatski2009testing}
Onatski, A. (2009).
\newblock Testing hypotheses about the number of factors in large factor
  models.
\newblock {\em Econometrica}, 77(5):1447--1479.

\bibitem[Priestley, 1981]{priestley1981spectral}
Priestley, M.~B. (1981).
\newblock {\em Spectral analysis and time series}.
\newblock Academic Press.

\bibitem[Sargent and Sims, 1977]{sargent1977business}
Sargent, T. and Sims, C. (1977).
\newblock Business cycle modeling without pretending to have too much a priori
  economic theory.
\newblock Technical report, Federal Reserve Bank of Minneapolis.

\bibitem[Stock and Watson, 1988]{stock1988testing}
Stock, J.~H. and Watson, M.~W. (1988).
\newblock Testing for common trends.
\newblock {\em Journal of the American statistical Association},
  83(404):1097--1107.

\bibitem[Velasco and Robinson, 2000]{velasco2000whittle}
Velasco, C. and Robinson, P.~M. (2000).
\newblock Whittle pseudo-maximum likelihood estimation for nonstationary time
  series.
\newblock {\em Journal of the American Statistical Association},
  95(452):1229--1243.

\bibitem[Wu and Zaffaroni, 2018]{wu2018asymptotic}
Wu, W.~B. and Zaffaroni, P. (2018).
\newblock Asymptotic theory for spectral density estimates of general
  multivariate time series.
\newblock {\em Econometric Theory}, 34:1--22.

\bibitem[Zhang and Wu, 2021]{zhang2021convergence}
Zhang, D. and Wu, W.~B. (2021).
\newblock Convergence of covariance and spectral density estimates for
  high-dimensional locally stationary processes.
\newblock {\em The Annals of Statistics}, 49(1):233--254.

\end{thebibliography}

\end{document}